\documentclass[titlepage,11pt]{article}
\usepackage{latexsym}
\usepackage{amsfonts}
\usepackage{amsmath}
\newcommand\blackslug{\hbox{\hskip 1pt \vrule width 4pt height 8pt depth 1.5pt
        \hskip 1pt}}
\newcommand\bbox{\hfill \quad \blackslug \medbreak}
\def\d{\hbox{-}}
\def\c{\hbox{-}\cdots\hbox{-}}
\def\l{,\ldots,}

\title{The Strong Perfect Graph Theorem}
\author{Maria Chudnovsky\\
Princeton University, Princeton NJ 08544
\and
Neil Robertson\thanks{Supported by ONR grant N00014-01-1-0608, NSF grant DMS-0071096, and AIM.}\\
Ohio State University, Columbus, Ohio 43210
\and
Paul Seymour\thanks{Supported by ONR grants N00014-97-1-0512 and N00014-01-1-0608, and NSF grant DMS-0070912.}\\
Princeton University, Princeton NJ 08544
\and
Robin Thomas\thanks{Supported by ONR grant N00014-01-1-0608, NSF grants DMS-9970514 and DMS-0200595, and AIM.}\\
Georgia Institute of Technology, Atlanta, GA 30332}

\date{June 20, 2002; revised \today}

\newtheorem{thm}{}[section]

\newcommand{\Proof}{\noindent{\bf Proof.}\ \ }

\begin{document}
\maketitle
\begin{abstract}
A graph $G$ is {\em perfect}
if for every induced subgraph $H$, the chromatic number of $H$ equals the
size of the largest complete subgraph of $H$, and $G$ is {\em Berge} if
no induced subgraph of $G$ is an odd cycle of length at least $5$ or the
complement of one.
The ``strong perfect graph conjecture'' (Berge, 1961) asserts that a graph
is perfect if and only if it is Berge.

A stronger conjecture was made recently by Conforti, Cornu\'{e}jols and  Vu\v{s}kovi\'{c} --- that every
Berge graph either falls into one of a few basic classes, or it has a kind
of separation that cannot occur in a minimal imperfect graph.

In this paper we prove both these conjectures.
\end{abstract}

\section{Introduction}

We begin with definitions of some of the terms we use which may be nonstandard.
All graphs in this paper are finite and simple.
The {\em complement} $\overline{G}$ of a graph $G$ has the same vertex set as $G$,
and distinct vertices $u,v$ are adjacent in $\overline{G}$ just when they are
not adjacent in $G$. A {\em hole} of $G$ is an induced subgraph of $G$ which
is a cycle of length $>3$. An {\em antihole} of $G$ is an induced subgraph of
$G$ whose complement is a hole in $\overline{G}$. A graph $G$ is {\em Berge}
if every hole and antihole of $G$ has even length.

A {\em clique} in $G$ is a subset $X$ of $V(G)$ so that every two members of $X$
are adjacent. A graph $G$ is {\em perfect}
if for every induced subgraph $H$ of $G$, the chromatic number of $H$ equals the
size of the largest clique of $H$. In 1961 Claude Berge~\cite{Berge} proposed the so-called
{\em Strong Perfect Graph Conjecture}, the main theorem of this paper:

\begin{thm}
\label{SPGC}
A graph is Berge if and only if it is perfect.
\end{thm}

It is easy to see that every perfect graph is Berge, and so to prove (\ref{SPGC})
it remains to prove the converse. This has received a great deal of attention over the
past 40 years, but so far has resisted solution. Most of the previous approaches on
(\ref{SPGC}) fall into two classes: proving that the theorem holds for graphs with
some particular graph excluded as an induced subgraph (there are a number of these for
different subgraphs, but such an approach obviously cannot do the whole thing), and using
linear programming methods to investigate the structure of a minimal counterexample.
(There are rich connections with linear and integer programming - see \cite{Reed} for
example.)

Our approach is different. Recently, Conforti, Cornu\'{e}jols and Vu\v{s}kovi\'{c}~\cite{CCV} conjectured that
every Berge graph either falls into one of four well-understood classes, or it admits one of
several kinds of decomposition. They pointed out that if this could be proved, and if also
it could be shown that no minimal counterexample to (\ref{SPGC}) admits any such decomposition,
then \ref{SPGC} would follow (for certainly no minimal counterexample to \ref{SPGC}
can fall into the four basic classes). We have been able to prove both (except we
need a fifth class).

Before we can be more precise we need more definitions. If $X \subseteq V(G)$ we denote
the subgraph of $G$ induced on $X$ by $G|X$. The {\em line graph} $L(G)$
of a graph $G$ has vertex set the set $E(G)$ of edges of $G$, and $e,f \in E(G)$ are adjacent
in $L(G)$ if they share an end in $G$.

We need one other class of graphs, defined as follows. Let $m,n \ge 2$ be integers, and let $\{a_1\l a_m\}$,
$\{b_1\l b_m\}$, $\{c_1\l c_n\}$, $\{d_1\l d_n\}$ be disjoint sets. Let $G$ have vertex
set their union, and edges as follows:
\begin{itemize}
\item $a_i$ is adjacent to $b_i$ for $1 \le i \le m$, and $c_j$ is nonadjacent to $d_j$ for
$1 \le j \le n$
\item there are no edges between $\{a_i,b_i\}$ and $\{a_{i'},b_{i'}\}$ for $1 \le i < i' \le m$, and all four
edges between  $\{c_j,d_j\}$ and $\{c_{j'},d_{j'}\}$ for $1 \le j < j' \le n$
\item there are exactly two edges between $\{a_i,b_i\}$ and $\{c_j,d_j\}$ for $1 \le i \le m$ and
$1 \le j \le n$, and these two edges are disjoint.
\end{itemize}
We call such a graph $G$ a {\em bicograph}. Let us say a graph $G$
is {\em basic} if either $G$ or $\overline{G}$ is bipartite or is the line graph of a
bipartite graph, or is a bicograph. (Note that if $G$ is a bicograph then so is $\overline{G}$.)
It is easy to see that all basic graphs are perfect.

Now we turn to the various kinds of decomposition that we need. First, a decomposition
due to Cornu\'{e}jols and Cunningham~\cite{2-join} ---
a {\em 2-join} in $G$ is a partition $(X_1,X_2)$ of $V(G)$ so that there exist
disjoint nonempty $A_i,B_i \subseteq X_i (i = 1,2)$ satisfying:
\begin{itemize}
\item every vertex of $A_1$ is
adjacent to every vertex of $A_2$, and every vertex of $B_1$ is
adjacent to every vertex of $B_2$,
\item there are no other edges between $X_1$ and $X_2$,
\item for $i = 1,2$, every component of $G|X_i$ meets both $A_i$ and $B_i$, and
\item for $i = 1,2$, if $|A_i| = |B_i| = 1$ and $G|X_i$ is a path joining the
members of $A_i$ and $B_i$, then it has length $\ge 3$.
\end{itemize}

If $X, Y \subseteq V(G)$ are disjoint, we say $X$ is {\em complete} to $Y$ (or the pair
$(X,Y)$ is {\em complete}) if every vertex
in $X$ is adjacent to every vertex in $Y$; and we say $X$ is {\em anticomplete} to $Y$
if there are no edges between $X$ and $Y$. Our second decomposition is a very slight variant
on the ``homogenous sets'' due to Chv\'{a}tal and Sbihi~\cite{M-join} ---
an {\em M-join} in $G$ is a partition of $V(G)$ into six nonempty sets, $(A,B,C,D,E,F)$,
so that:
\begin{itemize}
\item
every vertex in $A$ has a neighbour in $B$ and a nonneighbour in $B$, and vice versa
\item
the pairs $(C,A),(A,F),(F,B),(B,D)$ are complete, and
\item
the pairs $(D,A), (A,E), (E,B),(B,C)$ are anticomplete.
\end{itemize}

A {\em path} in $G$ is an {\em induced} subgraph of $G$ which is non-null, connected
and in which every vertex has degree $\le 2$
(this definition is highly nonstandard, and we apologise, but it avoids writing
``induced'' about 600 times), and an {\em antipath} is an induced subgraph
whose complement is a path. The {\em length} of a path is the number of edges
in it (and the length of an antipath is the number of edges in its complement).
We therefore recognize paths and antipaths of length $0$.
If $P$ is a path, $P^*$ denotes the set of internal vertices of $P$, called
the {\em interior} of $P$; and similarly for antipaths.
Let $A,B$ be disjoint subsets of $V(G)$. We say the pair $(A,B)$ is {\em balanced} if
there is no odd path between nonadjacent vertices in $B$ with interior in $A$, and
there is no odd antipath between adjacent vertices in $A$ with interior in $B$.
A set $X \subseteq V(G)$ is {\em connected} if $G|X$ is connected (so $\emptyset$ is connected); and
{\em anticonnected} if $\overline{G}|X$ is connected.

The third kind of decomposition we use is due to Chv\'{a}tal~\cite{starcut} ---
a {\em skew partition} in $G$ is a partition $(A,B)$ of $V(G)$ so that
$A$ is not connected and $B$ is not anticonnected. Skew partitions pose a difficulty that
the other two decompositions do not, for it has not been shown that a minimal counterexample
to \ref{SPGC} cannot admit a skew partition; indeed, this is a well-known open question,
first raised by Chv\'{a}tal~\cite{starcut}, the so-called ``skew partition conjecture''.
We get around it by confining ourselves to balanced skew partitions, which do not present
this difficulty.

We shall prove the following, a form of which was conjectured in \cite{CCV}:

\begin{thm}
\label{decomp}
For every Berge graph $G$, either $G$ is basic, or one of $G$, $\overline{G}$ admits
a 2-join, or $G$ admits an M-join, or $G$ admits a balanced skew partition.
\end{thm}
\noindent
{\bf Proof of \ref{SPGC}, assuming \ref{decomp} and \ref{evenskew}.}
\\
Suppose that \ref{SPGC} is false, and let $G$ be a counterexample with $|V(G)|$ as small
as possible. Since every perfect graph is Berge, it follows that $G$ is Berge and not
perfect. Every basic graph is perfect, and so $G$ is not basic. It is shown in \cite{2-join} that
$G$ does not admit a 2-join. Since Lov\'{a}sz~\cite{Lovasz} showed that the complement
of a perfect graph is perfect, it follows that $\overline{G}$ is also a counterexample
to \ref{SPGC} of minimum size, and therefore $\overline{G}$ also does not admit a
2-join. It is shown in \cite{M-join} that $G$ does not admit an M-join, and we shall
prove in \ref{evenskew}
that $G$ does not admit a balanced skew partition. It follows that $G$ violates \ref{decomp},
and therefore there is no such graph $G$. This proves \ref{SPGC}. \bbox

All nontrivial bicographs admit skew partitions, so if we delete ``balanced'' from \ref{decomp} then
we no longer need to consider bicographs as basic --- four basic classes suffice. Unfortunately,
nontrivial bicographs do not admit balanced skew partitions, and general skew partitions are not
good enough for the application to \ref{SPGC}, so we have to do it the way we did.

How can we prove a theorem of the form of \ref{decomp}? There are several other theorems
of this kind in graph theory - eg \cite{Tutte}, \cite{tu}, \cite{Pfaffian}, \cite{Wagner}
 and others. All these theorems say that ``every graph (or matroid) not containing an object
of type X either falls into one of a few basic classes or admits a decomposition''.
And for each of these theorems, the technique of the proof is the same: we judiciously choose
a small X-free graph $H$ ({\em X-free} means not containing an object of type X)
which does not fall into any of the basic classes;
check that it has a decomposition of the kind it is supposed to have; and show that this
decomposition extends to a decomposition of every bigger X-free graph containing $H$.
That proves that the theorem is true for all X-free graphs that contain $H$, so now we
focus on the X-free graphs that do not contain $H$. Repeat as often as necessary (with different
small graphs in place of $H$) until the graphs that remain are sufficiently restricted
that they can be handled by some other means.
We shall see that the same technique also works for decomposing Berge graphs.

The paper is organized as follows. The next three sections develop tools that will be needed all through the
paper. Section 2 concerns a fundamental lemma of Roussel and Rubio; we give several variations and
extensions of it, and more in section 3,
of a different kind. In section 4 we develop some features of skew partitions, to make them
easier to handle in the main proof, and in particular, prove that no minimum imperfect graph admits
a balanced skew partition. In section 5 we begin on the main proof. Sections 5-8 prove that every
Berge graph containing a ``substantial'' line graph as an induced subgraph, satisfies \ref{decomp}
(``substantial'' means a line graph of a bipartite subdivision of a 3-connected graph $J$, with some
more conditions if $J = K_4$). Section 9 proves the same thing for line graphs of subdivisions
of $K_4$ that are not ``substantial'' --- this is where bicographs come in. In section 10 we prove that
Berge graphs containing an ``even prism'' satisfy \ref{decomp}. (To prove this we may assume we are looking
at a Berge graph that does not contain the line graph of a subdivision of $K_4$, for otherwise we
could apply the results of the earlier sections. The same thing happens later - at each step we
may assume the current Berge graph does not contain any of the subgraphs that were handled in
earlier steps.) Sections 11-13 do the same for ``long odd prisms'', and section 14 does the same for a subgraph
we call the ``double diamond''. Section 15 is a break for resharpening the tools we proved in the first four
sections, and in particular, here we prove that no minimum
imperfect graph admits a skew partition. Section 16 proves that any Berge graph containing
what we call an ``odd wheel'' satisfies \ref{decomp}, in sections 17-23 we prove the same for wheels in
general, and finally in section 24 we handle Berge graphs not containing wheels.

These steps are summarized more precisely in the next theorem, which we include now
in the hope that it will be helpful to the reader, although some necessary definitions have not been
given yet --- for the missing definitions, the reader should see the appropriate section(s) later.
Let $\mathcal{F}_1\l \mathcal{F}_{11}$ be the classes of Berge graphs defined as follows
(each is a subclass of the previous class):
\begin{itemize}
\item $\mathcal{F}_1$ is the class of all Berge graphs $G$ such that for every bipartite subdivision
$H$ of $K_4$, every induced subgraph of $G$ isomorphic to $L(H)$ is degenerate
\item $\mathcal{F}_2$ is the class of all graphs $G$ such that $G, \overline{G} \in \mathcal{F}_1$ and
no induced subgraph of $G$ is isomorphic to $L(K_{3,3})$
\item $\mathcal{F}_3$ is the class of all Berge graphs $G$ so that for every bipartite subdivision
$H$ of $K_4$, no induced subgraph of $G$ or of $\overline{G}$ is isomorphic to $L(H)$
\item $\mathcal{F}_4$ is the class of all $G \in \mathcal{F}_3$ so that no induced subgraph
of $G$ is an even prism
\item $\mathcal{F}_5$ is the class of all $G \in \mathcal{F}_3$ so that no induced subgraph
of $G$ or of $\overline{G}$ is a long prism
\item $\mathcal{F}_6$ is the class of all $G \in \mathcal{F}_5$ such that no induced subgraph of $G$
is isomorphic to a double diamond
\item $\mathcal{F}_7$ is the class of all $G \in \mathcal{F}_6$ so that $G$ and $\overline{G}$
do not contain odd wheels
\item $\mathcal{F}_8$ is the class of all $G \in \mathcal{F}_7$ so that $G$ and $\overline{G}$
do not contain pseudowheels
\item $\mathcal{F}_9$ is the class of all $G \in \mathcal{F}_8$ such that $G$ and $\overline{G}$
do not contain wheels
\item $\mathcal{F}_{10}$ is the class of all $G \in \mathcal{F}_9$ such that, for every hole $C$ in $G$
of length $\ge 6$, no vertex of $G$ has three consecutive neighbours in $C$, and the same holds
in $\overline{G}$
\item $\mathcal{F}_{11}$ is the class of all $G \in \mathcal{F}_{10}$ such that every antihole in $G$ has
length $4$.
\end{itemize}

\begin{thm}\label{summary}
The following are the main steps of the proof of \ref{decomp}:
\begin{enumerate}
\item For every Berge graph $G$, either $G$ is a line graph of a bipartite graph, or $G$ admits a 2-join
or a balanced skew partition, or $G \in \mathcal{F}_1$; and (consequently), either one of $G,\overline{G}$
is a line graph of a bipartite graph, or one of $G,\overline{G}$
admits a 2-join, or $G$ admits a balanced skew partition, or $G, \overline{G} \in \mathcal{F}_1$
\item For every $G$ with $G,\overline{G} \in \mathcal{F}_1$, either $G = L(K_{3,3})$, or one of $G, \overline{G}$
admits a 2-join, or $G$ admits a balanced skew partition, or $G \in \mathcal{F}_2$
\item For every $G \in \mathcal{F}_2$,  either $G$ is a bicograph, or one of $G, \overline{G}$ admits a
2-join, or $G$ admits a balanced skew partition, or $G \in \mathcal{F}_3$
\item For every $G\in \mathcal{F}_1$,  either $G$ is an even prism with exactly $9$ vertices, or
$G$ admits a 2-join or a balanced skew partition, or $G \in \mathcal{F}_4$
\item For every $G \in \mathcal{F}_3$,  either one of $G, \overline{G}$
admits a 2-join, or $G$ admits an M-join, or $G$ admits a balanced skew partition, or $G \in \mathcal{F}_5$
\item For every $G \in \mathcal{F}_5$,  either one of $G, \overline{G}$
admits a 2-join, or $G$ admits a balanced skew partition, or $G \in \mathcal{F}_6$
\item For every $G \in \mathcal{F}_6$, either $G$ admits a balanced skew partition, or $G \in \mathcal{F}_7$
\item For every $G \in \mathcal{F}_7$, either $G$ admits a balanced skew partition, or $G \in \mathcal{F}_8$
\item For every $G \in \mathcal{F}_8$, either $G$ admits a balanced skew partition, or $G \in \mathcal{F}_9$
\item For every $G \in \mathcal{F}_9$, either $G$ admits a balanced skew partition, or $G \in \mathcal{F}_{10}$
\item For every graph $G \in \mathcal{F}_{10}$, either $G \in \mathcal{F}_{11}$ or
$\overline{G} \in \mathcal{F}_{11}$
\item For every graph $G \in \mathcal{F}_{11}$, either $G$ admits a balanced skew partition, or $G$ is complete or
bipartite.
\end{enumerate}
\end{thm}

The twelve statements of \ref{summary} are proved in \ref{linegraph}, \ref{linegraph2.5}, \ref{bicographs0},
 \ref{evenprism}, \ref{longprism}, \ref{cube}, \ref{oddwheel}, \ref{pseudowheel}, \ref{nowheel},
 \ref{notriad}, \ref{noantihole}, and \ref{bipartite} respectively.

\section{The Roussel-Rubio lemma}

There is a beautiful and very powerful theorem of \cite{R&R} which we use
many times throughout the paper. (We proved it independently, in joint work with Carsten
Thomassen, but Roussel and Rubio found it earlier.)
Its main use is to show that in some respects, the common neighbours
of an anticonnected set of vertices (in a Berge graph) act like or almost like the neighbours of
a single vertex.

If $X \subseteq V(G)$ and $v \in V(G)$, we
say $v$ is $X$-$complete$ if $v$ is adjacent to every vertex in $X$ (and
consequently $v \notin X$), and an edge $uv$ is $X$-complete if $u,v$ are both $X$-complete.
Let $P$ be a path in $G$ (we remind the reader
that this means $P$ is an induced subgraph which is a path), of length
$\ge 2$, and let the vertices of $P$ be $p_1,\ldots,p_n$ in order. A {\em leap}
for $P$ ({\em in} $G$) is a pair of nonadjacent vertices $a,b$ of $G$ so that there are exactly six edges of
$G$ between ${a,b}$ and $V(P)$, namely $ap_1,ap_2,ap_n,bp_1,bp_{n-1},bp_n$.

The Roussel-Rubio lemma (slightly reformulated for convenience) is the
following.
\begin{thm}
\label{RR}
Let $G$ be Berge, let $X$ be an anticonnected subset of $V(G)$, and
$P$ be a path in $G \setminus X$ with odd length, such that both ends of $P$ are
$X$-complete. Then either:
\begin{enumerate}
\item some edge of $P$ is $X$-complete, or
\item $P$ has length $\ge 5$ and $X$ contains a leap for $P$, or
\item $P$ has length $3$ and there is an odd antipath joining the internal vertices of $P$
with interior in $X$.
\end{enumerate}
\end{thm}

This has a number of corollaries that again we shall need throughout the paper,
and in this section we prove some of them.

\begin{thm} \label{greentouch}
Let $G$ be Berge, let $X$ be an anticonnected subset of $V(G)$, and
$P$ be a path in $G \setminus X$ with odd length, such that both ends of $P$ are
$X$-complete, and no edge of $P$ is $X$-complete.
Then every $X$-complete vertex has a neighbour in $P^*$.
\end{thm}
\Proof
Let $v$ be $X$-complete.
Certainly $P$ has length $>1$, since its ends are $X$-complete and therefore nonadjacent.
Suppose first it has length $>3$. Then by \ref{RR}, $X$ contains a leap, and so there is
a path $Q$ with ends in $X$ and with $Q^* = P^*$. Then $v$ is adjacent
to both ends of $Q$, and since $G|(V(Q)\cup\{v\})$ is not an odd hole, it follows that
$v$ has a neighbour in $Q^* = P^*$, as required. Now suppose $P$ has length $3$, and let
its vertices be $p_1\c p_4$ in order. By \ref{RR}, there is an odd antipath $Q$
between $p_2$ and $p_3$ with interior in $X$.  Since $Q$ cannot be completed
to an odd antihole via $p_3\d v\d p_2$, it follows that $v$ is adjacent to one of $p_2,p_3$, as
required. \bbox

Here is another easy lemma that gets used enough that it is worth stating separately.

\begin{thm}\label{evengap}
Let $G$ be Berge, let $X \subseteq V(G)$ be anticonnected, and let
$P$ be a path or hole in $G\setminus X$,
so that at least three vertices of $P$ are $X$-complete. Let $Q$ be a subpath of $P$ (and hence
of $G$) with both ends $X$-complete. Then the number of $X$-complete edges in $Q$
has the same parity as the length of $Q$. In particular, if $P$ is a hole, then
either there are an even number of $X$-complete edges in $P$, or there are exactly two
$X$-complete vertices and they are adjacent.
\end{thm}
\Proof The second assertion follows from the first. For the first, we use induction on the
length of $Q$. If some internal vertex of $Q$ is $X$-complete
then the result follows from the inductive hypothesis, so we may assume not. If $Q$ has
length  $1$ or even then the theorem holds, so we may assume its length is $\ge 3$ and odd.
By hypothesis there is an $X$-complete vertex $v$ say of $P$ that is not an end of $Q$,
and therefore does not belong to $Q$; and since $P$ is a path or hole, it follows that
$v$ has no neighbour in $Q^*$, contrary to \ref{greentouch}. This proves \ref{evengap}.\bbox

A {\em triangle} in $G$ is a set of three vertices, mutually adjacent. We say
a vertex $v$ can be {\em linked} onto a triangle $\{a_1,a_2,a_3\}$ (via paths $P_1,P_2,P_3$)
if:
\begin{itemize}
\item the three paths $P_1,P_2,P_3$ are mutually vertex-disjoint
\item for $i = 1,2,3$ $a_i$ is an end of $P_i$
\item for $1 \le i < j \le 3$, $a_ia_j$ is the unique edge of $G$ between $V(P_i)$ and $V(P_j)$
\item $v$ has a neighbour in each of $P_1,P_2$ and $P_3$.
\end{itemize}

The following is well-known and quite useful:
\begin{thm} \label{trianglev}
Let $G$ be Berge, and suppose $v$ can be linked onto a triangle $\{a_1,a_2,a_3\}$.
Then {v} is adjacent to at least two of $a_1,a_2,a_3$.
\end{thm}
\Proof Let $v$ be linked via paths $P_1,P_2,P_3$.
For $1 \le i \le 3$, $v$ has a neighbour in $P_i$;  let $P_i$ be the
path from $v$ to $a_i$ with interior in $V(Q_i)$.
At least two of $Q_1,Q_2,Q_3$ have lengths of the same parity,
say $Q_1,Q_2$; and since $G|(V(Q_1)\cup V(Q_2))$ is not an odd hole, it is a cycle of length $3$,
and the claim follows. \bbox

A variant of \ref{greentouch} is sometimes useful, the following:

\begin{thm}\label{greentouch2}
Let $G$ be Berge, let $X \subseteq V(G)$, and let $P$ be a path in $G\setminus X$ of odd length,
with vertices be $p_1\c p_n$, so that $p_1,p_n$ are $X$-complete, and no edge of $P$ is $X$-complete.
Let $v \in V(G)$ be $X$-complete. Then either $v$ is adjacent to one of $p_1$,$p_2$, or
the only neighbour of $v$ in $P^*$ is $p_{n-1}$.
\end{thm}
\Proof By \ref{greentouch}, $v$ has a neighbour in $P^*$, and we may assume that
$p_{n-1}$ is not its only such neighbour, so $v$ has a neighbour in
$\{p_2,\ldots,p_{n-2}\}$. If $P$ has length $\le 3$ then the result follows, so we may assume
its length is at least $5$. By \ref{RR}, there is a leap $a,b$ for $P$ in $X$; so there is a
path $a\d p_2\c p_{n-1}\d b$. Now $\{p_1,p_2,a\}$ is a triangle,
and $v$ can be linked onto it via the three paths $b\d p_1$, $P\setminus \{p_1,p_{n-1},p_n\}$, $a$;
and so $v$ has two neighbours in the triangle, by \ref{trianglev}, and the claim follows.\bbox

\begin{thm}\label{balancev}
If $G$ is Berge and $A,B \subseteq V(G)$ are disjoint, and $v \in V(G) \setminus (A \cup B)$,
and $v$ is complete to $B$ and anticomplete to $A$, then $(A,B)$ is balanced.
\end{thm}
The proof is clear.

\begin{thm}\label{shiftbalance}
Let $(A,B)$ be balanced in a Berge graph $G$. Let $C \subseteq V(G) \setminus (A \cup B)$. Then :
\begin{enumerate}
\item if $A$ is connected and every vertex in $B$ has a neighbour in $A$, and $A$ is anticomplete to
$C$, then $(C,B)$ is balanced
\item if $B$ is anticonnected and no vertex in $A$ is $B$-complete, and $B$ is complete to
$C$, then $(A,C)$ is balanced.
\end{enumerate}
\end{thm}
\Proof The first statement follows from the second by taking complements, so it suffices
to prove the second. Suppose $u,v \in A$ are adjacent and joined by an odd antipath $P$ with interior
in $C$. Since $B$ is anticonnected and $u,v$ both have non-neighbours in $B$, they are also joined
by an antipath $Q$ with interior in $B$, which is even since $(A,B)$ is balanced.
But then $u\d P\d v\d Q\d u$ is an odd antihole, a contradiction. Now suppose there are nonadjacent
$u,v \in C$, joined by an odd path $P$ with interior in $A$. Then $P$ has length $\ge 5$, since
otherwise its vertices could be reordered to be an odd antipath of the kind we already handled.
The ends of $P$ are $B$-complete, and no internal vertex is $B$-complete, and so $B$ contains
a leap for $P$ by \ref{RR}; and hence there is an odd path with ends in $B$ and interior in $A$, which
is impossible since $(A,B)$ is balanced. This proves \ref{shiftbalance}.\bbox

We already said what we mean by linking a vertex onto a triangle, but now we do the same
for an anticonnected set. We say
an anticonnected set $X$ can be {\em linked} onto a triangle $\{a_1,a_2,a_3\}$ (via paths $P_1,P_2,P_3$)
if:
\begin{itemize}
\item the three paths $P_1,P_2,P_3$ are mutually vertex-disjoint
\item for $i = 1,2,3$ $a_i$ is an end of $P_i$
\item for $1 \le i < j \le 3$, $a_ia_j$ is the unique edge of $G$ between $V(P_i)$ and $V(P_j)$
\item each of $P_1,P_2$ and $P_3$ contains an $X$-complete vertex.
\end{itemize}

There is a corresponding extension of \ref{trianglev}, the following:
\begin{thm} \label{triangleX}
Let $G$ be Berge, let $X$ be an anticonnected set, and suppose $X$ can be linked
onto a triangle $\{a_1,a_2,a_3\}$  via $P_1,P_2,P_3$. For $i = 1,2,3$ let $P_i$ have ends
$a_i$ and $b_i$, and let $b_i$ be the unique vertex of $P_i$ that is $X$-complete. Then
either at least two of $P_1,P_2,P_3$ have length $0$ (and hence two of $a_1,a_2,a_3$ are $X$-complete)
or one of $P_1,P_2,P_3$ has length $0$ and the other two have length $1$ (say $P_3$ has length $0$);
and in this case, every $X$-complete vertex in $G$ is adjacent to one of $a_1$,$a_2$.
\end{thm}
\Proof Some two of $P_1,P_2,P_3$ have lengths of the same parity, say $P_1$ and $P_2$. Hence the
path $Q = b_1\d P_1\d a_1\d a_2\d P_2\d b_2$ (with the obvious meaning - we shall feel free to specify paths by
whatever notation is most convenient) is odd, and its ends are $X$-complete, and none of its
internal vertices are $X$-complete. If $Q$ has length $1$ then the theorem holds, so we assume it has
length $\ge 3$. By \ref{greentouch}, every $X$-complete vertex has a neighbour in $Q^*$, and
since $b_3$ is $X$-complete, it follows that $b_3 = a_3$. Hence we may assume both $P_1$ and
$P_2$ have length $\ge 1$ for otherwise the claim holds. Suppose that $Q$ has length $3$. Then
$P_1$ and $P_2$ have length $1$, and the claim holds again. So we may assume (for a contradiction)
that $Q$ has length
$\ge 5$, and from the symmetry we may assume $P_1$ has length $\ge 2$.  Since $b_3$ is not adjacent
to the end $b_1$ of $Q$ or to its neighbour in $Q$, and yet it has at least two neighbours in $Q^*$
(namely $a_1$ and $a_2$), this contradicts \ref{greentouch2}. This proves \ref{triangleX}.\bbox

As we said earlier, the main use of \ref{RR} is to show that the common neighbours of an
anticonnected set behave in some respects like the neighbours of a single vertex. From this
point of view, \ref{RR} itself tells us something about when there can be an odd ``pseudohole'',
in which one ``vertex'' is actually an anticonnected set. We also need a version of this when there
are two such vertices, the following.

\begin{thm}\label{doubleRR}
Let $G$ be Berge, and let $X, Y$ be disjoint nonempty anticonnected subsets of $V(G)$,
complete to each other.  Let $P$ be a path in $G\setminus (X\cup Y)$ with even length $>0$, with vertices
$p_1,\ldots,p_n$ in order, so that $p_1$ is the unique
$X$-complete vertex of $P$ and $p_n$ is the unique $Y$-complete vertex of $P$.
Then either:
\begin{enumerate}
\item $P$ has length $\ge 4$ and there are nonadjacent $x_1,x_2 \in X$ so that $x_1\d p_2\d \cdots\d p_n\d x_2$ is a path, or
\item $P$ has length $\ge 4$ and there are nonadjacent $y_1,y_2 \in Y$ so that $y_1\d p_1\d \cdots\d p_{n-1}\d y_2$ is a path,  or
\item $P$ has length 2 and there is an antipath $Q$ between $p_2$ and $p_3$ with interior in $X$,
and an antipath $R$ between $p_1$ and $p_2$ with interior in $Y$, and exactly one of $Q$,$R$ has
odd length.
\end{enumerate}
In each case, either $(V(P\setminus p_1), X)$ or $(V(P\setminus p_n),Y)$ is not balanced.
\end{thm}
\Proof It follows from the hypotheses that $X$,$Y$ and $V(P)$ are mutually disjoint.
If $P$ has length $2$, choose an antipath $Q$ between $p_2$ and $p_3$ with interior in $X$,
and an antipath $R$ between $p_1$ and $p_2$ with interior in $Y$. Then $p_2\d Q\d p_3\d p_1\d R\d p_2$ is
an antihole, and so exactly one of $Q$,$R$ has odd length and the theorem holds. So we may assume
$P$ has length $\ge4$.
We may assume that $V(G) = V(P) \cup X \cup Y$, by deleting any other vertices. Let $G'$ be
obtained from $G\setminus Y$ by adding a new vertex {y} with neighbour set $X \cup \{p_n\}$.
Let $P'$ be the path $p_1\d \cdots\d p_n\d y$ of $G'$. Then $P'$ has odd length $\ge 5$. If $G'$ is
Berge then by \ref{RR} there is a leap for $P'$ in $X$, and the result follows. So we may
assume $G'$ is not Berge.

Assume first that there is an odd hole $C$ of length $\ge 7$ in $G'$. It
necessarily uses $y$, and the neighbours of $y$ in $C$ are $Y$-complete, and no other
vertices of $C\setminus y$ are $Y$-complete. Hence there is an odd path $Q$ in $G\setminus Y$
of length
$\ge 5$, with both ends $Y$-complete and no internal vertices $Y$-complete. So the ends
of $Q$ belong to $X\cup \{p_n\}$ and its interior to $V(P)\setminus p_n$. By \ref{RR}
$Y$ contains a leap for $Q$; so there is an odd path $R$ of length $\ge 5$
with ends ($y_1,y_2$ say) in $Y$ and with interior in $V(P)\setminus p_n$. Since $R$ cannot
be completed to a hole via $y_2\d p_n\d y_1$ it follows that $p_n$ has a neighbour in $R^*$, and
so $p_{n-1}$ belongs to $R$. If also $p_1$ belongs to $R$ then the theorem holds, so we may
assume it does not. Since $R$ is odd and $P$ is even it follows that $p_2$ also does not
belong to $R$, and so $p_1$ has no neighbour in $R^*$; yet the ends of $R$ are $X$-complete
and its internal vertices are not, contrary to \ref{greentouch}. This completes the case
when there is an odd hole in $G'$ of length $\ge 7$.

Since an odd hole of length $5$ is
also an odd antihole, we may assume that there is an odd antihole in $G'$, say $D$. Again
$D$ must use $y$, and uses exactly two nonneighbours of $y$; so in $G$ there is an odd
antipath $Q$ between adjacent vertices of $P \setminus p_n$ (say $u$ and $v$),
and with interior in $X\cup \{p_n\}$. Since $u$ and $v$ are not $Y$-complete, they are also
joined by an antipath $R$ with interior in $Y$, and $R$ must also be odd since its union
with $Q$ is an antihole. Since $R$ cannot be completed to an antihole via $v\d p_n\d u$ it
follows that $p_n$ is adjacent to one of $u$,$v$, and hence we may assume that
$u = p_{n-2}$ and $v = p_{n-1}$. Since $P$ has length $\ge 4$ it follows that $u$,$v$ are
also joined by an antipath with interior in $X$, say $S$, and again $S$ is odd since its
union with $R$ is an antihole. But $S$ can be completed to an antihole via $v\d p_1\d u$,
a contradiction. This proves \ref{doubleRR}. \bbox

Next we need a version of \ref{RR} for holes.
Let $C$ be a hole in $G$, and let $e = uv$ be an edge of it.  A {\em leap} for $C$
({\em in} $G$, {\em at} $uv$) is a leap for the path $C \setminus e$ in $G \setminus e$. A {\em hat}
for $C$ ({\em in} $G$, {\em at} $uv$) is a vertex of $G$ adjacent to $u$ and $v$
and to no other vertex of $C$.

\begin{thm} \label{RRC}

Let $G$ be Berge,  let $X \subseteq V(G)$ be anticonnected, let $C$ be a hole in $G\setminus X$
with length $>4$, and let $e = uv$ be an edge of $C$.
Assume that $u,v$ are $X$-complete and no other vertex of $C$ is $X$-complete. Then either
$X$ contains a hat for $C$ at $uv$, or $X$ contains a leap for $C$ at $uv$.
\end{thm}
\Proof
Let the vertices of $C$ be $p_1,\ldots,p_n$ in order, where $u = p_1$ and $v = p_n$.
Let $G_1 = G|(V(C) \cup X)$, and let $G_2 = G_1 \setminus e$. If $G_2$ is Berge, then from
\ref{RR} applied to the path $C \setminus e$ in $G_2$ it follows that $X$ contains a leap for $C$ at $uv$.
So we may assume that $G_2$ is not Berge. Consequently it has an odd hole or antihole $D$ say,
and since $D$ is not an odd hole or antihole in $G_1$ it must use both $p_1$ and
$p_n$. Suppose first that $D$ is an odd hole. Since every vertex in $X$ is adjacent to both
$p_1$ and $p_n$ it follows that at most one vertex of $X$ is in $D$; and since $G_2 \setminus X$
has no cycles, there is exactly one vertex of $X$ in $D$, say $x$. Hence $D \setminus x$ is a path
of $G_2 \setminus X$ between $p_1$ and $p_n$, and so $D \setminus x = C \setminus e$; and
since $D$ is a hole of $G_2$
it follows that $x$ has no neighbours in $\{p_2,\ldots,p_{n-1}\}$, and therefore is a hat
as required. Next assume that $D$ is an antihole. Since it uses both $p_1$ and $p_n$,
and they are nonadjacent in $G_2$, it follows that they are consecutive in $D$, so the
vertices of $D$ can be numbered $d_1,\ldots,d_m$ in order, where $d_1 = p_1$ and $d_m = p_n$,
and therefore $m \ge 5$. Consequently, both $d_2$ and $d_{m-1}$ are not in $X$, since
they are not complete to $\{p_1,p_n\}$, and therefore $d_1,d_2,d_{m-1},d_m$ are
vertices of $C$. Yet $d_1d_{m-1}, d_{m-1}d_2,d_2d_m$ are edges of $G_1$, which
is impossible since $n\ge 6$. This proves \ref{RRC}. \bbox

There is an analogous version of \ref{doubleRR}, as follows.

\begin{thm} \label{doubleRRC}
Let $G$ be Berge, and let $X, Y$ be disjoint nonempty anticonnected subsets of $V(G)$,
complete to each other.  Let $P$ be a path in $G\setminus (X\cup Y)$
with even length $\ge 4$, with vertices $p_1,\ldots,p_n$ in order,
so that $p_1$ is the unique
$X$-complete vertex of $P$, and $p_1,p_n$ are the only $Y$-complete vertices of $P$.
Then either:
\begin{enumerate}
\item there exists $x \in X$ non-adjacent to all of $p_2,\ldots,p_n$, or
\item there are nonadjacent $x_1,x_2 \in X$ so that $x_1\d p_2\d \cdots\d p_n\d x_2$ is a path.
\end{enumerate}
\end{thm}
\Proof The proof is similar to that of \ref{doubleRR}. We may assume $V(G) = V(P)\cup X\cup Y$.
Let $G'$ be
obtained from $G\setminus Y$ by adding a new vertex {y} with neighbour set $X \cup \{p_1,p_n\}$.
If $G'$ is Berge then the result follows from \ref{RRC}, so we may assume $G'$ is not Berge.
Assume first that there is an odd hole $C$ of length $\ge 7$ in $G'$.
Hence there is an odd path $Q$ in $G\setminus Y$ of length
$\ge 5$, with both ends $Y$-complete and no internal vertices $Y$-complete. So the ends
of $Q$ belong to $X\cup \{p_1,p_n\}$ and its interior to $V(P^*)$. By \ref{RR}
$Y$ contains a leap for $Q$; so there is an odd path $R$ of length $\ge 5$
with ends ($y_1,y_2$ say) in $Y$ and with interior in $V(P^*)$. Since $R$ is odd and $R^*$
is a subpath of the even path $P^*$, it follows that not both $p_2$ and $p_{n-1}$ belong
to $R$; but then $R$ can be completed to an odd hole via one of $y_2\d p_n\d y_1$ , $y_2\d p_1\d y_1$,
a contradiction.  This completes the case when there is an odd hole in $G'$ of length $\ge 7$,
so now we may assume that there is an odd antihole in $G'$, say $D$. Again
$D$ must use $y$, and uses exactly two nonneighbours of $y$; so in $G$ there is an odd
antipath $Q$ between adjacent vertices of $P^*$ (say $u$ and $v$),
and with interior in $X\cup \{p_n\}$. Since $u$ and $v$ are not $Y$-complete, they are also
joined by an antipath $R$ with interior in $Y$, and $R$ must also be odd since its union
with $Q$ is an antihole. Since one of $p_1$,$p_n$ is nonadjacent to both of $u,v$, we may
complete $R$ to an odd antihole via one of $u\d p_1\d v$,$u\d p_n\d v$, a contradiction.
This proves \ref{doubleRRC}. \bbox

\section{Paths and antipaths meeting}
Another class of applications of \ref{RR} is to the situation when a big path or hole
meets a big antipath or antihole. In this section we prove a collection of useful lemmas of this type.
First, a neat application of \ref{RR} (we include this only because it is
striking --- in fact we do not use it at all).

\begin{thm}\label{bighole&antihole}
Let $G$ be Berge, let $C$ be a hole in $G$, and $D$ an antihole in $G$, both of length
$\ge 8$. Then $|V(C) \cap V(D)| \le 3$.
\end{thm}
\Proof
It is easy to see that $|V(C) \cap V(D)| \le 4$, without using that $G$ is Berge. Suppose that
$|V(C) \cap V(D)| = 4$; then $V(C) \cap V(D)$ is the vertex set of a $3$-edge path. Let $C$ have vertices
$p_1,\ldots,p_m$ in order, and $D$ have vertices $q_1,\ldots,q_n$ in order, where $m,n \ge 8$
and $p_1 = q_2, p_2 = q_4, p_3 = q_1, p_4 = q_3$. Let $P$ be the path $p_4\d p_5\c p_m\d p_1$,
and $Q$ the antipath $q_4\d q_5 \c q_n\d q_1$. Let $X$ be the interior of $Q$.
Then $p_1$ and $p_4$ are $X$-complete (since $D$ is an antihole), and $P$
is a path with length odd and $\ge 5$ between these two vertices. If some vertex $p_i$ say in the
interior of $P$ is $X$-complete, then since $p_i$ is nonadjacent to both $p_2$ and $p_3$
we can complete $Q$ to an odd antihole via $q_1\d p_i\d q_4$, a contradiction. So by \ref{RR}
$X$ contains a leap for $P$; so there exists $i$ with $5 \le i < n$ and a path $P'$ joining
$q_i$ and $q_{i+1}$ with the same interior as $P$. Since $n \ge 8$, either $i >5$ or $i+1 <n$
and from the symmetry we may assume the first. But then $P'$ can be completed to an odd hole via
$q_{i+1}\d p_2\d q_i$, a contradiction. This proves \ref{bighole&antihole}.\bbox

The next two lemmas are results of the same kind:

\begin{thm}\label{evenantipath2}
Let $p_1\d \cdots\d p_m$ be a path in a Berge graph $G$.
Let $2 \le s\le m-2$, and let $p_s\d q_1\d \cdots\d q_n\d p_{s+1}$ be an antipath, where
$n \ge 2$ is odd. Assume that each of $q_1\l q_n$ has a neighbour in $\{p_1\l p_{s-1}\}$
and a neighbour in $\{p_{s+2}\l p_m\}$.  Then either:
\begin{itemize}
\item $s \ge 3$ and the only nonedges between $\{p_{s-2},p_{s-1},p_s,p_{s+1},p_{s+2}\}$
and $\{q_1,\ldots,q_n\}$ are $p_{s-1}q_n,p_sq_1,p_{s+1}q_n$, or
\item $s \le m-3$ and the only nonedges between $\{p_{s-1},p_s,p_{s+1},p_{s+2},p_{s+3}\}$
and $\{q_1,\ldots,q_n\}$ are $p_sq_1,p_{s+1}q_n,p_{s+2}q_1$.
\end{itemize}
\end{thm}
\Proof
The antipath $p_s\d q_1\d \cdots\d q_n\d p_{s+1}$ is even, of length $\ge 4$; all its vertices have
neighbours in $\{p_1\l p_{s-1}\}$ except $p_{s+1}$, and they all have neighbours in $\{p_{s+2}\l p_m\}$
except $p_s$. Since the sets $\{p_1\l p_{s-1}\}, \{p_{s+2}\l p_m\}$ are both connected and are
anticomplete to each other, it follows from \ref{doubleRR} applied in $\overline{G}$ and the symmetry
that we may assume that there are adjacent vertices $u,v \in \{p_1\l p_{s-1}\}$ so that
$u\d p_s \d q_1 \c q_n\d v$ is an antipath. Since $v$ is adjacent to $p_s$ and to $u$ it follows
that $s \ge 3$, $v = p_{s-1}$ and $u = p_{s-2}$. Since $p_{s-2}\d p_s \d q_1 \c q_n\d p_{s-1}$ is an odd
antipath of length $\ge 5$, and its ends are anticomplete to $\{p_{s+1}\l p_m\}$ and its internal vertices are
not, it follows from \ref{RR} applied in $\overline{G}$ that there are adjacent $w,x \in \{p_{s+1}\l p_m\}$
such that $w \d p_s \d q_1 \c q_n\d x$ is an antipath. Since $x$ is adjacent to $p_s$ and to $w$ it
follows that $x = p_{s+1}$ and $w = p_{s+2}$. But then the theorem holds. This proves \ref{evenantipath2}.\bbox

\begin{thm}\label{hole&antipath}
Let $G$ be Berge, let $C$ be a hole in $G$ of length $\ge 6$, with vertices $p_1,\ldots,p_m$ in
order, and let $Q$ be an antipath with vertices $p_1,q_1,\ldots,q_n,p_2$, with length $\ge 4$
and even. Let $z \in V(G)$, complete to $V(Q)$ and with no neighbours among $p_3,\ldots,p_m$.
There is at most one vertex in $\{p_3,\ldots,p_m\}$ complete to either $\{q_1,\ldots,q_{n-1}\}$
or $\{q_2,\ldots,q_n\}$, and any such vertex is one of $p_3,p_m$.
\end{thm}
\Proof
It follows that none of $q_1,\ldots,q_n$ belong to $C$, since they are all adjacent to $z$.
Let $X = \{q_1,\ldots,q_n\}$, and let $Y_1,Y_2$ be
the sets of vertices in $\{p_3,\ldots,p_m\}$ complete to $X\setminus q_n$,$X\setminus q_1$
respectively.
\\
\\
(1) $Y_1 \subseteq Y_2 \cup \{p_m\}$, {\em and} $Y_2 \subseteq Y_1 \cup \{p_3\}$.
\\
\\
For suppose some $p_i \in Y_1$, and is not in $Y_2$; then since the odd antipath
$Q \setminus p_2$ cannot be completed to an odd antihole via $q_n\d p_i\d p_1$, it follows that
$i = m$. This proves (1).
\\
\\
(2) {\em If $Y_1 \not \subseteq \{p_m\}$ then $p_3 \in Y_1 \cap Y_2$, and if
$Y_2 \not \subseteq \{p_3\}$ then $p_m \in Y_1 \cap Y_2$.}
\\
\\
For assume $Y_1 \not \subseteq \{p_m\}$, and choose $i$ with $3 \le i \le m-1$ minimum so
that $p_i \in Y_1$. By (1), $p_i \in Y_2$, so we may assume $i>3$, for otherwise the claim
holds. If $i$ is odd, then the path $p_2\d p_3\d \cdots\d p_i$ is odd and between
$X\setminus q_n$-complete vertices, and no internal vertex is $X\setminus q_n$-complete, and yet
the $X\setminus q_n$-complete vertex $z$ does not have a neighbour in its interior, contrary to
\ref{greentouch}. So $i$ is even. The path $p_i\d \cdots\d p_m\d p_1$
is therefore odd, and has length $\ge 3$, and its ends are $X\setminus q_1$-complete,
and the $X\setminus q_1$-complete vertex $z$ does not have a
neighbour in its interior; so by \ref{greentouch} some vertex $v$ of its interior is in $Y_2$,
and therefore in $Y_1 \cap Y_2$ by (1).  But the path $z\d p_2 \cdots\d p_i$ is odd, and between
$X$-complete
vertices, and has no more such vertices in its interior, and $v$ has no neighbour in its interior,
contrary to \ref{greentouch}. This proves (2).

\bigskip
Now not both $p_3,p_m$ are in $Y_1 \cap Y_2$, for otherwise $Q$ could be completed to an
odd antihole via $p_2\d p_m\d p_3\d p_1$. Hence we may assume $p_3  \notin Y_1 \cap Y_2$, and
so from (2), $Y_1 \subseteq \{p_m\}$. By (1), $Y_2 \subseteq \{p_3\} \cup Y_1$, and so
$Y_1 \cup Y_2 \subseteq \{p_3,p_m\}$. We may therefore assume that $Y_1 \cup Y_2 = \{p_3,p_m\}$,
for otherwise the theorem holds. In particular, $p_3 \in Y_2$. If also $p_m \in Y_2$, then
$p_3\d p_4\d \cdots\d p_m$ is an odd path between $X\setminus q_1$-complete vertices, and none of
its internal vertices are $X\setminus q_1$-complete, and yet the $X\setminus q_1$-complete
vertex $z$ does not have a neighbour in its interior, contrary to \ref{greentouch}. So
$p_m \notin Y_2$, and so $p_m \in Y_1$; but then $p_3\d q_1\d q_2\d \cdots\d q_n\d p_m\d p_3$ is an
odd antihole, a contradiction. This proves \ref{hole&antipath}. \bbox

\section{Skew partitions}

A maximal connected subset of a
nonempty set $A \subseteq V(G)$ is called a {\em component} of $A$, and a maximal
anticonnected subset is called an {\em anticomponent} of $A$. Let us say a skew partition $A,B$
of $G$ is {\em loose} if either some vertex in $B$ has no neighbour in some component of $A$,
or some vertex in $A$ is complete to some anticomponent of $B$.

In this section we investigate what skew partitions look like
in a ``minimum imperfect graph'', a Berge graph which is the smallest counterexample to
\ref{SPGC}; and in particular, we shall show that no skew partition in such a graph can be
either balanced (defined in the first section) or loose. So, in order to prove \ref{SPGC}, it
would be enough to prove a weaker form of our main theorem \ref{decomp}, in which
``admits a balanced skew
partition'' is repaced by ``admits a balanced or loose skew partition''. In a way this would
be easier, because many of the skew partitions we find later in the paper are loose. But
we might as well prove the stronger form, because as we shall show below, any Berge graph
admitting a loose skew partition also admits a balanced one.

\begin{thm}\label{singleton}
Let $G$ be Berge, and suppose that $G$ admits a skew partition $(A,B)$ so that either
some component of $A$ or some anticomponent of $B$ has only one vertex.
Then $G$ admits a balanced skew partition.
\end{thm}
\Proof Let $A_1,\ldots,A_m$ be the components of $A$, and $B_1,\ldots,B_n$ the anticomponents of $B$.
So the sets $A_1,\ldots,A_m$ are pairwise disjoint, all non-empty, and partition $A$, and
$m \ge 2$; and similarly for $B$.
By taking complements if necessary we may assume that $|A_1| = 1$, $A_1 = \{a_1\}$ say.
Let $N$ be the set of vertices of $G$ adjacent to $a_1$; so $N \subseteq B$. Assume first that
$N$ is not anticonnected. Then $(V(G)\setminus N,N)$ is a skew partition of $G$, and it is
easy to check that it is balanced, as required. So we may assume that $N$ is anticonnected.
Consequently $N$ is a subset of one of $B_1,\ldots,B_n$, say $B_1$. Choose $b_2 \in B_2$.
Then $N' = N \cup \{ b_2 \}$ is not anticonnected, and so $(V(G)\setminus N',N')$ is a skew
partition of $G$, and once again it is easily checked to be balanced. This proves \ref{singleton}.
\bbox

\begin{thm}\label{geteven}
If $G$ is Berge, and admits a loose skew partition, then it admits a balanced skew partition.
\end{thm}
\Proof
Let $(A,B)$ be a loose skew partition of $G$.
By taking complements if necessary, we may assume that some vertex in $B$ has no neighbour
in some component of $A$. With $G$ fixed, let us choose the skew partition $(A,B)$ and
a component $A_1$ of $A$ and an anticomponent $B_1$ of $B$ with $|B| - 2|B_1|$  minimum, such
that some vertex in $B_1$ (say $b_1$) has no neighbour in $A_1$. (We call this property the
``optimality'' of $(A,B)$.) Let the other components of
$A$ be $A_2,\ldots,A_m$, and the other anticomponents of $B$ be $B_2,\ldots,B_n$.
By \ref{singleton} we may assume that
no $|A_i|$ or $|B_j| = 1$, and in this case we shall show that the skew partition
$(A,B)$ is balanced.
\\
\\
(1) {\em For $2 \le j \le n$, no vertex in $A$ is $B_j$-complete and not $B_1$-complete, and
every vertex in $B\setminus B_1$ has a neighbour in $A_1$.}
\\
\\
For the first claim, assume some vertex $v \in A$ is $B_2$-complete and not $B_1$-complete, say.
Let $A_1' = A_1$ if $v \not \in A_1$, and let $A_1'$ be a maximal connected subset of
$A_1 \setminus v$ otherwise. (So $A_1'$ is nonempty since we assumed $|A_1| \ge 2$.) Let
$A' = A \setminus v$ and $B' = B \cup \{v\}$; then $B_2$ is still an anticomponent of $B'$, so
$(A',B')$ is a skew partition, violating the optimality of $(A,B)$ (for since $v$ is not
$B_1$-complete, there is an anticomponent of $B'$ including $B \cup \{v\}$). For the second claim,
assume that some vertex $v \in B_2$ say has no neighbour in $A_1$. Then since $B_2 \ge 2$,
it follows that $(A \cup \{v\},B \setminus v)$ is a skew partition of $G$, again violating
the optimality of $(A,B)$. This proves (1).

\bigskip
By \ref{balancev}, the pair $(A_1,B_j)$ is balanced, for $2 \le j \le n$, since $b_1$ is
complete to $B_j$ and has no neighbours in $A_1$. By (1) and \ref{shiftbalance}.1, it
follows that $(A_i,B_j)$ is balanced for $2 \le i \le m$ and $2 \le j \le n$. It remains
to check all the pairs $(A_i,B_1)$. Let $1 \le i \le m$, and let $A_i'$ be the set
of vertices in $A_i$ that are not $B_1$-complete. By (1), no vertex in $A_i'$ is $B_2$-complete,
and $(A_i', B_2)$ is balanced, and
hence by \ref{shiftbalance}.2, so is $(A_i',B_1)$, and consequently so is $(A_i,B_1)$.
This proves that $(A,B)$ is balanced, and so completes the proof of \ref{geteven}. \bbox

\begin{thm}\label{mixedpair}
Let $(A,B)$ be a skew partition of a Berge graph $G$. If either:
\begin{itemize}
\item there exist $u,v \in B$ joined by an odd path with interior in $A$, and joined by an even
path with interior in $A$, or
\item there exist $u,v \in A$ joined by an odd antipath with interior in $B$, and joined by an
even antipath with interior in $B$,
\end{itemize}
then $(A,B)$ is loose and therefore $G$ admits a balanced skew partition.
\end{thm}
\Proof
By taking complements we may assume that
the first case of the theorem applies. Let $A_1,\ldots,A_m$ be the components of $A$, and
$B_1,\ldots,B_n$ the anticomponents of $B$. Since $u,v \in B$, and they are joined by an even path,
they are therefore nonadjacent, and so belong to the same $B_j$, say $B_1$. There is an even path
$P_1$ and an odd path $P_2$ joining $u,v$, both with interior in $A$. We may assume that $P_1$
has interior in $A_1$. Since the union of $P_1$ and $P_2$ is not a hole, it follows that
$P_2$ also has interior in $A_1$. If $u,v$ are joined by a path with interior in $A_2$, then
its union with one of $P_1$,$P_2$ would be an odd hole, a contradiction; so there is no
such path. Hence one of $u,v$ has no neighbours in $A_2$, and hence $(A,B)$ is loose, and
the theorem follows from \ref{geteven}. This proves \ref{mixedpair}. \bbox

Let $(A,B)$ be a skew partition of $G$, and
let $A_1,\ldots,A_m$ be the components of $A$, and $B_1,\ldots,B_n$ the anticomponents of $B$.
For $1 \le i \le m$ and $1 \le j \le n$, we say $(i,j)$ is a {\em path pair} if there is an
odd path in $G$ with ends nonadjacent vertices of $B_j$ and with interior in $A_i$; and
$(i,j)$ is an {\em antipath pair} if there is an
odd antipath in $G$ with ends adjacent vertices of $A_i$ and with interior in $B_j$.

\begin{thm}\label{allpairs}
Let $(A,B)$ be a skew partition of a Berge graph $G$, and let
$A_1,\ldots,A_m$ be the components of $A$, and $B_1,\ldots,B_n$ the anticomponents of $B$.
Then either:
\begin{itemize}
\item $(A,B)$ is loose or balanced, or
\item $(i,j)$ is a path pair for all $i,j$ with $1 \le i \le m$ and $1 \le j \le n$, or
\item $(i,j)$ is an antipath pair for all $i,j$ with $1 \le i \le m$ and $1 \le j \le n$.
\end{itemize}
\end{thm}
\Proof
We may assume $(A,B)$ is not loose and not balanced.
\\
\\
(1) {\em If for some $i,j$ there is an odd path of length $\ge 5$ with ends in $B_j$ and interior
in $A_i$, then the theorem holds.}
\\
\\
For assume there is such a path for $i = j = 1$ say. Let this path, $P_1$ say, have vertices
$b_1\d p_1\d p_2\d \dots\d p_n\d b_1'$, where $b_1,b_1' \in B_1$ and $p_1,\ldots,p_n \in A_1$.
Let $2 \le j \le n$. Then $P_1$ is an odd path of length $\ge 5$ between common neighbours
of $B_j$, and no internal vertex of it is $B_j$-complete since $(A,B)$ is not loose. By
\ref{RR}, $B_j$ contains a leap; so there exist nonadjacent $b_j,b_j' \in B_j$ so that
$b_j\d p_1\d p_2\d \dots\d p_n\d b_j'$ is a path. Hence $(1,j)$ is a path pair. Now let
$2 \le i \le m$ and $1 \le j \le n$. Since $(A,B)$ is not loose, $b_j$ and $b_j'$ both have
neighbours in $A_i$, and so there is a path $P_2$ say joining them with interior in $A_i$; it
is odd by \ref{mixedpair}, and so $(i,j)$ is a path pair. This proves (1).

From (1) we may assume that for all $i,j$, every odd path of length $>1$
with ends in $B_j$ and interior in $A_i$ has length 3; and similarly every odd antipath of length
$>1$ with ends in $A_i$ and interior in $B_j$ has length 3. Consequently, every path pair is also
an antipath pair (because a path of length 3 can be reordered to be an antipath of length 3).
We may assume that $(1,1)$ is a path pair, and so there exist $b_1,b_1' \in B_1$ and
$a_1,a_1' \in A_1$ so that $b_1\d a_1\d a_1'\d b_1'$ is a path $P_1$ say. Let $2 \le i \le m$.
Since $b_1$ and
$b_1'$ both have neighbours in $A_i$, they are joined by a path with interior in $A_i$, odd by
\ref{mixedpair} ; and so by (1) it has length 3. Hence there exist
$a_i,a_i' \in A_i$ so that $b_1\d a_i\d a_i'\d b_1'$ is a path. By the same argument in the complement,
it follows that for all $1 \le i \le m$ and $2 \le j \le n$, there exist $b_j,b_j' \in B_j$ so that
$b_j\d a_i\d a_i'\d b_j'$ is a path. So every pair $(i,j)$ is both a path and antipath pair. This proves
\ref{allpairs}. \bbox

We can reformulate the previous result in a form that is easier to apply, as follows.
\begin{thm}\label{onepair}
Let $G$ be Berge.
Suppose that there is a partition of $V(G)$ into four nonempty sets $X,Y,L,R$, such that there
are no edges between $L$ and $R$, and $X$ is complete to $Y$. If either:
\begin {itemize}
\item some vertex in $X\cup Y$ has no neighbours in $L$ or no neighbours in $R$, or
\item some vertex in $L \cup R$ is complete to $X$ or complete to $Y$, or
\item $(L,Y)$ is balanced
\end{itemize}
then $G$ admits a balanced skew partition.
\end{thm}
\Proof Certainly $(L \cup R, X \cup Y)$ is a skew partition, so by \ref{geteven} we may assume it
is not loose, and therefore neither of the first two alternative hypotheses holds.
So we assume the third hypothesis holds.
Let $A_1,\ldots,A_m$ be the components of $L \cup R$, and let $B_1,\ldots,B_n$ be the
anticomponents of $X \cup Y$. Since $X,Y,L,R$ are all nonempty we may assume that
$A_1 \subseteq L$, and $B_1 \subseteq X$. By hypothesis, $(1,1)$ is not a path or antipath
pair, and so by \ref{allpairs} the skew partition is balanced. This proves \ref{onepair}. \bbox

In the main proof there will be several occasions when we need to show that a given skew partition
is either loose or balanced, and some of them can be handled by the following. Let $(A,B)$ be a
skew partition of $G$. We say that an anticonnected subset $W$ of $B$ is a {\em kernel} for the
skew partition if some component of $A$ contains no $W$-complete vertex.

\begin{thm}\label{kernel}
Let $(A,B)$ be a skew partition of a Berge graph $G$, and let $W$ be a kernel for it.
Let $A_1$ be a component of $A$, and suppose that $(A_1,W)$ is balanced.
Then $G$ admits a balanced skew partition.
\end{thm}
\Proof
By \ref{geteven} we may assume $(A,B)$ is not loose. Let the components of $A$ be $A_1,\ldots,A_m$, and the
anticomponents of $B$ be $B_1,\ldots,B_n$.
\\
\\
(1) {\em $(A_i,W)$ is balanced for $1 \le i \le m$.}
\\
\\
For this is true if $i = 1$, so assume $i >1$. From \ref{mixedpair} there is no odd path
between nonadjacent vertices of $W$ with interior in $A_i$.  Suppose there is an odd
antipath $Q$ of length $>1$, with ends in $A_i$ and interior in $W$. Then it has length $\ge 5$,
for otherwise it can be reordered to be an odd path that we have already shown impossible.
Now the ends of $Q$ have no neighbours in the connected set $A_1$, and its internal vertices all
have neighbours in $A_1$; and so by \ref{RR} in the complement, there is a leap in the complement;
that is, there is an antipath with ends in $A_1$ and with the same interior as $Q$, which is
impossible. This proves (1).

Since $W$ is anticonnected, we may assume that $W \subseteq B_1$. Since (1) restores the
symmetry between $A_1,\ldots,A_m$, we may assume that there is no $W$-complete vertex in $A_1$.
By \ref{allpairs} we may assume $(1,2)$ is a path or antipath pair. Suppose first that it is
an antipath pair. Then there is an odd antipath $Q_1$ of length $\ge 3$ with ends in $A_1$ and
interior in $B_2$. Since its ends both have nonneighbours in $W$, its ends are also joined by
an antipath $Q_2$ with interior in $W$, odd by \ref{mixedpair}, contrary to (1).
So there is no such $Q_1$. Hence there is an
odd path $P$ with ends in $B_2$ and interior in $A_1$, necessarily of length $\ge 5$ (since we
already did the antipath case). Since the interior of $P$ contains no $W$-complete
vertex, \ref{RR} implies that $W$ contains a leap; and so there is a path with ends in $W$ with
the same interior as $P$, a contradiction. This proves \ref{kernel}. \bbox

One can refine \ref{kernel} a little more, as follows.
\begin{thm}\label{kernel2}
Let $(A,B)$ be a skew partition of a Berge graph $G$, and let $W$ be a kernel for it.
Let $A_1$ be a component of $A$, and suppose that
\begin{itemize}
\item every pair of nonadjacent vertices of $W$ with neighbours in $A_1$ are joined by an even path
with interior in $A_1$
\item every pair of adjacent vertices of $A_1$ with nonneighbours in $W$ are joined by an even
antipath with interior in $W$.
\end{itemize}
Then $G$ admits a balanced skew partition.
\end{thm}
\Proof This is immediate from \ref{kernel} and \ref{mixedpair}. \bbox

By a {\em minimum imperfect graph} we mean a {\em Berge} graph $G$, not perfect, with $|V(G)|$
minimum.  Now let us investigate skew partitions in a minimum imperfect graph.

\begin{thm}\label{evenskew}
Let $(A,B)$ be a skew partition in a minimum imperfect graph $G$, and let $A_1,\ldots,A_m$ and
$B_1,\ldots,B_n$ be defined as usual.
For all $i$ with $1 \le i \le m$ there exists $j$ with $1 \le j \le n$ such that $(i,j)$ is a
path or antipath pair, and for all $j$ with $1 \le j \le n$ there exists $i$ with $1 \le i \le m$
such that $(i,j)$ is a path or antipath pair.
\end{thm}
\Proof The first statement is equivalent to the second by taking complements, since
$\overline{G}$ also
satisfies the hypotheses of the theorem and $(B,A)$ is a skew partition in it. It therefore
suffices to prove the second statement, and we may assume $j = 1$.
Let $G'$ be the graph obtained from $G$ by adding a new
vertex $z$ with neighbour set $B_1$.
\\
\\
(1){\em We may assume $G'$ is Berge.}
\\
\\
For suppose it is not. Then in $G'$ there is an odd hole or antihole using $z$. Suppose first
that there is an odd hole, $C$ say. Then the neighbours of $z$ in $C$ (say $x,y$)
belong to $B_1$, and
no other vertex of $B_1$ is in $C$. For $2 \le j\le n$ no vertex of $B_j$ is in $C$ since
it would be adjacent to $x,y$ and $C$ would have length $4$; so $C\setminus z$ is an odd path of
$G$, with ends in $B_1$ and with interior in $A$. Since the interior of this odd path is
connected, it is a subset of one of $A_1,\ldots,A_m$, say $A_i$; but then $(i,1)$ is a
path pair and the theorem holds. So we may assume there is no such $C$. Now assume there
is an odd antihole $D$ in $G'$, again using $z$. Then exactly two vertices of $D\setminus z$
are nonadjacent to $z$, so all the others belong to $B_1$. Hence in $G$ there is an odd
antipath $Q$ of length $\ge 3$, with ends $x,y \not \in B_1$ and with interior in $B_1$.
Since both $x$ and $y$ have nonneighbours in the interior of $Q$ it follows that $x,y \not \in B$;
and since they are adjacent they both belong to $A_i$ for some $i$. But then $(i,1)$ is an
antipath pair. This proves (1).

\bigskip
For a subset $X$ of $V(G)$, we denote the size of the largest clique in $X$ by $\omega(X)$.
Let $\omega(B_1) = s$, and $\omega(A \cup B) = t$. Since $G$ is minimum imperfect it cannot
be $t$-coloured.
\\
\\
(2) {\em For $1 \le i \le m$ there is a subset $C_i \subseteq A_i$ so that
$\omega(C_i \cup B_1) = s$ and \[\omega((A_i\setminus C_i) \cup (B \setminus B_1)) \le t-s.\] }
For let $H = G'|(B \cup A_i \cup \{ z\})$; then $H$ is Berge, by (1). Now by \cite{starcut},
 there are at least two vertices of $G$ not in $H$ (all the vertices in $A \setminus A_i$), and
since $H$ has only one new vertex it follows that $|V(H)| < |V(G)|$. From the minimality
of $|V(G)|$ we deduce that $H$ is perfect. Now a theorem of Lov\'{a}sz~\cite{Lovasz} shows that
replicating a vertex of a perfect graph makes another perfect graph; so if we replace $z$
by a set $Z$ of $t-s$ vertices all complete to $B_1$ and to each other, and with no other
neighbours in
$A_i \cup B$, then the graph we make is perfect. From the construction, the largest clique
in this graph has size $\le t$, and so it is $t$-colourable. Since $Z$ is a clique of size $t-s$,
we may assume that colours $1,\ldots,s$ do not occur in $Z$, and colours $s+1,\ldots,t$ do.
Since $B_1$ is complete to $Z$, colours $s+1,\ldots,t$ do not occur in $B_1$, and so only
colours $1,\ldots,s$ occur in $B_1$; and since $\omega(B_1) = s$, all these colours do occur
in $B_1$. Since $B_1$ is complete to $B \setminus B_1$, none of colours $1,\ldots,s$
occur in $B \setminus B_1$.  Let $C_i$ be the set of vertices $v \in A_i$ with colours
$1,\ldots,s$. Then $C_i \cup B_1$ has been coloured using only $s$ colours, and so
$\omega(C_i \cup B_1) = s$; and the remainder of $H \setminus z $ has been coloured using
only colours $s+1,\ldots,t$, and so
\[\omega((A_i\setminus C_i) \cup (B \setminus B_1)) \le t-s.\]
This proves (2).
\bigskip

Now let $C = B_1 \cup C_1 \cup \cdots \cup C_m$ and $D = V(G) \setminus C$. Since there are no edges
between different $A_i$'s, it follows from (2) that $\omega(C) = s$, and similarly
$\omega(D) \le t-s$. Since $|C|,|D| < |V(G)|$ it follows that $G|C, G|D$ are both perfect;
so they are $s$-colourable and $(t-s)$-colourable, respectively. But then $G$ is $t$-colourable,
a contradiction. This proves \ref{evenskew}.\bbox

\begin{thm}\label{tightskew}
Let $G$ be a minimum imperfect graph. Then
$G$ admits no balanced skew partition, and consequently no skew partition of $G$ is loose.
\end{thm}
\Proof
The first claim follows from \ref{evenskew}, and the second from \ref{geteven}. This proves
\ref{tightskew}. \bbox

At this stage we are not ready to prove that $G$ admits no skew partition at all; that will
be shown later in the paper. But let us sketch the route. We will need three special graphs:
\begin{itemize}
\item A {\em prism} means a graph consisting of two vertex-disjoint triangles $\{a_1,a_2,a_3\}$,
$\{b_1,b_2,b_3\}$, and three paths $P_1,P_2,P_3$, where each $P_i$ has ends $a_i,b_i$,
and for $1 \le i < j \le 3$ the only edges between $V(P_i)$ and $V(P_j)$ are $a_ia_j$ and
$b_ib_j$.  The prism is {\em long} if at least one of the three paths has length $>1$.
\item A {\em double diamond} means the graph with eight vertices $a_1,\ldots,a_4,b_1,\ldots,b_4$
and with the following adjacencies:
every two $a_i$'s are adjacent except $a_3a_4$, every two $b_i$'s are
adjacent except $b_3b_4$, and $a_ib_i$ is an edge for $1 \le i \le 4$.
\item The third graph is just $L(K_{3,3}\setminus e)$, the line graph of the graph obtained
from $L(K_{3,3})$ by deleting one edge.
\end{itemize}
Note that the second and third graphs in this list are isomorphic to their complements.
We shall eventually prove that every Berge graph containing
as an induced subgraph either a long prism or a double diamond or $L(K_{3,3}\setminus e)$
must satisfy the conclusion of \ref{decomp}, and consequently cannot be a minimum imperfect
graph (and therefore nor is its complement). So once that is established, it will follow
from the next theorem that no minimum imperfect graph admits a skew partition.

\begin{thm}\label{findprism}
If $G$ is Berge, and admits a skew partition, then either $G$ admits a balanced skew partition, or
one of $G,\overline{G}$ contains as
an induced subgraph either a long prism, or a double diamond, or $L(K_{3,3}\setminus e)$.
\end{thm}
\Proof
Let $(A,B)$ be a skew partition in $G$, which is not loose by \ref{geteven}.
Let $A_1,\ldots,A_m,\linebreak[0] B_1,\ldots,B_n$ be as
before. Suppose first that for some path pair $(i,j)$ there is an odd path $P$ of length $\ge 5$
with ends in $B_j$ and with interior in $A_i$; and we may assume $i = j = 1$.
Let the vertices of $P$ be $p_1,p_2,\ldots,p_n$ in order. Now the ends of $P$ are
$B_2$-complete, and its internal vertices are not, since the skew partition is not loose;
so by \ref{RR}, $B_2$ contains a leap $x,y$, where $x$ is adjacent to $p_2$. But then the subgraph
induced on $V(P)\cup \{x,y\}$ is a long prism, as required.
So we may assume that no such path has length $\ge 5$; and similarly no odd antipath with
ends in some $A_i$ and interior in some $B_j$ has length $\ge 5$. So every path pair is
also an antipath pair and vice versa (because all the corresponding odd paths and antipaths
have length $3$ and so each is both a path and an antipath). We may therefore assume that
$(1,1)$ is a path pair, and that there exist nonadjacent $b_1,b_1' \in B_1$ and adjacent
$a_1,a_1' \in A_1$ so that $b_1\d a_1\d a_1'\d b_1'$ is a path. Since the skew partition is not loose,
$a_1,a_1'$ both have non-neighbours in $B_2$, and hence are joined by an antipath with interior in
$B_2$; this antipath is odd, since its union with ${b_1,b_1'}$ induces an antihole, and
since all such antipaths have length $3$ it follows that there exist nonadjacent
$b_2,b_2' \in B_2$ so that $b_2\d a_1\d a_1'\d b_2'$ is a path. Now $b_1,b_1'$ both have neighbours
in $A_2$, since the skew partition is not loose, and hence are joined by a path with interior in
$A_2$, and it is
odd as usual, and hence has length 3; so there exist adjacent $a_2,a_2' \in A_2$ so that
$b_1\d a_2\d a_2'\d b_1'$ is a path. Since $b_2\d b_1\d a_2\d a_2'\d b_1'\d b_2$ is not an odd hole, $b_2$
is adjacent to one of $a_2,a_2'$, and similarly so is $b_2'$. But $b_2,b_2'$ have no
common neighbour in $A_2$, for if $v \in A_2$ were adjacent to them both then
$v\d b_2\d a_1\d a_1'\d b_2'\d v$ would be an odd hole. So there are exactly two edges between
$\{a_2,a_2'\}$ and $\{b_2,b_2'\}$, forming an induced $2$-edge matching. There are two
possible pairings; in one case the subgraph induced on these eight vertices is a double
diamond, and in the other it is $L(K_{3,3}\setminus e)$. This proves \ref{findprism}.\bbox

\section{Small attachments to a line graph}

We come now to the first of the major steps of the proof. Suppose that $G$ is Berge, and
contains as an induced subgraph a substantial line graph $L(H)$. Then in general,
$G$ itself can only be basic by being a line graph, so \ref{decomp} would imply that
either $G$ is a line graph, or it has a decomposition in accordance with \ref{decomp}.
Proving a result of this kind is our first main goal, but exactly how
it goes depends on what we mean by ``substantial''. To make the theorem as powerful
as possible, we need to weaken what we mean by ``substantial'' as much as we can; but when
$L(H)$ gets very small, all sorts of bad things start to happen. One is that the theorem
is not true any more. For instance, when $H = K_{3,3}$ or $K_{3,3}\setminus e$, then
$L(H)$ is not only a line graph but also the complement of a line graph (indeed, it is
isomorphic to its own complement). So $L(H)$ can live happily inside bigger graphs that are
complements of line graphs, without inducing any kind of decomposition. The best we can
hope for, when $L(H)$ is so small, is therefore to prove that either $G$ is a line graph or
the complement of a line graph, or has a decomposition of the kind we like.
This works for $L(K_{3,3})$, but for $L(K_{3,3}\setminus e)$ the situation is even
worse, because this graph is basic in {\em three} ways - it is a line graph, the complement
of a line graph, and a bicograph. So for Berge graphs that contain $L(K_{3,3}\setminus e)$,
the best we can hope is that either $G$ is a line graph or the complement of one or a bicograph,
or it has a decomposition.  And that turns out to be true, but it also explains why the small
cases will be something of a headache, as the reader will see.

Our route through these complications is as follows. If $H$ is a subdivision of $K_4$, we say that
$L(H)$  is
``degenerate'' if there is some cycle $C$ in $K_4$ of length 4 such that the corresponding cycle of
$H$ still has length 4 (so the four edges of $C$ were not subdivided at all in producing $H$),
and {\em nondegenerate} otherwise.
First we assume that $G$ contains some $L(H)$ where $H$ is a subdivision of a 3-connected graph,
such that if $H$ is a subdivision of $K_4$ then $L(H)$ is not degenerate. Then everything works
properly --- we
can show that either $G$ is a line graph, or $G$ has a decomposition, or $H = K_{3,3}$ and
$\overline{G}$ is a line graph. Next we shall show that if $G$ does not contain any such $L(H)$,
and it does contain some degenerate $L(H)$ where $H$ is a subdivision of $K_4$, then either $G$ is
a bicograph or it has a decomposition.
The third step is, now assume that neither $G$ nor its complement contains the line graph of any
bipartite subdivision of $K_4$, and that $G$ does contain a long prism (even or odd, though the odd case
is much more difficult); then it has a decomposition (except for one graph which is basic).

The goal of the next few sections is therefore to prove the following two theorems (the bicographs and
long prisms come later).

\begin{thm}\label{linegraph}
Let $G$ be Berge, and assume some nondegenerate $L(H)$ is an induced subgraph of $G$, where $H$ is a
bipartite subdivision of $K_4$. Then either $G$ is a line graph, or $G$
admits a 2-join, or $G$ admits a balanced skew partition. In particular, \ref{summary}.1 holds.
\end{thm}

\begin{thm}\label{linegraph2}
Let $G$ be Berge, and assume it contains an induced subgraph isomorphic to $L(K_{3,3})$. Then either
one of $G, \overline{G}$ is a line graph, or one of $G,\overline{G}$ admits a 2-join, or
$G$ admits a balanced skew partition.
\end{thm}

The proof (of both --- we prove them together) is roughly as follows. We choose a $3$-connected
graph $J$, as large as
possible so that $G$ contains $L(H)$ for some bipartite subdivision $H$ of $J$ (and when
$H = K_{3,3}$, we also assume that passing to the complement will not give us
a better choice of $J$). Now
we investigate how the remainder of $G$ can attach onto $L(H)$. The edges of $J$ correspond to
edge-disjoint paths of $H$, which in turn become vertex-disjoint paths of $L(H)$, which we
call ``rungs'' (we will do the definitions properly later). One thing
we find is that the remainder of $G$ can contain alternative rungs - paths that could replace
one of the rungs in $L(H)$ to give a new $L(H')$, for some other bipartite subdivision
$H'$ of the same graph $J$. We find it advantageous to assemble all these alternative rungs
in one ``strip'', for each edge of $J$, and to maximize the union of these strips (being
careful that there are no unexpected edges of $G$ between strips). Each strip
corresponds to an edge of $J$, and runs between two sets of vertices (called ``potatoes'' for
now) that correspond to vertices of $J$. Let the union of the strips be $Z$ say. Again we
ask, how does the remainder of $G$ attach onto this ``generalized line graph'' $Z$?
This turns out to be quite pretty. There are only two kinds of vertices in the remainder of
$G$, vertices with very few neighbours in $Z$, and vertices
with a lot of neighbours. For the first kind, all their neighbours lie either in one of the
strips, or in one of the potatoes; and we can show
that for any connected set of these ``minor'' vertices, the union of their neighbours in $Z$
has the same property (they all lie in one strip or in one potato). For the second kind of
vertex, they have so many neighbours in $Z$ that all their {\em non-neighbours}
in any one potato lie inside one strip incident with the potato; and the same is true
for the union of the nonneighbours of any anticonnected set of these ``major'' vertices. In other
words, every anticonnected set of these major vertices has a great many common neighbours
in $Z$, so many that they separate all the strips from one another, and that is where we
find skew partitions. If there are no major vertices, then we find that $G$ either admits
a 2-join, or $G$ is a line graph.

First, we assume that $G$ contains $L(H)$, and we shall study how the remaining vertices of
$G$ attach to $L(H)$. Any such vertex has a set of neighbours in $V(L(H))$, that we want
to investigate; but this set is more conveniently thought of as a subset of $E(H)$, and we
begin with some lemmas about subsets of edges of a graph $H$. (Our lemmas also apply to the
common neighbours of an anticonnected set.) Our goal in this section is to examine how
individual vertices attach to $L(H)$, and how connected sets of minor vertices attach.
In the next section we think about anticonnected sets of major vertices.

In this section and the next few, we have to pay for our convention that ``path'' means
``induced path'',
because here we need paths in the conventional sense. So to compound the confusion, let us
use a different word for them. A {\em track} $P$ is a non-null connected graph in which
every vertex has degree $\le 2$ ; and its {\em length} is the number of edges in it.
(Its ends and internal vertices are defined in the natural way.)
A {\em track} in a graph $H$ means a subgraph of $H$ (not necessarily induced) which is a
track.  Note that there is a correspondence between the tracks (with at least one edge)
in a graph $H$ and the paths
in $L(H)$; the edge-set of a track becomes the vertex-set of a path, and vice versa. And two
tracks are vertex-disjoint if and only if the corresponding paths are vertex-disjoint and there
is no edge of $L(H)$ between them. However, {\em the parity changes}; a track in $H$ and the
corresponding path in $L(H)$ have lengths of opposite parity.

A {\em branch-vertex} of a graph $H$ means a vertex with degree $\ge 3$; and a
{\em branch} of $H$ means a maximal track $P$ in $H$ such that no internal vertex of
$P$ is a branch-vertex.
{\em Subdividing} an edge $uv$ means deleting the edge $uv$, adding a new vertex $w$,
and adding two new edges $uw$ and $wv$. Starting with a graph $J$, the effect of repeatedly
subdividing edges is to replace each edge of $J$ by a track joining the same pair of
vertices, where these tracks are disjoint except for their ends. We call the graph we obtain
a {\em subdivision} of $J$. Note that $V(J) \subseteq V(H)$.
Let $J$ be a $3$-connected graph. (We use the convention that a $k$-connected graph
must have $>k$ vertices.) If $H$ is a subdivision of $J$ then $V(J)$ is the set of
branch-vertices of $H$, and the branches of $H$ are in 1-1 correspondence with
the edges of $J$.  We say $H$ is {\em cyclically $3$-connected} if it is a subdivision
of some $3$-connected graph $J$. (We remind the reader that in this paper, all graphs are simple
by definition.)

We need the following lemma:

\begin{thm}\label{nondegen}
Let $H$ be bipartite and cyclically 3-connected. Then either $H= K_{3,3}$, or $H$ is a subdivision
of $K_4$, or $H$ has a subgraph $H'$ such that $H'$ is a subdivision of $K_4$,
and there is no cycle of $H'$ with vertex set the set of branch-vertices of $H'$.
\end{thm}
\Proof
There is a subgraph of $H$ which is a subdivision of $K_4$, and we may assume that it does
not satisfy the theorem. Hence there are tracks $p_1\c p_m$ ($=P$ say) and $q_1\c q_n$ ($=Q$ say)
of $H$, vertex-disjoint, so that $p_1q_1,p_1q_n,p_mq_1,p_mq_n$ are edges, and $m,n \ge 3$ are odd.
Suppose every track in $H$ between $\{p_1\l p_m\}$ and $\{q_1\l q_n\}$ uses one of the edges
$p_1q_1,p_1q_n,p_mq_1,p_mq_n$. Then there are no edges between $P$ and $Q$ except the given four,
and for every component $F$ of $H \setminus (V(P) \cup V(Q))$, the set of attachments of $F$ in
$V(P) \cup V(Q)$ is a subset of one of $V(P),V(Q)$. Since $H$ is cyclically 3-connected, it follows that
$H$ is a subdivision of $K_4$ and the theorem holds. So we may assume that
there is a track $R$ of $H$, say $r_1\c r_t$, from $V(P)$ to $V(Q)$, not using any of
$p_1q_1,p_1q_n,p_mq_1,p_mq_n$. We may assume that $r_1 \in \{p_1\l p_{m-1}\}$, $r_t \in \{q_1\l q_{n-1}\}$,
and none of $r_2\l r_{t-1}$ belong to $V(P) \cup V(Q)$. The subgraph $H'$ formed by the edges
$E(P) \cup E(Q) \cup E(R) \cup \{p_1q_n,p_mq_1,p_mq_n\}$ (and the vertices of $H$ incident with them)
is a subdivision of $K_4$, and we may assume it does not satisfy the theorem. There is therefore a cycle
of $H'$ with vertex set $\{r_1,r_t,p_m,q_n\}$. Since $H$ is bipartite and $p_mq_n$ is an edge, it follows
that $t = 2$. Hence not both $r_1 = p_1$ and $r_2 = q_1$, and so $r_1 = p_{m-1}$ and $r_2 = q_{n-1}$.
By the same argument with $p_1,p_m$ exchanged, it follows that $r_1 = p_2$, and so $m = 3$, and similarly
$n = 3$. Hence there is a subgraph $J$ of $H$ isomorphic to $K_{3,3}$.

It is helpful now to change the notation. Let $J$ have vertex set $\{a_1,a_2,a_3,b_1,b_2,b_3\}$, where
$a_1,a_2,a_3$ are adjacent to $b_1,b_2,b_3$. Suppose that there is a component $F$ of $H \setminus V(J)$.
Since $H$ is cyclically 3-connected, at least two vertices of $J$
are attachments of $F$. If say $a_1,b_1$ are attachments, choose a track $P$ between $a_1,b_1$ with
interior in $F$; then the union of $P$ and $J \setminus \{a_1b_1,a_2b_2\}$ satisfies the theorem.
If say $a_1,a_2$ are attachments of $F$, choose a track $P$ between $a_1,a_2$ with interior in $F$; then
the union of $P$ and $J \setminus\{a_1b_1,a_2b_3\}$ satisfies the theorem. So we may assume there is no such $F$.
Since $H$ is bipartite, it follows that $H = J = K_{3,3}$ , and so the theorem holds.
This proves \ref{nondegen}.\bbox

This lemma can be used to reformulate \ref{linegraph2} in the following, perhaps more informative, way:
\begin{thm}\label{linegraph2.5}
Let $G$ be Berge, and assume it contains $L(K_{3,3})$ as an induced subgraph. Then either:
\begin{itemize}
\item $G = L(K_{3,3})$, or
\item one of $G, \overline{G}$ contains some nondegenerate $L(H)$ as an induced subgraph, where
$H$ is a bipartite subdivision of $K_4$, or
\item one of $G,\overline{G}$ admits a 2-join, or $G$ admits a balanced skew partition.
\end{itemize}
In particular, \ref{summary}.2 holds.
\end{thm}
\noindent{\bf Proof of \ref{linegraph2.5}, assuming \ref{linegraph2}.}
\newline
From \ref{linegraph2} we may assume that $G$ is a line graph, $G = L(H)$ say. If $H$ is not cyclically
3-connected, then $G = L(H)$ admits a 2-join or a balanced skew partition and we are done, so we may
assume that $H$ is cyclically
3-connected. It follows easily that $H$ is bipartite, and we may assume that $L(H')$ is degenerate for
every subgraph $H'$ of $H$
isomorphic to a subdivision of $K_4$. By \ref{nondegen}, it follows that
$H = K_{3,3}$, and so $G = L(K_{3,3})$. This proves \ref{linegraph2.5}.\bbox

We observe:

\begin{thm}\label{branchcut}
Let $H$ be cyclically 3-connected, and let $C,D$ be subgraphs with $C \cup D = H$,
$|V(C \cap D)| \le 2$, and $V(C), V(D) \neq V(H)$. Then one of $C,D$ is contained
in a branch of $H$.
\end{thm}
The proof is clear.

\begin{thm}\label{twisty}
Let $c_1,c_2$ be nonadjacent vertices of a graph $H$, so that $H\setminus \{c_1,c_2\}$ is
connected. For $i = 1,2$, let the edges incident with $c_i$ be partitioned into two sets $A_i,B_i$,
where $A_1,A_2$ are both nonempty and at least one of $B_1,B_2$ is nonempty.
Assume that for every edge $uv \in A_1 \cup A_2$, $H\setminus \{u,v\}$ is connected, and that
no vertex of $V(H)$ is incident with all edges in $A_1\cup A_2$.
Then one of the following holds:
\begin{enumerate}
\item
there is a track in $H$ with first edge in $A_1$, second edge in $B_1$ (and hence second vertex
$c_1$), last vertex $c_2$ and last edge in $A_2$, or
\item
there is a track in $H$ with first edge in $A_2$, second edge in $B_2$ (and hence second vertex
$c_2$), last vertex $c_1$ and last edge in $A_1$.
\end{enumerate}
\end{thm}
\Proof
For $i = 1,2$ let $X_i$ be the set of ends (different from $c_i$) of edges in $A_i$, and define
$Y_i$ similarly for $B_i$. So by hypothesis, $X_1,X_2$ are nonempty,
$|X_1\cup X_2| \ge 2$, and we may assume $Y_1$ is
nonempty.  Choose $x_1 \in X_1$ so that $X_2 \not \subseteq \{x_1\}$ (this is possible since
$|X_1\cup X_2| \ge 2$). Both $Y_1$ and $X_2$ meet the connected graph $H\setminus \{c_1,x_1\}$, and
so there is a track in $H\setminus \{c_1,x_1\}$ from $Y_1$ to $X_2\cup Y_2$, say $P$, with vertices
$p_1,\ldots,p_n$ say. We may assume that $p_1 \in Y_1$, and no other $p_i$ is in $Y_1$; and $p_n \in X_2 \cup Y_2$,
and no other $p_i$ is in $X_2\cup Y_2$. In particular it follows that $c_2  \not \in V(P)$.
Since $x_1\not \in V(P)$ we may assume that $p_n \not \in X_2$ (for otherwise the theorem holds),
so $p_n \in Y_2$. If any vertex of $X_1$ is in $P$ then again the theorem holds
(since $X_2$ is nonempty and none of its vertices are in $P$), so we may assume that $P$ is
disjoint from $X_1 \cup X_2$. Since $H\setminus \{c_1,c_2 \}$ is connected, there is a minimal
track $Q$ in $H\setminus \{c_1,c_2 \}$ from $X_1 \cup X_2$ to $V(P)$, and we may assume that
only its first vertex ($q$ say)
is in $X_1 \cup X_2$. If $q \in X_1\setminus X_2$, choose $x \in X_2$; if
$q \in X_2\setminus X_1$ choose $x \in X_1$; and if $q \in X_1 \cap X_2$ choose
$x \in X_1 \cup X_2$ different from $q$.
Thus we may assume that $q \in X_1$ and there exists $x \in X_2$ different from $q$ and
hence not in $Q$. So $P \cup Q$ contains a path from $q$ to $B_2$ not containing $x$, and
hence the theorem holds. This proves \ref{twisty}.\bbox

If $v$ is a vertex of $H$, the set of edges of $H$ incident with $v$ is denoted by
$\delta(v)$ or $\delta_H(v)$.
Let $H$ be bipartite and cyclically $3$-connected, and let $X$ be some set.
We say that $X$ {\em saturates} $L(H)$
if for every branch-vertex $v$ of $H$, at most one edge of $\delta_H(v)$ is
not in $X$ (or equivalently, for every $K_3$ subgraph of $L(H)$, at least two of its
vertices are in $X$). When $H$ is connected and bipartite, we speak of vertices having the
same or different {\em biparity} depending whether every track between them is even or odd
respectively. Two edges of $G$ are {\em disjoint} if they have no end in common.

\begin{thm}\label{greenedges}
Let $H$ be bipartite and cyclically $3$-connected. Let $X \subseteq E(H)$, satisfying:
\newline\noindent
(a) for every track $P$ of $H$ of length $\ge 4$ and even, with both end-edges in $X$ and
with no internal edge in $X$,  every edge in $X$ has an end in the interior of $P$
\newline\noindent
(b) there do not exist three tracks of $H$ with an end ($b$ say) in common
and otherwise vertex-disjoint, such that each contains an edge in $X$, and at least two of the
three edges of the tracks incident with $b$ do not belong to $X$.
\newline\noindent
Then either:
\begin{enumerate}
\item $X$ saturates $L(H)$, or
\item there is a branch $B$ of $H$ so that every edge in $X$ has an end in $B$, or
\item $|X| = 2$ and there is a track $P$ in $H$ of even length $\ge 4$ with both end-edges
in $X$ and no internal edge in $X$, such that there is a branch-vertex of $H$ in $P$ not
incident with either end-edge of $P$, or
\item $|X| = 4$, and the edges in $X$ form a $4$-cycle whose four vertices are all
branch-vertices, or
\item there are two vertices $c_1,c_2$ of $H$, of different biparity and not in the
same branch of $H$, so that $X = \delta(c_1) \cup \delta(c_2)$.
\end{enumerate}
\end{thm}
\Proof
\\
\\
(1) {\em There do not exist a connected subgraph $T$ of $H \setminus X$
and three mutually disjoint edges $x_1,x_2,x_3 \in X$ so that each $x_i$ has at least one
end in $T$.}
\\
\\
For suppose such $T, x_1,x_2,x_3$ exist. We may assume $T$ is a maximal connected subgraph
of $H \setminus X$.  Let $x_i$ have ends $a_i,b_i$ $(i = 1,2,3)$, where
$a_1,a_2,a_3$ have the same biparity. Make a graph $K$ with vertex set
$a_1,a_2,a_3,b_1,b_2,b_3$, where we say two vertices of $K$ are adjacent if there is a track
in $T$ joining them not using any other vertex of $K$. Since $T$ is connected and meets
all of $x_1,x_2,x_3$ it follows that there is a component of $K$ containing
an end of each of these three edges. Now if $a_1a_2$ is an edge of $K$, then the corresponding
track in $T$ is even, and hypothesis (a) is contradicted. So the only possible edges
in $K$ join some $a_i$ to some $b_j$. Also, if say $a_3$ is adjacent in $K$ to both
$b_1$ and $b_2$, then hypothesis (b) is contradicted. Since there is a component of $K$
containing an end of each of $x_1,x_2,x_3$, we may assume that $a_1b_3, b_2a_3, a_3b_3 \in E(K)$,
and the only other possible edges of $K$ are $a_1b_1, a_2b_2, a_2b_1$. In
particular, there are no more edges of $K$ incident with $a_3$ or $b_3$. Let the tracks
in $T$ corresponding to $a_1b_3, b_2a_3, a_3b_3 \in E(K)$ be $P_1,P_2,P_3$ respectively. Since
$P_3$ joins the adjacent vertices $a_3,b_3$ and does not use the edge $x_3$, it follows that
$P_3$ has nonempty interior. Choose a maximal connected subgraph $S$ of $T$ including the
interior of $P_3$ and not containing either of $a_3,b_3$. Since there are no more edges of
$K$ incident with $a_3$ or $b_3$, it follows that none of $a_1,b_1,a_2,b_2$ is in $V(S)$, and
therefore $S$ is vertex-disjoint from $P_1$ and $P_2$ as well. Consequently the only
edges of $T$ between $V(S)\cup \{a_3,b_3\}$ and the remainder of $H$ are incident with
$a_3$ or $b_3$. Since $H$ is cyclically $3$-connected and $a_3,b_3$ are adjacent, it
follows that $H \setminus  \{a_3,b_3\}$ is connected, and therefore there is an edge $sv$ of
$H$ such that $s \in V(S)$ and $v \in V(H) \setminus (V(S) \cup \{a_3,b_3\})$. Since
$T$ is maximal, it follows that $sv \in X$; and from the symmetry we may assume $v\notin
\{a_1,b_1\}$. Choose a minimal track in $S$  between $s$ and the interior of $P_3$; then it
can be extended via a subpath of $P_3$ and via $sv$ to become a track $P_4$ in $H$, of
length $\ge 2$,  from
$v$ to $b_3$, using none of $a_1,b_1,a_3$, and with only its first edge in $X$. But then the
tracks $b_1\d a_1 \d P_1 \d b_3$, $P_4$, and the one-edge track made by $x_3$, violate hypothesis (b). This
proves (1).
\\
\\
(2) {\em We may assume that $|X| \ge 3$.}
\\
\\
For if $|X| \le 1$ then statement 2 of the theorem holds; suppose $|X| = 2$, and
$X = \{a_1b_1,a_2b_2\}$ say,
where $a_1,a_2$ have the same biparity. If some branch of $H$ meets both these edges then
statement 2 of the theorem holds, so we may assume not, and in particular these four vertices are
distinct. Since $H$ is cyclically $3$-connected it follows from \ref{branchcut}
that $H \setminus \{b_1,b_2\}$ is connected, so there is a track $P$ of $H$ between
$b_1,b_2$, with even length
$\ge 4$, with first edge $b_1a_1$ and last edge $a_2b_2$. Since $a_1,a_2$ do not belong to
the same branch of $H$, there is a branch-vertex of
$H$ in $P$ not incident with either end-edge of $P$; and so statement 3 of the theorem holds.
This proves (2).
\\
\\
Now we may assume that $X$ does not saturate $L(H)$, and so there is a branch-vertex of $H$
incident with $\ge 2$ edges not in $X$. Hence there is a connected subgraph $A$ of
$H\setminus X$, containing a branch-vertex and at least two edges incident with it. Choose
such a subgraph $A$ maximal. It follows that $A$ is not contained in any branch of $H$.
By (1), there is no 3-edge matching
among the edges in $X$ that meet $A$; and since this set of edges makes a bipartite subgraph,
it follows from K\"{o}nig's theorem that there are two vertices $c_1,c_2$ so that every edge in
$X$ with an end in $A$ is incident with one of $c_1,c_2$.
\\
\\
(3) {\em We may assume that every edge in $X$ is incident with one of $c_1,c_2$.}
\\
\\
For suppose not; then there is an edge in $X$ vertex-disjoint from $V(A)\cup\{c_1,c_2\}$.
Let $B = H\setminus V(A)$.  From the maximality of $A$, every edge
of $H$ between $V(A)$ and $V(B) $ belongs to $X$ and therefore is incident with one of
$c_1,c_2$, and so there are two subgraphs $C,D$ of $H$ with $V(C) = V(A) \cup \{c_1,c_2 \}$,
$C \cup D = H$, $V(C \cap D) = \{c_1,c_2 \}$, $A \subseteq C$ and $B \subseteq D$.
In particular, $V(C),V(D) \neq V(G)$.  Since $H$ is cyclically $3$-connected it follows
from \ref{branchcut} that
one of $C$,$D$ is contained in a branch of $H$. Now $C$ is not, because it contains $A$
(and we already saw that $A$ is not contained in a branch); so $D$ is contained in a branch.
In particular, this branch contains $c_1$ and $c_2$, and also meets all edges in $X$
with no end in $V(A)$, and therefore meets all edges in $X$; but then statement 2 of the
theorem holds.  This proves (3).
\bigskip

We may assume that $c_1,c_2$ do not belong to the same branch, for otherwise statement 2 of
the theorem holds; and consequently $c_1,c_2$ are nonadjacent, and
$H\setminus \{c_1,c_2\}$ is connected, by \ref{branchcut}.

Assume that $c_1,c_2$ have the same biparity. Since $|X| \ge 3$, we may assume there
are at least two edges in $X$ incident with $c_1$, say $c_1a_1,c_1a_2$. If there is an edge
$c_2a_3$ incident with $c_2$ where $a_3 \neq a_1,a_2$, take a minimal track in
$H\setminus \{c_1,c_2\}$ between $a_3$ and one of $a_1,a_2$; it violates hypothesis (a)
of the theorem. So the only possible edges in $X$ incident with $c_2$ are $c_2a_1$ and
$c_2a_2$. If both are present, then by exchanging $c_1$ and $c_2$ it follows that there
are no more edges incident with $c_1$, and so either statement 2 or 4 of the theorem holds.
If exactly one is present, say $c_2a_1$, then the branch of $H$ containing $c_1a_1$
satisfies statment 2 of the theorem. If none are present, then any branch containing
$c_1$ satisfies statement 2. This completes the proof if $c_1,c_2$ have the same
biparity.

Now assume that  $c_1,c_2$ have different biparity. For $i = 1,2$ let $A_i =
\delta(c_i) \cap X$, and let $B_i = \delta(c_i) \setminus A_i$. We may assume that
$A_1,A_2$ are nonempty. Since $c_1,c_2$ have different
biparity, no vertex is incident with all the edges in $A_1 \cup A_2$.
If both $B_1,B_2$ are empty, then statement 5 of the theorem
holds, so we may assume at least one of them is nonempty. By \ref{twisty}, we may assume
there is a track in $H$ with first edge in $A_1$, second edge in $B_1$ (and hence second vertex
$c_1$), last vertex $c_2$ and last edge in $A_2$. By choosing such a track as short as
possible, it follows that only one edge in $A_2$ meets its interior. By hypothesis (a), all
edges in $X$ meet its interior, and hence in particular $|A_2| = 1$. But then we can replace
$c_2$ by the other end of the edge in $A_2$, and will be in the ``same biparity'' case that we
have already done. This proves \ref{greenedges}. \bbox

Now we apply what we just proved to the neighbours of a single vertex. If $G$,$J$ are
graphs, we say that $J$ {\em appears} in $G$ if there is a bipartite subdivision $H$ of $J$ so
that $L(H)$ is isomorphic to an induced subgraph of $G$. We call $L(H)$ an {\em appearance}
of $J$ in $G$. Note that if $L(H)$ is isomorphic to some induced subgraph $K$ of $G$, there
is another subdivision $H'$ isomorphic to $H$, made from $H$ by replacing each edge of $H$
by the corresponding vertex of $K$; and now $L(H') = K$ (rather than just being isomorphic to
it). So whenever it is convenient we shall assume that the isomorphism between $L(H)$ and $K$
is just equality, without further explanation. Note in particular that $E(H) = V(K)$, and so
some vertices of $G$ are edges of $H$.

When $J = K_4$, we already defined what we mean by a degenerate appearance of $J$. When $J \not = K_4$,
let us say that an appearance $L(H)$ of $J$ in $G$ is {\em degenerate} if $J = H = K_{3,3}$, and otherwise
it is {\em nondegenerate}.  So all appearances of any graph $J \not = K_4,K_{3,3}$ are nondegenerate.
If $J$ is 3-connected, we say a graph $J'$ is a
{\em $J$-enlargement} if $J'$ is 3-connected, and has a proper subgraph which is isomorphic
to a subdivision of $J$.

We remind the reader that we are currently trying to prove two statements,
\ref{linegraph} and \ref{linegraph2}.
To do so we shall prove the following:
\begin{thm} \label{hypHuse}
Let $G$ be Berge. Let $J$ be a 3-connected graph, such that either:
\begin{itemize}
\item there is a nondegenerate appearance $L(H)$ of $J$ in $G$, and there is no
$J$-enlargement with a nondegenerate appearance in $G$, or
\item $J = K_{3,3}$, there is an appearance $L(H)$ of $J$ in $G$, and
no $J$-enlargement appears in either $G$ or $\overline{G}$.
\end{itemize}
Then either $G = L(H)$, or $G$ admits a 2-join or a balanced skew partition.
\end{thm}

The proof of this will take several sections; but let us see now that \ref{hypHuse} implies
\ref{linegraph} and \ref{linegraph2}.
\\
\\
\noindent{\bf Proof of \ref{linegraph}, assuming \ref{hypHuse}.}

Let $G$ be Berge, and assume there is a nondegenerate appearance
of $K_4$ in $G$. Choose a 3-connected graph $J$ maximal (under $J$-enlargement) so that
there is a nondegenerate
appearance of $J$ in $G$; then the hypotheses of \ref{hypHuse} are satisfied, and the claim follows
from \ref{hypHuse}.
This proves \ref{linegraph}. \bbox

\noindent{\bf Proof of \ref{linegraph2}, assuming \ref{hypHuse}.}

Let $G$ be Berge, and assume it contains an induced subgraph isomorphic to $L(K_{3,3})$. We
may therefore choose a 3-connected graph $J$, either equal to $K_{3,3}$ or a $K_{3,3}$-enlargement,
maximal (in the sense of $J$-enlargement) so that there is an
appearance of $J$ in $G$. If there is a nondegenerate appearance of $J$ in $G$, then
the hypotheses of \ref{hypHuse} hold, and the claim follows from \ref{hypHuse}. So we may assume that
every appearance of $J$ in $G$ is degenerate, and in particular $J = K_{3,3}$.
If there is a $J$-enlargement which appears in $\overline{G}$,
choose it maximal; then the claim follows by applying \ref{hypHuse} to $\overline{G}$.
So we may assume that there is no
$J$-enlargement that appears in $\overline{G}$. But then again the claim follows from
\ref{hypHuse}. This proves \ref{linegraph2}. \bbox

\begin{thm}\label{vnbrs}
Let $G$ be Berge. Let $J$ be a 3-connected graph, and let $L(H)$ be an appearance of $J$ in $G$.
Let $y \in V(G) \setminus V(L(H))$, and let $X$ be the set of vertices
of $L(H)$ that are adjacent to $y$ in $G$. Then either:
\begin{enumerate}
\item $X$ saturates $L(H)$, or
\item there is a branch-vertex $v$ of $H$ with $X \subseteq \delta_H(v)$, or
\item there is a branch $B$ of $H$ with $X \subseteq E(B)$, or
\item there is a branch $B$ of $H$ with ends $b_1,b_2$ say, so that
$X \setminus E(B) = \delta_H(b_1) \setminus E(B)$, or
\item there is a branch $B$ of $H$ of odd length with ends $b_1,b_2$ say, so that
$X \setminus E(B) = (\delta_H(b_1) \cup \delta_H(b_2)) \setminus E(B)$, or
\item there are two vertices $c_1,c_2$ of $H$, of different biparity and not in the
same branch of $H$, so that $X = \delta(c_1) \cup \delta(c_2)$.
\end{enumerate}
In particular, either statements 1 or 6 hold, or there are at most two branch-vertices of $H$
incident with more than one edge in $X$; and exactly two only if statement 5 holds.
\end{thm}
\Proof

The second assertion (the final sentence) follows from the first, because if statements 2,3 or 4
hold then there is at most one branch-vertex incident with more than one edge in $X$; while if
$B,b_1,b_2$ are as
in statement 5, then since $B$ is odd, it follows that $b_1$, $b_2$ have no common neighbour, and
so no branch-vertex different from $b_1,b_2$ is incident with more than one edge in $X$.
So it remains to prove the first assertion.
\\
\\
(1) {\em Every track of $H$ with both end-edges in $X$ and no internal edge in
$X$ has length odd or 2.}
\\
\\
For since $X$ is the set of neighbours in $L(H)$ of a single vertex, it follows that every path
in $L(H)$ with both ends in $X$ and with no interior vertex in $X$ has length even or 1, and
this proves (1).
\bigskip

In particular, hypothesis (a) of \ref{greenedges} is satisfied, and so is hypothesis (b), by \ref{trianglev}.
Hence one of statements 1-5 of \ref{greenedges} applies. If \ref{greenedges}.1
applies then statement 1 of the theorem holds, and \ref{greenedges}.3
cannot apply, by (1). If \ref{greenedges}.4 applies, let the corresponding 4-cycle have vertices $b_1,b_2,b_3,b_4$
in order; then there is a track $T$ between $b_1,b_3$ not using $b_2,b_4$, since $b_1, b_3$ are branch-vertices,
and it is even since $b_1,b_3$ have the same biparity, and then the track $b_2\d b_1 \d T \d b_3 \d b_4$
contradicts (1). So \ref{greenedges}.4 cannot apply.
If \ref{greenedges}.5 holds then statement 6 of the theorem holds.  So we
may assume that \ref{greenedges}.2 applies.

Let $C$ be a branch of $H$ meeting all the
edges in $X$, and let $c_1,c_2$ be the ends of $C$. For $i = 1,2$ let $A_i$ be the set
of edges in $\delta(c_i)$ that are in $X$ and not in $C$; and let $B_i$ be the
set of edges in $\delta(c_i)$ that are not in $X$ and not in $C$. If one of $A_1,B_1$ is
empty and also $A_2$ is empty, then statement 3 or 4 of the theorem holds as required.
Suppose that $B_1,B_2$ are both empty. Choose $a_1\in A_1$ and $a_2 \in A_2$, disjoint, and
let $T$ be a track of $H$ from $c_1$ to $c_2$ with end-edges $a_1$ and $a_2$. Then no internal
edge of $T$ is in $X$, and its end-edges are in $X$, and so it cannot be even by (1); therefore $c_1$,$c_2$
have different biparity, and so $C$ is odd. But then statement 5 of the theorem holds.
So we may assume that $A_1,B_1$ are both nonempty.
\\
\\
(2) {\em We may assume $A_2$ is nonempty.}
\\
\\
For assume $A_2$ is empty. If $c_1$ meets every edge in $X$ then statement 2 of the theorem
holds, so we may assume not; and hence some edge of $C \setminus c_1$ is in $X$. Let $P$ be
a minimal subtrack of $C$ containing $c_2$ and some edge in $X$. Choose an edge $c_1a_1 \in A_1$.
Since $H \setminus \{c_1,a_1\}$ is connected, there is a track $Q$ in this graph starting from
$b_1$ say to $c_2$, where $b_1c_1 \in B_1$;
and we can extend it to a track from $a_1$ to $c_2$ with second vertex $c_1$,
first edge in $A_1$ and second edge in $B_1$. By combining the latter with $P$ we deduce
from (1) that the lengths of $P$ and $Q$ have opposite parity. On the other hand, there is
a track $R$ in $H\setminus c_1$ between $a_1$ and $c_2$; and by extending it via $a_1c_1$,
combining it with $P$ and applying (1) we deduce that the lengths of $R$ and $P$ have
the same parity. But since $H$ is bipartite, the lengths of $Q$ and $R$ have the same parity,
a contradiction. This proves (2).
\\
\\
(3) {\em $C$ has even length.}
\\
\\
For assume it is odd, and so $c_1,c_2$ have different biparity. Consequently no
vertex is incident with all the edges in $A_1 \cup A_2$, and hence we can apply \ref{twisty}
to the graph $H'$ obtained from $H$ by deleting the internal vertices and edges of $C$. We
deduce that (without loss of generality) there is a track in $H'$ with first edge in $A_1$,
second edge in $B_1$ (and hence second vertex $c_1$), last vertex $c_2$ and last edge in $A_2$.
But since $c_1,c_2$ have different biparity, this track has even length $\ge 4$, contrary
to (1). This proves (3).
\bigskip

Assume next that for $i = 1,2$ there are edges $c_ia_i \in A_i$ disjoint from each other.
There is a track in $H \setminus\{c_1,c_2\}$ between $a_1$ and $a_2$, and it is even
since $c_1,c_2$ have the same biparity, and it has length $\ge 2$. By extending it via
the edges $c_1a_1$ and $c_2a_2$ we obtain a track violating (1). So there do not exist
such edges. Hence there is a vertex $a \in V(H)$ so that $A_i = \{c_ia\}$ for $i = 1,2$.
Now there is only one branch of $H$ incident with $c_1$ and $c_2$, since $J$ is simple,
so $a$ is not in the interior of a branch, and so it is a branch-vertex.
Choose a branch-vertex $b$ of $H$ different from $c_1,c_2,a$, and choose three paths
$P_1,P_2,P_3$ between $b$ and $c_1,c_2,a$ respectively, pairwise disjoint except for $b$.
So $P_1$ and $P_2$ have lengths of the same parity, and $P_3$ has length of different parity.
We may assume that there is an edge in $X$ not incident with $a$, for otherwise statement 2
of the theorem holds, so for $i = 1,2$ there is a minimal subtrack $Q_i$ of $C$ containing
$c_i$ and an edge in $X$. If $Q_1 = C$ then (since $C$ has even length) $P_1 \cup P_2$ is the
interior of an even track with end-edges in $X$ and no internal edges in $X$, contrary to
(1). So $c_2$ is not a vertex of $Q_1$, and similarly $c_1$ is not in $Q_2$. From the
track $Q_1\d c_1\d P_1\d b\d P_2\d c_2\d a$ and (1) it follows that $Q_1$ is even; and from the
track $Q_1\d c_1\d P_1\d b\d P_3\d a\d c_2$ and (1) it follows that $Q_1$ is odd, a contradiction.
This proves \ref{vnbrs}. \bbox

We recall that $H$ is a subdivision of $J$, and $L(H)$ is an induced subgraph of $G$.
For each vertex $v$ of $J$, we denote the set of edges of $H$ incident with $v$ by $N_v$,
and for each edge $uv$ of $J$, we denote the set of edges of the branch of $H$ between
$u$ and $v$ by $R_{uv}$. So each $N_v$ and each $R_{uv}$ is a subset of $V(L(H))$.
We say a subset $X$ of $V(L(H))$ is {\em local} (with respect to $L(H)$) if either
$X \subseteq N_v$ for some
vertex $v$ of $J$, or $X \subseteq R_{uv}$ for some edge $uv$ of $J$.
In general, if $K$ is an induced subgraph of $G$,
and $F \subseteq V(G)$ is a connected set disjoint from $V(K)$, a vertex in
$V(K)$ is an {\em attachment} of $F$ if it has a neighbour in $F$.

\begin{thm}\label{smallcomp}
Let $G$ be Berge. Let $J$ be a 3-connected graph, let $L(H)$ be an appearance of $J$ in $G$,
and let $F$ be a connected set of vertices, disjoint from $V(L(H))$, such that the set
of attachments of $F$ in $L(H)$ is not local. Assume that for every $v \in F$ the set
of neighbours of $v$ in $L(H)$ does not saturate $L(H)$.  Then there is a path $P$ of $G$
with $V(P) \subseteq F$ and with ends $p_1$ and $p_2$, such that either:
\begin{enumerate}
\item there are vertices $c_1,c_2$ of $H$, not in the same branch of $H$, so that for $i = 1,2$
$p_i$ is complete in $G$ to $\delta_H(c_i)$, and there are no other edges between $V(P)$ and
$V(L(H))$, or
\item there is an edge $b_1b_2$ of $J$ (for $i = 1,2$, $r_i$ denotes
the unique vertex in $N_{b_i}\cap R_{b_1b_2}$) such that one of the following holds:
\begin{enumerate}
\item $p_1$ is adjacent in $G$ to all vertices in $N_{b_1} \setminus r_1$, and $p_2$ has a neighbour in
$R_{b_1b_2} \setminus r_1$, and every edge from $V(P)$ to
$V(L(H))\setminus r_1$ is either from $p_1$ to $N_{b_1} \setminus r_1$, or
from $p_2$ to  $R_{b_1b_2} \setminus r_1$, or
\item for $i = 1,2$, $p_i$ is adjacent in $G$ to all vertices in $N_{b_i} \setminus r_i$, and there are
no other edges between $V(P)$ and $V(L(H))$ except possibly $p_1r_1,p_2r_2$,
and $P$ has the same parity as $R_{b_1b_2}$, or
\item $p_1 = p_2$, and $p_1$ is adjacent to all vertices in $(N_{b_1} \cup N_{b_2})\setminus\{r_1, r_2\}$,
and all neighbours of $p_1$ in $V(L(H))$ belong to $N_{b_1} \cup N_{b_2} \cup R_{b_1b_2}$,
and $R_{b_1b_2}$ is even, or
\item $r_1 = r_2$, and  for $i = 1,2$, $p_i$ is adjacent in $G$ to all vertices in $N_{b_i} \setminus r_i$,
and there are
no other edges between $V(P)$ and $V(L(H))\setminus r_1$, and $P$ is even.
\end{enumerate}
\end{enumerate}
\end{thm}
\Proof
We may assume $F$ is minimal so that its set of attachments is not local. Let $X$ be the set
of attachments of $F$ in $L(H)$. Suppose first that $|F|$ = 1, $F$ = $\{y\}$ say.
Apply \ref{vnbrs} to $y$. Now \ref{vnbrs}.1 is false since by hypothesis $X$
does not saturate $L(H)$, and \ref{vnbrs}.2, and \ref{vnbrs}.3 are false since $X$ is not local.
So one of \ref{vnbrs}.4-6 holds, and the claim follows.
Consequently we may assume that $|F| \ge 2$.
\\
\\
(1) {\em There exist two attachments $x_1$,$x_2$ of $F$ so that $\{x_1,x_2\}$ is
not local. }
\\
\\
For $X \subseteq E(H)$. If there exists $x_1 \in X$
not incident in $H$ with a branch-vertex, and in some branch $B$, choose any $x_2 \in X$ not
in $B$, then  $\{x_1,x_2\}$ is not local. So we may assume that every edge in $X$ is incident
with a branch-vertex of $H$. Choose $x_1 \in X$, in some branch $B_1$ of $H$, and incident
with a branch-vertex $b_1$. There exists $x_2 \in X$ not incident with $b_1$, and we may assume
that $x_2 \in E(B_1)$, for otherwise $\{x_1,x_2\}$ is not local. Hence $x_2$ is incident
with the other end $b_2$ say of $B_1$. There exists $x_3 \in X$ not belonging to $E(B)$, and
it cannot share an end both with $x_1$ and with $x_2$, so we may assume $x_3$ is not incident
with $b_1$. But then $\{x_1,x_3\}$ is not local, as required. This proves (1).

\bigskip
From the minimality of $F$, it follows that $F$ is minimal such that $x_1$ and $x_2$ are both
attachments of $F$, and so (since $x_1$ and $x_2$ are nonadjacent), $F$ is the interior of
a path with vertices $x_1, p_1,\ldots,p_n,x_2$ in order. Let $X_1$ be the set of attachments in
$L(H)$ of $F\setminus p_n$, and let $X_2$ be the attachments of $F\setminus p_1$.
From the minimality of $F$, $X_1$ and $X_2$ are both local.
\\
\\
(2) {\em If there is an edge $uv$ of $J$ so that $X_1 \subseteq N_u$ and
$X_2 \subseteq R_{uv}$ then the theorem holds.}
\\
\\
For let the ends of $R_{uv}$ be $r_1,r_2$ where $r_1 \in N_u$. Since $X$ is not local, it
follows that $p_1$ has a neighbour in $N_u \setminus r_1$ and $p_n$ has a neighbour in
$R_{uv} \setminus r_1$. If $p_1$ is adjacent to every vertex in $N_u \setminus r_1$ then
statement 2.a of the theorem holds, so we may assume $p_1$ has a neighbour $s_1$ and a nonneighbour $s_2$
in $N_u \setminus r_1$. Let $Q$ be the path between $r_2$ and $s_1$ with interior in  $F \cup R_{uv}\setminus r_1$.
Choose $w \in V(J)$ so that $s_1 \in R_{uw}$.  Now $H$ is a
subdivision of a $3$-connected graph, so if we delete all edges of $H$ incident with $u$ except
$s_1$, the graph we produce is still connected. Consequently there is a track of $H$ from
$u$ to $v$ with first edge $s_1$; and hence there is a path $S_1$ of $L(H)$ from $s_1$ to $r_2$,
vertex-disjoint from $R_{uv}\cup N_u$ except for its ends. Indeed, if we delete from $H$
both the vertex $w$ and all edges incident with $u$ except $s_2$, the graph remains connected;
so there is a path $S_2$ of $L(H)$ between $s_2$ and $r_2$, vertex-disjoint from
$R_{uv}\cup N_u \cup R_{uw} \cup N_w$ except for its ends. Now $S_1$ and $S_2$ have the same
parity since $H$ is bipartite.
Yet $S_1$ can be completed via $r_2\d Q\d s_1$ and $S_2$ can be completed via
$r_2\d Q\d s_1\d s_2$, a contradiction. This proves (2).
\\
\\
(3) {\em If there are nonadjacent vertices $v_1,v_2 \in V(J)$ so that
$X_i \subseteq N_{v_i}$ for $i = 1,2$, then the theorem holds.}
\\
\\
Let $A_1$ be the set of vertices in $N_{v_1}$ adjacent to $p_1$, and
$B_1 = N_{v_1} \setminus A_1$; and let $A_2$ be the set of vertices in $N_{v_2}$ adjacent
to $p_n$, and $B_2 = N_{v_2} \setminus A_2$. So $X = A_1 \cup A_2$.  If both $B_1$ and $B_2$
are empty then statement 1 of the theorem holds, so we may assume that at least one of
$B_1$,$B_2$ is nonempty. Certainly $A_1$ and $A_2$ are both nonempty, so there is a track in
$H$ from $v_1$ to $v_2$ with end-edges in $A_1$ and $A_2$ respectively. Hence there is a path
$S_1$ in $L(H)$ from $A_1$ to $A_2$, vertex-disjoint from $N_{v_1}\cup N_{v_2}$
except for its ends. Since $X = A_1 \cup A_2$ is not local, there is no $w \in V(J)$ with
$A_1 \cup A_2 \subseteq N_w$.  Hence we can apply \ref{twisty}, and we deduce (possibly after
exchanging $v_1$ and $v_2$) that there is a path $S_2$ in $L(H)$ with first vertex in
$A_1$, second vertex in $B_1$, last vertex in $A_2$, and otherwise disjoint from
$N_{v_1}\cup N_{v_2}$. Since $H$ is bipartite, $S_1$ and $S_2$ have opposite
parity; but they can both be completed via $F$, a contradiction. This proves (3).
\\
\\
(4) {\em If there are adjacent vertices $v_1,v_2 \in V(J)$ so that
$X_i \subseteq N_{v_i}$ for $i = 1,2$, then the theorem holds.}
\\
\\
For $i = 1,2$ let $r_i$ be the end of $R_{v_1v_2}$ in $N_{v_i}$. Let $A_1$ be the set of
vertices in $N_{v_1}\setminus r_1$ adjacent to $p_1$, and
$B_1 = N_{v_1}\setminus (A_1 \cup r_1)$; and define $A_2,B_2$ similarly. Then
$X \subseteq A_1 \cup A_2 \cup \{r_1,r_2\}$. By (2) we may assume that
$A_1$ and $A_2$ are both nonempty. Suppose that both $B_1$ and $B_2$ are empty. Then
there is a cycle in $J$ of length $\ge 4$ using the edge $v_1v_2$, and so there
is a path in $L(H)$ of length $\ge 2$ from $A_1$ to $A_2$ with no internal vertex in
$N_{v_1} \cup V(R_{v_1v_2}) \cup N_{v_2}$. The union of this path with $R_{v_1v_2}$
induces a hole, and so does its union with $F$, and therefore these two paths have
lengths of the same parity. Consequently either statement 2.b or 2.d of the theorem
holds. So we may assume that at least one of $B_1,B_2$ is nonempty. There is a path $S_1$
from $A_1$ to $A_2$ with no vertex in $N_{v_1}\cup N_{v_2} \cup R_{v_1v_2}$ except for its ends.
Suppose that there is no vertex $w \in V(J)$ with $A_1 \cup A_2 \subseteq N_w$. Then we can
apply \ref{twisty} to the graph obtained from $H$ by deleting the edges and internal
vertices of the branch between $v_1$ and $v_2$. We deduce (possibly after exchanging $v_1$
and $v_2$) that there is a path $S_2$ of $L(H)$ with first vertex in $A_1$, second vertex in
$B_1$, last vertex in $A_2$, and otherwise disjoint from $N_{v_1}\cup N_{v_2} \cup R_{v_1v_2}$.
Since $H$ is bipartite, $S_1$ and $S_2$ have opposite
parity; but they can both be completed via $F$, a contradiction. Consequently there is a
vertex $w  \in V(J)$ with $A_1 \cup A_2 \subseteq N_w$. Since $H$ is bipartite, it follows
that $R_{v_1v_2}$ has odd length, and in particular $r_1 \not = r_2$. Since $|N_{v_i}\cap N_w|
\le 1$ (since $J$ is simple) it follows that $|A_i| = 1$, $A_i = \{a_i\}$ say, for $i = 1,2$.
Since $X$ is not local it is not a subset of $N_w$ and so there is a vertex of $R_{v_1v_2}$
in $X$. Since $X_i \subseteq N_{v_i}$ for $i = 1,2$, no internal vertex of $R_{v_1v_2}$
is in $X$, so we may assume that $r_1 \in X$. Since $r_1 \notin N_{v_2}$ it follows that
$r_1 \notin X_2$, and hence $p_1$ is the only vertex in $F$ adjacent to $r_1$. Now the
hole $p_1\d \cdots\d p_n\d a_2\d a_1\d p_1$ is even, and so $n$ is even. If we delete the vertex $v_2$
and the edge $a_1$ from $H$, what remains is still connected, and so contains a track from
$w$ to $v_1$. Hence there is a path $T$ in $L(H)$ from some $a_3 \in N(w)$ to $r_1$, disjoint
from $N_{v_2} \cup a_1$. But $T$ can be completed to a hole via $r_1\d R_{v_1v_2}\d r_2\d a_2\d a_3$
and via $r_1\d p_1\d \cdots\d p_n\d a_2\d a_3$, and these two completions have different parity,
a contradiction. This proves (4).
\\
\\
(5) {\em If $X_1 \cap X_2$ is nonempty, and in particular if one of $p_2,\ldots,p_{n-1}$ has a
neighbour in $L(H)$, then the theorem holds.}
\\
\\
For any neighbour in $L(H)$ of one of $p_2,\ldots,p_{n-1}$ belongs to $X_1 \cap X_2$, so assume
$x \in X_1 \cap X_2$.  Then $x \in R_{v_1v_2}$ for a
unique edge $v_1v_2$ of $J$, and $x \in N_v$ for at most two $v \in V(J)$, namely $v_1$ and
$v_2$. Since both $X_1$ and $X_2$ are local, each is a subset of one of $N_{v_1}, N_{v_2},
R_{v_1v_2}$, and they are not both subsets of the same one. So we may assume that
$X_1 \subseteq N_{v_1}$.  Hence either $X_2 \subseteq N_{v_2}$  or $X_2 \subseteq R_{v_1v_2}$,
and therefore the theorem holds by (5) or (2). This proves (5).
\\
\\
(6) {\em If there is a vertex $u$ and an edge $v_1v_2$ of $J$ so that $X_1 \subseteq N_u$
and $X_2 \subseteq R_{v_1v_2}$ then the theorem holds.}
\\
\\
For by (2) we may assume $u$ is different from $v_1$ and $v_2$. Choose a cycle $C_1$ of
$H$ using the branch between $v_1$ and $v_2$ and not using $u$, and choose a minimal track $S$
in $H\setminus \{v_1,v_2\}$ between $u$ and $V(C_1)$. Let the ends of $S$ be
$u$ and $w$ say. Hence in $L(H)$ there are three vertex-disjoint paths, from $N_{v_1}$,
$N_{v_2}$,$N_u$ respectively to $N_w$, and there are no edges between them except in the
triangle $T$ formed by their
ends in $N_w$. If $p_n$ has a unique neighbour (say $r$) in $R_{v_1v_2}$, then $r$ can be
linked onto the triangle $T$, contrary to \ref{trianglev}. If $p_n$ has two nonadjacent
neighbours in $R_{v_1v_2}$, then $p_n$ can be linked onto the triangle $T$, contrary to
\ref{trianglev}. So $p_n$ has exactly two neighbours in $R_{v_1v_2}$, and they are adjacent.
If $p_1$ is adjacent to all of $N_u$, then statement 1 of the theorem holds, so we may
assume that $p_1$ has a neighbour and a
non-neighbour in $N_u$. Let $A$ be the neighbours of $p_1$ in $N_u$ and $B = N_u \setminus A$.
In $H$ there is a cycle $C_2$ using the branch between $v_1$ and $v_2$, and using an edge
in $A$ and an edge in $B$. (To see this, divide $u$ into two adjacent vertices, one incident
with the edges in $A$ and the other with those in $B$, and use Menger's theorem to deduce that
there are two vertex-disjoint paths between these two vertices and $\{v_1,v_2\}$.) Hence
in $G$, there is a path between $N_{v_1}$ and $N_{v_2}$ using a unique edge of $N(u)$, and that
edge is between a vertex $a \in A$ say and some vertex in $B$. Hence $a$ can be linked onto
the triangle formed by $p_n$ and its two neighbours in $R_{v_1v_2}$, a contradiction. This proves (6).
\\
\\
(7) {\em If there are edges $u_1v_1$ and $u_2v_2$ of $J$ with $X_i \subseteq R_{u_iv_i}$ for
$i = 1,2$, then the theorem holds.}
\\
\\
For in this case it follows that the edges $u_1v_1$ and $u_2v_2$ are different, and hence
we may assume that $v_2$ is different from $u_1$ and $v_1$, and $v_1$ is different from
$u_2$ and $v_2$; possibly $u_1 = u_2$. If $p_1$ has exactly two neighbours in $R_{u_1v_1}$
and they are adjacent, and also $p_n$ has exactly two neighbours in $R_{u_2v_2}$
and they are adjacent, then statement 1 of the theorem holds; so we may assume that $p_1$ has either
only one neighbour, or two nonadjacent neighbours, in $R_{u_1v_1}$.  There is a cycle
in $H$ using the branch between $u_1$ and $v_1$, and using $u_2$ and not $v_2$ (since
$J \setminus v_2$ is $2$-connected). There correspond two paths in $L(H)$, say $P$ and $Q$,
from $N_{u_1}$ and $N_{v_1}$ respectively to $N_{u_2}$, disjoint from each other, and there is a
third path $R$ say from $p_1$ to $N_u$ via $F$ and a subpath of $R_{u_2v_2}$. There are no edges
between these paths except within the triangle $T$ formed by their ends in $N_u$.
If $p_1$ has only one neighbour $r \in R_{u_1v_1}$, then we may assume that $r$ is in the
interior of $R_{u_1v_1}$, by (6), and so $r$ can be linked onto $T$, contrary to
\ref{trianglev}. If $p_1$ has two nonadjacent neighbours in $R_{u_1v_1}$, then $p_1$
can be linked onto $T$, again a contradiction. This proves (7).

\bigskip
But (2)-(7) cover all the possibilities for the local sets $X_1$ and $X_2$, and so this proves
\ref{smallcomp}. \bbox

\section{Major attachments to a line graph}

In this section we study anticonnected sets of ``major'' vertices, and their common neighbours
in $L(H)$. Conveniently we can again apply \ref{greenedges}, because of the following.

\begin{thm}\label{hypothesisb}
Let $G$ be Berge, and let $L(H)$ be an appearance in $G$ of a 3-connected graph $J$,
let $Y$ be an anticonnected set of vertices in
$V(G) \setminus V(L(H))$, and let $X$ be the set of $Y$-complete vertices in $L(H)$.
Then either:
\begin{itemize}
\item $J = K_{3,3}$ or $K_4$, and $L(H)$ is degenerate, and there is a $J$-enlargement
that appears in $\overline{G}$, or
\item $X$ satisfies hypotheses (a) and (b) of \ref{greenedges}.
\end{itemize}
\end{thm}
\Proof  That $X$ satisfies hypothesis (a) of \ref{greenedges}
is immediate from \ref{greentouch}, so it remains to prove that if hypothesis (b) is false then
the first alternative of the theorem holds.  We begin with:
\\
\\
(1) {\em If $\{x_1,x_2,x_3\} \subseteq X$ is a matching in $H$, and $P$ is an even track of $H$
with end-edges $x_1$ and $x_2$ and with no edge of $X$ in its interior, then
$P$ has length $4$ and its middle vertex is incident with $x_3$.}
\\
\\
Let the vertices of $P$ be $p_1,\ldots,p_n$; so $n \ge 4$ and is even, $x_1$ is $p_1p_2$,
and $x_2$ is $p_{n-1}p_n$. Since ${x_1,x_2,x_3}$ is a matching, it follows that $x_3$
is not incident with any of $p_1,p_2,p_{n-1},p_n$. By \ref{greentouch}, $x_3$ in incident
with some vertex of $P$, and hence shares an end with at least two internal edges of $P$.
By \ref{greentouch2}, $x_3$ shares an end with one of the first two edges of $P$,
and with one of the last two. It is therefore incident with $p_2$ and with $p_{n-2}$. Since
$H$ is bipartite and $n$ is even, it follows that $p_2 = p_{n-2}$ and so $n = 4$. This
proves (1).

\bigskip
Now assume that hypothesis (b) of \ref{greenedges} fails to hold; and we may therefore
assume that $Y$ is minimal such that it is anticonnected and hypothesis (b) is false. Then
there are three tracks $T_1,T_2,T_3$ of $H$ with an end ($v$ say) in common
and otherwise vertex-disjoint, such that each contains an edge in $X$, and at least two of the
three edges of the tracks incident with $v$ do not belong to $X$. Let $T_i$ be from $a_i$
to $v$ for $ i = 1,2,3$. We may assume that the only edge of each $T_i$ in $X$ is the edge
($x_i$, say) incident with $a_i$. Now two of $T_1,T_2,T_3$ have lengths of the same parity,
say $T_1,T_2$; and yet ${x_1,x_2,x_3}$ is a matching in $H$, since at most one of these edges
is incident with $v$. Hence from (1), $T_1$ and $T_2$ both have length $2$, and $T_3$ has
length $1$. For $i = 1,2$, let $u_i$ be the middle vertex of $T_i$, let $e_i$ be the edge $vu_i$
of $H$, and let $A_i$ be the
set of all $a \in V(H)$ so that $au_i \in E(H) \cap X$. Let $C$ be the set of all
$c \in V(H)$ so that $cv \in E(H) \cap X$. Since $H$ is bipartite it follows that
$C$ is disjoint from $A_1 \cup A_2$. Let $Q$ be an antipath of $G$ between $e_1$ and $e_2$,
with interior in $Y$. Since $Q$ can be completed to an antihole via $e_2\d x_1\d x_2\d e_1$
it follows that $Q$ is odd. From the minimality of $Y$ it follows that $Y = V(Q^*)$.
Now $A_1 \cap A_2$ is empty; for if there exists
$a \in A_1 \cap A_2$, then $au_1 = v_1$ and $au_2 = v_2$ are vertices of $G$ in $X$, and $Q$ could be
completed to an odd antihole via $e_2\d v_1\d x_3\d v_2\d e_1$, a contradiction. So $A_1,A_2,C$ are
mutually disjoint. Note also that every edge of $H$ in $X$ is incident with one of $v,u_1,u_2$,
by \ref{greentouch}.
\\
\\
(2) {\em We may assume there is no track in $H\setminus \{u_1,u_2,v\}$ between $A_1$ and $A_2$.}
\\
\\
For assume $T$ is a minimal such track. It is even (since all the vertices in $A_1 \cup A_2$
have the same biparity), and no internal vertex is in $A_1 \cup A_2$; and none of its
edges are in $X$, since all edges in $X$ are incident with one of $v,u_1,u_2$. Hence
by (1) applied to the track formed by $T$, $x_1$ and $x_2$, it follows that $T$ has length 2
and its middle vertex is $c$ (for every $c \in C$.) So $C = \{c\}$,
say, where $ca_1$ and c$a_2$ are edges. Also, every edge in $X$ is incident with one of
the vertices of $T$, and so for $i = 1,2$, $|A_i| = 1$, and $A_i = \{a_i\}$. Suppose that
$|V(H)| \ge 7$, and let $Z$
be any component of $H \setminus \{v,u_1,u_2,a_1,a_2,c\}$. At most one of
$u_1,u_2,c$ has a neighbour in $Z$, for otherwise since they have the same biparity there
would be an even track in $H$ with end-edges two of $x_1,x_2,x_3$ disjoint from the
third, contrary to (1). Similarly at most one of $a_1,a_2,v$ has neighbours in $Z$.
Since $H$ is cyclically 3-connected, there is exactly one of each biparity with a neighbour
in $Z$, and they are nonadjacent, so from the symmetry we may assume that
$u_1$ and $a_2$ are the two vertices of these six with neighbours in $Z$. Let
$S$ be a track between $u_1$ and $a_2$ with interior in $Z$; then $c\d v\d u_1\d S\d a_2\d u_2$
is an even path with end-edges in $X$ and no internal edge in $X$, of length $\ge 6$,
contrary to (1). So there is no such $Z$, and hence $|V(H)| = 6$, and the only other
possible edges of $H$ are $a_1u_2$ and $a_2u_1$. Since $H$ is cyclically 3-connected,
at least one of these is present, so we may assume $a_1u_2$ is an edge. Thus $J = K_{3,3}$
or $K_4$, and $L(H)$ is degenerate. We claim that also there is an appearance in $\overline{G}$
of a 3-connected graph with more edges than $J$. For recall that
$Y = V(Q^*)$, where $Q$ is an odd antipath between $e_1$ and $e_2$. Let the end-vertices
of $Q^*$ be $y_1$ and $y_2$, where $y_i$ is nonadjacent to $e_i$. For $i = 1,2$, let
$Y_i = Y \setminus y_i$, and let $X_i$ be the set of $Y_i$-complete vertices in $V(H)$.
From the minimality of $Y$ it follows that each $X_i$ satisfies hypothesis (b) of
\ref{greenedges}, and so by \ref{greenedges}, $X_i$ saturates $L(H)$ (for the other
possibilities listed in \ref{greenedges} cannot hold). Since $e_1 \in X_1 \setminus X_2$ and
$e_2 \in X_2 \setminus X_1$, and $X_1 \cap X_2 = X$, it follows that every edge of $H$ is
in one of $X_1,X_2$, and every branch-vertex of $H$ is incident with exactly two edges
in $X_1$ and two in $X_2$. So $a_1u_2 \in X_1$, $a_1c \in X_2$, $a_2c \in X_1$,
and if the edge $u_1a_2$ exists it belongs to $X_2$. But then $\overline{G}|(V(L(H))\cup Y)$
is an appearance in $\overline{G}$ of a $J$-enlargement, and so the
theorem holds. This proves (2).

\bigskip
From (2), there is a partition $(F_1,F_2)$ of $V(H) \setminus \{u_1,u_2,v\}$ with
$A_i \subseteq F_i$ $(i = 1,2)$ so that there is no edge between $F_1$ and $F_2$.
We may assume that $C \cap F_2$ is nonempty, so (since $H \setminus \{u_2,v\}$ is connected,
by \ref{branchcut})
there is a track in $H \setminus \{u_2,v\}$ between $u_1$ and $C \cap F_2$. Since there are
no edges between $F_1$ and $F_2$, this track does not meet $F_1$, and so none of its edges
are in $X$. Since $u_1$ and the vertices in $C$ have the same biparity, this track is even; and
it is the interior of a track $P$ say from $a_1$ to $v$ with first and last edges in $X$ and
no other edges in $X$. From \ref{greentouch}, the edge $cv$ must have an end in the interior
of $P$, and so $c$ is a vertex of $P$, for every $c \in C$; and so (since we may assume
$P$ is minimal), $|C| = 1$, $C = \{c\}$ say. From (1), we deduce that $P$ has length 4,
and every vertex in $A_2$ is the middle vertex of $P$, and so $A_2 = \{a_2\}$, and $a_2$
is adjacent to both $u_1$ and $c$. The three tracks $a_1\d u_1\d a_2$ , $v\d c\d a_2$,$u_2\d a_2$
are another instance of three tracks violating hypothesis (b) of \ref{greenedges}, and so
by (2) applied to these three tracks, we deduce that there is no track in
$H \setminus \{c,u_1,a_2\}$ between $A_1$ and $v$. Since $A_1$ is nonempty and is contained
in $F_1$, there is a maximal connected subset $F \subseteq F_1$ with nonempty intersection with
$A_1$. So the only vertices of $L(H)$ not in $F$ which might have neighbours in $F$ are
$u_1,u_2,v$. On the other hand, neither $u_2$ nor $v$ has a neighbour in $F$, since there is
no track in $H \setminus {c,u_1,a_2}$ between $A_1$ and $v$. So $u_1$ is the only such vertex,
contradicting that $H$ is cyclically 3-connected.  This proves \ref{hypothesisb}.\bbox

The remainder of this section is concerned with analyzing the five possible outcomes
of \ref{greenedges}. Let us say a vertex $v \in V(G) \setminus V(L(H))$ is {\em major} (with
respect to $L(H)$) if the set of its neighbours in $L(H)$ saturates $L(H)$.
An appearance $L(H)$ of $J$ in $G$ is {\em overshadowed} if there is a branch $B$ of $H$ with
odd length $\ge 3$, with ends $b_1,b_2$, so that some vertex of $G$ is nonadjacent in $G$ to at
most one vertex in $\delta(b_1)$ and at most one in $\delta(b_2)$.

\begin{thm}\label{nonsat}
Let $G$ be Berge, let $L(H)$ be an appearance in $G$ of a 3-connected graph $J$,
let $Y$ be an anticonnected set of major vertices in
$V(G) \setminus V(L(H))$, and let $X$ be the set of $Y$-complete vertices in $L(H)$.
Assume that $X$ does not saturate $L(H)$. Then $J = K_{3,3}$ or $K_4$. Moreover:
\begin{itemize}
\item If $J = K_{3,3}$, then either:
\begin{itemize}
\item there is an overshadowed appearance of $J$ in $G$, or
\item $L(H)$ is degenerate and there is an overshadowed appearance of $J$ in $\overline{G}$, or
\item $L(H)$ is degenerate and there is a $J$-enlargement that appears in $\overline{G}$.
\end{itemize}
\item If $J = K_4$ and $L(H)$ is nondegenerate then there is an overshadowed appearance of $J$ in $G$.
\item If $J = K_4$ and $L(H)$ is degenerate, let $V(J) = \{1,2,3,4\}$, and for $1 \le i < j \le 4$
let $B_{ij} = B_{ji}$ be the branch of $J$ joining $i$ and $j$; let the end-edges of $B_{ij}$
be $r_{ij}$ (incident with $i$) and $r_{ji}$ (incident with $j$); and let $R_{ij}$
$R_{ij}$ be the path $L(B_{ij})$ in $G$ (so $r_{ij}$ and $r_{ji}$ are the end-vertices of this path).
Let $R_{1,3},R_{1,4},R_{2,3},R_{2,4}$ have length $0$. Then either:
\begin{enumerate}
\item there is an overshadowed appearance of $J$ in $G$, or
\item (up to symmetry) there exist
nonadjacent $y,y' \in Y$ so that the neighbours of $y$ in $L(H)$ are $r_{1,2},r_{1,4},r_{3,2},
r_{3,4},r_{2,4}$ and possibly $r_{1,3}$, and the neighbours of $y'$ are $r_{2,1},r_{2,3},
r_{4,1},r_{4,3},r_{1,3}$ and possibly $r_{2,4}$, or
\item $R_{1,2},R_{3,4}$ have length $1$.
\end{enumerate}
\end{itemize}
\end{thm}
\Proof
Suppose not; then we may assume that $Y$ is minimal such that it is anticonnected and its
common neighbours do not saturate $L(H)$. Choose two vertices of $L(H)$, both incident in
$H$ with the same branch-vertex of $H$, and both not in $X$. Then there is an antipath
joining them with interior in $Y$, and the common neighbours of the interior of this
antipath do not saturate $L(H)$. From the minimality of $Y$ it follows that this antipath
contains all vertices in $Y$. Consequently, $Y$ is the vertex set of an antipath $Q$ say, with
ends $y_1$,$y_2$ say. From the hypothesis, $|Y| \ge 2$, since the neighbours of any vertex in
$Y$ saturate $L(H)$, so $y_1,y_2$ are distinct. Now for $i = 1,2$,
$Y\setminus y_i$ (= $Y_i$ say) is anticonnected; let $X_i$ be the set of $Y_i$-complete
vertices in $L(H)$. From the minimality of $Y$, both $X_1$ and $X_2$ saturate $L(H)$.
\\
\\
(1) {\em For every branch-vertex $b$ of $H$, $X$ contains all edges of $H$ incident with
$b$ except at most two; and if there are two such edges incident with $b$ not in $X$, then
one is in $X_1 \setminus X_2$ and the other in $X_2\setminus X_1$.}
\\
\\
For both $X_1$ and $X_2$ saturate $L(H)$. Therefore, $X_1$
contains at least all except one of the edges of $H$ incident with $b$, and so does $X_2$.
Since $X_1 \cap X_2  = X$ this proves (1).

\bigskip
Now by \ref{hypothesisb}, we may assume that $X$ satisfies one of the five alternatives of
\ref{greenedges}, for otherwise the theorem holds.
Let us examine them. It is convenient to postpone the case $J = K_4$ until later.
\\
\\
(2) {\em If $J \not = K_4$ then the theorem holds.}
\\
\\
For \ref{greenedges}.1 does not hold by hypothesis. Since
$H$ has at least five branch-vertices, and each of them is incident with an edge in $X$, it
follows that \ref{greenedges}.3 and \ref{greenedges}.4 do not hold.
We may therefore assume that either \ref{greenedges}.2 or \ref{greenedges}.5 holds.
In the first case, let $B$ be a branch of $H$ with ends $b_1,b_2$ so that every edge in $X$
has an end in $V(B)$; and in the second case, let $b_1,b_2$ be vertices of $H$, of opposite
biparity and not in the same branch, so that $X = \delta(b_1) \cup \delta(b_2)$, and let $B$
be the subgraph of $H$ consisting just of $b_1$ and $b_2$ (and no edges). (Note that in
this second case, $b_1$ and $b_2$ need not be branch-vertices.)
Let $H' = H\setminus V(B)$. By \ref{branchcut}, $H'$ is connected,
but none of its edges are in $X$. By (1), all vertices of $H'$
have degree $\le 2$ , and so $H'$ is a path or hole in $H$. If $H'$
is a hole, let the branch-vertices of $H$ that lie in this hole be $p_1,\ldots,p_n$ in order
(there are at least three, since $J \not = K_4$); define $p_0 = p_n$, and let $B_i$ the branch
of $H$ with ends $p_{i-1}$ and $p_i$ , for $1 \le i \le n$. (So $H'$ is the union
of $B_1,\ldots,B_n$.) If $H'$ is a path, let its ends be $p_0$ and $p_n$, and let
the branch-vertices of $H$ in its interior be $p_1,\ldots,p_{n-1}$, so that $p_0,\ldots,p_n$
are in order (so $n \ge 2$); and for $1 \le i \le n$ let $B_i$ be the path of $H'$
between $p_{i-1}$ and $p_i$. In either case, for $1 \le i \le n$ let the end-edges of $B_i$ be
$e_i$ (incident with $p_{i-1}$) and $f_i$ (incident with $p_i$). So for $1 \le i < n$, one
of $e_{i+1},f_i$ is in $X_1$ and the other is in $X_2$; and the same holds for $e_1,f_n$
if $H'$ is a hole. We recall that $Q$ is an antipath in $Y$ between $y_1$ and $y_2$; there
are two cases depending whether $Q$ is odd or even.

Assume first that $Q$ is even. Then there do
not exist two disjoint edges of $H'$, one in $X_1$ and the other in $X_2$, for
if there were we could complete $Q$ to an odd antihole in $G$ using them. Since there is
a branch-vertex in $H'$ different from $p_0,p_n$, it follows that $H'$
has an edge in $X_1$ and one in $X_2$, and since all edges of the first type meet all those
of the second type, there are at most two of each. Suppose that there are two of each. Then
$H'$ is a hole of length 4, and its four edges are alternately in $X_1$ and in
$X_2$. If $b_1$ is not a branch-vertex then it is not adjacent to $b_2$, and so all its neighbours
in $H$ are in $H'$; while if $b_1$ is a branch-vertex then it has degree $\ge 3$, and at least all
except one of its neighbours are in $H'$. So in either case $b_1$ has at least two neighbours in
$H'$, and so does $b_2$. Since $H$ is bipartite, they are each adjacent to exactly two vertices
of $H'$, and these two vertices are nonadjacent; and $b_1$ and $b_2$ cannot
be adjacent to the same pair of vertices of $H'$, since $H$ is cyclically
3-connected.  Since $H$ has at least five branch-vertices, one of $b_1,b_2$ is a branch-vertex,
and so $B$ is a branch, and hence both $b_1$ and $b_2$ are branch-vertices.
All four vertices of $H'$ have degree $\ge 3$ in $H$, and so
$J = K_{3,3}$. Now $B$ has odd length. If it has length $\ge 3$, then $L(H)$ is overshadowed
since any vertex in $Q$ is complete in $G$ to the four vertices of $L(H)$ that in $H$ are the
edges joining $V(B)$ and $V(H')$. So we may assume $B$ has length 1; let its edge be $b$.
All branches of $H$ now have length 1, so $L(H)$ is degenerate. Then
$\overline{G} | (V(L(H)) \setminus b) \cup V(Q)$ is an appearance of $J = K_{3,3}$ in
$\overline{G}$; and
this appearance is overshadowed, because $Q$ has even length $\ge 2$, and because the vertex $b$
is adjacent in $\overline{G}$ to each of $e_1,e_2,e_3,e_4$. So in either case the theorem holds.
Hence we may assume that exactly one edge of
$H'$ is in $X_2$ say. If $H'$ is a hole, then it contains at
least 3 branch-vertices of $H$, each incident with an edge of $H'$ in $X_2$,
a contradiction. So $H'$ is a path, and by the same argument, $n \le 3$. Assume
that $n = 3$. Then it follows that $e_2 = f_2 \in X_2$, and $f_1,e_3 \in X_1$. Since
$p_1$ and $p_2$ are branch-vertices they are each adjacent to one of $b_1,b_2$, and not the
same one since $H$ is bipartite. Hence $b_1, b_2$ have opposite biparity, and hence have
no common neighbours in $H'$. Since they each have at least two neighbours in $V(H')$,
it follows that
they both have exactly two, and each of $p_0,p_1,p_2,p_3$ is adjacent to exactly one of
$b_1,b_2$. But then $p_0,p_3$ have degree 2 in $H$, and so $H$ has only four branch-vertices,
a contradiction. So $n = 2$. Since $H$ has $\ge 5$ branch-vertices, it follows that
$b_1,b_2,p_0,p_2$ are all branch-vertices, and so $B$ is a branch, and $p_0,p_2$ are both
adjacent to both $b_1,b_2$.
Since $p_1$ is adjacent to one of $b_1,b_2$, it follows that $p_0,p_1,p_2$ all have the
same biparity, and so $B_1,B_2$ both have even length, and so does $B$. We may assume that
$f_1 \in X_1$ and $e_2 \in X_2$. Consequently, $e_1 \not \in X_1$, and $f_2 \not \in X_2$;
and not both $e_1 \in X_2$, and $f_2 \in X_1$, for all edges in $X_1$ meet all those in $X_2$.
So we may assume that $e_1$ is not in $X_1 \cup X_2$. Since both $X_1$ and $X_2$ saturate
$L(H)$, it follows that the edges $p_0b_1$ and $p_0b_2$ are both in both $X_1$ and $X_2$,
and hence are both in $X$. Also, $p_1$ is adjacent to one of $b_1,b_2$, say $b_1$, and
the edge $p_1b_1 \in X$. So the track $b_2\d p_0\d B_1\d p_1\d b_1$ is even, has length $\ge 4$,
and both its end-edges are in $X$ and none of its other edges are in $X$. Consequently
every edge of $H$ in $X$ has an end in the interior of this track, by \ref{greentouch}, and
in particular no edge incident with $p_2$ is in $X$, a contradiction. This completes the
proof of (2) when $Q$ is even.

Now we prove (2) when $Q$ is odd. For this we shall use repeatedly the fact that
if $e_1 \in X_1 \setminus X_2$, and $e_2 \in X_2 \setminus X_1$, and $e_0 \in X$, and
$e_1$ meets $e_2$, then one of them meets $e_0$; for otherwise $Q$ could be completed to
an odd antihole in $G$ via $y_2\d e_1\d e_0\d e_2\d y_1$. If $H'$ is a hole, then
we may choose $p_i,p_j$ nonadjacent; and since one of $f_i,e_{i+1}$ is in $X_1$ and the
other in $X_2$, and there is an edge in $X$ between $p_j$ and $\{b_1,b_2\}$, this is
a contradiction. So $H'$ is a path. Since one of $f_1,e_2$ is in $X_1$ and the
other in $X_2$, every edge in $X$ meets one of these two edges, and since for $1 \le i <n$ there
is an edge in $X$ between $p_i$ and $\{b_1,b_2\}$, it follows that $n \le 3$. Assume that
$n = 3$; then there is no edge in $X$ incident with $p_3$, and so $p_3$ is not a branch-vertex,
and similarly $p_0$ is not a branch-vertex (by the same argument with $p_1$ replaced by $p_2$);
but then $H$ has only four branch-vertices, a contradiction. So $n = 2$, and therefore
$b_1,b_2, p_0,p_2$ are all branch-vertices, and $B$ is a branch, and $p_0,p_2$ are adjacent
to both of $b_1,b_2$.
Since there is an edge in $X$ between $p_0$ and $\{b_1,b_2\}$, which therefore meets either
$f_1$ or $e_2$, it follows that $B_1$ has
length 1, and similarly so does $B_2$; but then since $H$ is bipartite, $p_1$ is nonadjacent
to both $b_1$ and $b_2$, a contradiction since it is a branch-vertex. This proves (2).

\bigskip
Henceforth we may assume that $J = K_4$. To simplify notation, let $V(J) = \{1,2,3,4\}$, and
let the branch of $H$ joining $i$ and $j$ be $B_{i,j} = B_{j,i}$; let its end-edges be $r_{i,j}$
(incident with $i$) and $r_{j,i}$ (incident with $j$); and let $R_{i,j}= R_{j,i} =
L(B_{i,j})$ (that is, the path in $G$ between $r_{i,j}$ and $r_{j,i}$ formed by the vertices
of $G$ in $E(B_{i,j})$). For $k = 1,2,3,4$, we denote by $C_k$ the hole in $G$ induced on the
union of the three paths $R_{i,j}$ ($i,j \not = k$), and we denote by $T_k$ the triangle
$\{r_{k,i}: 1 \le i \le 4, i \not = k\}$. We observe first that
\\
\\
(3) {\em If $r_{3,1},r_{3,2} \not \in X$, then every vertex of $C_4$ in $X$ is adjacent to one
of $r_{3,1},r_{3,2}$.}
\\
\\
For this is trivial if $C_4$ has length 4, so assume it has length $\ge 6$.
Suppose that some vertex $x$ of $C_4$ nonadjacent to $r_{3,1},r_{3,2}$ belongs to $X$.
Now by (1), $r_{3,4} \in X$, and we may assume $r_{3,i} \in X_i$  for
$i = 1,2$. So $r_{3,1}\d y_1\d Q\d y_2\d r_{3,2}$ is an antipath, of length $\ge 3$. It is even,
since it can be completed to an antihole via $r_{3,2}\d x\d r_{3,1}$; and this contradicts
\ref{hole&antipath}. This proves (3).
\\
\\
(4)  {\em If one of $C_1,\ldots,C_4$ contains no member of $X$ then the theorem holds.}
\\
\\
For suppose that no vertex of $C_4$ say is in $X$. Then from (1),
$r_{1,4},r_{2,4},r_{3,4}$ are all in $X$. Now one of $R_{1,2},R_{1,3},R_{2,3}$ has odd length,
since $C_4$ has even length, say $R_{1,2}$. If $r_{1,4},r_{2,4}$ are not adjacent then
$r_{1,4}\d r_{1,2}\d R_{1,2}\d r_{2,1}\d r_{2,4}$
is an odd path between vertices in $X$, and none of its internal vertices are in $X$, and
$r_{3,4}$ has no neighbour in its interior, contrary to \ref{greentouch}. So they are adjacent,
that is, $R_{1,4},R_{2,4}$ both have length 0. Since we may assume that $r_{2,1},r_{2,3}$ are in
$X_1 \setminus X_2$ and $X_2 \setminus X_1$ respectively, it follows that $Q$ can be completed
to an antihole via $y_2\d r_{2,3}\d r_{1,4}\d r_{2,1}\d y_1$, and hence $Q$ has even length $\ge 2$.
By \ref{hole&antipath}, since $Q$ is the interior of an even antipath between $r_{1,3}$
and $r_{1,2}$, and $r_{2,1},r_{2,3}$ are both different from $r_{1,3},r_{1,2}$, and one of
$r_{2,1}$ and $r_{2,3}$ is in $X_1 \setminus X_2$ and the other in $X_2 \setminus X_1$, it follows
that $C_4$ has length 4, and therefore $R_{1,3},R_{2,3}$ have length 0, and $R_{1,2}$ has
length 1. Also, $r_{1,3} \in X_1 \setminus X_2$ by (1), and hence $r_{1,2} \in X_2 \setminus X_1$.
But then $\overline{G} | (T_1 \cup T_2 \cup T_3 \cup Y)$ is an appearance in
$\overline {G}$ of $K_4$; and since $Q$ has even length $\ge 2$,  it is overshadowed in
$\overline{G}$ because of the vertex $r_{4,3}$, and so statement 1 of the theorem holds.
This proves (4).
\\
\\
(5) {\em If one of $C_1\l C_4$ contains at most one member of $X$ then the theorem holds.}
\\
\\
By (4) we may assume at least one vertex of each $C_k$ is in $X$. We may therefore assume
that exactly one vertex (say $x$) of $C_4$ say
is in $X$; and so from the symmetry we may assume that $x$ is in $R_{1,2}$. So neither
of $r_{3,1},r_{3,2}$ is in $X$, and so from \ref{hypothesisb}, we may assume that
$Y$ cannot be linked onto
$T_3$. Consequently $x = r_{2,1}$ (for otherwise $r_{2,4} \in X$ and we could use the
paths $r_{3,4}$, $r_{2,4}\d r_{2,3}\d R_{2,3}\d r_{3,2}$, and a subpath of $C_4\setminus r_{3,2}$
from $x$ to
$r_{3,1}$ to link $Y$ onto $T_3$), and similarly $x = r_{1,2}$. So $R_{1,2}$ has length 0.
Also, from (3), $x$ is adjacent to one of $r_{3,1},r_{3,2}$, so we may assume that
$R_{2,3}$ has length 0 . Therefore $R_{1,3}$ is odd. Hence
$r_{1,2}\d r_{1,3}\d R_{1,3}\d r_{3,1}\d r_{3,4}$ is an odd path between
vertices in $X$, and its internal vertices are not in $X$. By \ref{greentouch}, every
vertex in $X$ has a neighbour in $R_{1,3}$. Hence $r_{3,4}$ is the only vertex of $C_1$
in $X$; and so (by exchanging $C_1$ and $C_4$) we deduce that $R_{3,4}$ has length 0.
So $R_{1,4}$ is even, and $R_{2,4}$ is odd, and the latter is the interior of a second odd
path between members of $X$.  We claim that $R_{1,4}$ has length 0; for if
$X = \{r_{1,2},r_{3,4}\}$, this follows from the symmetry between $R_{1,4}$ and $R_{2,3}$,
while if $|X| > 2$ then since every vertex in $X$ has a neighbour in both $R_{1,3}$ and
$R_{2,4}$, the third vertex in $X$ must be both $r_{1,4}$ and $r_{4,1}$, and hence
again $R_{1,4}$ has length 0. So $R_{1,2},R_{2,3},R_{3,4},R_{1,4}$ all have length 0, and
so $L(H)$ is degenerate, and $X = \{r_{1,2},r_{3,4}\}$ or $ \{r_{1,2},r_{3,4},r_{1,4}\}$.

In particular there is symmetry between $R_{1,3}$ and $R_{4,2}$. If both these paths have length
1 then statement 3 of the theorem holds, so from the symmetry we may assume that $R_{1,3}$ has length $>1$.
Now we recall that there is an antipath $Q$ in $Y$ between $y_1$ and $y_2$, with $V(Q) = Y$.
Let its vertices be $y_1 = q_1,q_2,\ldots,q_k = y_2$ in order. We may assume that
$r_{3,1} \in X_1$ and $r_{3,2} \in X_2$; and consequently $r_{2,4} \in X_1$.
Now $R_{1,3}$ has length $>1$ and is odd, so
$r_{1,2}\d r_{1,3}\d R_{1,3}\d r_{3,1}\d r_{3,4}$ is an odd path of length $\ge 5$,
with ends in $X$ and no internal vertex in $X$. By \ref{RR}, $Y$ contains a leap.
But the only vertex in $Y$ nonadjacent to $r_{3,1}$ is $y_1 = q_1$, and its only nonneighbour
in $Y$ is $q_2$; so $q_1,q_2$ is the leap. Hence $q_2$ is adjacent to $r_{1,4}$, since it
has two neighbours in $T_1$. Now the path $r_{1,2}\d r_{2,4}\d R_{2,4}\d r_{4,2}\d r_{4,3}$ is odd
and has length $\ge 3$, between common neighbours of $\{q_1,q_2\}$; and these two vertices
have no common neighbour in its interior, since we could complete the path
$q_1\d r_{1,3}\d R_{1,3}\d r_{3,1}\d q_2$ to an odd hole through any such common neighbour. So
$\{q_1,q_2\}$ is also a leap for this second path. But then the only neighbours of $q_1$
in $L(H)$ are $r_{1,2},r_{1,3},r_{4,2},r_{4,3},r_{2,3}$ and possibly $r_{1,4}$; and the
only neighbours of $q_2$ are $r_{2,1},r_{2,4},r_{3,1},r_{3,4},r_{1,4}$ and possibly
$r_{2,3}$. But then statement 2 of the theorem holds. This proves (5).
\\
\\
(6) {\em If $r_{3,1},r_{3,2} \not \in X$ then either the theorem holds, or $C_4$ has
length $4$ and its other two vertices are in $X$.}
\\
\\
For  by (5), we may assume that at least two vertices of $C_4$ are in $X$. Since $Y$
cannot be linked onto $T_3$, there are exactly two such vertices and they are adjacent. By
(3) they are both adjacent to one of $r_{3,1},r_{3,2}$. This proves (6).
\\
\\
(7) {\em If $r_{3,1},r_{3,2} \not \in X$ and $R_{1,2}$ has length $>0$ then the theorem holds.}
\\
\\
For then by (6) we may assume that $R_{1,2}$ has length 1, and $R_{1,3},R_{2,3}$ both have
length 0, and $r_{1,2},r_{2,1}$ are both in $X$. Since $Y$ cannot be linked onto $T_3$ it
follows that $r_{1,4} \not \in X$, and similarly $r_{2,4} \not \in X$. By (1) we may assume that
$r_{1,3} \in X_1 \setminus X_2$, and $r_{1,4} \in X_2 \setminus X_1$. But then
$r_{2,3} \in X_2 \setminus X_1$, and $r_{2,4} \in X_1 \setminus X_2$.
Hence $Q$ can be completed to an antihole via $y_2\d r_{1,3}\d r_{2,1}\d r_{1,4}\d y_2$, and
so $Q$ is even; and it therefore cannot be completed via $y_2\d r_{2,4}\d r_{1,4}\d y_1$,
and so $r_{1,4},r_{2,4}$ are adjacent. So $R_{1,4},R_{2,4}$ both have length 0, and therefore
$R_{3,4}$ has odd length.  Since $r_{1,3}$ and $r_{1,4}$ are not in $X$, it follows from (6) applied to
$C_2$ that we may assume $C_2$ has length 4, and hence
$R_{3,4}$ has length 1; but then statement 3 of the theorem holds. This proves (7).

\bigskip
Now, to complete the proof: we may assume that $r_{3,1},r_{3,2} \not \in X$.
By (6) and (7) we may assume that $R_{1,2}$ and $R_{2,3}$ have length 0, and $R_{1,3}$
has length 1, and $r_{1,3},r_{1,2} \in X$. Since $Y$ cannot be linked onto $T_3$ it follows that
$r_{2,4} \not \in X$; so by (6) and (7) applied to $C_1$, we may assume that $R_{3,4}$ has
length 0, $R_{2,4}$ has length 1, and $r_{4,2} \in X$. If $R_{1,4}$ has length 0 then statement 3
of the theorem holds, and otherwise $L(H)$ is overshadowed (since any
vertex in $Y$ is adjacent to all of $r_{1,2},r_{1,3},r_{4,2},r_{4,3}$), and so statement 1
of the theorem holds.  This proves \ref{nonsat}. \bbox

\section {Rung replacement}

Before we apply \ref{nonsat}, let us simplify it a little. We can effectively eliminate the
two cases of $L(H)$ being overshadowed. We need a few lemmas.

\begin{thm}\label{getprism}
Let $c_1,c_2$ be adjacent vertices of a $3$-connected graph $J$, and let $e,f$ be edges
of $J$ incident with $c_1$ and different from $c_1c_2$. There are three tracks of $J$ from
$c_1$ to $c_2$, pairwise vertex-disjoint except for their ends, and with first edges
$c_1c_2, e,f$ respectively.
\end{thm}
\Proof
Since $J$ is 3-connected, if we delete from $J$ all edges incident with $c_1$ except $e$ and $f$,
the graph we make is still 2-connected, and so it has a cycle containing $c_1$ and $c_2$.
This proves \ref{getprism}.\bbox

Let $P_1,P_2,P_3$ be paths of $G$. We say they {\em form} a prism $K$ if they are pairwise disjoint, each has
length $\ge 1$, and for each $i$ the ends of $P_i$ can be labelled $a_i,b_i$ so that for $1 \le i < j \le 3$
the only edges of $G$ between $V(P_i)$ and $V(P_j)$ are $a_ia_j$ and $b_ib_j$; and $K$ is the prism in $G$
induced on the union of these three paths.

\begin{thm}\label{prismparity}
Let $R_1,R_2,R_3$ be a prism in a Berge graph $G$; then $R_1,R_2,R_3$ all have the same parity.
\end{thm}
The proof is clear.

\begin{thm}\label{prismbig}
Let $G$ be Berge, let $Y \subseteq V(G)$ be anticonnected, and for $i = 1,2,3$ let $a_i\d P_i \d b_i$ be a
path in $G\setminus Y$, forming a prism with triangles $\{a_1,a_2,a_3\},\{b_1,b_2,b_3\}$.
Assume $P_1,P_2,P_3$ all have length $>1$, and that every vertex in $Y$
is adjacent to at least two of $\{a_1,a_2,a_3\}$ and to at least two of $\{b_1,b_2,b_3\}$.
Then at least two of $\{a_1,a_2,a_3\}$ and at least two of $\{b_1,b_2,b_3\}$ are $Y$-complete.
\end{thm}
\Proof
Suppose not; then there is an antipath with interior in $Y$, joining two vertices either both in
$\{a_1,a_2,a_3\}$ or both in $\{b_1,b_2,b_3\}$. Let $Q$ be the shortest such antipath. We may
assume $Q$ joins $a_1$ and $a_2$ say. Since every vertex in $Y$ is adjacent to either $a_1$ or
$a_2$ it follows that $Q$ has length $\ge 3$. From the minimality of $Q$, $a_3$ is $Q^*$-complete, and
so is at least one of $b_1,b_2,b_3$, say $b_i$. Since $Q$ can be completed to an antihole via
$a_1\d b_i\d a_2$ it follows that $Q$ is even. From \ref{hole&antipath} applied to the hole
formed by $P_1 \cup P_2$ and hat $a_3$, neither of $b_1,b_2$ is $Q^*$-complete, and so there is an
antipath between $b_1$ and $b_2$ with interior in $Q^*$. By the minimality of $Q$, the two
antipaths have the same interior; but this again contradicts \ref{hole&antipath}. This proves
\ref{prismbig}.\bbox

In fact it is easy to find strengthenings of \ref{prismbig} in which some of the paths $P_i$
have length 1, but for the moment \ref{prismbig} will suffice.

\begin{thm}\label{prismrung}
Let $G$ be Berge, and for $1 \le i \le 3$ let $P_i$ be a path of even length $\ge 2$, from $a_i$
to $b_i$, so that these three paths form a prism with triangles $A = \{a_1,a_2,a_3\}$ and
$B = \{b_1,b_2,b_3\}$. Let $P_1'$ be a path from $a_1'$ to $b_1$, so that $P_1',P_2,P_3$ also
form a prism. Let $y \in V(G)$ have at least two neighbours in $A$ and in $B$. Then it also
has at least two neighbours in $\{a_1',a_2,a_3\}$.
\end{thm}
\Proof
Suppose not. By \ref{prismparity} $P_1'$ has even length.
Let $X$ be the set of neighbours of $Y$ in $G$. Then $a_1' \not \in X$, and
$a_1 \in X$, and exactly one of $a_2,a_3 \in X$, say $a_2 \in X$. Also, $y$ cannot be linked
onto the triangle $A' = \{a_1',a_2,a_3\}$, by \ref{trianglev}, and since one of $b_2, b_3 \in X$
it follows that no internal vertex of $P_1'$ is in $X$. Hence $b_1 \not \in X$, for otherwise
$y \d a_2\d a_1'\d P_1'\d b_1 \d$ would be an odd hole.
So $b_2,b_3 \in X$. Since $y \d a_1\d a_3\d P_3\d b_3\d y$ is not an odd hole,
there is a member of $X$ in $P_3 \setminus b_3$. But then $y$ can be
linked onto $A'$, via $b_2\d b_1\d P_1'\d a_1$, the path $a_2$, and the path between $y$ and $a_3$ with interior
in $V(P_3)\setminus\{b_3\}$, contrary to \ref{trianglev}. This proves \ref{prismrung}. \bbox

We shall only need the following when $J = K_4$ or $K_{3,3}$, but we might as well prove it
in general.

\begin{thm}\label{overshadowed}
Let $G$ be Berge, and let $L(H)$ be an overshadowed appearance of $J$ in $G$, where $J$ is 3-connected.
Then either:
\begin{itemize}
\item there is a $J$-enlargement with a nondegenerate appearance in $G$, or
\item $G$ admits a balanced skew partition.
\end{itemize}
\end{thm}
\Proof For each edge $uv$ of $J$, let $B_{uv}$ be the branch of $H$
with ends $u,v$, and let $R_{uv}$ be the path $L(B_{uv})$ of $L(H)$. For each $v \in V(J)$
let $N_v$ be the clique of $L(H)$ with vertex set $\delta_H(v)$.
There is an edge $c_1c_2$ of $J$ so that $B_{c_1c_2}$ has odd length $\ge 3$, and some
vertex of $G$ is nonadjacent in $G$ to at most one vertex of $N_{c_1}$ and to at most one vertex
of $N_{c_2}$.  We say such a vertex $v$ is {\em $B_{c_1c_2}$-dominant with respect to L(H)}.
Let the ends of $R_{c_1c_2}$
(that is, the end-edges of $B_{c_1c_2}$) be $r_1,r_2$, where $r_i \in N_{c_i}$. Let
$Y$ be a maximal anticonnected set of vertices each with at most one non-neighbour in
$N_{c_1}$ and at most one non-neighbour in $N_{c_2}$. We shall prove that $Y$ and some of its
common neighbours separate the interior of $R_{c_1c_2}$ from the remainder of $L(H)$ in $G$, so
that will be the skew partition we are looking for. Let $X$ be the set of all $Y$-complete vertices
in $G$.
\\
\\
(1) {\em For $i = 1,2$, at most one vertex of $N_{c_i}$ is not in $X$.}
\\
\\
For let $a_1,a_2$ be any two distinct vertices in $N_{c_1}\setminus r_1$; we shall show that at
most one of $a_1,a_2, r_1$ is not in $X$. By \ref{getprism}, there are two paths $Q_1$,$Q_2$ of
$H$ between $c_1$ and $c_2$, so that $Q_1$,$Q_2$,$B_{c_1c_2}$ are vertex-disjoint except for their ends,
and for $i = 1,2$, $a_i$ is the first edge of $Q_i$. Let $b_i$ be the other end-edge of $Q_i$.
Both $Q_1$ and $Q_2$ have odd length, since $B_{c_1c_2}$ is odd and $H$ is bipartite; and they have
length $\ge 3$ since $b_1,b_2$ are nonadjacent (for they are the ends of a branch of length $>1$.)
Hence there are two paths $P_1$,$P_2$ of $L(H)$ from $N_{c_1}$ to $N_{c_2}$, so that
$P_1$,$P_2$,$R_{c_1c_2}$ are vertex-disjoint and form a prism, and $P_i$ is from $a_i$ to $b_i$.
Now $B_{c_1c_2}$ is odd and therefore $R_{c_1c_2}$ is even, and similarly $P_1$ and $P_2$ are even.
By hypothesis, each member of $Y$ is adjacent to at least two vertices of the triangle
$\{a_1,a_2,r_1\}$ and to two vertices of the triangle $\{b_1,b_2,r_2\}$. By \ref{prismbig}
it follows that $X$ contains at least two members of $\{a_1,a_2,r_1\}$. This proves (1).

\bigskip

Let
\begin{eqnarray*}
X_1 &=& X \cap (N_{c_1}\cup N_{c_2}) \\
X_2 &=& X \cap (V(L(H))\setminus (N_{c_1}\cup N_{c_2}))\\
X_0 &=& X \setminus V(L(H))\\
S &=& V(R_{c_1c_2})\setminus X_1\\
T &=& (V(L(H))\setminus V(R_{c_1c_2}))\setminus X_1.
\end{eqnarray*}

We observe first that
no vertex of $S$ is adjacent to any vertex in $T$; for such an edge would join two vertices both
in $N_{c_i}$ for some $i$, and therefore both not in $X$, contradicting (1).
\\
\\
(2) {\em If $F \subseteq V(G)$
is connected and some vertex of $S$ has a neighbour in $F$, and so does some vertex of $T$,
and $F\cap (X_0 \cup X_1 \cup Y) = \emptyset$, then the theorem holds.}
\\
\\
We shall prove this by induction on $|F|$; so, we assume it holds for all smaller choices of $F$
(even for different choices of $L(H)$). Hence
we may assume that $G|F$ is a path with vertices $f_1,\ldots,f_n$ say, where $f_1$
is the only vertex of $F$ with a neighbour in $S$, and $f_n$ is the only vertex with a
neighbour in $T$. From the minimality of $F$ it also follows that $F$ is disjoint from
$L(H)$; for any vertex of $F$ in $L(H)$ would be in $S$ or $T$, since it is not in $X_1$, and
then we could make $F$ shorter by omitting this vertex. Consequently $F \cap X = \emptyset$.
Suppose some vertex in $v \in F$ is
major with respect to $L(H)$. Then since $v \not \in X$ it follows that $v$ has a nonneighbour in $Y$,
and so $Y \cup v$ is anticonnected; the maximality of $Y$ therefore implies that $v \in Y$, and hence
$F \cap Y \not = \emptyset$ and the claim holds. So we may assume that no vertex
in $F$ is major. On the other hand, the set of attachments of $F$ in $L(H)$ is not local, because
it has an attachment in $R_{c_1c_2}$, and its attachments are not all contained in any of
$V(R_{c_1c_2})$, $N_{b_1}$,$N_{b_2}$. Let us apply \ref{smallcomp}. Suppose first that
\ref{smallcomp}.1 holds. Then we obtain an appearance $L(H')$ in $G$ of some $J$-enlargement, with
$L(H)$ an induced subgraph of $L(H')$. Since $R_{c_1c_2}$ has even nonzero length, it follows
that $L(H)$ is not degenerate, and therefore neither is $L(H')$, and hence the theorem holds.
So we may assume that \ref{smallcomp}.2 holds, and there is an edge $b_1b_2$ of $J$,
(for $i = 1,2$, $s_i$ denotes the unique vertex in $N_{b_i}\cap R_{b_1b_2}$)
and a path $P$ of $G$ with $V(P) \subseteq F$ and with ends
$p_1$ and $p_2$, such that one of the following holds:
\begin{enumerate}
\item $p_1$ is adjacent in $G$ to all vertices in $N_{b_1} \setminus s_1$, and $p_2$ has a neighbour in
$R_{b_1b_2} \setminus s_1$, and every edge from $V(P)$ to
$V(L(H))\setminus s_1$ is either from $p_1$ to $N_{b_1} \setminus s_1$, or
from $p_2$ to  $R_{b_1b_2} \setminus s_1$, or
\item for $i = 1,2$, $p_i$ is adjacent in $G$ to all vertices in $N_{b_i} \setminus s_i$, and there are
no other edges between $V(P)$ and $V(L(H))$ except possibly $p_1s_1,p_2s_2$,
and $P$ has the same parity as $R_{b_1b_2}$, or
\item $p_1 = p_2$, and $p_1$ is adjacent to all vertices in $(N_{b_1} \cup N_{b_2})\setminus\{s_1, s_2\}$,
and all neighbours of $p_1$ in $V(L(H))$ belong to $N_{b_1} \cup N_{b_2} \cup R_{b_1b_2}$,
and $R_{b_1b_2}$ is even, or
\item $s_1 = s_2$, and  for $i = 1,2$, $p_i$ is adjacent in $G$ to all vertices in $N_{b_i} \setminus s_i$,
and there are
no other edges between $V(P)$ and $V(L(H))\setminus s_1$, and $P$ is even.
\end{enumerate}
In case 1, let $R'$ be the (unique) path from $p_1$ to $s_2$ in $(V(P) \cup V(R_{b_1b_2})) \setminus s_1$,
and in the other cases
let $R'$ be $P$. So if in $L(H)$ we replace $R_{b_1b_2}$ by $R'$ we obtain another
appearance of $J$ in $G$, say $L(H')$, where $H'$ is obtained from $H$ by replacing the branch
$B_{b_1b_2}$ by some new branch $B'$ joining the same two vertices. For each $v \in V(J)$ let
$N'_v$ be the clique in $L(H')$ formed by the edges in $\delta_{H'}(v)$. So $N'_v = N_v$ for all
vertices $v$ of $J$ except for $b_1$ and $b_2$.  Let $R'$ be between $r_1'$
and $r_2'$, where $r_i' \in N_{b_i}'$ for $i = 1,2$.

Now suppose that $b_1b_2$ and $c_1c_2$ are different edges of $J$. Then $B_{c_1c_2}$ is still a
branch of $H'$, and we claim
that every $y \in Y$ is $B_{c_1c_2}$-dominant with respect to $L(H')$. For let $e,f$ be two edges of $J$
incident with $c_1$ and different from $c_1c_2$. By \ref{getprism} there are three tracks of $J$
from $c_1$ to $c_2$, vertex-disjoint except for their ends, and one of them is the edge $c_1c_2$,
and the first edges of the other two are $e$ and $f$. There are three tracks corresponding to
these in $H$, and their line graph is a prism in $L(H)$. There also correspond three tracks in $H'$,
yielding a prism in $L(H')$. Since  $R_{b_1b_2} \not = R_{c_1c_2}$, it follows that
$R_{b_1b_2}$ is incident with at most one of $c_1,c_2$, so these two prisms are related as in
\ref{prismrung}. Hence by \ref{prismrung}, since $y$ has two neighbours in both triangles
of the first prism, it also has two neighbours in the triangles of the second. This proves that
$y$ is $B_{c_1c_2}$-dominant with respect to $L(H')$. The same argument in the reverse direction shows that
$Y$ remains a maximal anticonnected set of $B_{c_1c_2}$-dominant vertices.
Since there is a proper subset $F'$ of $F$ with attachments in $S$ and in the new
set $T'$ in $V(H')$ corresponding to $T$ (for $T'$ contains all the vertices of $R'$
that are in $F$, and there is at least one such vertex), it follows that we may apply the inductive
hypothesis. So $F'$, and hence $F$, contains a vertex of $X$. This completes the argument
when $b_1b_2$ and $c_1c_2$ are distinct edges.

Now we assume that $b_i = c_i$ for $i = 1,2$.
There were four cases in the definition of $P$, listed above. Case 3 is impossible, since then the
vertex $p_1$ would be $B_{c_1c_2}$-dominant with respect to $L(H)$, and therefore would be in either
$X$ or $Y$, a contradiction. Also, case 1 is impossible, by applying \ref{prismrung} as before
to show that $Y$ remains a maximal anticonnected set of $B'$-dominant vertices, and applying the
inductive hypothesis. Case 4 is impossible since $B_{c_1c_2}$ has length $\ge 3$. So case 2 applies;
that is, $p_2$ is adjacent to all vertices in
$N_{c_2} \setminus r_2$, and to no vertex of $R_{c_1c_2}$ except possibly $r_2$.  So
$N'_{c_i} = (N_{c_i} \setminus \{r_i\}) \cup \{r_i'\}$ for $i = 1,2$. We recall that in this
case $R' = P$, and $P$ is a subpath of the path with vertices $f_1,\ldots,f_n$. Choose $h$ with
$1 \le h \le n$ minimum so that $f_h$ is a vertex of $R'$. Since both $R'$ and $G|F$ are paths
it follows that $f_h$ is one end of $R'$, say $r_1'$. (This is without loss of generality,
because in this case 2, there is symmetry between $b_1 = c_1$ and $b_2 = c_2$.) From the minimality
of $F$, $r_1'$ has no neighbour in $T$, and in particular every vertex in $N_{c_1}\setminus r_1$ is in $X$.
We claim also that every vertex of $N_{c_2}\setminus r_2$
is in $X$. For if not, then $r_2 \in X$, and by \ref{getprism} there is a prism $R_{c_1c_2},P_1,P_2$
say, in $L(H)$, where each $P_i$ has an end $a_i \in N_{c_1}$ and an end $b_i \in N_{c_2}$, and
$b_2 \not \in X$. (Consequently $r_2,b_1 \in X$.) Hence at most one vertex of the triangle
$\{r_2',b_1,b_2\}$ is in $X$, and some
vertex in $X$ (namely $a_1$) has no neighbour in this triangle, so by \ref{triangleX}, $Y$
cannot be linked onto this triangle. In particular, no vertex of $P_2$ is in $X$ except $a_2$.
But then $a_2\d P_2\d b_2\d r_2$ is an odd path between members of $X$, and none of its internal
vertices are in $X$, and $a_1$ has no neighbour in its interior, contrary to \ref{greentouch}.
This proves that every vertex of $N_{c_2}\setminus r_2$ is in $X$. Consequently all vertices of $Y$
are $B'$-dominant with respect to $L(H')$. We claim also that $Y$ is still maximal. For suppose
not, and let $Y \subset Y'$ for some larger anticonnected set $Y'$ of $B'$-dominant vertices.
Since $r_1',r_2'$ are not in $X$, they are certainly not $Y'$-complete, and since by (1) applied
to $Y'$, at most one vertex of $N_{c_i}'$ is not $Y'$-complete for $i = 1,2$, it follows that
every vertex of $N_{c_1}'\setminus r_1'$ and $N_{c_2}'\setminus r_2'$ are $Y'$-complete. But then all
the members of $Y'$ are $B_{c_1c_2}$-dominant with respect to $L(H)$, contrary to the maximality of $Y$.
This proves that $Y$ is a maximal anticonnected set of $B'$-dominant vertices with respect to $L(H')$.
Hence we can apply induction on $F$, and the result follows. This proves (2).

\bigskip

It follows from (2) that there is a partition of $V(G) \setminus (X_0 \cup X_1 \cup Y)$ into two sets
$L$ and $M$ say, where there is no edge between $L$ and $M$, and $S \subseteq L$ and $T \subseteq M$.
So $(L \cup M, X_0 \cup X_1 \cup Y)$ is a skew partition of $G$. By \ref{geteven} we may assume
it is not loose, and so $X_2$ is empty; and we shall show it is balanced.
\\
\\
(3) {\em For $i = 1,2$, all vertices of $N_{c_i} \setminus r_i$ belong to $X_1$.}
\\
\\
For suppose there is a vertex $n_1$ of $N_{c_1} \setminus r_1$ not in $X$. Therefore all other
vertices of $N_{c_1}$ belong to $X$, and in particular, $r_1 \in X$. Suppose no other vertex of
$R_{c_1c_2}$ is in $X$; then $r_2 \not \in X$, so $X$ includes $N_{c_2} \setminus r_2$.
Choose any $n_2 \in N_{c_2} \setminus r_2$, and any $n_1' \in N_{c_1} \setminus r_1$ different
from $n_1$. Then $r_1\d R_{c_1c_2}\d r_2\d n_2$ is an odd path between $Y$-complete vertices, and none
of its internal vertices are $Y$-complete, and yet $n_1'$ has no neighbour in its interior, contrary
to \ref{greentouch}. This proves that some vertex of $R_{c_1c_2}$ different from $r_1$ is in $X$;
yet $X_2$ is empty, so the interior of $R_{c_1c_2}$ contains no vertex in $X$. Consequently
$r_2 \in X$. Choose $n_2 \in N_{c_2} \setminus r_2$ so that $N_{c_2} \setminus n_2 \subseteq X$.
Since $J$ is 3-connected, there is a track of $H$ from $c_1$ to $c_2$ with first edge
$n_1$ and last edge different from $n_2$. This track is odd since $c_1$ and $c_2$ have opposite
biparity; and so in $G$ there is an even path, $P$ say, from $n_1$ to some
$n_2' \in N_{c_2} \setminus n_2$, with no vertex in $N_{c_1} \cup N_{c_2}$ except its ends. But then
$r_1\d n_1\d P\d n_2'$ is an odd path between $Y$-complete vertices, no vertex in its interior is
$Y$-complete, and the $Y$-complete vertex $r_2$ has no neighbour in its interior, contrary
to \ref{greentouch}. This proves (3).

\bigskip

Let $W = (N_{c_1} \setminus r_1) \cup (N_{c_2} \setminus r_2)$. Then $W \subseteq X_1$ by (3),
and since there are no edges between $N_{c_1}$ and $N_{c_2}$, it follows that $W$ has
exactly two components, both cliques. In particular, $W$ is anticonnected. Now every $W$-complete
vertex is $B_{c_1c_2}$-dominant, and so belongs to $X \cup Y$; and hence there are no $W$-complete
vertices in $L \cup M$. Consequently $W$ is a kernel for the skew partition. Suppose
$u_1,u_2 \in W$ are nonadjacent and joined by an odd path with interior in $L$. Then
one of $u_1,u_2$ is in $N_{c_1} \setminus r_1$ and the other in $N_{c_2} \setminus r_2$, and
therefore they are joined by a path in $L(H)$ using no more vertices in $N_{c_1} \cup N_{c_2}$,
which is even (since $H$ is bipartite). But this contradicts \ref{mixedpair}, and so there are no
such $u_1,u_2$. Finally, suppose there is a pair of vertices of $L$ joined by an odd antipath with
interior in $W$, necessarily of length $\ge 5$ (since we already did the odd path case). Then
$G|W$ contains an antipath of length 3, which is impossible since its components are cliques.
From \ref{kernel} it follows that the skew partition is balanced. This proves \ref{overshadowed}. \bbox

\section{Generalized line graphs}

As we said earlier, our strategy is to find the biggest line graph in $G$ that we can, and then
assemble all the alternative rungs for a given edge of $J$ into a ``strip''. In this section
we make that precise.

Let $J$ be 3-connected, and let $G$ be Berge. A {\em $J$-strip system $(S,N)$} in $G$ means:
\begin{itemize}
\item for each edge $uv$ of $J$, a subset $S_{uv}=S_{vu} \subseteq V(G)$
\item for each vertex $v$ of $J$, a subset $N_v \subseteq V(G)$
\end{itemize}
satisfying the following conditions (for $uv \in E(J)$, a {\em $uv$-rung } means a path $R$ of
$G$ with ends $s,t$ say, where $V(R) \subseteq S_{uv}$, and $s$ is the unique vertex of
$R$ in $N_u$, and $t$ is the unique vertex of $R$ in $N_v$):
\begin{itemize}
\item The sets $S_{uv} (uv \in E(J))$ are pairwise disjoint
\item For each $u \in V(J)$, $N_u \subseteq \bigcup (S_{uv}: v \in V(J)$ adjacent to $u)$
\item For each $uv \in E(J)$, every vertex of $S_{uv}$ is in a $uv$-rung
\item If $uv,wx \in E(J)$ with $u,v,w,x$ all distinct, then there are
no edges between $S_{uv}$ and $S_{wx}$
\item If $uv,uw \in E(J)$ with $v \not = w$, then $N_u \cap S_{uv}$ is complete to
$N_u \cap S_{uw}$, and there are no other edges between $S_{uv}$ and $S_{uw}$
\item For each $uv \in E(J)$ there is a $uv$-rung so that for every cycle $C$ of $J$, the
sum of the lengths of the $uv$-rungs for $uv \in E(C)$ has the same parity as $|V(C)|$.
\end{itemize}

It follows that for distinct $u,v \in V(J)$, $N_u \cap N_v$ is empty if $u$,$v$ are nonadjacent, and
otherwise $N_u \cap N_v \subseteq S_{uv}$;
and for $uv \in E(J)$ and $w \in V(J)$, if $w \not = u,v$ then $S_{uv} \cap N_w = \emptyset$.
The final axiom looks strange, but we shall show immediately that the same property holds
for {\em every} choice of $uv$-rungs.

\begin{thm}\label{runglength}
If $(S,N)$ is a $J$-strip system in a Berge graph $G$, where $J$ is 3-connected.
Then for every $uv \in E(J)$, all $uv$-rungs have lengths of the same parity.
\end{thm}
\Proof
Since $J$ is 3-connected, there is a cycle $C$ of $J$ with $|V(C)| \ge 4$ and with
$uv \in E(C)$. For each $xy \in E(C)$ different from $uv$, choose an $xy$-rung $R_{xy}$.
For every $uv$-rung $R$, the union of $V(R)$ and all the $V(R_{xy})$'s induces a cycle in $G$. This
has length $\ge4$ since $C$ has length $\ge 4$, so it is a hole and therefore even.
Hence all choices of $R$ have lengths of the same parity. This proves \ref{runglength}. \bbox

For each edge $uv$ of $J$, choose a $uv$-rung $R_{uv}$.
It follows from \ref{runglength} and the final axiom above that
the subgraph of $G$ induced on the union of the vertex sets of these rungs is a line graph of
a bipartite subdivision $H$ of $J$.  For brevity we say that this choice of rungs {\em forms} $L(H)$.

We need two easy observations:

\begin{thm}\label{mixedrungs}
Let $(S,N)$ be a $J$-strip system in a Berge graph $G$, where $J$ is 3-connected.
If there is an edge $uv$ of $J$ so that some $uv$-rung has length 0 and another $uv$-rung has
length $\ge 1$, then there is an overshadowed appearance of $J$ in $G$.
\end{thm}
\Proof For each edge $ij$ of $J$ choose an $ij$-rung $R_{ij}$, so that $R_{uv}$ has length
$\ge 1$ and otherwise arbitrarily; and let this choice of rungs form $L(H)$. Let $y$ be the
vertex of some $uv$-rung of length 0.  By \ref{runglength}, $R_{uv}$ has even length. Let
$B$ be the branch of $H$ between $u$ and $v$, so $E(B) = V(R_{uv})$. Then $B$ is odd and has
length $\ge 3$ and $y$ is nonadjacent in $G$ to at most one vertex of $G$ in $\delta_H(u)$
and at most one in $\delta_H(v)$. Hence $L(H)$ is overshadowed.  This proves \ref{mixedrungs}.\bbox

A $J$-strip system is {\em nondegenerate}
if there is some choice of rungs so that the line graph $L(H)$ they form is a nondegenerate
appearance of $J$.  \ref{mixedrungs} has the following corollary:

\begin{thm}\label{mixedrungs2}
Let $(S,N)$ be a nondegenerate $J$-strip system in a Berge graph $G$, where $J$ is 3-connected.
If there is no overshadowed appearance of $J$ in $G$, then
for every choice of rungs, the line graph they form is a nondegenerate appearance of $J$ in $G$.
\end{thm}
\Proof  Take a choice of rungs $R_{ij} (ij \in E(J))$, forming $L(H)$ say, where $L(H)$ is nondegenerate;
and suppose there is another choice, $R_{ij}' (ij \in E(J))$, forming $L(H')$ say, where
$L(H')$ is degenerate. Then for some $ij \in E(J)$, $R_{ij}$ has nonzero length and $R_{ij}'$ has length 0.
By \ref{mixedrungs} there is an overshadowed appearance of $J$ in $G$. This proves \ref{mixedrungs2}.\bbox

Given a $J$-strip system $(S,N)$, we define $V(S,N) = \bigcup(S_{uv}: uv \in E(J))$. Hence
every $N_v \subseteq V(S,N)$. If $u,v \in V(J)$ are adjacent, we define $N_{uv} = N_u \cap S_{uv}$.
So every vertex of $N_u$ belongs to $N_{uv}$ for exactly one $v$. Note that $N_{uv}$ is
in general different from $N_{vu}$, but $S_{uv}$ and $S_{vu}$ mean the same thing.  We say
$X \subseteq V(S,N)$ {\em saturates} the strip system if for every $u \in V(J)$, there is at most
one neighbour $v$ of $u$ in $J$ so that $N_{uv} \not \subseteq X$; and a vertex
$y \in V(G) \setminus V(S,N)$ is {\em major} (with respect to the strip system) if the set of its
neighbours in $V(S,N)$ saturates $(S,N)$. We say $X \subseteq V(S,N)$ is {\em local}
(with respect to the strip system) if either $X \subseteq N_v$ for some $v \in V(J)$, or
$X \subseteq S_{uv}$ for some edge $uv \in E(J)$.

\begin{thm}\label{bigtosmall}
Let $G$ be Berge, and let $J$ be a 3-connected graph. Let $(S,N)$ be a $J$-strip system in
$G$, nondegenerate if $J = K_4$.
Let $y \in V(G) \setminus V(S,N)$, and let $X$ be the set of neighbours of $y$. If there is a
choice of rungs, forming a line graph $L(H)$, so that $X$ saturates $L(H)$, then
either:
\begin{itemize}
\item  $X$ saturates the strip system, or
\item there is a $J$-enlargement with a nondegenerate appearance in $G$, or
\item $J = K_4$ and there is an overshadowed appearance of $J$ in $G$.
\end{itemize}
\end{thm}

\Proof We define the {\em fork number} of a choice of rungs to be the number of branch-vertices of $H$
incident in $H$ with $\ge 2$ edges in $X\cap E(H)$, where $L(H)$ is the line graph formed by this
choice of rungs.
Let us say that a choice of rungs $R_{ij}$ forming a line graph $L(H)$ is {\em saturated} if $X$
saturates $L(H)$,  and in this case its fork number is $|V(J)|$. If every choice of rungs is
saturated, then $X$ saturates the strip system as required, so we may therefore assume that there
is some choice of rungs that is not saturated. Let this choice of rungs form $L(H)$, and let us apply
\ref{vnbrs} to $L(H)$. Hence \ref{vnbrs}.1 is false; suppose that \ref{vnbrs}.6 holds.
Then $G | (V(L(H)) \cup \{y\}) = L(H')$, and $L(H')$ is an appearance in $G$ of a $J$-enlargement.
We may assume that $L(H')$ is degenerate, for otherwise the theorem holds. Hence $J = K_4$ and
$L(H)$ is degenerate. Since the strip system is nondegenerate, the result follows from
\ref{mixedrungs2}. So we may assume that one of \ref{vnbrs}.2-5 holds. Hence
for any nonsaturated choice of rungs, the fork number is $\le 2$. Now
there are two choices of rungs  $R_{ij} (ij \in E(J))$ and $R_{ij}' (ij \in E(J))$, so
that the first is saturated and the second is not, differing only on one edge of $J$; say
$R_{ij} = R_{ij}'$ for all edges $ij$ of $J$ except the edge $1\d 2$.  Since these two choices
of rungs differ only on one edge of $J$, their fork numbers differ by at most 2; and so
$|V(J)| = 4$, and so $J = K_4$.

Let $V(J) = \{1,2,3,4\}$, and $R_{ij} \not = R_{ij}'$ only for the edge $1\d 2$.
As usual, the ends of each $R_{ij}$ will be denoted by $r_{ij}$ and $r_{ji}$, and the ends of
each $R_{ij}'$ will be denoted by $r_{ij}'$ and $r_{ji}'$. For $i = 1,2,3,4$, we denote the
triangle $\{r_{ij}: j \in \{1,\ldots,4\}\setminus i\}$ by $T_i$, and the
triangle  $\{r_{ij}': j \in \{1,\ldots,4\}\setminus i\}$ by $T_i'$. The line graphs made
by $R_{ij}$ and $R'_{ij}$ are $L(H)$ and $L(H')$ respectively. Since $X$ saturates $L(H)$, it
has at least two members in each of $T_1,\ldots,T_4$; and since $X$ does not saturate $L(H')$,
there is some $T_i'$ containing at most one member of $X$. Since $T_3 = T_3'$ and $T_4 = T_4'$, we
may assume that $|X\cap T_1| = 2$ and $|X\cap T_1'| = 1$; and so $r_{1,2} \in X$,
$r_{1,2}' \not \in X$, and exactly one of $r_{1,3},r_{1,4} \in X$, say $r_{1,3} \in X$ and
$r_{1,4} \not \in X$.

Also, at least two vertices of $T_3$ and $T_4$ are in $X$, so there are at least two branch-vertices
of $H'$ incident in $H'$ with more than one edge in $X$. By \ref{vnbrs} applied to $H'$, we deduce
that \ref{vnbrs}.5 holds, and so there is an edge $ij$ of $J$ so that $R_{ij}'$ is even and
\[(X \cap V(L(H')))\setminus V(R_{ij}') = (T_i' \cup T_j') \setminus V(R_{ij}').\]
In particular,
$T_i'$ and $T_j'$ both contain at least two vertices in $X$, and so $i,j \ge 2$. Since
$r_{1,3} \in X$ it follows that one of $i,j = 3$, say $j = 3$, and $r_{1,3} \in T_3$; so
$R_{1,3}$ has length 0. Now there are two cases, $i = 2$ and $i = 4$. Suppose first that
$i = 2$. Then
\[(X \cap V(L(H'))) \setminus V(R_{2,3}) = \{r_{1,3}, r_{3,4}, r_{2,4},r_{2,1}'\},\]
 and since at least
two vertices of $T_4$ are in $X$ it follows that $R_{2,4},R_{3,4}$ both have length 0, a contradiction
since $R_{ij}' = R_{2,3}$ is even. So $i = 4$, and hence $R_{3,4}$ is even and
\[(X\cap V(L(H'))) \setminus V(R_{3,4}) = \{r_{3,1},r_{4,1},r_{3,2},r_{4,2}\}.\]
Since the path
$r_{3,2} \d R_{2,3} \d r_{2,3} \d r_{2,4} \d R_{2,4} \d r_{4,2}$ can be completed to a hole via
$r_{4,2} \d r_{4,3} \d R_{3,4} \d r_{3,4} \d r_{3,2}$, it follows that the first path is even, and
so exactly one of $R_{2,3},R_{2,4}$ is odd; and since the same path can be completed to a hole via
$r_{4,2} \d r_{4,1} \d R_{1,4} \d r_{1,4} \d r_{1,3}\d r_{3,2}$ it follows that $R_{1,4}$ is odd.
Since one of $R_{2,3},R_{2,4}$ is odd, they do not both have length 0, and hence at most
one of $r_{2,3},r_{2,4} \in X$. Since $X$ saturates $L(H)$, it follows that exactly one of
$r_{2,3},r_{2,4} \in X$ (and hence one of $R_{2,3},R_{2,4}$ has length 0), and also
that $r_{2,1} \in X$. Since no vertex of $R_{1,2}'$ is in $X$, this restores the symmetry between
$T_1'$ and $T_2'$.

Suppose that $R_{2,3}$ has length 0. Then $R_{2,4}$ and $R_{1,2}$ are odd, and in particular
$r_{2,1} \not = r_{1,2}$. If $r_{2,1}$ has no neighbour in $R_{1,2}'$, then
$y\d r_{2,1}\d r_{2,4}\d r_{2,1}'\d R_{1,2}'\d r_{1,2}'\d r_{1,4}\d R_{1,4}\d r_{4,1}\d y$ is
an odd hole, a contradiction. So
$r_{2,1}$ has a neighbour in $R_{1,2}'$; but then $y$ can be linked onto the triangle $T_1'$
via $R_{1,2}'$ and $R_{1,4}$, contrary to \ref{trianglev}. This proves that $R_{2,3}$ has
length $\ge 1$. Hence $R_{2,3}$ has odd length and $R_{2,4}$ has length 0, and consequently
$R_{1,2}$, $R_{3,4}$ have even length and $R_{1,4}$ is odd. If $R_{3,4}$ has positive length
then $L(H)$ is overshadowed (because of the vertex $y$), and so the theorem holds. We
may therefore assume that $R_{3,4}$ has length 0. If $r_{2,1} \not = r_{1,2}$ and $r_{2,1}$ has no
neighbour in $R_{1,2}'$, then $y \d r_{2,1}\d r_{2,4}\d r_{2,1}'\d R_{1,2}'\d r_{1,2}'\d r_{1,3}\d y$ is an
odd hole, a contradiction;
while if $r_{2,1} \not = r_{1,2}$ and $r_{2,1}$ has a neighbour in $R_{1,2}'$, then
then $y$ can be linked onto the triangle $T_1'$
via $R_{1,2}'$ and $R_{1,4}$, contrary to \ref{trianglev}. So $r_{2,1} = r_{1,2}$.
But then $L(H)$ is degenerate. Since the strip system is
nondegenerate, it follows from \ref{mixedrungs2} that there is an overshadowed appearance
of $K_4$ in $G$.  This proves \ref{bigtosmall}.\bbox

A $J$-strip system $(S,N)$ in $G$ is {\em maximal} if there is no $J$-strip system $(S',N')$ in $G$
such that $V(S,N) \subset V(S',N')$, and $S_{uv}' \cap V(S,N) = S_{uv}$ for every $uv \in E(J)$,
and $N_v \subseteq N_v'$ for every $v \in V(J)$.  We need to analyze maximal strip systems.
For an edge $uv \in E(J)$, we call the set $S_{uv}$ a {\em strip} of the strip system.

\begin{thm}\label{stripnbrs1}
Let $G$ be Berge, let $J$ be a 3-connected graph, and let $(S,N)$ be a maximal $J$-strip system in $G$.
Assume that either:
\begin{itemize}
\item $(S,N)$ is nondegenerate, and there is no $J$-enlargement with a nondegenerate appearance in $G$, or
\item $J = K_{3,3}$, and there is no $J$-enlargement that appears in either $G$ or $\overline{G}$.
\end{itemize}
Let $F \subseteq V(G) \setminus V(S,N)$ be connected and contain no vertices that are major with
respect to $(S,N)$. Then either $J = K_4$ and there is an overshadowed appearance of $J$ in $G$, or
the set of the attachments of $F$ in $V(S,N)$ is local.
\end{thm}
\Proof
We may assume that if $J = K_4$ then there is no overshadowed appearance of $J$ in $G$.
Let $X$ be the set of attachments of $F$ in $V(S,N)$, and we suppose for a contradiction that $X$
is not local. We may assume that $F$ is minimal (connected) with this property.
\\
\\
(1) {\em For every choice of rungs, forming $L(H)$ say:
\begin{itemize}
\item for each $y \in F$, the set of neighbours of $y$ does not saturate $L(H)$, and
\item if $J = K_4$ then $L(H)$ is not degenerate.
\end{itemize}}
For no $y \in F$ is major with respect to the strip system, and
no $J$-enlargement has a nondegenerate appearance in $G$, and there is no overshadowed
appearance of $J$ in $G$, so the first claim follows from \ref{bigtosmall}.
For the second claim, assume $J = K_4$; then by hypothesis, the strip system is not
degenerate, and the claim follows from \ref{mixedrungs2}.  This proves (1).
\\
\\
(2) {\em There is no $v \in V(J)$ so that $X \subseteq \bigcup (S_{uv}: uv \in E(J))$.}
\\
\\
For assume that $v$ is such a vertex. Since $X$ is not local, there exists $x \in
X \cap S_{uv} \setminus N_v$ for some edge $uv$ of $J$. Since $X \not \subseteq S_{uv}$,
there exists $x' \in X \cap S_{u'v}$ for some edge $u'v$ of $J$ with $u' \not = u$.
For $w \in V(J)$, $x$ belongs to $N_w$ only if $w = u$, and $x'$ belongs to $N_w$ only if
$w \in \{v,u'\}$; and since $x,x'$ do not belong to the same strip it follows that
$\{x,x'\}$ is not local with respect to the strip system.
Make a choice of rungs $R_{ij}$ ($ij \in E(J)$) so that $x \in V(R_{uv})$ and $x' \in V(R_{u'v})$,
forming $L(H)$. Then $\{x,x'\}$ is not local with respect to $L(H)$, so by (1) we can apply
\ref{smallcomp}. Suppose that \ref{smallcomp}.1 holds. Then there is an appearance $L(H')$ in $G$ of
some $J$-enlargement $J'$, with $L(H)$ an induced subgraph of $L(H')$. Moreover,
if $J' = K_{3,3}$ then $J = K_4$, and so $L(H)$ is nondegenerate and therefore so is $L(H')$. Since
$J' \not = K_4$ it follows that $L(H')$ is nondegenerate, contrary to hypothesis. So \ref{smallcomp}.1
does not hold, and hence
\ref{smallcomp}.2.a holds, and there is a branch $D$ of $H$ with an end $d$ so that
$\delta_H(d)\setminus E(D) = (X \cap E(H)) \setminus E(D)$.  Since $x$ and $x'$ are disjoint edges
in $X \cap E(H)$, they are not both incident with $d$, and so one of them is in $E(D\setminus d)$.
The branch containing $x'$ does not meet $x$, so $D$ is the branch between $u$ and $v$, and
$d = v$.  Hence $x'$ is incident
with $v$ in $H$, and $\delta_H(v) \subseteq X \cup E(D)$. Consequently, for all neighbours
$w \not = u$ of $v$ in $J$, $X$ contains the vertex of $R_{vw}$ that belongs to $N_v$, and contains
no other vertex of $R_{vw}$. This restores the symmetry between $u'$ and the other neighbours
of $v$ different from $u$; and since it holds for all choices of the rungs $R_{vw}$, we deduce that
$X \setminus S_{uv} = N_v \setminus S_{uv}$. The minimality of $F$ implies that there is a path $P$
with $V(P) = F$, with ends $p_1,p_2$ so that $p_1$ is complete to $N_v \setminus N_{vu}$, and
no other vertex of $P$ has any neighbours in $N_v \setminus N_{vu}$, and $p_2$ is adjacent to
$x$, and no other vertex of $P$ has any neighbours in $S_{uv} \setminus N_v$. But then we can
add $p_1$ to $N_v$ and $F$ to $S_{uv}$, contradicting the maximality of $(S,N)$. This proves (2).

\bigskip

Let $K = \{uv \in E(J): X \cap S_{uv} \not = \emptyset \}$.
\\
\\
(3) {\em There are two disjoint edges in $K$.}
\\
\\
For make a choice of rungs  $R_{uv}$ ($uv \in E(J)$) so that $X \cap V(R_{uv}) \not = \emptyset$
for each $uv \in K$, forming $L(H)$. If there are no two disjoint edges in $K$, then by (1)
and \ref{smallcomp}, it follows that either $X \cap V(L(H))$ is local (with respect to $L(H)$)
or \ref{smallcomp}.2.a holds, and in either case
there is a branch $D$ of $H$ with an end $d$ so that every edge of $X \cap E(H)$
either is in $E(D)$ or is incident with $d$. In particular, every branch containing an edge of
$X$ is incident with $d$, and so $d$ meets all edges of $J$ in $K$, contrary to (2). This proves (3).

\bigskip

From (3) it follows that there exists a 2-element subset of $X$ that is not local, and so
from the minimality of $F$ it follows that $F$ is the vertex set of a path, say $f_1,\ldots,f_n$.
Let us say a choice $R_{uv}$ ($uv \in E(J)$) of rungs is {\em broad} if
there are two disjoint edges $ij$ and $hk$ of $J$ so that $X$ meets both $R_{ij}$ and $R_{hk}$.
From (3) there is a broad choice.  We denote the ends of $R_{uv}$ by $r_{uv}$ and $r_{vu}$,
where $r_{uv} \in N_u$ and $r_{vu} \in N_v$.
\\
\\
(4) {\em For every broad choice of rungs $R_{uv} (uv \in E(J))$, there is a unique
pair $(i,j)$ of adjacent vertices of $J$ so that:
\begin{itemize}
\item for every $w\in V(J)$ different from $j$ and adjacent to $i$ in $J$, $r_{iw}f_1$ is the unique
edge of $G$ between $V(R_{iw})$ and $F$,
\item for every $w\in V(J)$ different from $i$ and adjacent to $j$ in $J$, $r_{jw}f_n$ is the unique
edge of $G$ between $V(R_{jw})$ and $F$,
\item for every edge $uv$ of $J$ disjoint from $ij$, there are no edges of $G$ between $V(R_{uv})$ and $F$.
\end{itemize}}
For by (1) we can apply \ref{smallcomp}, and since the choice of rungs is broad, the minimality of $F$
implies  that one of \ref{smallcomp}.2.b, \ref{smallcomp}.2.c, \ref{smallcomp}.2.d holds.
Hence there is an edge $ij$ as in (4).  Suppose there is another, say
$i'j'$. Since $i'j'$ meets all edges of $J$ that share exactly one end with $ij$, and $J$ is
3-connected, it follows that $J = K_4$ and the two edges $ij, i'j'$ are disjoint. Moreover,
the unique vertex of $R_{ii'}$ in $X$ is both $r_{ii'}$ and $r_{i'i}$, so $R_{ii'}$ has length 0.
Similarly $R_{ij'},R_{ji'},R_{jj'}$ all have length 0, and so $L(H)$ is degenerate, contrary to (1).
This proves (4).
\\
\\
(5) {\em Every choice of rungs is broad.}
\\
\\
For from (3), there is a broad choice, and from (4) in any broad choice $R_{uv}$ ($uv \in E(J)$)
there are four different edges $a_1b_1,\ldots,a_4b_4$ of $J$, so that $a_1b_1$ is disjoint from $a_2$,
and $a_3b_3$ is disjoint from $a_4b_4$, and $X$ meets $R_{a_ib_i}$ for $1 \le i \le 4$. Consequently,
if we take another choice of rungs, differing from this one on only one edge, then it too is
broad. It follows that every choice is broad. This proves (5).

\bigskip
For a given choice of rungs, let us call the edge $ij$ as in (4) the {\em traversal} for the choice.
\\
\\
(6) {\em There are two choices of rungs with different traversals.}
\\\
\\
Take a choice of rungs, and let $ij$ be its traversal; and suppose that all other
choices of rungs have the same traversal. Let $A_1 = N_i \setminus S_{ij}$,
and $A_2 = N_j \setminus S_{ij}$. From (4),(5), and the uniqueness of $ij$ it follows that
$X \cap (V(S,N) \setminus S_{ij}) = A_1 \cup A_2$. Hence $n \ge 2$, for if $n = 1$ then
we can add $f_1$ to $N_i,N_j$ and $S_{ij}$, contrary to the maximality of the strip system.
Choose $x_1 \in A_1$ and $x_2 \in A_2$ in disjoint strips. From (4), $x_1$ is adjacent
to exactly one of $f_1,f_n$, say $f_1$. For any other vertex $x_3 \in A_2$, let $R_{uv}$ ($uv \in E(J)$)
be a choice of rungs forming $L(H)$ say, so that $x_1,x_3 \in V(H)$. From (4) and (5)
it follows that $f_n$ is adjacent to $x_3$; and so $f_n$ is complete to $A_2$, and similarly
$f_1$ is complete to $A_1$. From the minimality of $F$, there are no other edges between
$F$ and $A_1 \cup A_2$; but then we can add $f_1$ to $N_i$, $f_n$ to $N_j$, and $F$ to
$S_{ij}$, contrary to the maximality of the strip system. This proves (6).

\bigskip

Let us say a choice $R_{uv}$ ($uv \in E(J)$) is {\em optimal} if $R_{uv}$ has a vertex in $X$ for all
edges $uv$ in $K$. For any choice of rungs, there is an optimal choice with the same traversal
(just replace rungs that miss $X$ by rungs that meet $X$ wherever possible); so (6)
implies that there are two optimal choices of rungs with different traversals.
Now for any optimal choice of rungs, if $hi$ is its traversal, then by (4) and the optimality
of the choice, it follows that $K$ consists precisely of the edges of $J$ with exactly one
end in common with $hi$, together possibly with $hi$ itself. In particular $hi$ meets all edges
in $K$. We may assume that some other edge $jk$ is the traversal for some other optimal
choice; and hence (since $J$ is 3-connected) it follows that $J = K_4$ and $jk$ is
disjoint from $hi$, and neither edge is in $K$. Hence $V(J) = \{h,i,j,k\}$. Now since the strip
system is not degenerate, there is one of the four edges $hj,hk,ij,ik$ whose strip contains
a rung of nonzero length; some $hj$-rung $R$ has length $>0$ say. From (4) it follows
that exactly one vertex of $R$ is in $X$, one of its ends; say the end in $N_h$.
Let $R_{uv}$ ($uv \in E(J)$) be any choice of rungs such that $R_{hj} = R$. Since
the end of $R$ in $N_j$ does not belong to $X$, it follows from (4) that for each of
$R_{hk},R_{ij},R_{ik}$, its unique vertex in $X$ is its end in $N_h \cup N_i$. Since the
choice of these rungs was arbitrary, it follows that $X \cap S_{hk} = N_{hk}$,
$X \cap S_{ij} = N_{ij}$, and $X \cap S_{ik} = N_{ik}$. If also $X \cap S_{hj} = N_{hj}$
then $hi$ is the traversal for every choice of rungs, contrary to (6), so
$X \cap S_{hj} \not = N_{hj}$. It follows that every $ij$-rung has length 0; for if one, $R'$
say, has length $>0$, then its unique vertex in $X$ is its end in $N_i$, and by exchanging
$h$ and $i$ it follows that $X \cap S_{hj} = N_{hj}$, a contradiction. Similarly all
$hk$ and $ik$-rungs have length 0, and therefore all $hj$-rungs have even length, since
$G$ is Berge. From (1), we may assume that $f_1$ is adjacent to $r_{hj}$ and complete to $S_{hk}$,
and $f_n$ is complete to $S_{ij}\cup S_{ik}$, and there are no other edges between $F$ and
$S_{hk} \cup S_{ij}\cup S_{ik}\cup \{r_{hj}\}$.
Let $R'$ be an $hj$-rung such that its vertex in $N_h$ ($r_{hj}'$, say) is not its unique vertex
in $X$. Consequently, its other end ($r_{jh}'$) is its unique vertex in $X$. By the same argument
with $hi$ and $jk$ exchanged, it follows that one of $f_1,f_n$ is complete to $S_{ij} \cup \{r_{jh}'\}$
and the other to $S_{hk} \cup S_{ik}$; and hence $n = 1$. But then
the path $f_1\d r_{hj}\d R_{hj}\d r_{jh}\d r_{ji}\d f_1$ is an odd hole, a contradiction.
This proves \ref{stripnbrs1}.\bbox

We are now ready to prove \ref{hypHuse}, which we restate:

\begin{thm}
Let $G$ be Berge. Let $J$ be a 3-connected graph, such that either:
\begin{itemize}
\item there is a nondegenerate appearance $L(H)$ of $J$ in $G$, and there is no
$J$-enlargement with a nondegenerate appearance in $G$, or
\item $J = K_{3,3}$, there is an appearance $L(H)$ of $J$ in $G$, and
no $J$-enlargement appears in either $G$ or $\overline{G}$.
\end{itemize}
Then either $G = L(H)$, or $G$ admits a 2-join or a balanced skew partition.
\end{thm}
\Proof If $J = K_4$ or $K_{3,3}$ and
some appearance of $J$ in $G$ or $\overline{G}$ is overshadowed, the theorem follows from
\ref{overshadowed}, so we assume not. Choose an appearance $L(H_0)$ of $J$ in $G$, nondegenerate if possible.
From the hypothesis, it follows that if $L(H_0)$ is degenerate, then $J = K_{3,3}$, and there is
no nondegenerate appearance of $K_{3,3}$ in $G$, and no $J$-enlargement appears in either $G$ or
$\overline{G}$.  Regard $L(H_0)$ as a $J$-strip system in the natural way, and enlarge it to
a maximal $J$-strip system $(S,N)$. If $L(H_0)$ is nondegenerate then so is the strip system.
Let $Y$ be the set of vertices in $V(G) \setminus V(S,N)$
that are major with respect to the strip system, and let $Z = V(G) \setminus (V(S,N) \cup Y)$.
By \ref{stripnbrs1}, for each component of $Z$, its set of attachments in $V(S,N)$ is local.
\\
\\
(1) {\em If $Y \not = \emptyset$ then $G$ admits a balanced skew partition.}
\\
\\
For let $Y'$ be an anticomponent of $Y$, and let $X$ be the set of all $Y'$-complete vertices
in $V(G)$. For every choice of rungs, forming $L(H)$ say, and for
every $y \in Y'$, the set of neighbours of $y$ in $L(H)$ saturates $L(H)$. Since we may assume
that $L(H)$ is nondegenerate if $J = K_4$ (because of \ref{mixedrungs2}), it follows from
\ref{nonsat} that $X$ saturates $L(H)$. Since this holds for every choice of rungs, it follows
that $X$ saturates the strip system. Let $b_1b_2$ be an edge of $J$, chosen if possible so that
$S_{b_1b_2} \not \subseteq X$. Now the sets $(N_{b_1v}$: $b_1v \in E(J))$ form a partition of
$N_{b_1}$ into say $m$ sets, and at least $m-1$ of them are subsets of $X$.
Choose $m-1$ of them that are subsets of $X$, not using $N_{b_1b_2}$ if possible (that is, if the
other $m-1$ sets are all subsets of $X$), and let their union be $X_1$. Define
$X_2 \subseteq N_{b_2}$ similarly. We note that $S_{b_1b_2} \not \subseteq X_1 \cup X_2$; for if
some vertex of $S_{b_1b_2}$ is not in $X$ then this is clear, while if $S_{b_1b_2} \subseteq X$
then $V(S,N) \subseteq X$ from our choice of $b_1b_2$, and then from the way we chose $X_1$ it
follows that $X_1 \cap S_{b_1b_2} = \emptyset$, and similarly $X_2 \cap S_{b_1b_2} = \emptyset$,
and again our claim holds.  This proves that $S_{b_1b_2} \not \subseteq X_1 \cup X_2$.
Define $X_3$ to be the set of vertices in $X \cap V(S,N)$
that are not in $X_1 \cup X_2$, and let $X_0$ be the set of vertices of $X$ that are not in $V(S,N)$.
So $X_0,X_1,X_2,X_3$ are four disjoint subsets of $X$, with union $X$. Note that
$Y \setminus Y' \subseteq X_0$. Let $B$ be the union of all components of
$G \setminus (Y' \cup X_0 \cup X_1 \cup X_2)$ that have nonempty intersection with
$V(L(H)) \setminus S_{b_1b_2}$, and let $A$ be the union of all the other components. We claim
that $B$ is nonempty; for there is an edge $c_1c_2$ of $J$ disjoint from $b_1b_2$, and no vertex
of $S_{c_1c_2}$ is in $N_{b_1} \cup N_{b_2} \cup S_{b_1b_2}$, and therefore no vertex of
$S_{c_1c_2}$ is in $Y' \cup X_0 \cup X_1 \cup X_2$. Suppose that $A$ is also nonempty. Then
$(A \cup B, Y' \cup X_0 \cup X_1 \cup X_2)$ is a skew partition of $G$. By \ref{geteven}
we may assume it is not loose; and so $X_3$ is empty (since any vertex of $X_3$ is in $A \cup B$
and yet is complete to $Y'$). In particular, $X \cap V(L(H)) \subseteq N_{b_1} \cup N_{b_2}$.
Since $X \cap V(L(H))$ saturates the strip system, it follows that for every vertex $w$ of $J$
different from $b_1,b_2$, $w$ has at most one neighbour in $J$ different from $b_1,b_2$, and
$w$ is adjacent in $J$ to both $b_1$ and $b_2$, and all $wb_1$ and $wb_2$-rungs have length 0. Since
$J$ is 3-connected it follows that $J = K_4$; but then the strip system is not nondegenerate,
a contradiction. So in this case the theorem
holds. We may therefore assume that $A$ is empty. Now we already saw that
$S_{b_1b_2} \not \subseteq X_1 \cup X_2$. Since $A$ is empty, it follows that there is
a path of $G$ between $S_{b_1b_2}$ and $V(L(H)) \setminus S_{b_1b_2}$, disjoint from
$Y' \cup X_0 \cup X_1 \cup X_2$. Choose such a path, minimal. From the choice of $X_1$ and $X_2$ this
path has a nonempty interior; from its minimality, none of its internal vertices belong to
$L(H)$; since all major vertices are in $Y' \cup X_0$, its interior contains no major vertices;
by \ref{stripnbrs1}, the set of attachments of its interior is local; yet its ends are
both attachments of its interior, so there exist $u \in S_{b_1b_2}$ and $v \in V(L(H)) \setminus
S_{b_1b_2}$, so that $u,v \not \in X_1 \cup X_2$, and yet
$\{u,v\}$ is local. Now $u,v$ do not lie in the same strip, and therefore there is some $N_{w}$
containing them both; and the only $w \in V(J)$ with $u \in N_w$ are $b_1,b_2$, so we may assume
that $u,v \in N_{b_1}$. Since they are not in $X_1$, and not in the same strip, this is impossible.
This proves (1).

\bigskip

We may therefore assume that $Y$ is empty.
\\
\\
(2) {\em If there is a component $F$ of $Z$ so that for some $v \in V(J)$, all attachments of $F$
in $L(H)$ belong to $N_v$, then $G$ admits a balanced skew partition.}
\\
\\
For let $F' = V(G) \setminus (F \cup N_v)$; then $F' \not = \emptyset$, and every path in $G$
from $F$ to $F'$ meets $N_v$. Since $N_v$ is not anticonnected, it follows that $F \cup F', N_v)$
is a skew partition. By \ref{geteven} we may assume it is not loose; and we will prove that
it is balanced. Let the neighbours of $v$ in $J$ be $u_1,\ldots,u_k$; then every anticomponent of
$N_v$ is a subset of one of $N_{vu_1},\ldots,N_{vu_k}$. Choose a neighbour $w$
of $u_1$ in $J$ different from $v,u_2$, choose $n_1 \in N_{u_1w}$, and choose $n_2 \in N_{vu_2}$.
Then $n_1,n_2$ belong to strips $S_{u_1w},S_{vu_2}$, where $u_1w, vu_2$ are disjoint edges of $J$;
and so $n_1,n_2$ are not adjacent in $G$.
Let $K = \{n_1\} \cup S_{vu_1} \setminus N_{vu_1}$. Then $K$ is connected (since every vertex of
$S_{vu_1}$ is in a $vu_1$-rung and $n_1$ is complete to $N_{u_1v}$), every vertex in $N_{vu_1}$
has a neighbour in $K$ (for the same reason), and $n_2$ is not in $K$ and has no neighbour in $K$.
(For the last claim, $n_2$ is not in $K$ since it is in only one strip; and it has no neighbour
in $S_{vu_1} \setminus N_{vu_1}$ from the definition of a strip system; and it is not adjacent
to $n_1$ as we already saw.) By \ref{balancev}, $(K,N_{u_1v})$ is balanced, and so by
\ref{shiftbalance}.1, so is $(K,F)$.  By \ref{onepair}, $G$ admits a balanced skew partition.
This proves (2).

\bigskip

We assume therefore that there are no such components $F$ of $Z$. Consequently, for every component
$F$ of $Z$, there is an edge $b_1b_2$ of $J$ so that all the attachments of $F$ in $L(H)$ are in
$S_{b_1b_2}$. If $Z$ is empty and for all $b_1b_2$ there is only one $b_1b_2$-rung,
then $G = L(H)$ and the theorem holds. So we may assume that there is an edge
$b_1b_2$ of $J$ such that either there is more than one $b_1b_2$-rung in $S_{b_1b_2}$ or there
is a component $F$ of $Z$ with all its attachments in $S_{b_1b_2}$. Let $A$ be the union of
$S_{b_1b_2}$ and any components of $Z$ that have an attachment in $S_{b_1b_2}$ (and which
therefore have attachments only in $S_{b_1b_2}$), and $B = V(G) \setminus A$. Let
$A_1 = N_{b_1b_2}$, $A_2 = N_{b_2b_1}$, $B_1 = N_{b_1} \setminus N_{b_1b_2}$,
and $B_2 = N_{b_2} \setminus N_{b_2b_1}$. Then $A_1,A_2 \subseteq A$, and $B_1,B_2$ are
disjoint subsets of $B$, and for $i = 1,2$ $A_i$ is complete to $B_i$, and there are no other
edges between $A$ and $B$. Also $|B_1| \ge 2$, and we chose $b_1b_2$ so that if $A_1,A_2$ both have
only one vertex then $A$ is not the vertex set of a path joining them.
Since we may assume that $G$ does not admit a 2-join, it follows that $A_1 \cap A_2 \not = \emptyset$.
Choose $a \in A_1 \cap A_2$. Then $a$ is complete to $B_1 \cup B_2$, and since $|A| \ge 2$, it
follows that $((B\setminus (B_1 \cup B_2)) \cup (A\setminus a), B_1 \cup B_2 \cup \{a\})$ is a skew
partition of $G$.  Since $\{a\}$ is an anticomponent of $B_1 \cup B_2 \cup \{a\}$, \ref{singleton}
implies that $G$ admits a balanced skew partition. This proves \ref{hypHuse}. \bbox

\section{Bicographs}

In this section we handle degenerate appearances of $K_4$. There is another way to view
them, not as line graphs but as sets of paths and antipaths with certain properties, as
we shall see.

Let $P_1,P_2$ be paths in a graph $G$, and let $Q_1,Q_2$ be antipaths. Suppose that
$P_1,P_2,Q_1,Q_2$ are pairwise disjoint, and we can label the ends of each $P_i$ as $a_i,b_i$,
and label the ends of each $Q_j$ as $x_j,y_j$, so that:
\begin{itemize}
\item $P_1,P_2,Q_1,Q_2$ all have length $\ge 1$
\item there are no edges between $P_1$ and $P_2$, and $Q_1$ is complete to $Q_2$
\item for $(i,j) = (1,1),(1,2)$ or $(2,1)$, the only edges between $V(P_i)$ and $\{x_j,y_j\}$
are $a_ix_j$ and $b_iy_j$, and the only edges between $V(P_2)$ and $\{x_2,y_2\}$
are $a_2y_2$ and $b_2x_2$,
\item  for $(i,j) = (1,1),(1,2)$ or $(2,1)$, the only nonedges between $V(Q_j)$ and $\{a_i,b_i\}$
are $a_iy_j$ and $b_ix_j$, and the only nonedges between $V(Q_2)$ and $\{a_2,b_2\}$
are $a_2x_2$ and $b_2y_2$.
\end{itemize}
In these circumstances we call the quadruple $(P_1,P_2,Q_1,Q_2)$ a {\em knot} in $G$. Note that
if $(P_1,P_2,Q_1,Q_2)$ is a knot then so is $(P_2,P_1,Q_1,Q_2)$, with a suitable relabelling
of the ends of the paths and antipaths.

If $L(H)$ is a degenerate appearance of $K_4$ in $G$, it can be viewed as a knot. For, in
our usual notation, let $R_{1,3},R_{1,4}, R_{2,3},R_{2,4}$ have length 0; let $P_1 = R_{1,2}$,
$P_2 = R_{3,4}$, let $Q_1$ be the antipath $r_{1,3}\d r_{2,4}$, and $Q_2$ the antipath
$r_{1,4}\d r_{2,3}$. It is easy to check that this is a knot.
In fact, this and its complement are the only knots in Berge graphs, as the next theorem shows.

\begin{thm}\label{knottype}
Let $(P_1,P_2,Q_1,Q_2)$ be a knot in a Berge graph $G$. Then all four of $P_1,P_2,Q_1,Q_2$
have odd length; and either both $P_1,P_2$ have length 1, or both $Q_1,Q_2$ have length 1.
\end{thm}
\Proof Certainly $P_1$ is odd since $x_1\d a_1\d P_1\d b_1\d y_2\d x_1$ is a hole, and similarly
the other three are odd. Suppose one of $P_1,P_2$ has length $>1$ and one of $Q_1,Q_2$ has
length $>1$. By exchanging
$P_1,P_2$ or $Q_1,Q_2$ we may therefore assume that $P_1,Q_1$ both have length $>1$. Let
$Y$ be the interior of $Q_1$. Then $a_1,b_1,a_2,b_2$ are all $Y$-complete, from the last
condition in the definition of a knot, and since $a_2$ has no neighbours in the interior
of $P_1$ it follows from \ref{greentouch} that there is a $Y$-complete vertex ($v$ say)
in the interior
of $P_1$. But $x_1,y_1$ are not $Y$-complete, and they are adjacent, so $a_1\d x_1\d y_1\d b_1$ is an odd
path between $Y$-complete vertices and $v$ has no neighbour in its interior, contrary to
\ref{greentouch}. This proves \ref{knottype}.\bbox

Nevertheless, it turns out to be advantageous to make only limited use of \ref{knottype}; it is
better to preserve the symmetry between the paths and the antipaths.

Let $(P_1,P_2,Q_1,Q_2)$ be a knot in a Berge graph $G$; we define $K$ to be the subgraph
of $G$ induced on $V(P_1) \cup V(P_2)\cup V(Q_1) \cup V(Q_2)$. (For brevity we say that
the knot {\em induces} $K$.)
We say a subset $X \subseteq V(K)$ is {\em local} (with respect to the knot) if
$X$ is disjoint from one of $V(P_1),V(P_2)$, and $X$ includes neither of $V(Q_1),V(Q_2)$,
and $X \cap (V(P_1) \cup V(P_2))$ is complete
to $X \cap (V(Q_1) \cup V(Q_2))$. We say $X$ {\em resolves} the knot
if $V(K) \setminus X$ is
local with respect to the knot $(Q_1,Q_2,P_1,P_2)$ in $\overline{G}$; that is, if $X$ includes
one of $V(Q_1),V(Q_2)$, and $X$ meets both $P_1$ and $P_2$, and $X$ contains at least one end
of every edge between $V(P_1) \cup V(P_2)$ and $V(Q_1) \cup V(Q_2)$.
Conveniently, these definitions almost agree with what we did for line graphs, because of the
following.

\begin{thm}\label{bigness}
Let $(P_1,P_2,Q_1,Q_2)$ be a knot in a graph $G$, inducing $K$, where $Q_1,Q_2$ both have
length 1, and so $K = L(H)$ is an appearance of $K_4$. Let $X \subseteq V(K)$.
Then:
\begin{itemize}
\item $X$ is local with respect to the knot if and only if it is local with respect
to $L(H)$
\item $X$ resolves the knot if and only if $X$ saturates $L(H)$ and
$X$ meets both $V(P_1)$ and $V(P_2)$.
\end{itemize}
\end{thm}
The proof is obvious and we omit it.
This allows us to unify some portions of \ref{smallcomp} and \ref{nonsat}, as follows.

\begin{thm}\label{smallbico}
Let $(P_1,P_2,Q_1,Q_2)$ be a knot in a Berge graph $G$, inducing $K$.
Assume that there is no appearance in $G$ or in $\overline{G}$ of any $K_4$-enlargement,
and there is no overshadowed appearance of $K_4$ in $G$ or in $\overline{G}$.
Let $F$ be a connected subset of $V(G) \setminus V(K)$, such that its set of attachments
in $K$ is not local. Then either:
\begin{enumerate}
\item there is a vertex in $F$ such that its neighbour set in $K$ resolves the knot, or
\item (up to symmetry) there is a path $R$ in $F$ with ends $r_1,r_2$ such that
$r_1,a_1$ have the same neighbours in $V(P_2)\cup V(Q_1)\cup V(Q_2)$, and there are no edges
between $R \setminus r_1$ and $V(P_2)\cup V(Q_1)\cup V(Q_2)$, and $r_2$ has a neighbour
in $P_1 \setminus a_1$, and there are no edges between $R \setminus r_2$ and $P_1 \setminus a_1$, or
\item (up to symmetry) there is an odd path $R$ in $F$ with ends $r_1,r_2$ such that
$r_1,a_1$ have the same neighbours in $V(P_2)\cup V(Q_1)\cup V(Q_2)$, and
$r_2,b_1$ have the same neighbours in $V(P_2)\cup V(Q_1)\cup V(Q_2)$,
and there are no edges between $V(R^*)$ and $V(P_2)\cup V(Q_1)\cup V(Q_2)$, and no edges between
$R$ and $P_1$ except possibly $r_1a_1$ and $r_2b_1$, or
\item there is a vertex $f\in F$  so that (up to symmetry)
$f,x_1$ have the same neighbours in $V(P_1) \cup V(P_2) \cup V(Q_2)$ and $f$ is not adjacent to $y_1$.
\end{enumerate}
\end{thm}
\Proof
By \ref{knottype} there are two cases, depending whether $Q_1$ and $Q_2$ have length 1 or
$P_1,P_2$ have length 1.
\\
\\
(1) {\em If $Q_1,Q_2$ have length 1 then the theorem holds.}
\\
\\
For assume $Q_1,Q_2$ have length 1. Then $K$  is a degenerate appearance of $K_4$ in $G$, say
$K = L(H)$. Suppose that the neighbour set of some $f\in F$ saturates $L(H)$.
If $f$ has a neighbour in both $V(P_1)$ and $V(P_2)$ then statement
1 of the theorem holds, so we assume it has no neighbour in $V(P_1)$. But then $f$ is adjacent to all four
of $x_1,x_2,y_1,y_2$, since it has two neighbours in every triangle of $K$, and then
$f\d x_1\d a_1\d P_1\d b_1\d y_1\d f$ is an odd hole, a contradiction. So we assume there is no such
$f$, and hence we may apply \ref{smallcomp}. If \ref{smallcomp}.1 holds then there is an appearance
in $G$ of some $K_4$-enlargement, a contradiction. So \ref{smallcomp}.2 holds. In the notation
of \ref{smallcomp}.2, the edge $b_1b_2$ of $J$ is of one of two types; either $N_{b_1}$
meets $N_{b_2}$ or it does not. In the first case, we may assume from the symmetry that
those two sets are $\{x_1,y_2,a_2\}$ and $\{x_1,x_2,a_1\}$, and there is a path $R$ of $G$ with
$V(R) \subseteq F$ and with ends $r_1$ and $r_2$, so that $r_1$ is adjacent to $a_1,x_2$, and
$r_2$ is adjacent to $a_2,y_2$, and there are no other edges between $V(P)$ and $K\setminus x_1$. If
$R$ has length 0 then statement 4 of the theorem holds, while if $R$ has length $>0$ then
it is even and there is an overshadowed appearance of $K_4$ in $G$, a contradiction. In the
second case, one of statements 2 and 3 of the theorem hold. This proves (1).

\bigskip

Henceforth we may therefore assume that one of $Q_1,Q_2$ has length $>1$, and therefore
by \ref{knottype}, both $P_1$ and $P_2$ have length 1. Hence $\overline{K} = L(H)$,
where $L(H)$ is a degenerate appearance of $K_4$ in $\overline{G}$.
\\
\\
(2) {\em If there exists $f \in F$ such that $f$ is not major with respect to $L(H)$ in
$\overline{G}$, then the theorem holds.}
\\
\\
For let $f\in F$ have this property.  If
the set of neighbours of $f$ in $K$ resolves the knot $(P_1,P_2,Q_1,Q_2)$, then statement 1 of the
theorem holds, so we assume not. Therefore,
in $\overline{G}$, the set of neighbours of $f$ in $\overline{K}$ is not local with respect
to the knot $(\overline{Q_1},\overline{Q_2},\overline{P_1},\overline{P_2})$. But
this set does not saturate $L(H)$; so we can apply \ref{smallcomp} (or, indeed, \ref{vnbrs})
to $\overline{G}$,
and deduce, as before, that either there is a $K_4$-enlargement that appears in $\overline{G}$
(a contradiction), or (up to symmetry) $f,a_1$ have the same neighbours in $K \setminus a_1$
(but then statement 2 of the theorem holds), or (up to symmetry)
$f,x_1$ have the same neighbours in $V(P_1)\cup V(P_2) \cup V(Q_2)$ (but then either
statement 1 or statement 4 of the theorem holds). This proves (2).

\bigskip

We may therefore assume that every $f \in F$ is major with respect to $L(H)$ in $\overline{G}$.
Let $X$ be the set of vertices of $K$ which, in $G$, have no neighbours in $F$.  By hypothesis,
$V(K) \setminus X$ is not local with respect to the knot $(P_1,P_2,Q_1,Q_2)$ in $G$, and hence $X$
does not resolve the knot $(\overline{Q_1},\overline{Q_2},\overline{P_1},\overline{P_2})$ in $\overline{G}$.
If $X$ does not saturate  $L(H)$ in $\overline{G}$, then by (2) we may apply \ref{nonsat}. Since
$Q_1$ has length $>1$ it follows that \ref{nonsat}.2 holds, and hence
statement 3 of the theorem holds. We may therefore assume that $X$ saturates $L(H)$ in $\overline{G}$.
By \ref{bigness}, $X$ is disjoint from one of $V(Q_1),V(Q_2)$, say $X \cap V(Q_1) = \emptyset$. Hence
$a_1,a_2,b_1,b_2  \in X$. Since $a_1\d y_1 \d Q_1 \d x_1 \d b_1$ is an odd antipath in $G$, and its internal
vertices all have neighbours in $F$, and its ends do not, it follows from \ref{greentouch} applied in
$\overline{G}$ that every vertex in $X$ has a non-neighbour in $V(Q_1)$; and hence no vertex of $Q_2$
belongs to $X$. This restores the symmetry between $Q_1,Q_2$.
Now one of $Q_1,Q_2$ has length $>1$, say $Q_1$ without loss of
generality. Hence, in $\overline{G}$, the path $a_1\d y_1\d Q_1\d x_1\d b_1$ is odd and has length $\ge 5$; its
ends are complete to $F$, and its internal vertices are not. By \ref{RR}, $F$ contains a leap;
so there exist nonadjacent $f_1,f_2 \in F$ such that $Q_1$ is the interior of a path $R$ between them.
(All this is in $\overline{G}$ - we will tell the reader when we switch back to $G$.) Now $f_1,f_2$
have no common neighbour in $Q_2$ (because $R$ could be completed to an odd hole through any such
common neighbour), so by \ref{RR}, $f_1,f_2$  is also a leap for the path $a_1\d y_2\d Q_2\d x_2\d b_1$
(this path might have length 3, but still we get a leap by \ref{RR}.3,
since $\{f_1,f_2\}$ cannot include the interior of any longer antipath between $x_2$ and $y_2$).
Hence from the symmetry we may assume that $f_1$ is adjacent to $y_1,y_2$, and $f_2$ to $x_1,x_2$,
and there are no other edges between $\{f_1,f_2\}$ and $V(Q_1) \cup V(Q_2)$. Therefore, back in $G$,
we see that $a_1,f_1$ have the same neighbours in $V(P_2) \cup V(Q_1) \cup V(Q_2)$, and so do
$b_1,f_2$, and therefore statement 3 of the theorem holds.  This proves \ref{smallbico}.\bbox

\ref{smallbico} suggests that we should attempt to combine paths into strips, as in the section
on ``Generalized line graphs'', and combine antipaths into ``antistrips''. Let us make that precise.

Let $A,B,C$ be disjoint subsets of $V(G)$. We call $S = (A,C,B)$ a {\em strip } if $A,B$ are
nonempty, and every vertex of $A \cup B \cup C$ belongs to a path between $A$ and $B$ with
only its first vertex in $A$, only its last vertex in $B$, and interior in $C$.
Such a path is called a {\em rung} of the strip $S$, or an $S$-rung. When
$S = (A,C,B)$ is a strip, $V(S)$ means $A\cup B\cup C$.
The {\em reverse} of a strip $(A,C,B)$ is the strip $(B,C,A)$.
An {\em antistrip} is a
triple that is a strip in $\overline{G}$, and the corresponding antipaths are called {\em antirungs}.
If $P$ is a rung with ends $a \in A$ and $b \in B$, we speak of the ``rung $a\d P\d b$'' for brevity;
the reader can deduce which end is in which set from the names of the ends, because we shall
always use $a,a',a_1$ etc. for ends in a set called something like $A$, and so on.

Let $S = (A,C,B)$ be a strip and $T = (X,Z,Y)$ an antistrip, with $V(S) \cap V(T) = \emptyset$.
We say $S,T$ are {\em parallel} if:
\begin{itemize}
\item $A$ is complete to $X$,
\item $B$ is complete to $Y$,
\item $A$ is anticomplete to $Y$,
\item $B$ is anticomplete to $X$,
\item $C$ is anticomplete to $X \cup Y$, and
\item $Z$ is complete to $A \cup B$.
\end{itemize}
We say $S,T$ are {\em co-parallel} if $S,T'$ are parallel, where $T'$ is the reverse of $T$.

Now let $S_1,S_2$ be strips and $T$ an antistrip, where $S_1,S_2,T$ are pairwise
disjoint. We say that $S_1,S_2$ {\em agree on} $T$ if either $S_1,T$ are parallel and
$S_2,T$ are parallel, or both pairs are co-parallel; and they {\em disagree} if one pair
is parallel and the other pair is co-parallel.
If $S$ is a strip and $T_1,T_2$ are antistrips, pairwise disjoint, we define
whether $T_1,T_2$ agree or disagree on $S$ similarly.

Now let $S_1,S_2$ be strips, and let  $T_1,T_2$ be antistrips, all pairwise disjoint. We call
the quadruple $(S_1,S_2,T_1,T_2)$ a {\em twist} if $S_1,S_2$ agree on one of $T_1,T_2$
and disagree on the other. (Equivalently, if $T_1,T_2$ agree on one of $S_1,S_2$, and disagree
on the other.) Note that if $(S_1,S_2,T_1,T_2)$ is a twist, then so is $(S_1',S_2,T_1,T_2)$, where
$S_1'$ is the reverse of $S_1$.

A {\em striation} in a graph $G$
is a family of strips $S_i = (A_i,C_i,B_i) (1 \le i \le m)$ together with a family
of antistrips $T_j = (X_j,Z_j,Y_j) (1 \le j \le n)$, satisfying the following conditions:
\begin{itemize}
\item all the strips and antistrips are pairwise disjoint, and all their rungs and antirungs
have odd length
\item $m,n \ge 2$
\item for $1 \le i < i' \le m$, $S_i$ is anticomplete to
$S_{i'}$, and for $1 \le j < j' \le n$, $T_j$ is complete to $T_{j'}$
\item for $1 \le i \le m$ and $1 \le j \le n$, $S_i$ and $T_j$ are either parallel or co-parallel
\item for $1 \le i < i' \le m$ there exist distinct $j,j'$ with $1 \le j,j' \le n$ such that
$(S_i,S_{i'},T_j,T_{j'})$ is a twist
\item for $1 \le j < j' \le n$ there exist distinct $i,i'$ with $1 \le i,i' \le m$ such that
$(S_i,S_{i'},T_j,T_{j'})$ is a twist.
\end{itemize}

(Note that if we replace some $(A_i, C_i,B_i)$ by its reverse, we obtain another striation.)
We denote the striation by $L$, and the union of the vertex sets of all its strips and antistrips
by $V(L)$.  By analogy with what we did for knots, let us say that a subset $X \subseteq V(L)$
is {\em local} with respect to $L$ if
\begin{itemize}
\item at most one of $X \cap V(S_1),\ldots,X\cap V(S_m)$ is nonempty,
\item for $1 \le j \le n$, every $T_j$-antirung has a vertex not in $X$, and
\item $X \cap (V(S_1) \cup \cdots \cup V(S_m))$ is complete to $X \cap (V(T_1) \cup \cdots \cup V(T_n))$.
\end{itemize}
We say $X$ {\em resolves} $L$ if $V(L) \setminus X$ is local with respect to the striation
in $\overline{G}$ obtained from $L$ by exchanging the strips and antistrips; that is, if
\begin{itemize}
\item there is at most one of $T_1,\ldots,T_n$ that is not a subset of $X$,
\item for $1 \le i \le m$, every $S_i$-rung meets $X$, and
\item $X$ contains at least one end of every edge between $V(S_1) \cup \cdots \cup V(S_m)$ and
$V(T_1) \cup \cdots \cup V(T_n)$.
\end{itemize}
A striation $L$ in $G$ is {\em maximal} if there is no striation $L'$ in $G$ with
$V(L) \subset V(L')$.

\begin{thm}\label{lacvnbrs}
Let $G$ be Berge, such that there is no appearance in $G$ or in $\overline{G}$ of any
$K_4$-enlargement, and there is no overshadowed appearance of $K_4$ in $G$ or in $\overline{G}$.
Let $L$ be a maximal striation in $G$.
Let $f \in V(G) \setminus V(L)$, and let $X$ be the set of neighbours of $f$ in $V(L)$.
Then either $X$ is local with respect to $L$, or $X$ resolves $L$.
\end{thm}
\Proof
Let $L$ have strips $S_i = (A_i,C_i,B_i) (1 \le i \le m)$ and antistrips
$T_j = (X_j,Z_j,Y_j) (1 \le j \le n)$.
\\
\\
(1) {\em Let $1 \le i \le m$, and $1 \le j \le n$; let $a_i\d P_i\d b_i$ be an $S_i$-rung,
and $x_j\d Q_j\d y_j$ a $T_j$-antirung.
Then either $X \cap V(P_i) \not = \emptyset$, or $V(Q_j) \not \subseteq X$.}
\\
\\
For suppose that $X$ includes $V(Q_1)$ and is disjoint from $V(P_1)$ say. By reversing $S_2$
we may assume that $S_1$ and $S_2$ agree on $T_1$; and we may assume they disagree on $T_2$.
Let $a_2\d P_2\d b_2$ be any $S_2$-rung, and $x_2\d Q_2\d y_2$ any $T_2$-antirung. Then
$(P_1,P_2,Q_1,Q_2)$ is a knot, so by \ref{knottype}, we may assume (taking complements
if necessary) that $Q_1$ has length 1. But then $f\d x_1\d a_1\d P_1\d b_1\d y_1\d f$ is an odd hole,
a contradiction. This proves (1).

\bigskip
From (1), taking complements if necessary, we may assume that for all $1 \le j \le n$, and
for all $T_j$-antirungs $Q_j$, $V(Q_j) \not \subseteq X$.
\\
\\
(2) {\em $X$ meets at most one of  $V(S_1),\ldots,V(S_m)$.}
\\
\\
For suppose that $X$ meets both $S_1$ and $S_2$ say. We may
assume that $(S_1,S_2,T_1,T_2)$ is a twist. For $i = 1,2$ choose an $S_i$-rung $P_i$ so that
$X \cap V(P_i) \not = \emptyset$, and for
$j = 1,2$ choose any $T_j$-antirung $Q_j$. By our assumption above, $f$ has nonneighbours in both $Q_1,Q_2$.
But then
$(P_1,P_2,Q_1,Q_2)$ is a knot, and setting $F = \{f\}$ violates \ref{smallbico}, a contradiction.
This proves (2).

\bigskip

We may assume that $X$ is not local with respect to $L$, and so we may assume that there is
an $S_1$-rung $a_1\d P_1\d b_1$ and a $T_1$-antirung $x_1\d Q_1\d y_1$ containing nonadjacent members
of $X$. By reversing each $T_j$ if necessary, we may assume that $S_1$ is parallel to each
$T_j$.  In particular, $a_1x_1$ is an edge, and so is $b_1y_1$.
Since the interior of $Q_1$ is complete to $V(P_1)$, we may assume that $x_1 \in X$, and
$X\cap (V(P_1)\setminus a_1) \not = \emptyset$. Let $2 \le j \le n$, and let $x_j\d Q_j\d y_j$ be any
$T_j$-antirung. For definiteness we assume $j = 2$. Now $T_1,T_2$ agree on $S_1$, and so there
is some $S_i$ on which they disagree, say $S_2$. Let $a_2\d P_2\d b_2$ be any $S_2$-rung.
Then $(P_1,P_2,Q_1,Q_2)$ is a knot, with union $K$ say, and $X \cap V(K)$ is not local
with respect to $K$ (since $x_1 \in X$, and $X\cap (V(P_1)\setminus a_1) \not = \emptyset$).
By \ref{smallbico}, it follows that \ref{smallbico}.2 holds, and hence $f, a_1$ have the
same neighbours in $V(Q_1) \cup V(Q_2)$. In particular, $V(Q_2) \setminus \{y_2\} \subseteq X$.
Since $V(Q_2) \not \subseteq X$, it follows that $y_2 \not \in X$;
since this holds for all $Q_2$, we deduce that $X \cap V(T_2) = X_2 \cup Z_2$; and since the
same holds for all antistrips of $L$ except $T_1$, we deduce that $X \cap V(T_j) = X_j \cup Z_j$
for $2 \le j \le n$. Since our only assumption about $T_1$ was that $X \cap X_1 \not = \emptyset$,
and since we have shown that the same is true for all $T_j$, we can replace $T_1$ by $T_2$ say,
and deduce similarly that $X \cap V(T_1) = X_1 \cup Z_1$. But then we can add $f$ to $A_1$,
contrary to the maximality of the striation. This proves \ref{lacvnbrs}.\bbox

\begin{thm}\label{lacFnbrs}
Let $G$ be Berge, such that there is no appearance in $G$ or in $\overline{G}$ of any
$K_4$-enlargement,
and there is no overshadowed appearance of $K_4$ in $G$ or in $\overline{G}$.
Let $L$ be a maximal striation in $G$.
Let $F \subseteq V(G) \setminus V(L)$ be connected, such that for each $f \in F$, the set of its
neighbours in $V(L)$ is local with respect to $L$. Then the set of attachments of $F$ in $V(L)$
is local with respect to $L$.
\end{thm}
\Proof
Let $L$ have strips $S_i = (A_i,C_i,B_i) (1 \le i \le m)$ and antistrips
$T_j = (X_j,Z_j,Y_j) (1 \le j \le n)$.
Suppose not, and choose a counterexample $F$ with $F$ minimal. Let $X$ be its set of attachments
in $V(L)$.
\\
\\
(1) {\em $X \not \subseteq V(T_1) \cup\c \cup V(T_n)$.}
\\
\\
For suppose it is. Since $X$ is not local, we may assume that
$X$ includes $V(Q_1)$ for some $T_1$-antirung $x_1\d Q_1\d y_1$.
Let $2 \le j \le n$, and let $x_j\d Q_j\d y_j$ be a $T_j$-antirung.
Then we can choose some $S_i,S_{i'}$ to make a twist, and if we choose an $S_i$-rung and $S_{i'}$-rung
and apply \ref{smallbico} to the resultant knot, we deduce (since no vertices of $S_i$ and $S_{i'}$
are in $X$) that \ref{smallbico}.3 holds. This has several consequences. First, it implies
that there is an odd path in $F$ with vertices $f_1,\ldots,f_k$ say, which is either parallel or
co-parallel to $Q_1$, and either parallel or co-parallel to $Q_j$; and there are no edges between
$\{f_2,\ldots,f_{k-1}\}$ and $Q_1 \cup Q_j$. Hence the set of attachments of $\{f_1,\ldots,f_k\}$
is not local with respect to $L$, and so $F = \{f_1,\ldots,f_k\}$ from the minimality of $F$.
Second, every vertex
of $Q_j$ is in $X$, and since this holds for all $Q_j$ it follows that $V(T_j) \subseteq X$.
By exchanging $T_1$ and $T_j$ it follows that $V(T_1) \subseteq X$. Moreover, since this
holds for all $j$ we deduce that $X = V(T_1) \cup \cdots \cup V(T_n)$. This restores the symmetry between
$T_1$ and $T_2,\ldots,T_n$.  Third, this shows that there are no edges
between $\{f_2,\ldots,f_{k-1}\}$ and $V(T_1) \cup \cdots \cup V(T_n)$. Fourth, for $1 \le j \le n$
every vertex in $Z_j$ is adjacent to both $f_1,f_k$. Since $k$ is even, this proves that either
$k = 2$ or $Z_1 \cup \cdots \cup Z_n = \emptyset$. Fifth, every vertex in
$X_1\cup Y_1 \cdots \cup X_n\cup Y_n$
is adjacent to exactly one of $f_1,f_n$; let $U$ be the set of those adjacent to $f_1$, and $V$ those
adjacent to $f_n$. For the moment fix $j$ with $1 \le j \le n$. Every $T_j$-antirung has one end in $U$
and the other in $V$; let $M_j$ be the union of the vertex sets of all $T_j$-antirungs $x_j\d Q_j\d y_j$
such that $x_j \in U$, and $N_j$ the union of all those with $x_j \in V$. Since there is no $T_j$-antirung
with both ends in $M_j$ or both ends in $N_j$, it follows that $M_j \cap N_j = \emptyset$, and there are
no nonedges betwen $M_j$ and $N_j$ except possibly between $M_j \cap X_j$ and $N_j \cap X_j$, or between
$M_j \cap Y_j$ and $N_j \cap Y_j$. Suppose there is such a nonedge; and choose $T_j$-antirungs
$x_j\d Q_j\d y_j$,$x_j'\d Q_j'\d y_j'$ where $x_j\in U$ is nonadjacent to $x_j' \in V$, say.
Now $X_j,x_j'$ have a common neighbour $d_1 \in A_1 \cup B_1$,
and then $d_1\d x_j\d f_1\d \cdots\d f_k\d x_j'\d d_1$ is an
odd hole. This proves that $M_j$ is complete to $N_j$. Now if $M_j$ is nonempty, then
$(M_j \cap X_j, M_j \cap Z_j, M_j \cap Y_j)$ is an antistrip, and similarly if $N_j$ is nonempty
it also induces an antistrip. We call these the {\em offspring} of $T_j$. (If one of $M_j,N_j$ is
empty, then the other equals $V(T_j)$, and so the only offspring of $T_j$ is $T_j$ itself;
and otherwise it has two.) Also, there is a new strip
$S_0 = (\{f_1\},\{f_2,\ldots,f_{k-1}\},\{f_k\})$. Note that
\begin{itemize}
\item for all $j$ with $1 \le j \le n$, $S_0$ is parallel or antiparallel with the offspring of $T_j$
\item for all $i$ with $1 \le i \le m$, there exists $j$ with $1 \le j \le n$ such that $S_0,S_i$
disagree on one of the offspring of $T_j$, and there exists $j$ so that $S_0,S_i$ agree on one
of the offspring of $T_j$. For if the first were false, say, then each of the $T_j$'s has
only one offspring, and we could add $f_1$ to $A_i$,
$\{f_2,\ldots,f_{k-1}\}$ to $C_i$, and $f_k$ to $B_i$, contradicting the maximality of the
striation; while if the second were false we could do the same with $f_1,f_k$ exchanged.
\item if $T_1',T_2'$ are each offspring of one of $T_1,\ldots,T_n$, then there exists
$i$ with $0 \le i \le m$ such that $T_1',T_2'$ agree on $S_i$; and there exists $i$ such
that they disagree. For this is clear if they are offspring of different parents, since
their parents were in a twist together; while if they are both offspring of the same $T_j$, then
they disagree on $S_0$ and agree on all of $S_1,\ldots,S_m$.
\end{itemize}
It follows from these observations that the set of strips $S_0,\ldots,S_m$, together with the
set of offspring of $T_1,\ldots,T_n$, forms a new striation, contrary to the maximality of $L$.
This proves (1).
\\
\\
(2) {\em $X$ meets exactly one of $S_1,\ldots,S_m$.}
\\
\\
For by (1) it meets at least one of these sets; suppose it meets two, say $S_1$ and $S_2$.
We may assume that $(S_1,S_2,T_1,T_2)$ is a twist.
For $i = 1,2$ choose an $S_i$-rung $a_i\d P_i\d b_i$ so that $X$ meets $P_i$,
and for $j = 1,2$ let $x_j\d Q_j\d y_j$ be a $Q_j$-antirung.
Then $(P_1,P_2,Q_1,Q_2)$ is a knot $K$ say, and $X \cap V(K)$ is not local with respect to
$K$. From the minimality of $F$, $F$ is minimal such that $X \cap V(K)$ is not local with respect to
$K$. It follows from \ref{smallbico}
that one of \ref{smallbico}.1,  \ref{smallbico}.4 holds;
and in either case there is a vertex $f \in F$ with neighbours in $P_1$ and in $P_2$. Hence
the set of neighbours of $f$ in $V(L)$ is not local with respect to $L$. But this contradicts
a hypothesis of the theorem, and hence proves (2).
\\
\\
(3) {\em $V(Q_j) \not \subseteq X$, for $1 \le j \le n$, and for every $T_j$-antirung $Q_j$. }
\\
\\
For suppose that $V(Q_1)\subseteq X$ for some $T_1$-antirung $x_1\d Q_1\d y_1$. By (2) we may assume that
$X$ meets $S_1$ and none of $S_2,\ldots,S_m$. Let $2 \le j \le n$, and choose $i$ with
$2 \le i \le m$ so that $(S_1,S_i,T_1,T_j)$
is a twist. Let $Q_j$ be an $x_j\d T_j\d y_j$-antirung, let $a_1\d P_1\d b_1$ be an $S_1$-rung so that
$X$ meets $P_1$, and let $a_i\d P_i\d b_i$ be an
$S_i$-rung. Hence $(P_1,P_i,Q_1,Q_j)$ is a knot. Let us apply
\ref{smallbico}. By (2) and the minimality of $F$ it follows that \ref{smallbico}.3 holds. This has several
consequences. First, from the minimality
of $F$, $G|F$ is an odd path $f_1\d \cdots\d f_k$ such that $f_1,a_1$ have the same neighbours in
$V(Q_1 \cup Q_j)$, and so do $f_k,b_1$, and there are no edges between $F$ and $V(P_1)$
except possibly $f_1a_1$ and $f_kb_1$. Since $X$ meets $P_1$, it follows that
at least one of these two edges is present; and therefore they both are, since  $f_1\d \cdots\d f_k$
is an odd path and so is $P_1$ (for otherwise the union of these two paths, with one of $x_1,y_1$,
would induce an odd hole). So $f_1$ is adjacent to $a_1$ and to no other vertex of $P_1$,
and $f_n$ to $b_1$ and to no other vertex of $P_1$. Second, $V(Q_j) \subseteq X$. Since this holds
for all $Q_j$ it follows that $V(T_j) \subseteq X$; and by exchanging $T_1$ and $T_j$ we
deduce that $V(T_1) \cup \cdots\cup V(T_n) \subseteq X$. Moreover $\{f_2,\ldots,f_{k-1}\}$
is anticomplete to  $V(T_1) \cup \cdots\cup V(T_n)$. Third, let $x_j'\d Q_j'\d y_j'$ be some other
$T_j$-antirung. By the same argument applied to the knot $(P_1,P_i,Q_1,Q_j')$, we deduce that
again \ref{smallbico}.3 holds, and so one of $f_1,f_k$ is adjacent to $x_j'$ and the other
to $y_j'$. Furthermore, the one adjacent to $x_j'$ is also adjacent to $a_1$; and so in
fact $f_1$ is adjacent to $x_j'$. Since this holds for all choices of $Q_j$ and of $j$,
it follows that $f_1,a_1$ have the same neighbours in $V(T_1) \cup \cdots\cup V(T_n)$, and
so do $f_k,b_1$. Hence we can add $f_1$ to $A_1$, $\{f_2,\ldots,f_{k-1}\}$ to $C_1$ and
$f_k$ to $B_1$, contrary to the maximality of the striation. This proves (3).

\bigskip

Since $X$ is not local with respect to $L$, we may assume from (2) and (3) that there
exist a vertex of $X \cap V(S_1)$ and a vertex of $X \cap V(T_1)$ that are nonadjacent. By
reversing $T_1,\ldots,T_n$ we may assume that $S_1$ is parallel to each $T_j$. Since every
vertex of $Z_1$ is complete to $V(S_1)$, we may assume that there is an $S_1$-rung
$a_1\d P_1\d b_1$ and a $T_1$-antirung $x_1\d Q_1\d y_1$ so that $x_1 \in X$ and $X \cap V(P_1\setminus a_1)
\not = \emptyset$. Let $2 \le j \le n$, and choose $i$
with $2 \le i \le m$ so that $(S_1,S_i,T_1,T_j)$ is a twist. Let $P_i$ be an $S_i$-rung, and
let $Q_j$ be a $T_j$-antirung. So $(P_1,P_i, Q_1,Q_j)$ is a knot $K$ say, and $X \cap V(K)$ is
not local with respect to $K$. By \ref{smallbico} and (3) and the minimality of $F$,
\ref{smallbico}.2 holds;
let $R$ be the path in $F$ satisfying \ref{smallbico}.2. Then the set of
attachments of $R$ in $V(L)$ is not local with respect to $L$, and so $V(R) = F$ from the
minimality of $F$. Hence $F$ is a path with vertices $f_1\d \cdots\d f_k$ say. Since $x_1 \in X$,
it follows that one of $f_1,f_k$ is adjacent to $x_1$, and we may assume that $f_1$ is
adjacent to $x_1$. By \ref{smallbico}.2, $f_1$ is also adjacent to $x_j$ and to all
internal vertices of $Q_1,Q_j$, and to neither of $y_1,y_j$, and none of $f_2,\ldots,f_{k-1}$
have neighbours in $V(Q_1 \cup Q_j)$, and $f_k$ has a neighbour
in $P_1 \setminus a_1$, and $f_k$ has no neighbours in $V(Q_1 \cup Q_j)$. For any other
choice of $Q_j$ the same happens, and $f_1,f_k$ cannot become exchanged since $f_1$ has
neighbours in $Q_1$ and $f_k$ has none. We deduce that $f_1$ is complete to $X_j \cup Z_j$
and anticomplete to $Y_j$; and $\{f_2,\ldots,f_k\}$ is anticomplete to $V(T_j)$. In particular
there is a vertex of $X \cap V(S_1)$ and a vertex of $X \cap V(T_j)$ that are nonadjacent,
and so by exchanging $T_1$ and $T_j$ in the above argument, we deduce that
$f_1$ is complete to $X_1 \cup Z_1$
and anticomplete to $Y_1$; and $\{f_2,\ldots,f_k\}$ is anticomplete to $V(T_1)$. Since this holds
for all $j$, it follows that $a_1,f_1$ have the same neighbours in $V(T_1) \cup \cdots\cup V(T_n)$,
and there are no edges between $\{f_2,\ldots,f_k\}$ and $V(T_1) \cup \cdots\cup V(T_n)$.
But then we can add $f_1$ to $A_1$ and  $\{f_2,\ldots,f_k\}$ to $C_1$, contrary to the maximality
of the striation. This proves \ref{lacFnbrs}.\bbox

Now we can prove \ref{summary}.3, which we restate.

\begin{thm}\label{bicographs0}
Let $G$ be a Berge graph, such that every appearance of $K_4$ in $G$ and in $\overline{G}$ is degenerate,
and there is no induced subgraph of $G$ isomorphic to $L(K_{3,3})$. Then either
$G$ is a bicograph, or $G$ admits a balanced skew partition, or
one of $G, \overline{G}$ admits a 2-join, or there is no appearance of $K_4$ in either $G$ or $\overline{G}$.
\end{thm}
\Proof If there is an appearance in $G$ of some $K_4$-enlargement, say $L(H')$, then by \ref{nondegen},
either $H' = K_{3,3}$, which is impossible by hypothesis, or
there is a subgraph $H''$ of $H'$ which is a bipartite subdivision of $K_4$, such that $L(H'')$ is nondegenerate,
and again this is impossible by hypothesis. So there is no appearance in $G$
of a $K_4$-enlargement, and similarly there is none in $\overline{G}$. Moreover, by \ref{overshadowed}, we
may assume that there is no overshadowed appearance of $K_4$ in $G$ or in $\overline{G}$.
We may assume that there is an appearance of $K_4$ in one of $G,\overline{G}$, and by taking complements
if necessary we may assume that $L(H)$ is an appearance of $K_4$ in $G$. By hypothesis it is degenerate,
and hence there is a striation in $G$; choose a maximal striation $L$.
Let $L$ have strips $S_i = (A_i,C_i,B_i) (1 \le i \le m)$ and antistrips
$T_j = (X_j,Z_j,Y_j) (1 \le j \le n)$. By \ref{lacvnbrs} we can partition $V(G) \setminus V(L)$
into two sets $M,N$, where for every vertex in $M$ its set of neighbours in $V(L)$ is local
with respect to $L$, and for every vertex in $N$, its set of neighbours in $V(L)$ resolves $L$.
\\
\\
(1) {\em If there exists $f \in N$ with a nonneighbour in $V(S_1) \cup \cdots \cup V(S_m)$
then the theorem holds.}
\\
\\
For let $f$ have a nonneighbour in $S_1$ say. Let $N_1$ be the anticomponent of $N$
containing $f$, and let $X$ be the set of all $N_1$-complete
vertices in $V(G)$. From \ref{lacFnbrs} applied in the complement, it follows that
$X$ resolves $L$. Since $f$ has a nonneighbour in $V(S_1)$, there is a vertex $u$ of $S_1$
not in $X$. Let $U$ be the component of $V(G) \setminus (X \cup N)$ containing $u$.
We claim that $U$ is disjoint from $V(L) \setminus V(S_1)$, and no vertex in
$V(S_2) \cup \cdots \cup V(S_m)$ has a neighbour in $U$. For suppose not; then there
is a path $P$ say in $G$, from $V(S_1)$ to $V(L) \setminus V(S_1)$, with
$(X \cup N) \cap V(P) \subseteq V(S_2) \cup \cdots \cup V(S_m) $; choose such a path minimal.
It follows that no internal vertex of $P$ is in $V(L)$ or in $X \cup N$; and since $X$ meets
every edge between $V(S_1)$ and $V(L) \setminus V(S_1)$,
and there are no edges between  $V(S_1)$ and $V(S_2) \cup \cdots \cup V(S_m)$,
it follows that $P^*$ is nonempty. Now no vertex of $P^*$ is in $N$,
since $N \subseteq N_1 \cup X$; and so there is a component $M_1$ of $M$ including $P^*$.
From \ref{lacFnbrs}, the set of attachments of $M_1$ in $V(L)$ is local with respect to
$L$. Since it has an attachment in $V(S_1)$ it therefore has none in
$V(S_2) \cup \cdots \cup V(S_m)$. But the ends of $P$ are attachments of $M_1$, they are
nonadjacent, and one is in $V(S_1)$ and the other is not, a contradiction. This proves
that $U$ is disjoint from $V(L) \setminus V(S_1)$. Let $X'$ be the set of vertices in $X$ with
neighbours in $U$, and let $V = V(G) \setminus (U \cup N_1 \cup X')$. Then $V$ is nonempty
because $V(S_2) \subseteq V$; and so $U \cup V, N_1 \cup X'$ is a skew partition of $G$.
Since there is a vertex of $S_2$ in $X$ (because $X$ resolves $L$), and this vertex is in
$V$, we deduce that the skew partition is loose, and hence by \ref{geteven}
$G$ admits a balanced skew partition.  This proves (1).

\bigskip

From (1) we may assume that $N$ is complete to $V(S_1) \cup \cdots \cup V(S_m)$, and by
taking complements, that $M$ is anticomplete to $V(T_1) \cup \cdots \cup V(T_n)$.
\\
\\
(2) {\em If $M,N$ are both nonempty then the theorem holds.}
\\
\\
For let $M_1$ be a component of $M$, and $N_1$ an anticomponent of $N$. By taking complements
we may assume that there is a nonedge between $M_1$ and $N_1$. Let $X$ be the set
of all $N_1$-complete vertices in $G$. Since the set of attachments of $M_1$ in $V(L)$ is local
by \ref{lacFnbrs}, and since it has no attachments in $V(T_1) \cup \cdots \cup V(T_n)$, we
may assume that all its attachments are in $V(S_1)$. Let
$V = V(G) \setminus (M_1 \cup N_1 \cup V(S_1))$. Since $V(S_1) \subseteq X$, it follows
that $(M_1 \cup V, N_1 \cup V(S_1))$ is a skew partition of $G$, and since there are
$N_1$-complete vertices with no neighbours in $M_1$ (for instance, any vertex of $V(S_2)$), the
skew partition is loose, and by \ref{geteven} $G$ admits a balanced skew partition. This proves (2).
\\
\\
(3) {\em If $M, N$ are both empty then the theorem holds}.
\\
\\
For then by \ref{knottype}, we may assume that for $1 \le j \le n$ all $Q_j$-antirungs
have length 1. If $|V(S_1)| >2$, then $(V(S_1), V(L) \setminus V(S_1))$ is a 2-join of $G$; for
every vertex in $V(T_1) \cup \cdots \cup V(T_n)$ is either complete to $A_1$ and
anticomplete to $B_1 \cup C_1$, or complete to $B_1$ and anticomplete to $A_1 \cup C_1$ (since
all the antirungs have length 1). So we may assume that each $S_i$ has only two vertices.
In particular, every $S_i$-rung has length 1, so by taking complements the same argument shows
that we may assume every $V(T_j)$ has only two vertices. But then $G$ is a bicograph and
the theorem holds. This proves (3).

\bigskip

From (2) and (3), and taking complements if necessary,  we may assume that $N$ is empty
and $M$ is nonempty. For
$1 \le i \le m$ let $M_i$ be the union of the components of $M$ that have an attachment in
$V(S_i)$, and let $M_0$ be the union of the components of $M$ that have no attachments in
$V(L)$. Then $M_0,M_1,\ldots,M_n$ are pairwise disjoint and have union $M$. If $M_0$ is
nonempty then $G$ is not connected, and therefore either it admits a balanced skew partition, or
$|V(G)| \le 4$ and $G$ is bipartite, so we may assume that $M_0$ is empty. Since
$M$ is nonempty we may assume that $M_1$ is nonempty. We recall that $T_1 = (X_1,Z_1,Y_1)$;
suppose that $z \in Z_1$.  Then $z$ is complete to $V(S_1)$, and hence if we define
$V = V(G) \setminus M_1 \cup V(S_1) \cup \{z\}$, then $(M_1 \cup V, V(S_1) \cup \{z\})$ is
a skew partition of $G$, and by \ref{singleton} $G$ admits a balanced skew partition. So we
may assume that $Z_1$ is empty, and similarly every $Z_j$ is empty. Then
$(M_1 \cup V(S_1), V(G) \setminus (M_1 \cup V(S_1))$ is a 2-join of $G$.
This proves \ref{bicographs0}. \bbox

It is convenient to combine three earlier results as follows.

\begin{thm}\label{bicographs}
Let $G$ be a Berge graph, such that there is an appearance of $K_4$ in $G$.
Then either one of $G,\overline{G}$ is a line graph, or $G$ is a bicograph,
or one of $G,\overline{G}$ admits a 2-join, or $G$ admits a balanced skew partition.
\end{thm}
\Proof
This is immediate from \ref{bicographs0},\ref{linegraph} and \ref{linegraph2}.\bbox

\section{The even prism}

We proved that if $G$ is Berge and contains a nondegenerate appearance of $K_4$, then
either $G$ is a line graph, or it admits a 2-join or a balanced skew partition.
Now we want to show that the same conclusion holds for Berge
graphs that contain an ``even'' prism.  (Incidentally, the results of this section are
independent of those in the previous one; these two sections could be in either order.)
We begin with some results about prisms in general.

For $i = 1,2,3$ let $a_i \d R_i\d b_i$ be a path in $G$, so that these three paths form a prism $K$
with triangles $\{a_1,a_2,a_3\}$ and $\{b_1,b_2,b_3\}$.
A subset $X \subseteq V(G)$ {\em saturates} the prism if
at least two vertices of each triangle belong to $X$; and a vertex is {\em major} with
respect to the prism if its neighbour set saturates it. A subset $X \subseteq V(K)$ is
{\em local} with respect to the prism if either $X \subseteq V(R_i)$ for some $i$, or
$X$ is a subset of one of the triangles.
By \ref{prismparity}, the three paths $R_1,R_2,R_3$ all have lengths of the same parity.
A prism is {\em even} if the three paths $R_1,R_2,R_3$ have even length, and {\em odd} otherwise.

For our purposes, even prisms are easier than odd ones. All the prisms contained in a
degenerate appearance of $K_4$ are odd, so if we succeed in growing an even prism to become
an appearance of $K_4$, this appearance is guaranteed to be nondegenerate.

\begin{thm} \label{prismjumps}
Let $R_1,R_2,R_3$ form a prism $K$ in a Berge graph $G$, with triangles $\{a_1,a_2,a_3\}$ and
$\{b_1,b_2,b_3\}$, where each $R_i$ has ends $a_i$ and $b_i$. Let
$F \subseteq V(G) \setminus V(K)$ be connected, such that its set of attachments in $K$
is not local. Assume no vertex in $F$ is major with respect to $K$. Then there is a path
$f_1 \d \cdots\d  f_n$ in $F$ with $n \ge 1$, such that (up to symmetry) either:
\begin{enumerate}
\item $f_1$ has two adjacent neighbours in $R_1$, and $f_n$ has two adjacent neighbours in
$R_2$, and there are no other edges between $\{f_1,\ldots,f_n\}$ and $V(K)$, and (therefore)
$G$ has an induced subgraph which is the line graph of a bipartite subdivision of $K_4$, or
\item $n \ge 2$, $f_1$ is adjacent to $a_1,a_2,a_3$, and $f_n$ is adjacent to $b_1,b_2,b_3$, and there are
no other edges between $\{f_1,\ldots,f_n\}$ and $V(K)$, or
\item $n \ge 2$, $f_1$ is adjacent to $a_1,a_2$, and $f_n$ is adjacent to $b_1,b_2$, and there are
no other edges between $\{f_1,\ldots,f_n\}$ and $V(K)$, or
\item $f_1$ is adjacent to $a_1,a_2$, and $f_n$ has at least one neighbour in $R_3 \setminus a_3$,
and there are no other edges between $\{f_1,\ldots,f_n\}$ and $V(K) \setminus a_3.$
\end{enumerate}
\end{thm}
\Proof
We may assume that $F$ is minimal such that it is connected and its set of attachments in $K$
is not local. Let $X$ be the set of attachments of $F$ in $K$. For $1 \le i \le 3$,
if $X \cap V(R_i) \not = \emptyset$, let $c_i$ and $d_i$ be the vertices of $R_i$ in $X$
closest (in $R_i$) to $a_i$ and to $b_i$ respectively, and let $C_i, D_i$ be the subpaths
of $R_i$ between $a_i$ and $c_i$, and between $d_i$ and $b_i$ respectively. Let
$A =  \{a_1,a_2,a_3\}$ and $B = \{b_1,b_2,b_3\}$.

We claim that some two-element subset of $X$ is not local. For since $X \not \subseteq B$
we may assume that $c_1$ exists and $c_1 \not = b_1$. Since $X \not \subseteq V(R_1)$, we may
assume $d_2$ exists. If $d_2 \not = a_2$ then $\{c_1,d_2\}$ is the desired subset; so we may
assume $d_2 = a_2$, and similarly $d_3 = a_3$ if $d_3$ exists. Since $X  \not \subseteq A$,
it follws that $d_1 \not = a_1$, and then $\{a_2,d_1\}$ is the desired subset.
So some two-element subset $\{x_1,x_2\}$ of $X$ is not local.
Consequently $x_1,x_2$ are not adjacent.  From the minimality of $F$, there is a path with
vertices $x_1,f_1,\ldots,f_n,x_2$ so that $F = \{f_1,\ldots,f_n\}$.
\\
\\
(1) {\em If $n = 1$ then the theorem holds.}
\\
\\
For assume $n = 1 $; then $F = \{f_1\}$.
Since $X$ is not local it meets at least two of the paths; suppose it only meets
$R_1$ and $R_2$. Suppose that $c_1 = d_1$. Then we may assume that $c_1 \not \in A$
and $c_2 \not = b_2$, by exchanging
$A$ and $B$ if necessary; but then $c_1$ can be linked onto the triangle $A$, via the paths
$c_1\d C_1\d a_1$, $c_1\d f_1\d c_2\d C_2\d a_2$, and $c_1\d D_1\d b_1\d b_3\d R_3\d a_3$, contrary to \ref{trianglev}, since
$f$ has at most one neighbour in $A$. So $c_1$ is different from $d_1$, and similarly
$c_2$ is different from $d_2$ (and in particular, $c_2 \not = b_2$). Suppose that
$c_1$ is nonadjacent to $d_1$. Then since $f_1$ is not major, we may assume it has at most one
neighbour in $A$, by exchanging $A$ and $B$ if necessary; but it can be linked onto $A$, via
$f_1\d c_1\d C_1\d a_1$, $f_1\d c_2\d C_2\d a_2$ and $f_1\d d_1\d D_1\d b_1\d b_3\d R_3\d a_3$, contrary to \ref{trianglev}.
So $c_1,d_1$ are adjacent, and similarly so are $c_2,d_2$, but then statement 1 of the theorem holds.
So we may assume that $X$ meets all three of $R_1,R_2,R_3$. Since $f_1$ is not major, we may
assume that it has at most one neighbour in $A$, by exchanging $A$ and $B$ if necessary, and
therefore cannot be linked onto $A$. Since it has neighbours in all three of $R_1,R_2,R_3$,
it follows that for at least two of these paths, the only neighbour of $f_1$ in this path
is in $B$. We may assume therefore that $c_1 = b_1$ and $c_2 = b_2$. Since $X$ is not local,
$c_3 \not = b_3$; but then statement 4 of the theorem holds. This proves (1).

\bigskip

We may therefore assume that $n \ge 2$. Let $X_1$ be the set of attachments of $F\setminus f_1$,
and $X_2$ the set of attachments of $F\setminus f_n$. From the minimality of $F$, both $X_1$ and $X_2$
are local. Moreover, $X = X_1 \cup X_2$, and for $2 \le i \le n-1$, every neighbour of
$f_i$ in $K$ belongs to $X_1 \cap X_2$.
\\
\\
(2) {\em If $X_1 \subseteq A$ and $X_2 \subseteq V(R_1)$ then the theorem holds.}
\\
\\
For then $f_1$ has at least one neighbour in $R_1\setminus a_1$, and $f_n$ is adjacent to at least one
of $a_2,a_3$, and there are no other edges between $F$ and $V(K)\setminus a_1$. If $f_n$ is
adjacent to both $a_2,a_3$ then statement 4 of the theorem holds, so we assume it is not adjacent
to $a_3$. But then $a_2$ can be linked onto the triangle $B$, via $a_2\d f_n \d f_{n-1} \c f_1\d d_1\d D_1\d b_1$,
$a_2\d R_2\d b_2$, $a_2\d a_3\d R_3\d b_3$, contrary to \ref{trianglev}. This proves (2).

\bigskip

From (2), since both $X_1$ and $X_2$ are local, we may assume that either
$X_1 \subseteq A$ and $X_2 \subseteq B$, or $X_1 \subseteq V(R_2)$ and $X_2 \subseteq V(R_1)$.
In either case $X_1 \cap X_2 = \emptyset$, so none of $f_2,\ldots,f_{n-1}$ has any neighbours
in $V(K)$. Therefore $X_1$ is the set of neighbours of $f_n$ in $V(K)$, and $X_2$ is the set of neighbours
of $f_1$ in $V(K)$.
\\
\\
(3) {\em If $X_1 \subseteq A$ and $X_2 \subseteq B$ then the theorem holds.}
\\
\\
For then we may assume that $f_n$ is adjacent to $a_1$ and $f_1$ to $b_2$. Suppose first that
$n$ has the same parity as the length of $R_1$. Since $a_2\d R_2\d b_2\d f_1\d \cdots\d f_n\d a_2$ is not an
odd hole, it follows that $f_n$ is not adjacent to $a_2$, and similarly $f_1$ is not adjacent
to $b_1$. Since $a_3\d R_3\d b_3\d b_2\d f_1\c f_n\d a_1\d a_3$ is not an odd hole, either $f_n$ is
adjacent to $a_3$ or $f_1$ to $b_3$, and not both, as we saw before. But then statement 4 of the
theorem holds. Now suppose that $n$ has different parity from the length of $R_1$. Since
$a_1\d a_2\d R_2\d b_2\d f_1\c f_n\d a_1$ is not an odd hole, $f_n$ is adjacent to $a_2$, and similarly
$f_1$ to $b_1$.  If there are no more edges between $F$ and $V(K)$ then statement 3 of the theorem
holds, so we may assume that $f_n$ is adjacent to $a_3$. By the same argument as before it follows
that $f_1$ is adjacent to $b_3$, and then statement 2 of the theorem holds. This proves (3).

\bigskip

From (2) and (3) we may assume that $X_1 \subseteq V(R_2)$ and $X_2 \subseteq V(R_1)$.
So $f_1$ is adjacent to the vertices of $R_1$ that are in $X$, and $f_n$ to those of $R_2$ in $X$.
If $c_1 = d_1$, then from the symmetry we may assume that $c_1 \not = a_1$, and
$c_2 \not = b_2$; but then $c_1$ can be linked onto $A$, via $c_1\d C_1\d a_1$,
$c_1\d f_1 \c f_n\d c_2\d C_2\d a_2$, $c_1\d D_1\d b_1\d b_3\d R_3\d a_3$, contrary to \ref{trianglev}.
So $c_1 \not = d_1$ and similarly $c_2 \not = d_2$; and in particular $c_2 \not = b_2$. If
$c_1,d_1$ are nonadjacent, then $f_1$ can be linked onto $A$ via
$f_1\d c_1\d C_1\d a_1$, $f_1 \c f_n\d c_2\d C_2\d a_2$, $f_1\d d_1\d D_1\d b_1\d b_3\d R_3\d a_3$;
but $f_1$ has at most one neighbour in $A$ (because $n \ge 2$), contrary to \ref{trianglev}.
So $c_1,d_1$ are adjacent, and similarly so are $c_2,d_2$; but then statement 1
of the theorem holds. This proves \ref{prismjumps}.\bbox

\begin{thm}\label{K4jump}
Let $R_1,R_2,R_3, K, F$ be as in \ref{prismjumps}, and suppose that
\ref{prismjumps}.1 holds. Then either $R_1$ and $R_2$ both have length 1,
or there is a nondegenerate appearance of $K_4$ in $G$.
\end{thm}
\Proof
For let $f_1 \c f_n$ be a path in $F$ so that
$f_1$ has two adjacent neighbours in $R_1$, and $f_n$ has two adjacent neighbours in
$P_2$, and there are no other edges between $\{f_1,\ldots,f_n\}$ and $V(K)$. Then
$G |(V(K) \cup \{f_1,\ldots,f_n\}$ is a line graph of a bipartite subdivision
of $K_4$. We may assume it is degenerate. Hence the prism is odd,
for all prisms contained in a degenerate appearance of $K_4$ are odd. So $R_3$ is odd,
and therefore so is the path $f_1\d \cdots\d f_n$, and the other four ``rungs'' of this
line graph have length 0. In particular, $R_1$ and $R_2$ both have length 1.
This proves \ref{K4jump}. \bbox

There is also a tighter version of \ref{prismjumps}, the following.

\begin{thm}\label{jumpfromv}
Let $G$ be a Berge graph, such that there is no nondegenerate appearance of $K_4$ in $G$.
Let $R_1,R_2,R_3$ form a prism $K$ in $G$, with triangles $\{a_1,a_2,a_3\}$ and
$\{b_1,b_2,b_3\}$, where each $R_i$ has ends $a_i$ and $b_i$. Let
$F \subseteq V(G) \setminus V(K)$ be connected, such that no vertex in $F$ is major with respect to $K$.
Let $x_1$ be an attachment of $F$ in the interior of $R_1$, and assume that there is another attachment
$x_2$ of $F$ not in $R_1$. Then there is a path $f_1 \c f_n$ in $F$ such that (up to the symmetry between
$A$ and $B$)
$f_1$ is adjacent to $a_2,a_3$, and $f_n$ has at least one neighbour in $R_1 \setminus a_1$,
and there are no other edges between $\{f_1,\ldots,f_n\}$ and $V(K) \setminus a_1.$

\end{thm}
\Proof
We may assume $F$ is minimal such that it is connected, $x_1$ is one of its attachments, and it has
some attachment $x_2$ in $R_2 \cup R_3$.
Hence there is a path $x_2\d v_1 \c v_m\d x_1$ where $F = \{v_1\l v_m\}$. By \ref{prismjumps},
there is a subpath $f_1 \c f_n$ of $v_1 \c v_m$ such that one of \ref{prismjumps}.1-4 holds.
From the minimality of $F$, $v_1$ is the only vertex of $F$ with a neighbour in $V(R_2)\cup V(R_3)$,
and in particular, at most one vertex of
$f_1 \c f_n$ has a neighbour in $V(R_2)\cup V(R_3)$. We deduce that $f_1 \c f_n$ does not
satisfy \ref{prismjumps}.2 or \ref{prismjumps}.3. Suppose it satisfies \ref{prismjumps}.1. By
\ref{K4jump} the path $f_1 \c f_n$ joins two of $R_1,R_2,R_3$ that are both of length 1, and therefore $n$ is even.
Since $R_1$ has length $\ge 2$ (because $x_1$ is in its interior) it follows that
$f_1,f_n$ are distinct vertices of $F$ both with neighbours in $V(R_2)\cup V(R_3)$, a contradiction. So
$f_1 \c f_n$ satisfies \ref{prismjumps}.4, and therefore we may assume that for some $i$ with $1 \le i \le 3$,
$f_1$ is adjacent to the two vertices in $A \setminus \{a_i\}$, and $f_n$ has at least one neighbour in
$R_i \setminus a_i$, and there are no other edges between $\{f_1,\ldots,f_n\}$ and $V(K) \setminus a_i.$
Suppose first that $i >1$, $i = 2$ say.
Then both $f_1,f_n$ have neighbours in $V(R_2) \cup V(R_3)$, and so from the minimality of
$F$ it follows that $n = 1$ and $f_1 = v_1$.  But then $f_1$ can be linked onto the triangle $B$, via
the path between $f_1$ and $b_1$ with interior in $\{v_2\l v_m\} \cup (V(R_1)\setminus a_1)$,
the path between $f_1$ and $b_2$ with interior in $V(R_2) \setminus a_2$, and the path $f_1\d a_3\d R_3\d b_3$,
contrary to \ref{trianglev}. Hence $i = 1$, and the theorem is satisfied. This proves \ref{jumpfromv}.\bbox

Another useful corollary of \ref{prismjumps} is the following.

\begin{thm}\label{legaljump}
Let $G$ be Berge, such that there is no nondegenerate appearance of $K_4$ in $G$.
Let $R_1,R_2,R_3$ form a prism $K$ in a Berge graph $G$, with triangles $\{a_1,a_2,a_3\}$ and
$\{b_1,b_2,b_3\}$, where each $R_i$ has ends $a_i$ and $b_i$. Let
$F \subseteq V(G) \setminus V(K)$ be connected, such that if the prism is even
then no vertex in $F$ is major with respect to $K$.
Assume that the set of attachments of $F$ in $K$ is not local, but none are in $V(R_3)$.
Then $|F| \ge 2$, and the set of attachments of $F$ in $K$ is precisely $\{a_1,b_1,a_2,b_2\}$.
\end{thm}
\Proof If there is a major vertex $v \in F$, then since it has no neighbours in $R_3$, it is adjacent
to $a_1$ and $b_2$, and since $v\d a_1\d a_3\d R_3\d b_3\d b_2\d v$ is a hole, it follows that the prism
is even, contrary to the hypothesis. So there is no major vertex in $F$.
By \ref{jumpfromv} no internal vertex of $R_1$ or $R_2$ is an attachment of $F$.
By \ref{prismjumps}, there is a path in $F$ satisfying one of  \ref{prismjumps}.1-4; and since it has
no attachments in $R_3$, it must satisfy \ref{prismjumps}.1 or \ref{prismjumps}.3, and
in either case $a_1,b_1,a_2,b_2$ are all attachments of $F$. Since no vertex in $F$ is major it
follows that $|F| \ge 2$. This proves \ref{legaljump}. \bbox

The next result is a close relative of \ref{overshadowed}.

\begin{thm}\label{evenbig}
Let $G$ be Berge, such that there is no nondegenerate appearance of $K_4$ in $G$.
If there is an even prism $K$ in a Berge graph $G$, such that some vertex of $G$ is major with
respect to $K$, then $G$ admits a balanced skew partition.
\end{thm}
\Proof
Any prism has six vertices of degree 3, called {\em triangle-vertices}; choose a prism $K$ and a
nonempty co-connected set $Y \subseteq V(G) \setminus V(K)$, such that
every vertex in $Y$ is major with respect to the prism, and as few triangle-vertices of $K$ are
$Y$-complete as possible.
Let the paths $a_i\d R_i \d b_i (i = 1,2,3)$ form $K$, where $\{a_1,a_2,a_3\},\{b_1,b_2,b_3\}$
are its triangles. We may assume that $Y$ is maximal with the given property.
Let $X$ be the set of all $Y$-complete vertices in $G$. By \ref{prismbig}, $X$ saturates $K$.
Consequently there is one of $R_1,R_2,R_3$ with both ends in $X$, say $R_1$.  Let
$X_0 = X \setminus V(K)$ and $X_1 = \{a_1,b_1\}.$
\\
\\
(1) {\em If $F \subseteq V(G)$
is connected and some vertex of $V(R_1^*)$ has a neighbour in $F$, and so does some vertex of
$V(R_2) \cup V(R_3)$, then $F\cap (X_0 \cup X_1 \cup Y) $ is nonempty.}
\\
\\
Suppose for a contradiction that some $F$ exists not satisfying (1), and choose it minimal.
Hence $G|F$ is a path, disjoint from $K$. Consequently $F \cap X = \emptyset$. Suppose some vertex in $v \in F$ is
major with respect to $K$. Then since $v \not \in X$ it follows that $v$ has a nonneighbour in $Y$,
and so $Y \cup \{v\}$ is anticonnected; the maximality of $Y$ therefore implies that $v \in Y$, and hence
$F \cap Y \not = \emptyset$ and the claim holds. So we may assume that no vertex
in $F$ is major. Let $x_1$ be an attachment of $F$ in $R_1^*$. By \ref{jumpfromv} we may assume that there is a path
$f_1\c f_n$ in $F$ such that $f_1$ is adjacent to $a_2,a_3$, and $f_n$ has neighbours in $R_1 \setminus a_1$,
and $f_1a_2,f_1a_3$ are the only edges between $\{f_1\l f_n\}$ and $V(R_2) \cup V(R_3)$
Now there is a path $R$ from $f_1$ to
$b_1$ with interior in $\{f_2,\ldots,f_n\} \cup V(R_1\setminus a_1)$, and hence $R,R_2,R_3$ form a prism $K'$
say.  By \ref{prismrung}, every vertex in $Y$ is major with respect to $K'$, and since $a_1$
is $Y$-complete and $f_1$ is not, it follows that the number of $Y$-complete triangle-vertices in $K'$
is smaller than the number in $K$, a contradiction.  This proves (1).

\bigskip

It follows from (1) that there is a partition of $V(G) \setminus (X_0 \cup X_1 \cup Y)$ into two sets
$L$ and $M$ say, where there is no edge between $L$ and $M$, and $V(R_1^*) \subseteq L$ and
$V(R_2) \cup V(R_3) \subseteq M$.
So $(L \cup M, X_0 \cup X_1 \cup Y)$ is a skew partition of $G$. Since at least two vertices of $A$ are in $X$
and only one is in $X_1$, there is a vertex of $X$ in $M$, and so the skew partition is loose.
By \ref{geteven} the result follows. This proves \ref{evenbig}.\bbox

The main result of this section is \ref{summary}.4, which we restate.

\begin{thm}\label{evenprism}
Let $G$ be a Berge graph, such that there is no nondegenerate appearance of $K_4$ in $G$.
If $G$ contains an even prism, then either $G$ is an even prism with $|V(G)| = 9$, or
$G$ admits a 2-join or a balanced skew partition.
\end{thm}
\Proof
Since $G$ contains an even prism, we can choose in $G$ a collection of nine sets
\\
\[\begin{array}{ccc}
A_1 & C_1 & B_1\\
A_2 & C_2 & B_2\\
A_3 & C_3 & B_3
\end{array}\]
\newline
with the following properties:
\begin{itemize}
\item all these sets are nonempty and pairwise disjoint
\item for $1 \le i < j \le 3$ $A_i$ is complete to $A_j$ and $B_i$ is complete to $B_j$, and there are
no other edges between $A_i \cup B_i \cup C_i$ and $A_j \cup B_j \cup C_j$
\item for $1 \le i \le 3$, every vertex of $A_i \cup B_i \cup C_i$ belongs to a path between $A_i$ and
$B_i$ with interior in $C_i$
\item some path between $A_1$ and $B_1$ with interior in $C_1$ is even.
\end{itemize}

We call this collection of nine sets a {\em hyperprism}. Let $H$ be the subgraph of $G$ induced on the
union of the nine sets. Choose the hyperprism with $V(H)$ maximal. For $1 \le i \le 3$, a path from $A_i$
to $B_i$ with interior in $C_i$ is called an $i$-rung. Let us write $S_i = A_i \cup B_i \cup C_i$ for
$1 \le i \le 3$, and $A = A_1 \cup A_2 \cup A_3$, and $B = B_1 \cup B_2 \cup B_3$.
\\
\\
(1) {\em For $1 \le i \le 3$, all $i$-rungs have even length.}
\\
\\
For we are given that some 1-rung $R_1$ say has even length. Let $R_2$ be an 2-rung; then the union of
$R_1$ and $R_2$ induces a hole, and so $R_2$ is even. Hence every 2- or 3-rung is even, and hence so is
every 1-rung. This proves (1).

\bigskip

A subset $X \subseteq V(H)$ is {\em local} (with respect to the hyperprism) if $X$ is a subset of
one of $S_1,S_2,S_3,A$ or $B$.
\\
\\
(2) {\em We may assume that
for every connected subset $F$ of $V(G) \setminus V(H)$, its set of attachments in $H$ is local.}
\\
\\
For suppose not. Choose $F$ minimal, and let $X$ be the set of attachments of $F$ in $H$.
Suppose first that there exists $x_1 \in X \cap C_1$. Since $X$ is not local, we may assume that
there exists $x_2 \in X \cap S_2$. For $i = 1,2,3$ choose an $i$-rung
$R_i$ with ends $a_i \in A_i$ and $b_i \in B_i$, so that for $i = 1,2$, $x_i \in V(R_i)$. Then
$R_1,R_2,R_3$ form an even prism $K$ say. By \ref{evenbig} we may assume no vertex in $F$ is major with
respect to $K$; so by \ref{jumpfromv}, we may assume that there is a path $f_1\c f_n$ in $F$
such that $f_1$ is adjacent to $a_2,a_3$, and $f_n$ has at least one neighbour in $R_1 \setminus a_1$,
and there are no other edges between $\{f_1,\ldots,f_n\}$ and $V(K) \setminus a_1$.
From the minimality of $F$ it follows that $F = \{f_1\l f_n\}$.
Since this holds for all choices of $R_3$ it follows that $f_1$ is complete to $A_3$ and there are
no edges between $\{f_1,\ldots,f_n\}$ and $B_3 \cup C_3$. Since $a_3 \in X$ the same conclusion
follows for all choices of $R_2$, and so $f_1$ is complete to $A_2$ and there are
no edges between $\{f_1,\ldots,f_n\}$ and $B_2 \cup C_2$. But then we can add $f_1$ to $A_1$
and $\{f_2,\ldots,f_n\}$ to $C_1$, contradicting the maximality of the hyperprism.
It follows that $X \cap C_1 = \emptyset$, and similarly $X \cap C_2, X \cap C_3 = \emptyset$.
We claim there is a 2-element subset of $X$ which is also not local.
For we may assume $X \cap A_1 \not = \emptyset$; and hence if $X$ meets $B_2$ or $B_3$
our claim holds. If not, then it meets $B_1$ (since it is not a subset of $A$) and meets $A_2 \cup A_3$
(since it is not a subset of $S_1$), and again the claim holds. So there is a subset $\{x_1,x_2\}$ of $X$
which is not local. We may assume that $x_1 \in A_1$ and $x_2 \in B_2$. From the minimality of $F$,
there is a path $x_1\d f_1 \c f_n\d x_2$ with $F = \{f_1,\ldots,f_n\}$.

Suppose first that
$n$ is even. For any 3-rung $R_3$ with ends $a_3 \in A_3$ and $b_3 \in B_3$,
$x_1\d f_1\c f_n\d x_2\d b_3\d R_3\d a_3\d x_1$ is not an odd hole, and so some vertex of $R_3$ is in $X$.
Since $X \cap C_3 = \emptyset$, and $a_3$ has no neighbour in $\{f_2,\ldots,f_n\}$ from the minimality
of $F$, and similarly $b_3$ has no neighbour in $\{f_1,\ldots,f_{n-1}\}$, it follows that either
$f_1$ is adjacent to $a_3$, or $f_n$ to $b_3$ (and not both, since otherwise
$f_1\c f_n\d b_3\d R_3\d a_3$ is an odd hole). From the symmetry we may assume that $f_n$ is adjacent to $b_3$.
By exchanging $S_2$ and $S_3$ it follows that for every 2-rung with ends $a_2 \in A_2$ and $b_2 \in B_2$,
either $f_1$ is adjacent to $a_2$ or $f_n$ to $b_2$, and not both. Suppose that $f_n$ is complete to
$B_2 \cup B_3$; then $f_1$ has no neighbours in $S_2 \cup S_3$, and we can add $f_n$ to $B_1$
 and $f_1,\ldots,f_{n-1}$ to $C_1$, contrary to the maximality of the hyperprism. So $f_n$ is not
complete to $B_2 \cup B_3$, and hence $f_1$ has a neighbour in one of $A_2,A_3$, say $A_3$; and by
exchanging $S_1$ and $S_2$ it follows that for every 1-rung with ends $a_1 \in A_1$ and $b_1 \in B_1$,
either $f_1$ is adjacent to $a_1$ or $f_n$ to $b_1$ and not both. In particular, $f_1$ has no neighbours in
$B$ and $f_n$ has none in $A$. For $i = 1,2,3$ let $A_i'$ be the set of neighbours of $f_1$ in $A_i$,
and let $A_i'' = A_i \setminus A_i'$; let $B_i''$ be the set of neighbours of $f_n$ in $B_i$, and let
$B_i' = B_i \setminus B_i''$. We have shown so far that every $i$-rung is either between $A_i'$ and $B_i'$ or
between $A_i''$ and $B_i''$. Let $C_i'$ be the union of the interiors of the $i$-rungs between $A_i'$ and $B_i'$,
and $C_i''$ the union of the interiors of the $i$-rungs between $A_i''$ and $B_i''$. We observe
that $C_i = C_i' \cup C_i''$. Moreover, $C_i' \cap C_i'' = \emptyset$, for otherwise there would be an $i$-rung
between $A_i'$ and $B_i''$. For the same reason there are no edges between $A_i' \cup C_i'$ and $C_i'' \cup B_i''$,
and no edges between $A_i'' \cup C_i''$ and $C_i' \cup B_i'$. We claim that $A_i'$ is complete to $A_i''$.
For if not, let $R''$ be an $i$-rung with ends $a'' \in A_i''$ and $b'' \in B_i''$, and let $a' \in A_i'$ be
nonadjacent to $a''$. Since we have seen that $f_1$ has neighbours in at least two of $A_1,A_2,A_3$, we
may choose $a \in A_j'$ for some $j \not = i$. Then $a\d f_1\d \cdots\d f_n\d b''\d R''\d a''\d a$ is an odd hole,
a contradiction. So $A_i'$ is complete to $A_i''$ for each $i$, and similarly $B_i'$ is complete to
$B_i''$ for each $i$. We showed already that we may assume that $A_1',A_2'',A_3',A_3''$ are all nonempty.
But then the nine sets

\[\begin{array}{ccc} A_1' & C_1' & B_1'\\ A_2'\cup A_3'& C_2'\cup C_3' & B_2'\cup B_3' \\
A_1''\cup A_2''\cup A_3''\cup \{f_1\} & C_1'' \cup C_2''\cup C_3''\cup \{f_2,\ldots,f_n\} &  B_1''\cup B_2''\cup B_3''
\end{array}\]
form a hyperprism, contrary to the maximality of $V(H)$. This completes the argument when
$n$ is even.

Now assume $n$ is odd. $f_1$ has a neighbour $a_1$ say in $A_1$; let $R_1$ be a 1-rung with ends
$a_1$ and $b_1$ say. Similarly let $R_2$ be a 2-rung with ends $a_2$ and $b_2$, where $b_2 \in B_2$
is adjacent to $f_n$. Since $a_1\d f_1\d \cdots\d f_n\d b_2\d b_1\d R_1\d a_1$ is not an odd hole, it follows
that $b_1 \in X$, and similarly $a_2 \in X$. From the minimality of $F$, one of $b_1,a_2$
is adjacent to $f_1$ and the other to $f_n$, and neither has any more neighbours in $F$.
Suppose that $f_n$ is not adjacent to $b_1$; so $f_1$ is adjacent to $b_1$, and $n \ge 2$,
and $f_n$ is adjacent to $a_2$. But then $b_1\d f_1\d \cdots\d f_n\d b_2\d b_1$ is an odd hole, a
contradiction. This proves that $f_n$ is adjacent to $b_1$ and $f_1$ to $a_2$. Hence for all
$1 \le i \le 3$, and for every $i$-rung with ends $a \in A$ and $b \in B$, $a \in X$ if and only
if $b \in X$, and if so then $f_1$ is adjacent to $a$ and $f_n$ to $b$. Consequently, for
every vertex in $X \cap A$, $f_1$ is its unique neighbour in $F$, and for every vertex in $X\cap B$,
$f_n$ is its unique neighbour in $F$. For $1 \le i \le 3$, let
\begin{eqnarray*}
A_i' = A_i \cap X \\
B_i' = B_i \cap X \\
A_i'' = A_i \setminus X \\
B_i'' = B_i \setminus X.
\end{eqnarray*}
Let $C_i'$ be the union of the interior of the $i$-rungs between $A_i'$ and $B_i'$, and
$C_i''$ the union of the interior of the $i$-rungs between $A_i''$ and $B_i''$. We have seen
that every $i$-rung is of one of these two types, and so $C_i = C_i' \cup C_i''$. Moreover,
since there is no rung between $A_i'$ and $B_i''$, it follows that $C_i' \cap C_i'' = \emptyset$,
and there are no edges between $A_i' \cup C_i'$ and $C_i'' \cup B_i''$,
and similarly no edges between $A_i'' \cup C_i''$ and $C_i' \cup B_i'$.
We have seen that $f_1$ has neighbours in at least two of $A_1,A_2,A_3$, and $f_n$ has neighbours
in at least two of $B_1,B_2,B_3$. We claim that also $f_1$ has nonneighbours in at least two
of $A_1,A_2,A_3$, and the same for $f_n$. For suppose not, and $f_1$ is complete to $A_1 \cup A_2$
say. Then $f_n$ is complete to $B_1 \cup B_2$; by \ref{evenbig} we may assume that $n>1$, and so we
can add $f_1$ to $A_3$, $f_n$ to $B_3$ and $f_2,\ldots,f_{n-1}$ to $C_3$, contrary to the maximality
of $V(H)$. This proves that $f_1$ has nonneighbours in at least two
of $A_1,A_2,A_3$, and similarly $f_n$ has nonneighbours in at least two of $B_1,B_2,B_3$.
Let $1 \le i \le 3$; we claim that $A_i'$ is complete to $A_i''$. For we may assume that
$i = 1$; suppose that $a' \in A_1'$ and $a'' \in A_1''$ are nonadjacent, and let $R''$ be a $1$-rung
with ends $a'',b''$. Choose $a \in A_2''\cup A_3''$ and $b \in B_2' \cup B_3'$; then
$a,b$ are not adjacent since all rungs have even length, and so
$a\d a'\d f_1 \c f_n\d b\d b''\d R''\d a''\d a$
is an odd hole, a contradiction. This proves that $A_i'$ is complete to $A_i''$ for $i = 1,2,3$, and
similarly $B_i'$ is complete to $B_i''$. We have seen that we may assume that $A_1',A_2'$ are
nonempty. But then
\[\begin{array}{ccc}
A_1' & C_1' & B_1' \\
A_2'\cup A_3' & C_2'\cup C_3' & B_2'\cup B_3' \\
A_1'' \cup A_2'' \cup A_3''\cup \{f_1\}& C_1'' \cup C_2'' \cup C_3''\cup \{f_2,\ldots,f_{n-1}\}&B_1'' \cup B_2'' \cup B_3'' \cup \{f_n\}
\end{array}\]
is a hyperprism, contrary to the maximality of $V(H)$. This proves (2).

\bigskip

Suppose $F$ is a component of $V(G) \setminus V(H)$, and all its attachments are in $A$. Then
$(V(G) \setminus A, A)$ is a skew partition of $G$. We must show that $G$ admits a balanced skew partition.
Choose $b_2 \in B_2$ and $a_3 \in A_3$. Then $B_1 \cup C_1 \cup \{b_2\}$ is connected, and all vertices
in $A_1$ have neighbours in it. By \ref{balancev}, $(B_1 \cup C_1 \cup \{b_2\}, A_1)$ is balanced,
and so by \ref{shiftbalance}.1, so is $(A_1,F)$. By \ref{onepair}, $G$ admits a balanced skew partition.
So we may assume there is no such $F$, and the same for $B$.

From (2) it follows that for every component of  $V(G) \setminus V(H)$, all its attachments in $H$ are
a subset of one of $S_1,S_2,S_3$. Let $X$ be the union of $S_1$ and all components of $V(G) \setminus V(H)$
whose attachment set is a subset of $S_1$, and let $Y = V(G) \setminus X$. Then $|Y| \ge 4$, and so
either $(X,Y)$ is a 2-join in $G$, or $|X| = 3$ and both $A_1,B_1$ have one element and $X$ is
the vertex set of a path between these two vertices. We may assume the latter, and the same for
$S_2$ and $S_3$; and so $|V(G)| = 9$, and $G$ is an even prism. This proves \ref{evenprism}.\bbox

\section{Step-connected strips}

Our next target is the statement analogous to \ref{evenprism} for long odd prisms, but we need
to creep up on it in stages. (A warning: we shall not prove the exact analogue, and we don't
know if it is true. We need to permit more types of decomposition, namely
2-joins in $\overline{G}$, and M-joins.) The key idea is to start with a prism of three
paths, $R_0,R_1,R_2$, where $R_0$ has length $\ge 3$, and to grow the union of the other two
paths into a kind of strip ({\em one} strip, not two)
with a richer internal structure than we have seen hitherto,
what we call being ``step-connected''. If we expand the union of these two paths into a
maximal step-connected strip, then the remainder of the graph attaches to this structure in
ways that we can exploit. In this section we introduce step-connected strips,
and prove some preliminary lemmas about them.

Let $(A,C,B)$ be a strip in $G$.  A {\em step} is
a pair $a_1\d R_1\d b_1,a_2\d R_2\d b_2$ of rungs so that
\begin{itemize}
\item $V(R_1) \cap V(R_2) = \emptyset$
\item $a_1$ is adjacent to $a_2$, and $b_1$ to $b_2$, and there are no other edges between
$V(R_1)$ and $V(R_2)$.
\end{itemize}
The edges $a_1a_2$ and $b_1b_2$ such that there exists a step as above are called
{\em stepped} edges.
We say that the strip is {\em step-connected} if every vertex of $A \cup B \cup C$ is in a
step, and for every partition $(X,Y)$ of $A$ or of $B$
into two nonempty sets, there is a step $R_1,R_2$ so that $R_1$ has an end in $X$ and
$R_2$ has an end in $Y$. (This second condition is equivalent to requiring that the subgraph
of $G$ with vertex set $A$ and edges the stepped edges within $A$ be connected, and the same for $B$.)

Let $(A,C,B)$ be a step-connected strip in a Berge graph $G$.
A vertex $v \in  V(G) \setminus (A \cup B \cup C) $ is a {\em left-star} for the strip if
it is complete to $A$ and anticomplete to $B \cup C$, and it is a {\em right-star} if
it is complete to $B$ and anticomplete to $A \cup C$. A {\em banister} (with respect to the strip)
is a path $a\d R\d b$ of $G \setminus (A \cup B \cup C)$, such that $a$ is a left-star,
$b$ is a right-star, and there are no edges
between the interior of $R$ and $V(S)$. (Here we distinguish between $a\d R\d b$ and $b\d R\d a$;
we follow the convention that when describing a banister relative to a strip, the end which
is the left-star is listed first.) A banister can have length 1.

\begin{thm}\label{leftcon}
Let $G$ be a Berge graph, such that there is no nondegenerate appearance of $K_4$ in $G$.
Let $S = (A,C,B)$ be a step-connected strip in $G$, and let $a_0\d R_0\d b_0$ be a banister.
Suppose that $v \in V(G) \setminus V(S)$ has a neighbour in $A \cup C$, and has no neighbour
in $B$; and that $P$ is a path in $G \setminus (V(S) \cup \{a_0\})$ from $v$ to $b_0$, and there
are no edges between $P^*$ and $V(S)$. Then $v$ is a left-star.
\end{thm}
\Proof Let $F$ be a connected subset of $V(P)$,
containing $v$ and disjoint from $V(R_0)$, and with an attachment in $R_0 \setminus a_0$.
\\
\\
(1) {\em For every step $a_1\d R_1\d b_1,a_2\d R_2\d b_2$ , if $v$ has a neighbour in $R_1 \cup R_2$
then $v$ is adjacent to $a_1,a_2$ and to no other vertices of $R_1 \cup R_2$.}
\\
\\
For assume $v$ has a neighbour in $R_1$ say, and hence in $R_1\setminus b_1$.
Now $R_0,R_1,R_2$ form a prism $K$ say, and no vertex in $F$ is major with respect to $K$ since
no vertex in $F$ is adjacent to $b_1$ or $b_2$. Yet $F$ has an attachment in $R_0\setminus a_0$
and one in $R_1\setminus b_1$, so its set of attachments is not local. Since $b_1$ is not an
attachment of $F$, it follows from \ref{legaljump} that $F$ has an attachment in $R_2$; and
therefore $v$ has a neighbour in $R_2\setminus b_2$. If $v$ has any neighbours in $R_1 \cup R_2$
different from $a_1,a_2$, say a neighbour in the interior of $R_1$, then
$v$ can be linked onto the triangle $b_0,b_1,b_2$, via the paths $v\d P\d b_0$,
from $v$ to $b_1$ with interior in $R_1\setminus a_1$, and from
$v$ to $b_2$ with interior in $R_2$; but this contradicts \ref{trianglev}. This proves (1).

\bigskip

From (1) it follows that $v$ has no neighbour in $C$ (since every vertex is in a step), and
therefore $v$ has at least one neighbour in $A$; and from (1) again, $v$ has no nonneighbour
in $A$ (for otherwise we could choose the step in (1) with $v$ adjacent to $a_1$ and not
to $a_2$, since the strip is step-connected.) This proves \ref{leftcon}.\bbox

\begin{thm}\label{leftcon2}
Let $G$ be Berge, such that there is no appearance of $K_4$ in $G$.
Let $S = (A,C,B)$ be a step-connected strip in $G$, and let $a_0\d R_0\d b_0$ be a banister.
Let $v \in V(G) \setminus V(S)$ have a neighbour in $V(S)$, and be nonadjacent to $b_0$. Let $P$ be a path
in $G \setminus (V(S) \cup \{a_0\})$ from $v$ to $b_0$, and let $Q$ be a path in
$G \setminus (V(S) \cup \{b_0\})$ from $v$ to $a_0$, such there are no edges from
$P^* \cup Q^*$ to $V(S)$. Then either $v$ is $B$-complete, or $v$ is a left-star.
\end{thm}
\Proof
If $v$ has no neighbours in $B$, then by \ref{leftcon} $v$ is a left-star, so we may assume $v$
has a neighbour in $B$. Since we may assume it is not $B$-complete, there is a step
$a_1\d R_1\d b_1,a_2\d R_2\d b_2$ such that $v$ is adjacent to $b_1$ and not to $b_2$.
Let $F\subseteq V(Q)$ be connected, containing $v$ and  disjoint
from $V(R_0)$, with an attachment in $R_0 \setminus b_0$. Now $R_0,R_1,R_2$ form a prism
$K$ say, and no vertex of $F$ is major with respect to $K$ since none of them has two neighbours
in $\{b_0,b_1,b_2\}$. But there is an attachment of $F$ in $R_0 \setminus b_0$, and $b_1$ is
also an attachment of $F$, so its set of attachments is not local with respect to the prism.
By \ref{prismjumps}, one of \ref{prismjumps}.1-4 holds. Since there is no appearance
of $K_4$ in $G$,
\ref{prismjumps}.1 does not hold. Also \ref{prismjumps}.2, \ref{prismjumps}.3 do not hold,
since $v$ is the only vertex in $F$ with neighbours in $A \cup B$.  So  \ref{prismjumps}.4
holds, and therefore $F$ has an attachment in $R_2$, and so $v$ has a neighbour in $R_2$.
But then $v$ can be linked onto the triangle $\{b_0,b_1,b_2\}$, via $v\d P\d b_0$, $v\d b_1$,
and the path from $v$ to $b_2$ with interior in $R_2$, contrary to \ref{trianglev}.
This proves \ref{leftcon2}.\bbox

We remark:
\begin{thm}\label{prismrunglength}
Let $G$ be Berge, containing no even prism,
let $S = (A,C,B)$ be a step-connected strip in $G$, and let $a_0\d R_0\d b_0$ be a banister.
Then every rung of the strip has odd length.
\end{thm}
\Proof
Let $a_1\d R_1\d b_1,a_2\d R_2\d b_2$ be a step. Then these three paths form a prism, and it
is not an even prism by hypothesis. In particular $R_0$ has odd length, by \ref{prismparity}.
For any rung $a\d R\d b$, the hole $a_0\d R_0\d b_0\d b\d R\d a\d a_0$ has even length, and so $R$ is odd.
This proves \ref{prismrunglength}. \bbox

\begin{thm}\label{stepsantipath1}
Let $G$ be a Berge graph, such that there is no appearance of $K_4$ in $G$ and no even prism in $G$.
Let $S = (A,C,B)$ be a step-connected strip in $G$.
Let $F \subseteq V(G) \setminus (A \cup B \cup C)$ be connected, such that there are no edges
between $F$ and $A \cup B \cup C$.
There is no anticonnected set $Q \subseteq V(G) \setminus (A \cup B\cup C\cup F)$ such that:
\begin{itemize}
\item some right-star has a neighbour in $F$ and a nonneighbour in $Q$,
\item some vertex in $B$ has a nonneighbour in $Q$,
\item some left-star with a neighbour in $F$ is $Q$-complete,
\item every vertex in $Q$ has a neighbour in $F$,
\item every vertex in $Q$ has a neighbour in $A \cup B \cup C$, and
\item no vertex in $Q$ is a left-star.
\end{itemize}
\end{thm}
\Proof
Suppose that such a set $Q$ exists. Let $a_0$ be a left-star with a neighbour in $F$
complete to $Q$, and let $b_0$ be a right-star with a neighbour in $F$ and a
nonneighbour in $Q$. Let $R_0$ be
a path between $a_0$ and $b_0$ with interior in $F$. Hence $a_0\d R_0\d b_0$ is a banister.
By \ref{prismrunglength} $R_0$ and every rung has odd length. Since
some vertex in $B$ has a nonneighbour in $F$, there is an antipath $q_1\d \cdots\d q_n$ in $Q$
such that $q_1$ is not adjacent to $b_0$ and $q_n$ is not adjacent to some vertex in $B$.
Choose such an antipath with $n$ minimum. Let $B_1$ be the set of neighbours of $q_n$ in
$B$, and $B_2 = B \setminus B_1$. So $B_2 \not = \emptyset$, and by \ref{leftcon} it
follows that $B_1 \not = \emptyset$.
Choose a step $a_1\d R_1\d b_1,a_2\d R_2\d b_2$ with $b_1 \in B_1$ and $b_2 \in B_2$.
\\
\\
(1) {\em $n \ge 2$.}
\\
\\
For suppose $n = 1$. Then $q_1$ is adjacent to $a_0$ and to $b_1$, and not to $b_0$, so by
\ref{legaljump}, $q_1$ has a neighbour in $R_2 \setminus b_2$. Since $q_1$ also has
a neighbour in $F$, it can be linked onto the
triangle $\{b_0,b_1,b_2\}$, via a path from $q_1$ to $b_0$ with interior in $F$, the path $q_1\d b_1$,
and the path from $q_1$ to $b_2$ with interior in $R_2$, contrary to \ref{trianglev}. This proves (1).
\\
\\
(2) {\em $(A\cup B \cup C,\{b_0,q_1,\ldots,q_n\})$ is balanced.}
\\
\\
For $b_1 \in B_1$ is complete to $\{b_0,q_1,\ldots,q_n\}$ from the minimality
of $n$. But $b_1$ has no neighbour in $F$, so by \ref{balancev}, $(F,\{b_0,q_1,\ldots,q_n\})$ is
balanced. Since $F$ is connected and every vertex in $\{b_0,q_1,\ldots,q_n\}$ has a neighbour
in $F$, the claim follows from \ref{shiftbalance}.1. This proves (2).

\bigskip

Now the path $a_0\d a_2\d R_2\d b_2\d b_1$ is odd, and its ends are complete to $\{q_1,\ldots,q_n\}$; so
by (2) and \ref{RR}, there are two adjacent vertices $u,v$ in this path, both complete to
$\{q_1,\ldots,q_n\}$. Since $b_2$ is not adjacent to $q_n$, it follows that
$u,v \in \{a_0\} \cup V(R_2\setminus b_2)$.
Suppose that the hole $a_0\d R_0\d b_0\d b_2\d R_2\d a_2\d a_0$ has length $\ge 6$.
Then one of $u,v$ is nonadjacent to both $b_0,b_2$, say $v$, and hence $n$ is odd, since
$v\d b_0\d q_1\d \cdots\d q_n\d b_2\d v$ is an antihole; but
$b_1$ is adjacent to $b_0$ and $b_2$, and has no other neighbours in this hole, and is complete
to $\{q_1,\ldots,q_n\}$, contrary to \ref{hole&antipath}. So the hole has length 4, and
in particular $a_2$ is adjacent to $b_2$ and is complete to $\{q_1,\ldots,q_n\}$, and
$a_0$ is adjacent to $b_0$. Hence
$n$ is odd, because $b_1\d a_2\d b_0\d q_1\d \cdots\d q_n\d b_2\d a_0\d b_1$ is an antihole, and so
$a_2\d b_0\d q_1\d \cdots\d q_n\d b_2$ is
an odd antipath, contrary to (2). This proves \ref{stepsantipath1}.\bbox

A triple $(S,F,Q)$ is called a {\em 1-breaker} in $G$ if it satisfies the following.
\begin{itemize}
\item $S = (A,C,B)$ is a step-connected strip in $G$,
\item $F \subseteq V(G) \setminus V(S)$ is connected, such that there are no edges
between $F$ and $V(S)$, and there is a left- and right-star, both with neighbours in $F$,
\item $Q \subseteq V(G) \setminus (V(S)\cup F)$ is anticonnected,
\item some vertex in $A$ has a nonneighbour in $Q$, and so does some vertex in $B$,
\item every vertex in $Q$ has a neighbour in $F$ and a neighbour in $A \cup B \cup C$,
\item some left-star with a neighbour in $F$ is $Q$-complete,
\item no vertex in $Q$ is a left-star.
\end{itemize}

\begin{thm}\label{1-hat}
Let $G$ be a Berge graph, such that there is no appearance of $K_4$ in $G$ and no even prism in $G$.
If there is a 1-breaker in $G$ then $G$ admits a balanced skew partition.
\end{thm}
\Proof
Suppose that some 1-breaker $(S,F,Q)$ exists, and
for fixed $G$ and $S$, choose $F$ and $Q$ with $|F| + |Q|$ maximum so that all the hypotheses of the
theorem remain satisfied (possibly exchanging ``left'' and ``right'').
Let $N$ be the set of vertices of $G$ not in $F$ but with a neighbour
in $F$. Hence $Q \subseteq N$, and every left- or right-star with a neighbour in $F$ is in $N$.
Let $S = (A,C,B)$.
\\
\\
(1) {\em Every vertex in $N$ has a neighbour in $A \cup B \cup C$.}
\\
\\
For suppose $v \in V(G) \setminus F$ has a neighbour in $F$ and has none in $A \cup B \cup C$. Let
$F' = F \cup \{v\}$. Certainly $F'$ is connected and disjoint from $A \cup B \cup C$, and there
are no edges between $F'$ and $A \cup B \cup C$; and $F'$ is disjoint from $Q$ since every
vertex in $Q$ has a neighbour in $A \cup B \cup C$.  It follows that
the hypotheses of the theorem remain true, contrary to the maximality of $|F|+|Q|$. This proves (1).
\\
\\
(2) {\em There is no left- or right-star in $Q$, and every left- and right-star
with a neighbour in $F$ is $Q$-complete.}
\\
\\
For we are given that there is no left-star in $Q$. Suppose there is a right-star
with a neighbour in $F$, either in $Q$ or with a nonneighbour in $Q$. Then there is an antipath
with interior in $Q$, between $B$ and some right-star with a neighbour in $F$;
but the set of vertices in such an antipath contradicts
\ref{stepsantipath1}. So there is no right-star in $Q$, and every right-star with
a neighbour in $F$ is
$Q$-complete. We are given that there is a right-star with a neighbour in $F$,
and so all hypotheses
of the theorem are true with ``left'' and ``right'' exchanged. It follows by the
same argument, therefore, that every left-star with a neighbour in $F$ is $Q$-complete.
This proves (2).

\bigskip

Since $Q \subseteq N$ is anticonnected, it is contained in some anticonnected component of $N$,
say $N_1$.
\\
\\
(3) {\em There is a left- or right-star in $N_1$.}
\\
\\
For let $N_2$ be the union of all the anticomponents of $N$ different from $N_1$.  Assume
that no left- and right-star is in $N_1$.
Let $Y = V(G) \setminus (F \cup N)$; then there are no edges between
$F$ and $Y$, from definition of $N$. Also, $A \cup B \cup C \subseteq Y$, so in particular
$Y \not = \emptyset$, and also $N_2 \not = \emptyset$ since by hypothesis there is a left-star in $N$.
Hence $(F \cup Y, N)$ is a skew partition of $G$. By (1), every
vertex in $N$ has a neighbour in $A \cup B\cup C$ and in $F$, and so every vertex in $N_1$
has a neighbour in $B$ (since otherwise it would be a left-star by \ref{leftcon} and
therefore belong to $N_2$). Now $(B \cup C, N_1)$ is balanced, by \ref{balancev}, since
any left-star is complete to $N_1$ and anticomplete to $B \cup C$. Since
$B\cup C$ is connected (because every vertex  of $B \cup C$ is in a step and the strip is
step-connected), it follows from \ref{shiftbalance}.1 that $(F,N_1)$ is balanced. From
\ref{onepair}, $G$ admits a balanced skew partition, a contradiction. This proves (3).

\bigskip

From (3), $N_1 \not = Q$; and hence there is a vertex $v \in N\setminus Q$
with a nonneighbour in $Q$. From the maximality of $|F| + |Q|$, replacing $Q$ by $Q \cup \{v\}$
violates one of the hypotheses of the theorem. But $v$ has a neighbour in $A \cup B \cup C$ by (1);
$v \not \in F$ since it belongs to $N$; $v$ is not a left-star since all left-stars in $N$ are
$Q$-complete by (2); and so no left-star in $N$ is $Q \cup \{v\}$-complete. Since they are all
$Q$-complete, it follows that $v$ is nonadjacent to every left-star in $N$ .
Similarly $v$ is nonadjacent to every right-star in $N$.
\\
\\
(4) {\em $v$ is complete to $A \cup B$.}
\\
\\
For suppose not; then from the symmetry we may assume that $v$ has a nonneighbour in $B$.
By \ref{leftcon2}, $v$ is a left-star, a contradiction. This proves (4).

\bigskip

Choose an antipath $v\d q_1\d \cdots\d q_k$ in $Q$, such that $q_k$ has a nonneighbour
in $A \cup B$, with $k$ minimum. From (4), $k \ge 1$. From the minimality of $k$,
$\{v, q_1,\ldots,q_{k-1}\}$ is complete to $A \cup B$. Let $A_1$ be the set of neighbours of $q_k$
in $A$, and $A_2 = A \setminus A_1$, and define $B_1,B_2 \subseteq B$ similarly. So
$A_2 \cup B_2$ is nonempty.
\\
\\
(5) {\em $k$ is odd.}
\\
\\
For $A_2 \cup B_2$ is nonempty. If there exists $a_2 \in A_2$,  let $b_0 \in N$ be a
right-star; then \[b_0\d v\d q_1\d \cdots\d q_k\d a_2\d b_0\] is
an antihole, so it follows that $k$ is odd. The result follows similarly if $B_2$ is nonempty.
\\
\\
(6) {\em $A_1$ is complete to $B_2$, and $A_2$ is complete to $B_1$.}
\\
\\
For suppose that $a_1 \in A_1$ and $b_2 \in B_2$ are nonadjacent. Let $b_0 \in N$ be a right-star;
then by (5), \[b_0\d v \d q_1\d \cdots\d q_k\d b_2\d a_1\d b_0\] is an odd antihole, a contradiction. So
$A_1$ is complete to $B_2$ and similarly $A_2$ is complete to $B_1$. This proves (6).
\\
\\
(7) {\em $A_1,B_1,A_2,B_2$ are all nonempty.}
\\
\\
For we may assume that $A_2$ is nonempty. Since the strip is step-connected, every vertex in $A$
has a nonneighbour in $B$, and so by (6), $B_1 \not = B$. Hence $B_2$ is also nonempty.
Since $q_k$ has a neighbour in $A \cup B \cup C$ it follows that it has a neighbour in $B$,
by \ref{leftcon}, and similarly it has a neighbour in $A$. This proves (7).

\bigskip

Now the strip is step-connected, and so there is a step $a_1\d R\d b_2, a_2\d R'\d b_1$
with $a_1 \in A_1$ and $a_2 \in A_2$. Since $a_1$ is not adjacent to $b_1$ it follows that
$b_1 \in B_1$ by (6), and similarly $b_2 \in B_2$. Also by (6), $R$ and $R'$ both have length 1.
Let $a_0 \in N$ be a left-star and $b_0 \in N$ a right-star. Since $v\d a_1\d a_0\d b_0\d b_2\d v$
is not an odd hole, it follows that $a_0$ is not adjacent to $b_0$.

For every vertex $u \in V(G) \setminus F$, let $F_u$ be the set of vertices in $F$ adjacent to $u$.
\\
\\
(8) {\em $F_{a_0} \cap F_{b_0} = \emptyset$, and every path in $F$ between $F_{a_0}$ and
$F_{b_0}$ meets both $F_{v}$ and $F_{q_k}$.}
\\
\\
For if $f \in F_{a_0} \cap F_{b_0}$, then $f\d a_0\d a_1\d b_2\d b_0\d f$ is an odd hole, so
$F_{a_0} \cap F_{b_0} = \emptyset$. Let $p_1\d P\d p_2$ be a path in $F$ between $F_{a_0}$ and
$F_{b_0}$, with $V(P)$ minimal, where $p_1 \in F_{a_0}$ and $p_2 \in F_{b_0}$. Hence
\[a_0\d p_1\d P\d p_2\d b_0\d b_1\d a_2\d a_0\] is a hole, and so $P$ is odd. If $P$ does not meet
$F_{v}$ then \[v\d a_1\d a_0\d p_1\d P\d p_2\d b_0\d b_1\d v\] is an odd hole, while if $P$ does not
meet $F_{q_k}$ then \[q_k\d a_0\d p_1\d P\d p_2\d b_0\d q_k\] is an odd hole, in both cases a contradiction.
This proves (8).
\\
\\
(9) {\em Every path in $F$ between $F_{v}$ and $F_{q_k}$ meets both $F_{a_0}$ and $F_{b_0}$.}
\\
\\
For suppose not; then since $F$ is connected and $F_{a_0} \cap F_{b_0} = \emptyset$, there is
a connected subset $F'$ of $F$ meeting both $F_{v}$, $F_{q_k}$ and meeting exactly one of
$F_{a_0},F_{b_0}$. From the symmetry we may assume $F'$ meets $F_{a_0}$ and not $F_{b_0}$.
Define $q_{k+1} = a_2$; then $q_{k+1}$ has no neighbour in $F'$, so we may choose
$i$ with $1 \le i \le k+1$ minimum so that $q_i$ has no neighbour in $F'$. Note that
$v$ has a neighbour in $F'$ (because $F'$ meets $F_{v}$). If $i$ is even, then
$b_0\d v\d q_1\c q_i$ is an odd antipath; its internal vertices have neighbours in $F'$,
and its ends do not, and $a_1$ is complete to its interior and has no neighbours in $F'$,
contrary to \ref{greentouch} in the complement. If $i$ is odd, then
$b_1\d a_0\d v \d q_1\c q_i$ is an odd antipath, and its internal vertices have neighbours in $F'$
and its ends do not, and again $a_1$ is complete to its interior and has no neighbours in $F'$,
contrary to \ref{greentouch} in the complement. This proves (9).

\bigskip

Let $f_1\d f_2\d \cdots\d f_n$ be a minimal path in $F$ between $F_{a_0}$ and
$F_{b_0}$, where $f_1 \in F_{a_0}$ and $f_n \in F_{b_0}$. Then $n \ge 2$ by (8), and by (8)
and (9) it follows that $f_1\d f_2\d \cdots\d f_n$ is also a minimal path between
$F_{v}$ and $F_{q_k}$, so we may assume that $f_1 \in F_{v}$, $f_n \in F_{q_k}$, and no
other vertex of the path is in either set. Then $f_1\d f_2 \c f_n\d q_k\d a_0\d f_1$ and
$f_1\d f_2 \c f_n\d b_0\d b_1\d v\d f_1$ are both holes, of different parity, a contradiction.
This proves \ref{1-hat}.\bbox

\section{Attachments in a staircase}

For the next step of our approach towards the long odd prism, let us fix a little more than
just the strip. Let $S = (A,C,B)$ be a step-connected strip in $G$, and let
$a_0\d R_0\d b_0$ be a banister of length $\ge 3$.
We call the pair $K = (S, R_0)$ a {\em staircase}, and define $V(K) = V(R_0) \cup V(S)$.
(For brevity we often speak of the staircase $K = (S = (A,C,B), a_0\d R_0\d b_0)$, meaning
that $K = (S, R_0)$ is a staircase, and $S = (A,C,B)$, and $R_0$ has ends $a_0,b_0$, where
$a_0$ is a left-star and $b_0$ is a right-star.)
The staircase is {\em maximal} if there is no staircase $(S' = (A',C',B'),a_0'\d R_0'\d b_0')$ such that
$A \subseteq A', B \subseteq B', C \subseteq C', V(S) \subset V(S')$.

Let $K = (S = (A,C,B), a_0\d R_0\d b_0)$ be a staircase in $G$. Some definitions (all with respect
to $K$):
\begin{itemize}
\item A subset $X \subseteq V(K)$ is {\em local} if $X$ is a subset of one of
$V(S), V(R_0), A\cup \{a_0\},B\cup\{b_0\}$
\item $v \in V(G) \setminus V(K)$ is {\em minor} if its set
of neighbours in $V(K)$ is local
\item $v \in V(G) \setminus V(K)$ is {\em major} if it has neighbours in
all of $A$,$B$ and $V(R_0)$
\item $v \in V(G) \setminus V(K)$ is {\em left-diagonal} if $v$ is $(A\cup \{b_0\})$-complete,
and {\em right-diagonal} if it is $(B \cup \{a_0\})$-complete
\item $v \in V(G) \setminus V(K)$ is {\em central} if it is $(A \cup B)$-complete, and is nonadjacent
to both $a_0$ and $b_0$
\end{itemize}

First let us examine the possible types of vertices outside the staircase.

\begin{thm}\label{stairvnbrs}
Let $G$ be a Berge graph, such that there is no appearance of $K_4$ in $G$, no even prism in $G$,
and no 1-breaker in $G$.
Let $K = (S = (A,C,B), a_0\d R_0\d b_0)$ be a maximal staircase in $G$, and let
$v \in V(G) \setminus V(K)$. Then exactly one of the following holds:
\begin{enumerate}
\item $v$ is minor.
\item $v$ is major; and in that case, it is either left- or right-diagonal or central.
\item $v$ is a left-star with a neighbour in $R_0 \setminus a_0$, or a right-star with a neighbour
in $R_0 \setminus b_0$.
\end{enumerate}
\end{thm}
\Proof
\\
\\
(1) {\em If $v$ is left- or right-diagonal then the theorem holds.}
\\
\\
For assume $v$ is right-diagonal say. If it has no neighbours in $A \cup C$ then statement 3
of the theorem holds, so we assume there is a step $a_1\d R_1\d b_1,a_2\d R_2\d b_2$ such that
$v$ has a neighbour in $R_1\setminus b_1$. Hence it can be linked onto the triangle
$\{a_0,a_1,a_2\}$, via $v\d a_0$, the path from $v$ to $a_1$ with interior in $R_1\setminus b_1$, and
the path from $v$ to $a_2$ with interior in $R_2$, and so by \ref{trianglev}, $v$ has
neighbours in $A$. So it is major, and therefore statement 2 holds. This proves (1).
\\
\\
(2) {\em If $v$ is adjacent to both $a_0,b_0$ then the theorem holds.}
\\
\\
For then it has a neighbour in $R_0^*$, since $R_0$ is odd and has length $\ge 3$ and $v$
is adjacent to both its ends; and we may assume that $v$ has a neighbour in
$V(S)$, for otherwise statement 1 of the theorem holds. If $v$ has no neighbour in $B$ then
it is a left-star by \ref{leftcon}, and statement
3 of the theorem holds, so we may assume it has neighbours in $B$ and similarly in $A$.
Hence it is major. Since $(S, V(R_0^*),\{v\})$ is not a 1-breaker, $v$ does not have
nonneighbours in both $A$ and $B$, so it is either left- or right-diagonal and the claim
follows from (1). This proves (2).
\\
\\
(3) {\em If $v$ is adjacent to $a_0$ and not to $b_0$ then the theorem holds.}
\\
\\
For we may assume $v$ has a neighbour in $V(S)$. If $v$ has a neighbour in $R_0^*$, then
by \ref{leftcon2} it is either $B$-complete (when it is right-diagonal and the claim follows
from (1)) or a left-star (when statement 3 holds). So we may assume it has no neighbour
in $R_0^*$. We may assume it has a neighbour in $B \cup C$, for otherwise it is minor; let
$a_1\d R_1\d b_1,a_2\d R_2\d b_2$ be a step so that $v$ has a neighbour in $R_1 \setminus a_1$, and
in addition so that $v$ is not adjacent to $b_2$ if possible. By \ref{legaljump}, $v$ has a
neighbour in $R_2$. If $a_2$ is its only neighbour in $R_2$, then the strip
$S' = (A \cup \{v\},C,B)$ is step-connected,
since $v\d R\d b_1,a_2\d R_2\d b_2$ is an $S'$-step where $R$ is the path from $v$ to $b_1$ with interior
in $R_1 \setminus a_1$; and this is contrary to the maximality of the staircase. So $v$ has
a neighbour in $R_2 \setminus a_2$; and hence $v$ can be linked onto the triangle $\{b_0,b_1,b_2\}$
via $v\d a_0\d R_0\d b_0$, and for $i = 1,2$, the path from $v$ to $b_i$ with interior in
$R_i\setminus a_i$. By \ref{trianglev} it follows that $v$ is adjacent to both $b_1,b_2$; and hence
from our choice of the step $R_1,R_2$, and since the strip is step-connected, it follows that
$v$ is right-diagonal, and the claim follows from (1). This proves (3).
\\
\\
(4) {\em If $v$ is nonadjacent to both $a_0,b_0$ then the theorem holds.}
\\
\\
For then we may assume that $v$ has a neighbour in $R_0^*$ and in $V(S)$, since otherwise
it is minor. If it is a left-star then statement 3 holds, so we assume not; and then by
\ref{leftcon2}, it is $B$-complete. Similarly it is $A$-complete and therefore central, and
statement 2 holds. This proves (4).

\bigskip

But (2)-(4) cover all the possibilities, up to symmetry, and this completes the proof of
\ref{stairvnbrs}. \bbox

Now let us do the same thing for connected sets.
\begin{thm}\label{stairFnbrs}
Let $G$ be a Berge graph, such that there is no appearance of $K_4$ in $G$, no even prism in $G$,
and no 1-breaker in $G$.
Let $K = (S = (A,C,B), a_0\d R_0\d b_0)$ be a maximal staircase in $G$, and let
$F \subseteq V(G) \setminus V(K)$ be connected, so that its set of attachments in $V(K)$
is not local with respect to $K$. Then $F$ contains either:
\begin{enumerate}
\item a major vertex, or
\item a banister $u\d R\d v$, so that there are no edges between $V(R)$ and $V(R_0)$, or
\item (up to symmetry) a path $u\d R\d v$, where $u$ is a left-star, $v$ has a neighbour in
$R_0\setminus a_0$, and there are no edges between $V(R\setminus u)$ and $V(S)$.
\end{enumerate}
\end{thm}
\Proof
Let $X$ be the set of attachments of $F$ in $V(K)$. We may assume that $F$ is minimal (connected)
so that $X$ is not local. Now a subset of $V(K)$ is local if and
only if it does not meet both $A\cup C$ and $V(R_0 \setminus a_0)$ and does not meet both
$B\cup C$ and $V(R_0 \setminus b_0)$; so we may assume that $X$ meets both $A\cup C$ and
$V(R_0 \setminus a_0)$, and therefore from the minimality of $F$, there is a path
$f_1\d \cdots\d f_k$ where $F = \{f_1,\ldots,f_k\}$ and $f_1$ is the unique vertex of $F$ with
a neighbour in $A \cup C$, and $f_k$ is the unique vertex of $F$ with a neighbour in
$V(R_0 \setminus a_0)$. If $k = 1$ then the claim follows from \ref{stairvnbrs}, so we may
assume that $k \ge 2$.
\\
\\
(1) {\em If $f_1$ is $A$-complete then the theorem holds.}
\\
\\
For assume $f_1$ is $A$-complete. If there is no edge between $F$ and
$B \cup C$, then statement 3 of the theorem holds, so we assume that there is such an edge.
Choose $i$ with $1 \le i \le k$ minimum so that $f_i$ has a neighbour in $B\cup C$.
Suppose first that there is no edge between $\{f_1,\ldots,f_i\}$ and $V(R_0)$.
Let $a_1\d R_1\d b_1,a_2\d R_2\d b_2$ be a step so that $f_i$ has a neighbour in $R_1 \setminus a_1$,
and in addition so that $f_i$ is nonadjacent to $b_2$ if possible.
With respect to the prism formed by $R_0,R_1,R_2$, the set of attachments of
$\{f_1,\ldots,f_i\}$ is not local, and so by \ref{legaljump}, $i \ge 2$ and its attachments in
the prism are $a_1,a_2,b_1,b_2$. Hence the only edges between $\{f_1,\ldots,f_i\}$ and
$V(R_1 \cup R_2)$ are $f_1a_1,f_1a_2,f_ib_1,f_ib_2$.  From our choice of the step it follows
that $f_i$ is $B$-complete. Consequently any step satisfies
the condition we imposed on $R_1,R_2$, and so the same conclusion follows for every step; that is,
statement 2 of the theorem holds. Now assume that there is an edge between
$\{f_1,\ldots,f_i\}$ and $V(R_0)$. Suppose that $i < k$; then there is no edge between
$\{f_1,\ldots,f_i\}$ and $R_0 \setminus a_0$, from the minimality of $F$, and so $a_0$ is an
attachment of $\{f_1,\ldots,f_i\}$. But this set also has an attachment in $B\cup C$, so its set of
attachments is not local, contrary to the minimality of $F$. This proves that $i = k$.
Since $k \ge 2$, the minimality of $F$ implies that there are no edges between $\{f_2,\ldots,f_k\}$ and
$V(R_0 \setminus b_0)$; and so $b_0$ is the unique neighbour of $f_k$ in $R_0$. Hence there
are no edges between $\{f_2,\ldots,f_k\}$ and $A \cup C$, from the minimality of $F$. Also,
there are no edges between $\{f_1,\ldots,f_{k-1}\}$ and $B \cup C$, from the minimality of $i$.
Choose a step $a_1\d R_1\d b_1,a_2\d R_2\d b_2$ so that $f_k$ is adjacent to $b_1$,
and in addition so that $f_k$ is nonadjacent to $b_2$ if possible.
Since $R_1$ is odd and $a_1\d f_1\d \cdots\d f_k\d b_1\d R_1\d a_1$ is a hole, it follows that $k$ is even.
Since $a_2\d f_1\d \ldots\d f_k\d b_0\d b_2\d R_2\d a_2$ is not an odd hole, $f_k$ is adjacent to
$b_2$, and therefore to all $B$ from our choice of the step. Since
$a_1\d f_1\d \cdots\d f_k\d b_0\d R_0\d a_0\d a_1$ is not an odd hole and $R_0$ is odd, it follows that
$f_1$ is adjacent to $a_0$. But then we can add $f_1$ to $A$, $f_k$ to $B$, and
$\{f_2,\ldots,f_{k-1}\}$ to $C$, contrary to the maximality of the staircase. This proves (1).

\bigskip
By (1), we may assume there is a step $a_1\d R_1\d b_1,a_2\d R_2\d b_2$  such that $f_1$ has a
neighbour in $R_1 \setminus b_1$, and $a_2$ is not adjacent to $f_1$.
Then $R_0,R_1,R_2$ form a prism $K'$ say, and the set
of attachments of $F$ in $V(K')$ is not local with respect to $K'$. Suppose that some vertex $v$
in $F$ is major with respect to $K'$. Then we claim $v$ is major with respect to $K$. For it
has a neighbour in $A$ and in $B$, and if it has none in $R_0$ then it is adjacent to all
of $a_1,a_2,b_1,b_2$, in which case $v\d a_1\d a_0\d R_0\d b_0\d b_2\d v$ is an odd hole. So $v$ is
major with respect to $K$, and hence the theorem holds. Hence we may assume that no vertex
in $F$ is major with
respect to $K'$, and so we may apply \ref{prismjumps}. By hypothesis, \ref{prismjumps}.1
does not hold. Since no vertex of $F$ is adjacent to $a_2$, \ref{prismjumps}.2 does not hold.

Suppose that  \ref{prismjumps}.3 holds. Since $f_1$ is not adjacent to $a_2$, it follows that
$f_1$ is adjacent to $a_0, a_1$, and there exists $i$ with $2 \le i \le k$
such that $f_i$ is adjacent to $b_0,b_1$, and there are
no other edges between $\{f_1,\ldots,f_i\}$ and $V(K')$. Then
we can add $f_1$ to $A$, $f_i$ to $B$ and $\{f_2,\ldots,f_{i-1}\}$
to $C$, contrary to the maximality of the staircase.  So \ref{prismjumps}.3 does not hold.

Hence \ref{prismjumps}.4 holds, that is, there is a path $p_1\d P\d p_2$ in $F$, such that for some
$j$ with $ 0 \le j \le 2$, either:
\begin{itemize}
\item $p_1$ is adjacent to the two vertices in $\{a_0,a_1,a_2\} \setminus \{a_j\}$, and $p_2$
has neighbours in $R_j \setminus a_j$, and there are no other edges between $V(P)$ and
$V(K')\setminus \{a_j\}$, or
\item $p_1$ is adjacent to the two vertices in $\{b_0,b_1,b_2\} \setminus \{b_j\}$, and $p_2$
has neighbours in $R_j \setminus b_j$, and there are no other edges between $V(P)$ and
$V(K')\setminus \{b_j\}$
\end{itemize}
From the minimality of $F$, $F = V(P)$.
If $j>0$ then in the first case we can add $p_1$ to $A$ and $V(P\setminus p_1)$ to $C$, contrary to the
maximality of the staircase; and in the second case we do the same with $A$ and $B$ exchanged.
So $j = 0$. The first case is impossible
since no vertex in $F$ is adjacent to $a_2$; and the second case is impossible since
$f_1 \in F = V(P)$ and $f_1$ has a neighbour in $R_1 \setminus b_1$. This proves \ref{stairFnbrs}.
\bbox

The previous result can be strengthened as follows.
\begin{thm}\label{prisminstair2}
Let $G$ be a Berge graph,  such that there is no appearance of $K_4$ in $G$, no even prism in $G$,
and no 1-breaker in $G$.
Let $K = (S = (A,C,B), a_0\d R_0\d b_0)$ be a maximal staircase in $G$, and let
$F \subseteq V(G) \setminus V(S)$ be connected, containing a left-star and with an attachment in
$B \cup C$. (Note that $F$ may intersect $V(R_0)$.) Then $F$ contains either a major vertex or a banister.
\end{thm}
\Proof
We may assume $F$ is minimal (possibly exchanging $A$ and $B$);
so $F$ is the vertex set of a path $f_1\d \cdots\d f_k$, where $f_1$
is the unique left-star in $F$, and $f_k$ is the only vertex in $F$ with a neighbour in $B \cup C$.
Since $f_1$ is a left-star and $f_k$ has a neighbour in $B \cup C$ it follows that $k \ge 2$.
We may assume there is no major vertex in $F$.
\\
\\
(1) {\em We may assume that none of $f_1\l f_k$ is a right-star, and that $f_k$ is not $B$-complete.}
\\
\\
For if there is a right-star in $F$, then
it must be $f_k$; and then from the minimality of $F$ (exchanging $A$ and $B$), no vertex of $F$
different from $f_1$ has a neighbour in $A \cup C$, and so $f_1\c f_k$ is a banister. So we may
assume that there is no right-star in $F$. Since $f_k$ is
neither major nor a right-star, by \ref{stairvnbrs} it is not $B$-complete. This proves (1).
\\
\\
(2) {\em $F \cap V(R_0) = \emptyset$, and there are no edges between $\{f_2\l f_k\}$ and
$V(R_0)\setminus \{b_0\}$.}
\\
\\
For by (1), $b_0 \notin F$. Suppose that either $\{f_2\l f_k\}$ intersects $V(R_0) \setminus \{b_0\}$, or
there is an edge joining these two sets. Choose $i$ with $2 \le i \le k$ maximum so that either
$f_i \in V(R_0)\setminus \{b_0\}$
or $f_i$ has a neighbour in $V(R_0) \setminus \{b_0\}$. We claim that $f_i \notin V(R_0)$. For if
$i = k$ this is true, since $f_k$ has neighbours in $B\cup C$; and if $i< k$ then $f_{i+1}$ has
no neighbour in $V(R_0)\setminus \{b_0\}$ from the maximality of $i$, and therefore again $f_i \notin V(R_0)$.
So none of $f_i \l f_k$ belong to $V(R_0)$. Since $\{f_i\l f_k\}$ has attachments in $V(R_0)\setminus \{b_0\}$
and in $B\cup C$, and contains no major vertex or left- or right-star, this contradicts \ref{stairFnbrs}.
So $\{f_2\l f_k\}$ is disjoint from $V(R_0) \setminus \{b_0\}$ and hence from $V(R_0)$, and there are
no edges between $\{f_2\l f_k\}$ and $V(R_0) \setminus \{b_0\}$. Since there is an edge between
$\{f_2\l f_k\}$ and $f_1$ it follows that $f_1 \notin V(R_0)$, and so $F \cap V(R_0) = \emptyset$.
This proves (2).

\bigskip
Let $a_1\d R_1\d b_1,a_2\d R_2\d b_2$ be a step such that $f_k$ has a neighbour in $R_1 \setminus \{a_1\}$
and $f_k$ is nonadjacent to $b_2$. (This exists since $f_k$ is not $B$-complete.)
\\
\\
(3) {\em $f_1a_2$ is the only edge between $F$ and $R_2$.}
\\
\\
For if $f_k$ has a neighbour in $R_2$, then its neighbour set in the prism formed by $R_0,R_1,R_2$
is not local with respect to that prism, and therefore by \ref{legaljump}, $f_k$ has a neighbour
in $R_0$; and then by \ref{stairvnbrs} it is major, a contradiction. So $f_k$ has no neighbours
in $R_2$. From the minimality of $F$, there are no edges between $F$ and $R_2\setminus a_2$.
Suppose that $a_2$ has a neighbour in $\{f_2\l f_k\}$, and choose $i$ maximum so that $a_2$ is
adjacent to $f_i$. Since $f_k$ has a neighbour in $V(R_1)\setminus \{a_1\}$, the
set of attachments of $\{f_i\l f_k\}$ is not local with respect to the prism formed by
$R_0,R_1,R_2$; and since $b_2$ is not an attachment, it follows from \ref{legaljump} that there
is an attachment of $\{f_i\l f_k\}$ in $V(R_0)$. By (2), $b_0$ has a neighbour in  $\{f_i\l f_k\}$;
but then  $\{f_i\l f_k\}$ violates \ref{stairFnbrs}. This proves (3).
\\
\\
(4) {\em $b_0$ has neighbours in $\{f_1\l f_{k-1}\}$.}
\\
\\
For first suppose that $b_0$ has no neighbour in $F$. Since $b_2$ is not an attachment of $F$, it
follows from \ref{legaljump} (applied to $F$ and the prism formed by $R_0,R_1,R_2$)
that there is an edge between $F$ and $V(R_0)$, and so $f_1$ has a neighbour in $R_0$.
But then $f_1$ can be linked onto the triangle $\{b_0,b_1,b_2\}$, via the path between $f_1$
and $b_0$ with interior in $V(R_0)$, the path between $f_1$ and $b_1$ with interior in
$\{f_2\l f_k\} \cup (V(R_1)\setminus \{a_1,b_1\})$, and the path $f_1\d a_2 \d R_2 \d b_2$. This
contradicts \ref{trianglev}, and therefore proves that $b_0$ has a neighbour in $F$. Suppose that
$f_k$ is the only neighbour of $b_0$ in $F$. Then since $f_k$ is not major, its unique neighbour in
$R_1$ is $b_1$. From \ref{prismrunglength}, $R_1,R_2$ are odd, and from the hole
$f_1\c f_k \d b_0 \d b_2 \d R_2 \d a_2 \d f_1$ it follows that $k$ is odd. If $a_1$ has no neighbour
in $\{f_2\l f_k\}$ then $f_1\c f_k \d b_1 \d R_1 \d a_1 \d f_1$ is an odd hole, and if
$a_1$ has a neighbour in  $\{f_2\l f_k\}$ then  $\{f_2\l f_k\}$ violates \ref{stairFnbrs}.
So $f_k$ is not the unique neighbour of $b_0$ in $F$. This proves (4).

\bigskip
Choose $i$ with $1 \le i < k$ minimum so that $b_0$ is adjacent to $f_i$, and let $R_0'$ be the path
$f_1\c f_i \d b_0$.
There are no edges between $\{f_1\l f_i\}$ and $B \cup C$ from the minimality of $F$, and from
\ref{stairFnbrs} there are no edges between $\{f_2\l f_i,b_0\}$ and $A \cup C$. Hence
$f_1\d R_0' \d b_0$ is a banister, and in particular the three paths $R_0',R_1,R_2$ form a prism, $K'$ say.
Let $F' = \{f_{i+1},\ldots,f_k\}$. Then $F'$ is connected and disjoint from $V(K')$,
and $F'$ has attachments in $R_1\setminus a_1$, and in $R_0' \setminus b_0$, and by (3) it has no
attachments in $R_2$.  By \ref{legaljump}
applied to $K'$, it follows that $F'$ contains a path with one end adjacent to
$a_1,f_1$, the other end adjacent to $b_0,b_1$, and with no more edges between this path and
$V(R_0') \cup V(R_1)$. Since the only vertex of $F'$ adjacent to $f_1$ is $f_2$, and that only
if $i = 1$, and the only vertex in $F'$ adjacent to $b_1$ is $f_k$, it follows that $i = 1$,
and the only edges between $\{f_2,\ldots,f_k\}$ and $V(R_0') \cup V(R_1)$ are $f_kb_1,f_kb_0,
f_2a_1,f_2f_1$. But then by (2), $a_1$ can be linked onto the triangle $\{b_0,b_1,f_k\}$, via
$a_1\d a_0\d R_0\d b_0$, $a_1\d R_1\d b_1$, $a_1\d f_2\c f_k$, contrary to \ref{trianglev}.
This proves \ref{prisminstair2}.\bbox

Now we turn to anticonnected sets of major vertices. We have already defined what it is for a staircase to
be maximal in $G$. We say a staircase $K = (S = (A,C,B), a_0\d R_0\d b_0)$ is {\em strongly maximal} if
it is maximal, and in addition, either $C \not = \emptyset$, or there is no staircase
$(R',S')$ in $\overline{G}$ with $V(S) \subset V(S')$.
A {\em 2-breaker} in $G$ is a pair $(K,Q)$ such that
\begin{itemize}
\item $K = (S = (A,C,B), a_0\d R_0\d b_0)$ is a strongly maximal staircase in $G$,
\item $Q \subseteq V(G) \setminus V(K)$ is anticonnected,
\item some vertex of $A$ is $Q$-complete, and some vertex of $B$ is $Q$-complete
\item $a_0,b_0$ are not $Q$-complete, and
\item some vertex of $R_0$ is $Q$-complete.
\end{itemize}

We observe that if $q$ is a central vertex with respect to a strongly maximal staircase $K$,
then $(K,\{q\})$ is a 2-breaker, so it follows from the next result that we no longer have to
worry about central vertices.

\begin{thm}\label{central}
Let $G$ be a Berge graph, containing no appearance of $K_4$, no even
prism, and no 1-breaker. If there is a 2-breaker in $G$ then $G$ admits a balanced skew partition.
\end{thm}
\Proof Choose a 2-breaker $(K,Q)$ in $G$, with notation as above, such that for fixed $K$ the set $Q$ is maximal.
Let $a_0\d S\d s$ and $b_0\d T\d t$ be the subpaths of $R_0$ such that $s$ is the unique $Q$-complete
vertex of $S$, and $t$ is the unique $Q$-complete vertex of $T$.
\\
\\
(1) {\em $S,T$ both have odd length, and therefore $s,t$ are different.}
\\
\\
For choose $a \in A$ and $b \in B$, both $Q$-complete; then $a\d a_0\d S\d s$ has length $>1$, and its ends are
$Q$-complete and its internal vertices are not, and $b$ is also $Q$-complete and has no
neighbours in the interior of $a\d a_0\d S\d s$. By \ref{greentouch}, this path is even, and
so $S$ is odd, and similarly $T$ is odd. Since $R_0$ is odd it follows that $s,t$ are different.
This proves (1).
\\
\\
(2) {\em Every vertex in $A \cup B$ is $Q$-complete.}
\\
\\
For suppose some vertex in $A$ say is not $Q$-complete. Choose a step $a_1\d R_1\d b_1,a_2\d R_2 \d b_2$
so that $a_1$ is $Q$-complete and $a_2$ is not. Since $s,t$ are different it follows that
$t$ is nonadjacent to both $a_0,a_2$; and so by \ref{triangleX}, $Q$ cannot be linked onto the
triangle $\{a_0,a_1,a_2\}$. Hence there is no $Q$-complete vertex in $R_2$. Assume $s,t$ are
nonadjacent; then the subpath of $R_0$ between them is odd, and $a_1$ has no neighbour in
its interior, so by \ref{greentouch} it contains another $Q$-complete vertex $u$ say; and then
$s\d S\d a_0\d a_2\d R_2\d b_2\d b_0\d T\d t$ is an odd path, its ends are $Q$-complete and its internal vertices
are not, and $u$ has no neighbour in its interior, contrary to \ref{greentouch}. So $s,t$
are adjacent. Hence the hole $a_0\d R_0\d b_0\d b_2\d R_2\d a_2\d a_0$ has length $\ge 6$, and the only
$Q$-complete vertices in it are the adjacent vertices $s,t$. By \ref{RRC} $Q$ contains a hat
or a leap; and in either case there is a vertex $q \in Q$ with no neighbours in $R_2$. But
$q$ is adjacent to $s$ and $a_1$, contrary to \ref{legaljump} applied to the prism
formed by $R_0,R_1,R_2$. This proves (2).
\\
\\
(3) {\em Every major vertex is either in $Q$ or complete to $Q$.}
\\
\\
For let $v$ be a major vertex, and suppose $v \not \in Q$, and $Q'$ is anticonnected, where
$Q' = Q \cup \{v\}$.
From \ref{stairvnbrs}, $v$ is either left- or right-diagonal, or central; and in either case
it has neighbours $a_1 \in A$ and $b_1 \in B$ that are nonadjacent. It follows that
$a_1\d a_0\d R_0\d b_0\d b_1$ is an odd path of length $\ge 5$, and its ends are $Q'$-complete. From the
maximality of $Q$, none of its internal vertices are $Q'$-complete, and so by \ref{RR}, $Q'$
contains a leap $q_1,q_2$ say. So neither of $q_1,q_2$ has neighbours in the interior of $R_0$; but this
is impossible since one of them is in $Q$ and is therefore adjacent to $s$. This proves (3).
\\
\\
(4) {\em There is no edge $uv$ of $G \setminus V(S)$ such that $u$ is a left-star, $v$ is a
right-star, and $u,v$ are not $Q$-complete.}
\\
\\
For suppose $uv$ is such an edge. Since $u,v$ have neighbours in $A \cup B$, they do not belong to
$R_0^*$. Since $u,v$ have nonneighbours in $Q$ and $Q$ is anticonnected,
there is an antipath $u\d q_1 \c q_k\d v$ with $q_1\l q_k \in Q$. Choose a step
$a_1\d R_1\d b_1,a_2\d R_2\d b_2$. Then $a_1\d b_2\d u\d q_1\c q_k\d v\d a_1$ is an antihole, so $k$ is even.
Hence every $Q$-complete vertex $w$ say is adjacent to one of $u,v$, for otherwise
$w\d u\d q_1 \c q_k\d v\d w$ would be an odd antihole. In particular, there are no $Q$-complete vertices
in $C$; and therefore $a_1\d R_1\d b_1$ is an odd path with both ends $Q$-complete and no internal vertex
$Q$-complete. Since $a_2$ is $Q$-complete and has no neighbour in the interior of $R_1$, it follows
from \ref{greentouch} that $R_1$ has length 1, and similarly $R_2$ has length 1. Since this step
was arbitrary, and every vertex is in a step, it follows that $C = \emptyset$.
Suppose that $u$ has no neighbour in $R_0^*$. Then all $Q$-complete vertices in $R_0^*$ are adjacent
to $v$. In particular,
$v$ is adjacent to $s,t$ and hence does not belong to $R_0$ (because $v$ is a right-star); and
$s\d S\d a_0\d a_1\d b_1$ is an odd path, its ends are $(Q \cup \{v\})$-complete, its internal
vertices are not, and the $(Q \cup \{v\})$-complete $t$ has no neighbour in its interior,
contrary to \ref{greentouch}. So $u$ has a neighbour in $R_0^*$, and similarly so does $v$. Now
$b_1\d u \d Q \d v \d a_1$
is an odd antipath, all its internal vertices have neighbours in the connected set $R_0^*$,
and its ends do not. By \ref{RR} applied in $\overline{G}$, there is a leap; that is, there exist
adjacent $a,b \in R_0^*$, both $Q$-complete, such that $b\d u \d Q \d v \d a$ is an antipath.
Define $A'= A \cup \{a\}$ and $B' = B\cup \{b\}$; then $(A',\emptyset,B')$
is a strip ($S'$ say) in $\overline{G}$. For every edge $a_1b_1$ of $G$ with $a_1\in A$ and $b_1 \in B$, the
pair $a\d b_1,a_1\d b$ is a step of $S'$ (in $\overline{G}$), and every
vertex of $A \cup B$ is in such an edge, and so $S'$ is step-connected.
Hence $((A',\emptyset,B'), v\d q_k\c q_1\d u)$ is a staircase in $\overline{G}$, contrary to
the hypothesis that $K$ is strongly maximal. This proves (4).
\\
\\
(5) {\em Every path in $G$ from an $A$-complete vertex to a vertex with a neighbour in $B \cup C$
contains either a vertex in $Q$ or a $Q$-complete vertex.}
\\
\\
For suppose not, and choose a path $p_1\d \cdots\d p_k$ say, with $k$ minimum such that $p_1$ is
$A$-complete and $p_k$ has a neighbour in $B \cup C$, and none of $p_1,\ldots,p_k$ is in $Q$
or $Q$-complete. Since $A \cup B$ is complete to $Q$ it follows that none of $p_1,\ldots,p_k$
is in $A \cup B$. Now $p_1$ is not in $C$ since no vertex in $C$ is $A$-complete (because they are
all in steps), and if some $p_i \in C$ for $i >1$, then $p_1\d \cdots\d p_{i-1}$ is a shorter path
with the same properties, contrary to the minimality of $k$. So none of $p_1,\ldots,p_k$
is in $V(S)$. (Some may be in $R_0$, however.) Since none of $p_1,\ldots,p_k$ is major by (3), it
follows from \ref{prisminstair2} and the minimality of $k$ that $p_1\d \cdots\d p_k$ is a banister.
From (4), since none of $p_1,\ldots,p_k$ is $Q$-complete, it follows that $k > 2$.
Let $a_1\d R_1\d b_1,a_2 \d R_2\d b_2$ be a step. From the hole $a_1\d p_1\c p_k\d b_1\d R_1\d a_1$
it follows that $k$ is even; and so  $a_1\d p_1 \c p_k\d b_2$ is an odd path of length
$\ge 5$; its ends are $Q$-complete, and its internal vertices are not. By \ref{RR}, $Q$ contains
a leap $a,b$; so $a\d p_1\d \cdots\d p_k\d b$ is a path. But then $(A \cup \{a\},C, B \cup \{b\})$
is a step-connected strip $S'$ say (since for every nonadjacent $a' \in A$ and $b' \in B$, the two
paths $a\d b',a'\d b$ make a step in this strip), and so $(S', p_1\d \cdots\d p_k)$ is a staircase,
contrary to the maximality of $(S,R_0)$. This proves (5).

\bigskip

Let $X$ be the set of all $Q$-complete vertices in $G$; let $M$ be the component of
$G \setminus (Q \cup X)$ that contains $a_0$, and $N$ the union of all the other components. By
(5), $b_0 \in N$, so $N$ is nonempty, and hence $(M \cup N, Q \cup X)$ is a skew partition
of $G$. Choose $b \in B$; then $b \in X$, and it has no neighbour in $M$ by (5). Hence the
skew partition is loose, and so $G$ admits a balanced skew partition, by \ref{geteven}.
This proves \ref{central}. \bbox

\begin{thm}\label{diagonalchain}
Let $G$ be a Berge graph, containing no appearance of $K_4$, no even prism, no 1-breaker and no 2-breaker.
Let $K = (S = (A,C,B), a_0\d R_0\d b_0)$ be a strongly maximal staircase in $G$.
Let $q_1\c q_k$ be an antipath such that $q_2,\ldots,q_{k-1}$ are both left- and right-diagonal,
and $q_1$ is left- and not right-diagonal, and $q_k$ is right- and not left-diagonal. Then
$q_1$ is a left-star and $q_k$ is a right-star.
\end{thm}
\Proof
First, obviously $k \ge 2$. Let $Q = \{q_1,\ldots,q_k\}$.
\\
\\
(1) {\em If $q_1$ is adjacent to $a_0$ and $q_k$ to $b_0$ then the theorem holds.}
\\
\\
For then both $a_0,b_0$ are $Q$-complete, and $q_1$ has a nonneighbour in $B$
(for otherwise it would be right-diagonal), and $q_k$ has a nonneighbour in $A$.
Since $R_0$ has odd length $\ge 3$, it follows that each of $q_1,\ldots,q_k$ has a neighbour
in $R_0^*$. Since $(S,R_0^*,Q)$ is not a 1-breaker, it follows
that $Q$ contains a left-star, which must be $q_1$; and similarly $q_k$ is a right-star. Then
the theorem holds. This proves (1).
\\
\\
(2) {\em If $q_1$ is adjacent to $a_0$ and $q_k$ is nonadjacent to $b_0$ then the theorem holds.}
\\
\\
For in this case, $q_1$ has a nonneighbour in $B$, say $b$. From the antihole
$a_0\d b\d q_1\d \cdots\d q_k\d b_0\d a_0$ we deduce that $k$ is odd. Now $R_0$ is odd, of length $\ge 3$,
and its ends are complete to $Q \setminus q_k$, and so is every $a \in A$, and $a$ has no neighbour
in the interior of $R_0$, so by \ref{greentouch}, there is a $(Q \setminus q_k)$-complete vertex in
the interior of $R_0$, say $t$. Let $T$ be the subpath of $t$ to $b_0$, and let us choose
$t$ with $T$ of minimum length, that is, so that $t$ is the unique $(Q \setminus q_k)$-complete
vertex of $T$ .  If $t$ is nonadjacent to $q_k$ then
$t\d b\d q_1\c q_k\d t$ is an odd antihole (since $k \ge 2$) , a contradiction. Hence
$t$ is $Q$-complete, and in particular, all of $q_1,\ldots,q_k$ have neighbours in
the interior of $R_0$. By \ref{stepsantipath1} it follows that $Q$ contains a left-star,
which must be $q_1$. We may assume that $q_k$ is not a right-star, for otherwise the theorem holds.
Since $q_k$ is right-diagonal, from \ref{stairvnbrs} it follows that $q_k$ is major and therefore
has a neighbour in $A$. Choose a step $a_1\d R_1\d b_1,a_2\d R_2\d b_2$ so that $q_k$ is adjacent to $a_1$,
and if possible nonadjacent to $a_2$.  Then $t\d T\d b_0\d b_1\d R_1\d a_1$ is a path, and both its ends are
$Q$-complete, and none of its internal vertices are $Q$-complete (since $q_1$ is a left-star).
By \ref{evenantipath2} applied to $t\d T\d b_0\d b_1\d R_1\d a_1$ and the antipath
$b_1\d q_1 \c q_k \d b_0$,
it follows that $t\d T\d b_0\d b_1\d R_1\d a_1$ has length 4, and so $R_1$ has length 1 and $T$ has
length 2; let its middle vertex be $u$ say. Also from \ref{evenantipath2}, $u$ is
$Q \setminus q_1$-complete, and nonadjacent to $q_1$. Suppose that $q_k$ is nonadjacent to $a_2$.
Then there is no $Q$-complete vertex in $R_2$. If $t$ is nonadjacent to $a_0$ then
$a_0\d a_2\d R_2\d b_2\d b_0\d u\d t$ is an odd path of length $\ge 5$; its ends are $Q$-complete and its
internal vertices are not, so by $\ref{RR}$, $Q$ contains a leap, which is impossible since
every vertex in $Q$ is adjacent to one of $b_0,b_2$. If $t$ is adjacent to $a_0$, then
$a_0\d a_2\d R_2\d b_2\d b_0\d R_0\d a_0$ is a hole of length $\ge 6$, and the only $Q$-complete vertices in
it are $a_0,t$, and these are adjacent; so by \ref{RRC} there is a hat or a leap in $Q$; and again
this is impossible since every vertex in $Q$ is adjacent to one of $b_0,b_2$. This proves that
$q_ka_2$ is an edge. From our choice of the step, it follows that $q_k$ is $A$-complete.
But therefore any step satisfies the condition we imposed on the step $R_1,R_2$; and therefore
every path in every step has length 1, that is $C = \emptyset$. Then
$S = (A \cup \{t\},\emptyset,B\cup \{u\})$ is a step-connected strip in $\overline{G}$, and
$(S', b_0\d q_k\d \cdots\d q_1)$ is a staircase in $\overline{G}$, contradicting that $(S,R_0)$ is
strongly maximal. This proves (2).
\\
\\
(3) {\em If $q_1$ is nonadjacent to $a_0$ and $q_k$ is nonadjacent to $b_0$ then the theorem holds.}
\\
\\
For then $a_0\d q_1\d \cdots\d q_k\d b_0\d a_0$ is an antihole, so $k$ is even. Let $A_1$ be the set of
vertices in $A$ adjacent to $q_k$, and $A_2 = A \setminus A_1$; and let $B_1$ be the set of
vertices in $B$ adjacent to $q_1$, and $B_2 = B \setminus B_1$. If $a_1 \in A_1$ and
$b_2 \in B_2$, then $a_1\d b_2\d q_1\d \cdots\d q_k\d b_0\d a_1$ is not an odd antihole, and so
$a_1$ is adjacent to $b_2$; and hence $A_1$ is complete to $B_2$, and similarly $A_2$ is complete
to $B_1$. If $A_1,B_1$ are both empty then by \ref{stairvnbrs}, the theorem holds; so we may
assume that $A_1$ is nonempty. Choose $a_1 \in A_1$. Since $a_1$ is in a step, it has a nonneighbour
in $B$, say $b_1$. Since $a_1$ is $B_2$-complete it follows that $b_1 \in B_1$.
Then $a_1,b_1$ are both $Q$-complete, and since $(K,Q)$ is not a 2-breaker, no internal vertex of $R_0$
is $Q$-complete. So $a_1\d a_0\d R_0\d b_0\d b_1$ is an odd path of length $\ge 5$, and its ends are
$Q$-complete, and its internal vertices are not. By \ref{RR}, $Q$ contains a leap. Since
every vertex of $Q$ except $q_1,q_k$ has $\ge 2$ neighbours in $R_0$, it follows that $k = 2$ and
$q_1,q_2$ both have no neighbours in the interior of $R_0$. Then
$S' = (A \cup \{q_2\}, C, B \cup \{q_1\})$ is a step-connected strip (since $a_1\d q_1,q_2\d b_1$ is a step
of it), and $(S', R_0)$ is a staircase, contrary to the maximality of $(S,R_0)$. This proves (3).

\bigskip

From (1),(2),(3), the theorem follows. This proves \ref{diagonalchain}.\bbox

\section{The long odd prism}

In this section we apply the results of the previous section to prove that a Berge graph
containing a long odd prism has a decomposition unless it is a line graph.

Let $K = ((A,C,B), a_0\d R_0\d b_0)$ be a strongly maximal staircase in a Berge graph $G$. From
\ref{stairvnbrs} there
are three possible kinds of $B$-complete vertices; right-stars, vertices complete to both $A$ and $B$,
and $B$-complete vertices adjacent to some but not all of $A$. The most difficult step in handling the
long odd prism is when there is a vertex of the third kind. In that case, we shall construct
a subset of $B$-complete vertices, including all these ``mixed'' vertices and some of the others,
such that they and their common neighbours form a cutset of the graph, and thereby give us a skew
partition. We define the set recursively as follows: initially let $X$ be the set of all
$B$-complete vertices adjacent to some but not all of $A$. Then enlarge $X$ by repeatedly
applying the following two rules, in any order:
\begin{enumerate}
\item if there is a $A\cup B$-complete vertex $v$ that is not in $X$ and not $X$-complete,
add $v$ to $X$
\item if there is a banister $a\d R\d b$ such that $a$ is not $X$-complete and $b$ is not in $X$,
add $b$ to $X$.
\end{enumerate}
The process eventually stops with some set $X$. We shall prove that $X$ and its common neighbours (say $Y$)
separate $A$ (or at least the part of $A$ that is not $X$-complete) from $b_0$, and this will provide
a balanced skew partition. To prove that $X\cup Y$ separates $G$ as described, we have
to show that every path from $A$ to $b_0$ meets $X \cup Y$, and it turns out that there are only two
kinds of paths to worry about; banisters, and 1-vertex paths consisting of a major vertex. Any banister $a\d R\d b$
is automatically hit, because of the rule above; if $a \not \in Y$ then $b \in X$. The 1-vertex paths are
trickier. Let $v$ be a major vertex. If it is $B$-complete, then it is in either $Y$ or $X$ by the rule above,
so assume it is not $B$-complete. By \ref{stairvnbrs}, it is left- and not right-diagonal, and now we have
to show it belongs to $Y$. If only we knew that every vertex in $X$ was adjacent to $a_0$, then it follows
easily that $v \in Y$, because of \ref{diagonalchain}. So that is what we need to do --- to prove that every
vertex in $X$ is adjacent to $a_0$.

Let us start again, more formally. Let $K = ((A,C,B), a_0\d R_0\d b_0)$ be a staircase in a Berge graph $G$.
We define a {\em right-sequence} to be a sequence $x_1,\ldots,x_t$, with the following properties
(which we refer to as the {\em right-sequence axioms}):
\begin{enumerate}
\item $x_1,\ldots,x_t$ are distinct and $B$-complete
\item for $1 \le i \le t$, if $x_i$ is $A$-complete then there exists $h$ with $1 \le h < i$ such that
$x_h$ is nonadjacent to $x_i$
\item for $1 \le i \le t$, if $x_i$ is $A$-anticomplete then there is a banister $r\d R\d x_i$ such that
$r$ has a nonneighbour in $\{x_1,\ldots,x_{i-1}\}$.
\end{enumerate}

Any initial subsequence of a right-sequence is therefore another right-sequence.
We say $x_i$ is {\em earlier} than $x_j$ if $i < j$. Let $X =  \{x_1,\ldots,x_t\}$.
For each $x_i \in X$ that has an earlier nonneighbour, we define its {\em predecessor}
to be $x_h$, where $h$ is minimum such that $1 \le h < i$ and $x_h$ is nonadjacent to $x_i$.
From the second axiom, every $x_i$ either has a nonneighbour in $A$ or a predecessor, so we can follow
the sequence of predecessors until we get to some vertex that is not $A$-complete.
For each $x_i$ we therefore define the {\em trajectory} of $x_i$ to be the sequence $w_1\d \cdots\d w_n$ with the
following properties:
\begin{itemize}
\item $n \ge 1$, and $w_1 = x_i$
\item $w_n$ has a nonneighbour in $A$
\item for $1 \le j < n$, $w_j$ is $A$-complete, and $w_{j+1}$ is the predecessor of $w_j$.
\end{itemize}
Clearly the trajectory is unique, and is an antipath. If $v \in V(G)$ is $A$-complete, not in
$X$ and not $X$-complete, we define the {\em trajectory} of $v$ to
be the antipath $v\d w_1\d \cdots\d w_n$, where $w_1$ is the earliest nonneighbour of $v$ in $X$,
and $w_1\d \cdots\d w_n$ is the trajectory of $w_1$.

Let $a$ be a left-star. If it is not $X$-complete, we define the {\em birth}
of $a$ to be the earliest nonneighbour of $a$ in $X$. Now let $b$ be a right-star. A
banister $a\d R\d b$ is said to be {\em $b$-optimal} if $a$ is not $X$-complete, and there is no
banister $a'\d R'\d b$ such that $a'$ is not $X$-complete and the birth of $a'$ is earlier than
the birth of $a$.

\begin{thm}\label{trajectory}
Let $G$ be Berge, containing no appearance of $K_4$, no even prism, no 1-breaker and no 2-breaker.
Let $K = (S = (A,C,B), a_0\d R_0\d b_0)$ be a strongly maximal staircase in $G$,
and let $x_1,\ldots,x_t$ be a right-sequence. Let $b$ be a right-star, and
let $a\d R\d b$ be a $b$-optimal banister.
Let $a\d w_1\d \cdots\d w_n$ be the trajectory of $a$. Then $n$ is odd, and either:
\begin{itemize}
\item $b$ is the unique vertex of $R$ which is $\{w_1,\ldots,w_n\}$-complete, or
\item $R$ has length 1, and there exists some even $m$ with $1 \le m < n$ such
that $a\d w_1\d \cdots\d w_m\d b$ is an antipath.
\end{itemize}
\end{thm}
\Proof We proceed by induction on $t$, and assume the result holds for all smaller values of $t$.
Hence we may assume that $w_1 = x_t$, for otherwise the result follows by induction.
Let  $W$ = $\{w_1,\ldots,w_n\}$; then every vertex in $B$ is $W$-complete.
\\
\\
(1) {\em $n$ is odd.}
\\
\\
For choose $a_2 \in A$ nonadjacent to $w_n$, and $b_1 \in B$ nonadjacent to $a_2$; then
$b_1\d a\d w_1\d \cdots\d w_n\d a_2\d b_1$ is an antihole, so $n$ is odd. This proves (1).
\\
\\
(2) {\em If $w_n$ has a neighbour in $A$ then the theorem holds.}
\\
\\
For choose a step $a_1\d R_1\d b_1,a_2\d R_2\d b_2$ such that $w_n$ is adjacent to $a_1$ and not to $a_2$.
Then $a_1$,$b_2$ are $W$-complete.
Suppose first that there are no $W$-complete vertices in $R$. Then $a_1\d a\d R\d b\d b_2$ is an odd path between
$W$-complete vertices. If $R$ has length 1 then there is an antipath $Q$ joining $a,b$ with
interior in $W$, and since it can be completed to an antihole via $b\d a_1\d b_2\d a$, it has odd length
and the theorem holds. So we may assume $R$ has length $>1$, and hence by \ref{RR} $W$ contains a leap.
Since all vertices of $W$ except $w_1$ are adjacent to $a$, the leap is $w_1,w_2$; and hence the only edges
between $w_1,w_2$ and $R$ are $w_1b$ and $w_2a$. Since $n$ is odd it follows that $n >2$ and so $w_1,w_2$
are both $A \cup B$-complete. But then $S' = (A \cup \{w_2\},C,B\cup \{w_1\})$ is a step-connected
strip, and $(S', a\d R\d b)$ is a staircase, contrary to the maximality of $(S,R_0)$. So we may assume there
are $W$-complete vertices in $R$. If $b$ is the only one then the theorem holds, so assume there is another.
But then $W$ can be linked onto the triangle $\{a,a_1,a_2\}$, via a subpath of $R \setminus b$, the
1-vertex path $a_1$, and a subpath of $R_2$. Since $b_1$ is $W$-complete and nonadjacent to both
$a,a_2$, this contradicts \ref{triangleX}. This proves (2).

\bigskip
From (2) we may assume that $w_n$ has no neighbour in $A$. Let $w_n = x_s$ say.
From the third axiom, there is a banister $r'\d R'\d w_n$, such that $r'$ has a nonneighbour in
$\{x_1,\ldots,x_{s-1}\}$, and therefore we may choose it to be $w_n$-optimal.
\\
\\
(3) {\em $R'$ is disjoint from $R$, and there are no edges between $V(R)\setminus a$ and $V(R')\setminus w_n$.}
\\
\\
Suppose that $(R\setminus a) \cup (R' \setminus w_n)$ is connected. Then it contains a path
between $r'$ and $b$, with interior in the union of the interior of $R$ and $R'$, and therefore
this path is a banister. But $R$ is $b$-optimal, and the birth of $r'$ is earlier than the birth
of $a$, a contradiction. So $R \setminus a$ is disjoint from $R' \setminus w_n$, and there are no
edges between them. Since $a \not = r'$ (because their births are different), and $b \not = w_n$
(because $R$ is optimal for $b$) it follows that $R$ is disjoint from $R'$. This proves (3).

\bigskip

Let $v_1\d \cdots\d v_m$ be the trajectory of $r'$, let
$V = \{v_1,\ldots,v_m\}$, and let $W' = \{a,w_1,\ldots,w_{n-1}\}$. Since each of
$v_1,\ldots,v_m$ is earlier than $w_n$, it follows from the definition of trajectory that
$v_1,\ldots,v_m$ are all $W'$-complete. By induction on $t$, it follows that either
$w_n$ is the unique $V$-complete vertex in $R'$, or $R'$ has length 1
and there is an odd antipath between $r'$ and $w_n$ with interior in $V$.
\\
\\
(4) {\em If $n = 1$ then the theorem holds.}
\\
\\
For let $n = 1$, and choose a step $a_1\d R_1\d b_1,a_2\d R_2\d b_2$. Suppose first that $a$ has no neighbour
in $R'$. Now $a$ is $V$-complete, and either $w_1$ is the unique $V$-complete vertex in $R'$, or $R'$
has length 1 and there is an odd antipath $Q$ between $r'$ and $w_1$ with interior in $V$. In the first case,
$a\d a_1\d r'\d R'\d w_1$ is an odd path, its ends are $V$-complete, its internal vertices are not, and the
$V$-complete vertex $b_2$ has no neighbour in its interior, contrary to \ref{greentouch}. In the second case,
$a\d r'\d Q\d w_1\d a$ is an odd antihole. This proves that $a$ has a neighbour in $R'$. Now suppose it has
a neighbour different from $r'$; then $R'$ has length $>1$, and so
 $w_1$ is the unique $V$-complete vertex in $R'$; and there is a path $P'$ say from $a$ to $w_1$ with
interior in $R' \setminus r'$. Since the ends of this path are $V$-complete and its internal vertices
are not, and the $V$-complete vertex $b_1$ has no neighbour in its interior, it is even by \ref{greentouch}.
But it can be completed to an odd hole via $w_1\d b_1\d R_1\d a_1\d a$, a contradiction. This proves that
$r'$ is the unique neighbour of $a$ in $R'$. Since $a\d r'\d R'\d w_1\d b_1\d b\d R\d a$ is not an odd hole,
it follows
from (3) that $w_1$ has neighbours in $R$. If $b$ is its unique neighbour in $R$ then the theorem holds, so
we assume not. Then there is a path $P$ say from $w_1$ to $a$ with interior in $R \setminus b$. Since
$w_1\d P\d a \d r' \d R' \d w_1$ is a hole it follows that $P$ is even; but $P$ can be completed via
$a\d a_1\d R_1\d b_1\d w_1$, a contradiction. This proves (4).

\bigskip

We may therefore assume that $n \ge 3$ (since it is odd.)
\\
\\
(5) {\em $C = \emptyset$.}
\\
\\
For suppose not, and choose a step $a_1\d R_1\d b_1,a_2\d R_2\d b_2$ where $R_1$ has length $>1$. Since $R_1$
is odd, and its ends are $(W\setminus w_n)$-complete, and the $(W\setminus w_n)$-complete vertex $b_2$
has no neighbour in its interior, there is a $(W\setminus w_n)$-complete vertex $v$ in the
interior of $R_1$, by \ref{greentouch}. But then $v$ is nonadjacent to both $a$ and $w_n$, since
they are left- and right-stars respectively, and so $v\d a\d w_1\d \cdots\d w_n\d v$ is an odd antihole, a
contradiction. This proves (5).
\\
\\
(6) {\em If $b$ is not $(W\setminus w_n)$-complete and no edge of $R$ is
$(W\setminus w_n)$-complete then the theorem holds.}
\\
\\
For choose a step $a_1\d R_1\d b_1,a_2\d R_2\d b_2$. Then $a_1\d a\d R\d b\d b_2$ is an odd path, its ends are
$(W\setminus w_n)$-complete, and none of its edges are $(W\setminus w_n)$-complete.
Suppose first that $R$ has length
$\ge 3$. Then by \ref{RR} there is a leap in $W\setminus w_n$; and so there are nonadjacent vertices
$x,y \in W\setminus w_n$ such that $x\d a\d R\d b\d y$ is a path. But then
$((A \cup \{x\},\emptyset,B\cup \{y\}),a\d R\d b)$
is a staircase, contrary to the maximality of $(S,R_0)$. So $R$ has length 1, and there exists $i$
with $1 \le i < n$ such that $a\d w_1\d \cdots\d w_i\d b$ is an odd antipath. But then the theorem holds.
This proves (6).
\\
\\
(7) {\em If no vertex in $R$ is $W$-complete then the theorem holds.}
\\
\\
For by (6) we may assume that there is a vertex $v$ of $R$ which is $(W\setminus w_n)$-complete.
Hence $v$ is nonadjacent to $w_n$. Since $n\ge 3$ and is odd, and $a\d w_1\d \cdots\d w_n\d v\d a$
is not an odd antihole, it follows that $v$ is adjacent to $a$. Consequently $v$ is the unique
$(W\setminus w_n)$-complete vertex in $R$. From (6) we may assume that $v= b$, and $R$ has length 1.
Choose a step $a_1\d R_1\d b_1,a_2\d R_2\d b_2$. Then $b_1\d a\d w_1\d \cdots\d w_n\d b$ is an odd antipath, of
length $\ge 5$. All its internal vertices have neighbours in the connected set $V(R'\setminus w_n)\cup \{a_2\}$,
and its ends do not. By \ref{RR} applied in $\overline{G}$, there are adjacent vertices $x,y$ in
$V(R'\setminus w_n)\cup \{a_2\}$, such that $x\d a\d w_1\d \cdots\d w_n\d y$ is an odd antipath. Since $x$ is
adjacent to $w_n$, it follows that $x$ is the neighbour of $w_n$ in $R'$, and therefore either
$y$ is the second neighbour of $x$ in $R'$, or $R'$ has length 1 and $y = a_2$. Assume first that
$R'$ has length $>1$, and so both $x,y$ belong to the interior of $R'$. Hence $x,y$ are both
anticomplete to $A \cup B$, and so $((B\cup \{x\}, \emptyset, A \cup \{y\}),a\d w_1\d \cdots\d w_n)$ is a
staircase in $\overline{G}$, contradicting that $(S,R_0)$ is strongly maximal. Now assume that
$R'$ has length 1. Then $x = r'$ and $y = a_2$, and
$((B\cup \{r'\}, \emptyset, A \cup \{b\}),a\d w_1\d \cdots\d w_n)$ is a
staircase in $\overline{G}$, a contradiction as before. This proves (7).

\bigskip

We may therefore assume that some vertex of $R\setminus b$ is $W$-complete, for otherwise the theorem
holds by (7).
Let $a\d P\d p$ be the subpath of $R\setminus b$ such that $p$ is the unique $W$-complete vertex of $P$.
Choose $a_1 \in A$ and $b_1 \in B$, adjacent (this is possible by (5)).
Let us apply \ref{evenantipath2} to the path $p\d P\d a\d a_1\d b_1$, and the even antipath
$a\d w_1\d \cdots\d w_n\d a_1$. Both ends of the path are complete to the interior of the antipath, so
by \ref{evenantipath2} it follows that $P$ has length 2, and if $q$ denotes its middle vertex then
$q$ is nonadjacent to $w_n$ and adjacent to $w_1,\ldots,w_{n-1}$. But then
$((B \cup \{p\},\emptyset,A\cup \{q\}),a\d w_1\d \cdots\d w_n)$ is a staircase in $\overline{G}$,
a contradiction. This completes the proof of \ref{trajectory}. \bbox

\begin{thm}\label{hitbanister}
Let $G$ be Berge, containing no appearance of $K_4$, no even prism, no 1-breaker and no 2-breaker.
Let $K = (S = (A,C,B), a_0\d R_0\d b_0)$ be a strongly maximal staircase in $G$,
and let $x_1,\ldots,x_t$ be a right-sequence. Then $x_1,\ldots,x_t$ are all adjacent to $a_0$.
\end{thm}
\Proof Suppose the theorem is false, and choose $t$ is small as possible so that the statement of the
theorem does not hold. So $t \ge 1$, and $x_1,\ldots,x_{t-1}$ are all adjacent to $a_0$, and $x_t$ is
not.
\\
\\
(1) {\em $a_0\d R_0\d b_0$ is not an optimal banister for $b_0$.}
\\
\\
For suppose it is, and let $a_0\d w_1\d \cdots\d w_n$ be the trajectory of $a_0$. Since $R_0$ has length $>1$
it follows from \ref{trajectory} that $b_0$ is the unique $W$-complete vertex of $R_0$, where
$W = \{w_1,\ldots,w_n\}$. Suppose first that $n = 1$. Then $w_1$ is nonadjacent to $a_0$ and has
a nonneighbour in $A$, and so by \ref{stairvnbrs} it is a right-star. By axiom 3 there is a banister
$r\d R\d w_1$ such that the birth of $r$ is earlier than $w_1$. Since $a_0\d R_0\d b_0$ is optimal for
$b_0$, it follows as in the proof of \ref{trajectory} that $R$ is disjoint from $R_0$, and there
are no edges between $R_0 \setminus a_0$ and $R \setminus w_1$. The only edge from $w_1$ to $R_0$
is $w_1b_0$, by \ref{trajectory}. Choose an $S$-rung $a_1\d R_1\d b_1$.
Since $a_1\d a_0\d R_0\d b_0\d w_1\d R\d r\d a_1$ is not an
odd hole it follows that $a_0$ has neighbours in $R$. If it has a neighbour different from $r$, then
the path from $a_0$ to $w_1$ with interior in $R\setminus r$ can be completed via $w_1\d b_0\d R_0\d a_0$
and via $w_1\d b_1\d R_1\d a_1\d a_0$, and one of the resulting holes is odd, a contradiction.
So the unique neighbour
of $a_0$ in $R$ is $r$. But then we can add $r$ to $A$, $w_1$ to $B$ and the interior of $R$ to $C$,
contradicting the maximality of $(S,R_0)$.
So $n \ge 2$. Now all of $w_1,\ldots,w_{n-1}$ are left-diagonals, and all of $w_2,\ldots,w_n$ are
right-diagonals. But $w_1$ is not a right-diagonal, and $w_n$ is not a
left-diagonal, and $w_1$ is not a right-star, contrary to \ref{diagonalchain}. This proves (1).

\bigskip

Now since $a_0$ has a nonneighbour in $\{x_1,\ldots,x_t\}$, it follows that there is an optimal
banister $r\d R\d b_0$ for $b_0$.
From (1), $r$ has a nonneighbour in $\{x_1,\ldots,x_{t-1}\}$. From the minimality of $t$
(replacing $R_0$ by $R$) it follows that $R$ has length 1, and so $rb_0$ is an edge.
Let $r\d w_1\d \cdots\d w_n$ be the trajectory of $r$; so $w_1$ is earlier than $x_t$.
Let $W = \{w_1,\ldots,w_n\}$; hence $a_0$ is $W$-complete. By \ref{trajectory}, $n$ is odd.
\\
\\
(2) {\em $b_0$ is $W$-complete.}
\\
\\
For suppose not. Then by \ref{trajectory}, there exists $i$ with $1 \le i < n$ so that
$r\d w_1\d \cdots\d w_i\d b_0$ is an
odd antipath. Now $r,w_1,\ldots,w_{i-1}$ are all left-diagonals; $w_1,\ldots,w_i$ are all right-diagonals;
$r$ is not a right-diagonal (since it is a left-star); and $w_i$ is not a left-diagonal (since it is
nonadjacent to $b_0$) and not a right- or left-star (since it is $A \cup B$-complete, because $i < n$).
This contradicts \ref{diagonalchain}, and so proves (2).
\\
\\
(3) {\em $a_0$ is adjacent to $r$, and $w_n$ is a right-star.}
\\
\\
Let $a_1\d R_1\d b_1$ be an $S$-rung with $w_n$ nonadjacent to $a_1$.
Since $a_0\d r\d w_1\d \cdots\d w_n\d a_1\d b_0\d a_0$ is not an odd antihole it follows that $a_0$ is adjacent to $r$.
So each of $r,w_1,\ldots,w_{n-1}$ is left-diagonal, each of $w_1,\ldots,w_n$ is right-diagonal,
$r$ is not right-diagonal, $w_n$ is not left-diagonal, and the claim follows from \ref{diagonalchain}.
This proves (3).
\\
\\
(4) {\em There is no $(W \cup \{r\})$-complete vertex in the interior of $R_0$.}
\\
\\
For suppose there is, $v$ say. Let $a_1\d R_1\d b_1$ be an $S$-rung. Then $a_0\d a_1\d R_1\d b_1\d b_0$ is an odd
path; both its ends are $(W \cup \{r\})$-complete; and the $(W \cup \{r\})$-complete vertex $v$ has
no neighbour in its interior, so by \ref{greentouch} there is a $(W \cup \{r\})$-complete vertex in
$R_1$. But by (3), $w_n$ is a right-star and $r$ is a left-star, so they have no
common neighbour in $R_1$, a contradiction. This proves (4).
\\
\\
(5) {\em $n = 1$.}
\\
\\
For assume $n > 1$.  Now $R_0$ is odd, and both its ends are $(W \cup \{r\})$-complete.
Suppose first that $R_0$ has length $\ge 5$. By \ref{RR} and (4) there is a leap; that is, there are
two nonadjacent vertices $x,y \in W \cup \{r\}$ joined by an odd path $P$ whose interior
is the interior of $R_0$. Choose $b_1 \in B$; then $b_1\d x\d P\d y\d b_1$ is not an odd hole, and so one of
$x,y$ is nonadjacent to $b_1$. Since $b_1$ is $W$-complete, we may assume $y = r$; and hence $x = w_1$ since
that is the only vertex in $W$ nonadjacent to $r$. Choose $a_1 \in A$; then since
$a_1\d r\d P\d w_1\d a_1$ is not
an odd hole it follows that $a_1$ is not adjacent to $w_1$ and so $n = 1$. Now assume that $R_0$
has length 3, and let its internal vertices be $x,y$ (in some order). By \ref{RR} there exists
an odd antipath $Q$ joining $x,y$ with interior in $W \cup \{r\}$. If $r \not \in V(Q)$ then
$b_1\d x\d Q\d y\d b_1$ is an odd antihole, where $b_1 \in B$; and if $w_n \not \in V(Q)$ then
$a_1\d x\d Q\d y\d a_1$ is an odd antihole, where $a_1 \in A$. Hence we may assume that
$x\d r\d w_1\d \cdots\d w_n\d y$ is an antipath.
We claim that $C = \emptyset$. For suppose there is an $S$-rung $a_1\d R_1\d b_1$ say of length $>1$.
Then $a_1\d R_1\d b_1\d b_0\d r\d a_1$ is a hole of length $\ge 6$; and $r\d w_1\d \cdots\d w_n\d a_1$ is an even antipath of
length $\ge 4$; and $a_0$ is complete to the antipath, and has no other neighbours on the
hole; and at least two vertices of the hole are complete to the interior of the antipath,
namely $b_0$ and $b_1$. This contradicts \ref{hole&antipath}. So $C = \emptyset$. Hence
$((B \cup \{x\}, \emptyset, A \cup \{y\}),r\d w_1\d \cdots\d w_n)$ is a staircase in $\overline{G}$,
a contradiction. This proves (5).

\bigskip

From (4), (5) we may apply \ref{RR} to $R_0$ and the anticonnected set $\{r,w_1\}$, and since the latter
has only two members, \ref{RR} implies that there is an odd path $P$ joining $r$ and $w_1$ with interior
equal to the interior
of $R_0$. From (3), $w_1$ is a right-star, and from axiom 3 there is a banister $r'\d R'\d w_1$ (and we may
choose it optimal for $w_1$) such that the birth of $r'$ is earlier than $w_1$. Now $R'$ is
disjoint from $R_0$, and there are no edges between $R_0 \setminus a_0$ and $R' \setminus w_1$; for otherwise
there would be a banister from $r'$ to $b_0$, contradicting that $r\d b_0$ is optimal for $b_0$.
Suppose that $r$ has a neighbour in $R'$; then the path between $r$ and $w_1$ with interior in $R'$
can be completed to holes via $w_1\d b_0\d r$ and via $w_1\d P\d r$, a contradiction since one of these holes is odd.
So $r$ has no neighbour in $R'$. Let $r'\d v_1\d \cdots\d v_m$ be the trajectory of $r'$. Since
$v_1,\ldots,v_m$ are earlier than $w_1$, and $w_1$ is the earliest nonneighbour of $r$, it follows
that $r$ is adjacent to all of $v_1,\ldots,v_m$. Now by \ref{trajectory}, either
\begin{itemize}
\item $w_1$ is the unique $\{v_1,\ldots,v_m\}$-complete vertex in $R'$; but then $w_1\d R'\d r'\d a_1\d r$
(where $a_1 \in A$ is nonadjacent to $v_m$) is an odd path; its ends are $\{v_1,\ldots,v_m\}$-complete
and its internal vertices are not; and the $\{v_1,\ldots,v_m\}$-complete vertex $b_1$ (for any $b_1 \in B$
nonadjacent to $a_1$) has no neighbour in its interior, contrary to \ref{greentouch}.
\item $R'$ has length 1, and there is an odd antipath $Q$ between $r'$ and $w_1$ with interior in
$\{v_1,\ldots,v_m\}$; but then $r\d r'\d Q\d w_1\d r$ is an odd antihole, a contradiction.
\end{itemize}
This completes the proof of \ref{hitbanister}.\bbox

Now we are ready to apply \ref{hitbanister} to produce a skew partition.
Let us say a {\em 3-breaker} in $G$ is a pair $(K,x)$ such that $K = (S = (A,C,B), a_0\d R_0\d b_0)$
is a strongly maximal staircase in $G$, and $x\in V(G) \setminus V(K)$ is $B$-complete, and not
$A$-complete, and not $A$-anticomplete.

\begin{thm}\label{mixedprism}
Let $G$ be Berge, containing no appearance of $K_4$, no even prism, no 1-breaker and no 2-breaker.
Suppose that there is a 3-breaker in $G$; then $G$ admits a balanced skew partition.
\end{thm}
\Proof Let $(K,x_1)$ be a 3-breaker, where $K = (S = (A,C,B), a_0\d R_0\d b_0)$.
The 1-vertex sequence $x_1$ is a right-sequence; so there exists a right-sequence $x_1,\ldots,x_t$
of maximum length, with $t \ge 1$. Let $X = \{x_1,\ldots,x_t\}$, and let $Y$ be the set of all
$A\cup X$-complete vertices in $V(G)\setminus V(S)$. So $a_0 \in Y$ by \ref{hitbanister}.
\\
\\
(1) {\em $X \cup Y \cup B$ meets the interior of every path in $G$ from $A\cup C$ to $b_0$.}
\\
\\
For suppose $P$ is a path from $A\cup C$ to $b_0$ with no internal vertex in $X\cup Y\cup B$. Note
that $b_0 \not \in X$ by \ref{hitbanister}, and so $b_0 \notin X \cup Y\cup B$ (since it is not $A$-complete).
We may assume $P$ is minimal, and therefore no internal vertex of $P$ is in $V(S)$. Let $P$ be from
$p \in A \cup C$ to $b_0$.  By \ref{prisminstair2},
$P\setminus p$ contains either a major vertex or a banister. Suppose first that it contains a banister $a\d R\d b$
say. Hence $a,b \notin  X \cup Y\cup B$. Since $a$ is $A$-complete it is therefore not $X$-complete
(because it is not in $Y$), and then we can set $x_{t+1} = b$, contradicting the maximality
of the right-sequence.  So $P\setminus p$ contains no banister. Now assume it contains a major
vertex $v$ say. Since $v \notin X\cup Y\cup B$, it follows that $v$ is not
$X\cup A$-complete. Suppose $v$ is $B$-complete. Since it is major it has a neighbour in $A$. If it is
not $A$-complete we can set $x_{t+1} = v$ and obtain a longer right-sequence, a contradiction; and if
$v$ is $A$-complete then since it is not $X \cup A$-complete, it is not $X$-complete and so again we
can set $x_{t+1} = v$ and obtain a longer right-sequence, a contradiction. So $v$ is not $B$-complete.
By \ref{stairvnbrs} and since there is no 2-breaker in $G$ and therefore no central vertex,
$v$ is left-diagonal, and not right-diagonal; and since it
is not $X\cup A$-complete, it is not $X$-complete. Let $v\d w_1\d \cdots\d w_n$
be the trajectory of $v$. Then each of $w_1,\ldots,w_n$ is right-diagonal, since they are all
$B \cup \{a_0\}$-complete. Since $w_n$ has a nonneighbour in $A$, it is not left-diagonal; and so there
is a minimum $i$ with $1 \le i\le n$ such that $w_i$ is not left-diagonal. By \ref{diagonalchain}
applied to the sequence $v,w_1,\ldots,w_i$, we deduce that $v$ is a left-star,
contradicting that $v$ is major. This proves (1).

\bigskip

Now since $S$ is step-connected, it follows that $A \cup C$ is connected; and therefore belongs to
a component $A_1$ of $G \setminus (X \cup Y \cup B)$. Let $A_2$ be the union of all
the other components. So by (1), $b_0 \in A_2$ , and $(A_1 \cup A_2, X \cup Y \cup B)$
is a skew partition of $G$ (since $Y \cup B$ is complete to $X$, and $X$ is nonempty). We need to find
a balanced skew partition. By \ref{geteven} we may assume
this skew partition is not loose; so every $X$-complete vertex in $G$ either belongs to $B$ or is also
$A$-complete. Every vertex in $Y \cup B$ has a neighbour in $A \cup C$, so $A \cup C$ is a kernel
for this skew partition, in $\overline{G}$. By \ref{kernel2} it suffices to show that in $G$,
any two nonadjacent
vertices in $Y \cup B$ are joined by an even path with interior in $A \cup C$, and any two adjacent
vertices of $A \cup C$ are joined by an even antipath with interior in $Y \cup B$. Now let $u,v \in Y \cup B$
be nonadjacent. If they are both adjacent to $b_0$, then any path joining them with interior in $A \cup C$
(and there is one) is even, since it can be completed to a hole via $v\d b_0\d u$. So we may assume that $u$ is
nonadjacent to $b_0$, and hence $u \notin B$, so $u \in Y$. If they are both in $Y$, then they are joined by an
even path $u\d a_1\d v$ for any $a_1 \in A$. So we may assume that $v \in B$. Since $u$ is nonadjacent to
$b_0$ and to $v$, it is neither left- nor right-diagonal, and it is not central since there is no
2-breaker; so from \ref{stairvnbrs} it is a left-star. Let $a_1\d R_1\d v$ be an $S$-rung; then
$u\d a_1\d R_1\d v$ is the desired even path between
$u$ and $v$. Now for antipaths, let $uv$ be an edge with $u,v \in A\cup C$. They both therefore have
nonneighbours in $B$, and since $B \cup \{a_0\}$ is anticonnected, they are joined by an antipath $Q$ with
interior in $B \cup \{a_0\}$. It suffices to show that $Q$ is even, since $Q^* \subseteq Y \cup B$. If
$a_0 \notin Q^*$, then $Q$ is even since $b_0\d u\d Q\d v\d b_0$ is an antihole. So $a_0$ is in $Q^*$. But there
are no edges between $a_0$ and $B$, and so $a_0$ is nonadjacent to every other vertex in the interior of $Q$;
and since $Q$ is an antipath, it therefore has at most 3 internal vertices, so its length is $\le 4$.
If it is odd, then it has length 3, that is, there are nonadjacent vertices $u' \in Y$ and $v' \in B$,
joined by an odd path with interior in $A \cup C$. But we have already shown that they are joined by an
even path, and the result follows from \ref{mixedpair}. This proves \ref{mixedprism}.\bbox

Now we can prove \ref{summary}.5, the main result of this section. We restate it (M-joins were defined
in section 1.)

\begin{thm}\label{longprism}
Let $G$ be Berge, such that there is no appearance of $K_4$ in either $G$ or $\overline{G}$.
Suppose that $G$ contains a long odd prism as an induced subgraph. Then either one of $G, \overline{G}$
admits a 2-join, or $G$ admits a balanced skew partition, or $G$ admits an M-join.
\end{thm}
\Proof
We assume that $G$ does not admit a balanced skew partition, and $G, \overline{G}$ do not admit 2-joins.
Since $G$ contains a long odd prism, and therefore $G,\overline{G}$ are not even prisms,
it follows from \ref{evenprism} that $G,\overline{G}$ contain no even prism. By \ref{1-hat},
\ref{central} and \ref{mixedprism}, $G,\overline{G}$ contain no 1-, 2- or 3-breaker.

Since $G$ contains a long odd prism, it contains a staircase; and therefore (possibly by replacing
$G$ by its complement) there is a strongly maximal staircase $K = (S = (A,C,B), a_0\d R_0\d b_0)$ say in $G$.
Let $A_0$ be the set of all left-stars, $B_0$ the set of all right-stars, and $N$ the set of all
vertices that are $A \cup B$-complete. By \ref{stairvnbrs}, every non-major $A$-complete vertex is in $A_0$,
and since there is no 3-breaker, every major $A$-complete vertex is in $N$, so every $A$-complete
vertex is in $A_0 \cup N$; and similarly
every $B$-complete vertex is in $B_0 \cup N$. Let $H = G \setminus (V(S) \cup A_0 \cup B_0 \cup N)$.
\\
\\
(1) {\em Let $F$ be a component of $H$, and let $X$ be the set of attachments of $F$ in $V(S)\cup A_0
\cup B_0$. Then either
$X \cap V(S) = \emptyset$, or $X \subseteq V(S)$ and $X$ meets both $A \cup C$ and $B \cup C$.}
\\
\\
We may assume that $X$ meets $V(S)$, and therefore from the symmetry we may assume that $X$ meets $A \cup C$.
Since no vertex in $F$ is $A$- or $B$-complete, and therefore no vertex in $F$ is major or a left- or
right-star, it follows from \ref{prisminstair2} that $X$ is disjoint from $B_0$. If
$X$ meets $B \cup C$ then similarly $X$ is disjoint from $A_0$, and so $X \subseteq V(S)$
and the claim holds. We assume therefore that $X \subseteq A \cup A_0$. Suppose there is a vertex
$v$ of $G$, not in $A \cup A_0 \cup N \cup F$ but with a neighbour in $F$. Then $v \notin V(H) \cup X$, and so
$v \in B_0$. By \ref{prisminstair2} applied to $F \cup \{v\}$, it follows that $F \cup \{v\}$ contains
a major vertex or a left-star, a contradiction. So there is no such $v$. Hence
$(V(G)\setminus (A \cup A_0 \cup N), A \cup A_0\cup N)$ is a skew partition of $G$, since $F$ is a
component of $V(G)\setminus (A \cup A_0 \cup N)$ and $b_0$ is in a different component, and $A,A_0\cup N$ are
both nonempty and complete to each other. Now by \ref{balancev}, $(B\cup C,A)$ is balanced, since
$a_0$ is complete to $A$ and anticomplete to $B \cup C$; and therefore from \ref{shiftbalance},
$(F,A)$ is balanced (since $B \cup C$ is connected and all vertices in $A$ have neighbours in it).
Hence from \ref{onepair}, $G$ admits a balanced skew partition, a contradiction. This proves (1).

\bigskip

Let $M$ be the union of all components of $H$ with no attachment in $V(S)$. Then $M$ is nonempty,
since it contains the interior of $R_0$.
Let $D$ be the union of all the other components of $H$.
Hence $V(G)$ is partitioned into $A,B,C,D,A_0,B_0,N,M$, where possibly $C,D$ or $N$ may be empty.
\\
\\
(2) {\em $N \not = \emptyset$.}
\\
\\
For assume that $N = \emptyset$. Then the only edges between
$V(S) \cup D$ and $A_0 \cup B_0 \cup M$ are the edges from $A$ to $A_0$ and those from $B$ to $B_0$;
and since $A_0 \cup B_0 \cup M$ contains at least 4 vertices (the vertices of $R_0$) and both $A$ and $B$
contain at least two, this is a 2-join in $G$, a contradiction. This proves (2).
\\
\\
(3) {\em $C \cup D= \emptyset$.}
\\
\\
For assume that $C \cup D$ is nonempty. By (1) there are no edges between $C\cup D $ and $A_0 \cup B_0 \cup M$.
Since $N$ is complete to $A \cup B$, it follows that
$(C\cup D \cup A_0 \cup B_0 \cup M, N \cup A \cup B)$ is a skew partition of $G$. By \ref{geteven},
it is not loose, and so there is no $N'$-complete vertex in $R_0$, where $N'$ is an anticomponent
of $N$. Let $a_1\d R_1\d b_1,a_2\d R_2\d b_2$ be a step; then $a_1\d a_0\d R_0\d b_0\d b_2$ is an odd path of length
$\ge 5$; its ends are $N'$-complete, and its internal vertices are not. By \ref{RR}, there is a leap
in $N'$, and so there exist nonadjacent $x,y$ in $N$ so that $x\d a_0\d R_0\d b_0\d y$ is a path. But then
$((A \cup \{x\},C,B\cup \{y\}), a_0\d R_0\d b_0)$ is a staircase, contradicting the maximality of
$(S,R_0)$. This proves (3).

\bigskip

But then the six sets $A,B,A_0,B_0,M,N$ form an M-join in $G$. This proves \ref{longprism}.\bbox

This is the only place in the entire paper where we use M-joins. It is natural to ask whether M-joins are
really necessary, or whether \ref{decomp} remains true if we omit them. It appears that the second holds;
one of us (Chudnovsky) has what seems to be a proof, which if correct will appear in her PhD thesis. But
in this paper we accept M-joins.

Let us say a graph $G$ is {\em recalcitrant} if:
\begin{itemize}
\item $G$ is Berge
\item $G$ and $\overline{G}$ are not line graphs, and $G$ is not a bicograph
\item $G$ and $\overline{G}$ do not admit 2-joins, and
\item $G$ does not admit an M-join or balanced skew partition.
\end{itemize}

The remainder of the paper is basically a proof of the following.
\begin{thm}\label{recalcitrant}
If $G$ is recalcitrant then either $G$ or $\overline{G}$ is bipartite.
\end{thm}
Clearly any counterexample to \ref{decomp} is recalcitrant, so \ref{recalcitrant} will imply
\ref{decomp}.

On the other hand, for some purposes (eliminating M-joins, and a possible recognition algorithm for perfection) it
is desirable to keep closer track of which results hold under which hypotheses, instead
of just using the blanket ``recalcitrant'' hypothesis. But at least, for the remainder of the paper
we shall only be concerned with Berge graphs $G$ such that in both $G,\overline{G}$ there is no
appearance of $K_4$ and no long prism; that is, with the members of the class $\mathcal{F}_5$
introduced in section 1. Certainly every recalcitrant graph belongs to  $\mathcal{F}_5$, by
\ref{evenprism} and \ref{bicographs}.

It turns out that for such graphs, there is a useful
strengthening of \ref{RR} --- the second alternative of that theorem can no longer hold.
\begin{thm} \label{RRR}
Let $G\in \mathcal{F}_5$, and let $P$ be a path in $G$ with odd length.
Let $X \subseteq V(G)$ be anticonnected, so that both ends of $P$ are $X$-complete. Then either:
\begin{enumerate}
\item some edge of $P$ is $X$-complete, or
\item $P$ has length $3$ and there is an odd antipath joining the internal vertices of $P$
with interior in $X$.
\end{enumerate}
\end{thm}
\Proof Let $P$ be $p_1\d \cdots\d p_n$.  By \ref{RR}, we may assume that $P$ has length $\ge 5$ and
$X$ contains a leap $u,v$ say; so $u\d p_2\d \cdots\d p_{n-1}\d v$ is a path. But then the three paths
$p_1\d v,u\d p_n, p_2\d \cdots\d p_{n-1}$ form a long prism, contrary to $G\in \mathcal{F}_5$.
This proves \ref{RRR}.\bbox

There is an analogous strengthening of \ref{doubleRR}, as follows.

\begin{thm}\label{doubleRR2}
Let $G\in \mathcal{F}_5$, and let $X, Y$ be disjoint nonempty anticonnected subsets of $V(G)$,
complete to each other.  Let $P$ be a path in $G$ with even length $>0$, with vertices
$p_1,\ldots,p_n$ in order, so that $p_1$ is the unique
$X$-complete vertex of $P$ and $p_n$ is the unique $Y$-complete vertex of $P$.  Then
$P$ has length $2$ and there is an antipath $Q$ between $p_2$ and $p_3$ with interior in $X$,
and an antipath $R$ between $p_1$ and $p_2$ with interior in $Y$, and exactly one of $Q$,$R$ has
odd length.
\end{thm}
\Proof
Let us apply \ref{doubleRR}. We may therefore assume that
 $P$ has length $\ge 4$ and there are nonadjacent $x_1,x_2 \in X$ so that
$x_1\d p_2\c p_n\d x_2$ is a path $P'$ say, of odd length $\ge 5$. But the ends of $P'$ are
$Y \cup \{p_1\}$-complete, and its internal vertices are not, contrary to \ref{RRR}.
This proves \ref{doubleRR2}.\bbox

\section{The double diamond}

We are finished with prisms --- we cannot dispose of the prism where all three paths have length 1 (yet), and
we have disposed of all others. Now we turn to a different type of subgraph, the double diamond.

Let $G$ be Berge. If $A,B$ are disjoint subsets of $V(G)$, we say a {\em square} in $(A,B)$ is a hole
$a_1\d b_1\d b_2\d a_2\d a_1$ of length 4, where $a_1,a_2 \in A$ and $b_1,b_2 \in B$. The pair $(A,B)$ is
{\em square-connected} if:
\begin{itemize}
\item  $|A|,|B| \ge 2$,
\item for every partition $(X,Y)$ of $A$ with $X,Y$ nonempty, there is a square $a_1\d b_1\d b_2\d a_2\d a_1$
with $a_1 \in X$ and $a_2 \in Y$
\item for every partition $(X,Y)$ of $B$ with $X,Y$ nonempty, there is a square $a_1\d b_1\d b_2\d a_2\d a_1$
with $b_1 \in X$ and $b_2 \in Y$.
\end{itemize}
It follows that if $(A,B)$ is square-connected then every vertex of $A\cup B$ is in a square. An
{\em antisquare} is a square in $\overline{G}$; that is, an antihole $a_1\d b_1\d b_2\d a_2\d a_1$ with
$a_1,a_2 \in A$ and $b_1,b_2 \in B$; and $(A,B)$ is {\em antisquare-connected} if it is square-connected
in  $\overline{G}$.
For strips in which every rung has length 1 (and from now on, those are the only kind of strips we
shall need), being square-connected is the same as being
step-connected. We have renamed the concepts because we wanted to improve our notation for a
step.

We say a quadruple $(A,B,C,D)$ of subsets of $V(G)$ is a {\em cube} in $G$ if it
satisfies the following conditions:
\begin{itemize}
\item $A,B,C,D$ are pairwise disjoint and nonempty
\item $A$ is complete to $C$, and $B$ to $D$, and $A$ is anticomplete to $D$, and $B$ to $C$
\item $(A,B)$ is square-connected, and $(C,D)$ is antisquare-connected.
\end{itemize}

If $G$ contains a double diamond, then it contains a cube in which $A,B,C,D$ all have two elements,
and that turns out to be the right approach to the double diamond --- grow the cube until it is maximal, and
analyze how the remainder of $G$ attaches to it. That is our goal in this section. A cube $(A,B,C,D)$
is {\em maximal} if there is no cube $(A',B',C',D')$ with $A \subseteq A', B \subseteq B', C \subseteq C',
D \subseteq D'$ such that $(A,B,C,D) \not = (A',B',C',D')$. The subgraph $G | (A \cup B \cup C \cup D)$ is
called the graph {\em formed} by the cube. Note that if  $(A,B,C,D)$ is a cube in $G$, then
$(C,D,B,A)$ is a cube in $\overline{G}$. (This is very convenient, because it reduces our work by half -
we are going to have the usual minor vertices and major vertices, and whatever we can prove about minor
vertices is also true in the complement for major ones.)

\begin{thm}\label{cubevnbrs}
Let $G\in \mathcal{F}_5$. Let $(A,B,C,D)$ be a maximal cube in $G$, forming $K$,
let $v \in V(G) \setminus V(K)$, and let $X$ be the set of neighbours of $v$ in $V(K)$. Then either
\begin{itemize}
\item $X$ is a subset of one of $A \cup B, C\cup D, A \cup C, B \cup D$, and
$X \cap (A \cup C)$  is complete to $X \cap (B \cup D)$, or
\item $X$ includes one of $A \cup B, C\cup D, A \cup D, B \cup C$, and
$(A \cup D)\setminus X$ is anticomplete to $(B \cup C)\setminus X$.
\end{itemize}
\end{thm}
\Proof Note that under taking complements the two outcomes become exchanged.
If $X \subseteq A\cup B$, and there exists $a \in X \cap A$ and $b \in X \cap B$, nonadjacent, then
choose $c \in C$ and $d \in D$, adjacent, and $v\d a\d c\d d\d b\d v$ is an odd hole. So if $X \subseteq A\cup B$
then the theorem holds. Similarly it holds if $X \subseteq C \cup D$; and trivially it holds if
$X$ is a subset of one of $A \cup C, B \cup D$. So we may assume that $X$ meets both $A$ and $D$.
From the same argument in $\overline{G}$, we may also assume that none of
$A \cup B, C\cup D, A \cup D, B \cup C$ is a
subset of $X$, that is, either $X$ includes neither of $A$,$C$ or it includes neither of $B$,$D$.
These two cases are exchanged when we pass to the complement; so we may assume by taking complements
that $X$ includes neither of $B,D$. Let $A_1 = A \cap X$, and $A_2 = A \setminus A_1$; and define $B_1,B_2$
etc. similarly. We have shown so far that $A_1,B_2,D_1,D_2$ are nonempty.
Choose an antisquare $c_2\d d_1\d d_2\d c_1\d c_2$ such that $d_1 \in D_1$ and $d_2 \in D_2$, and choose
$b_2 \in B_2$. Since $v\d c_2\d d_2\d b_2\d d_1\d v$ is not an odd hole, it follows that $c_2 \in C_2$.
Hence $A_1$ is complete to $B_1$; for if $a_1 \in A_1$ and $b_1 \in B_1$ are nonadjacent
then $v\d a_1\d c_2\d d_2\d b_1\d v$ is an odd hole. If $A_1 = A$, then since $(A,B)$ is square-connected
and $A_1$ is complete to $B_1$ it follows
that $B_1$ is empty; but then we can add $v$ to $C$ (because $v\d d_2\d d_1\d c_2\d v$ becomes a
new antisquare), contrary to the maximality of the cube. So $A_2$ is nonempty. Hence there is a square
$a_1\d b_1\d b_2\d a_2\d a_1$
with $a_1 \in A_1$ and $a_2 \in A_2$. Since $a_1$ is nonadjacent to $b_2$ and complete to $B_1$,
it follows that $b_2 \in B_2$; but then $v\d a_1\d a_2\d b_2\d d_1\d v$ is an odd hole, a contradiction.
This proves \ref{cubevnbrs}.\bbox

Say a vertex $v \in V(G) \setminus V(K)$ is {\em minor} if the first case of \ref{cubevnbrs} applies
to it, and {\em major} if the second case applies. Then every such vertex is either minor or major and
not both; and by taking complements, the minor and major vertices are exchanged.

\begin{thm}\label{cubeFnbrs}
Let $G\in \mathcal{F}_5$. Let $(A,B,C,D)$ be a maximal cube in $G$, forming $K$,
let $F \subseteq V(G) \setminus V(K)$ be a connected set of minor vertices, and let $X$ be the set
of attachments of $F$ in $V(K)$. Then $X$ is a subset of one of $A \cup B, C\cup D, A \cup C, B \cup D$.
Moreover, $X \cap (A \cup C)$ is complete to $X \cap (B \cup D)$.
\end{thm}
\Proof Suppose the first assertion is false,
and choose $F$ minimal with this property. We may assume that $X$ meets both of $A,D$.
Since all vertices in $F$ are minor, it follows that $F$ is a path $f_1\d f_2\d \cdots\d f_k$ of length $\ge 1$.
We may assume $f_1$ is the unique vertex of $F$ with a neighbour in $A$, and $f_k$ is the unique vertex of $F$
with a neighbour in $D$. Let $X_1$,$X_2$ be the sets of attachments in $V(K)$ of $F \setminus f_k$,
$F \setminus f_1$ respectively. From the minimality of $F$ it follows that $X_1$ is a subset of one of
$A\cup B, A \cup C$, and $X_2$ is a subset of one of $B\cup D, C \cup D$.
\\
\\
(1) {\em Not both $X_1 \subseteq A \cup B$ and $X_2 \subseteq B\cup D$ .}
\\
\\
For suppose that both these hold. If $k$ is even, choose $a\in A$ is adjacent to $f_1$, and $d \in D$
is adjacent to $f_k$, and $c\in C$ is adjacent to $d$; then $a\d f_1\d \cdots\d f_k\d d\d c\d a$
is an odd hole , a contradiction. So $k$ is odd.
Suppose first that $f_1$ is complete to $A$. Since it is minor, it has
no neighbours in $B$ (for no vertex in $B$ is $A$-complete). If there are no edges between $B$ and
$F$, let $a_1\d b_1\d b_2\d a_2\d a_1$ be a square, and let $d \in D$ be adjacent to $f_k$; then
$a_1\d b_1,a_2\d b_2,f_1\c f_k\d d$ form a long prism, a contradiction.
So there are edges between $B$ and $F$. Choose $i$ with $1 \le i \le k$ minimum so that
$f_i$ has a neighbour in $B$. If $f_i$ is not complete to $B$, choose a square $a_1\d b_1\d b_2\d a_2\d a_1$
so that $f_i$ is adjacent to $b_1$ and not to $b_2$; then $b_1$ can be linked onto the triangle
$\{f_1,a_1,a_2\}$, via $b_1\d f_i\c f_1$, $b_1\d a_1$, $b_1\d b_2\d a_2$, contrary to \ref{trianglev}.
So $f_i$ is complete to $B$. Let $a_1\d b_1\d b_2\d a_2\d a_1$ be a square; then since
$a_1\d b_1,a_2\d b_2, f_1 \c f_i$ do not form
a long prism (because $G\in \mathcal{F}_5$), it follows that $i = 2$. But $k >2$ since $k$ is odd; so we can add
$f_1$ to $C$ and $f_2$ to $D$, contrary to the maximality of the cube. This proves (1) if
$f_1$ is $A$-complete. Now assume $f_1$ is not $A$-complete, and choose a square
$a_1\d b_1\d b_2\d a_2\d a_1$ so that $f_1$ is adjacent to $a_1$ and not to $a_2$. Since
$a_1\d f_1\d \cdots\d f_k\d d\d b_2\d a_2\d a_1$ is not an odd hole (where $d \in D$ is adjacent to $f_k$), it follows
that $b_2$ has a neighbour in $F$. Choose $i$ minimum so that $b_2$ is adjacent to $f_i$. Let
$c \in C$ and $d \in D$ be any adjacent pair of vertices. Then the three paths $a_1\d b_1,a_2\d b_2,c\d d$
form a prism, and since the set of attachments of $\{f_1,\ldots,f_i\}$ in this
prism is not local, and does not include $a_2$, it has an attachment in the third path $c\d d$,
by \ref{legaljump}; and hence
$i = k$, and $f_k$ is $D$-complete. Again, let $c \in C$ and $d \in D$ be adjacent. Then the prism
formed by $a_1\d f_1\d \cdots\d f_k$,$a_2\d b_2$,$c\d d$ is long, contrary to $G\in \mathcal{F}_5$. This proves (1).
\\
\\
(2)  {\em Not both $X_1 \subseteq A \cup C$ and $X_2 \subseteq C\cup D$ .}
\\
\\
For assume these both hold. Choose a square $a_1\d b_1\d b_2\d a_2\d a_1$ such that
$f_1$ is adjacent to $a_1$, and choose $d \in D$ adjacent to $f_k$. If $a_2$ is adjacent to $f_1$
then $a_1\d b_1,a_2\d b_2,f_1\c f_k \d d$ form a long odd prism, a contradiction. If $a_2$ is not adjacent
to $f_1$ then $a_1$ can be linked onto the triangle $\{b_1,b_2,d\}$, via $a_1\d b_1$,$a_1\d a_2\d b_2$,
$a_1\d f_1\d \cdots\d f_k\d d$, a contradiction. This proves (2).
\\
\\
(3) {\em Not both  $X_1 \subseteq A \cup B$ and $X_2 \subseteq C\cup D$ .}
\\
\\
For assume these both hold. Then $X_1 \cap X_2 = \emptyset$, and so $f_1$ is the unique neighbour
in $F$ of the vertices in $X_1$, and $f_k$ is the unique neighbour of those in $X_2$.
From (1), $X_2\not \subseteq B\cup D$ and so $X_2 \cap C \not = \emptyset$; and similarly from (2),
$X_1 \cap B  \not = \emptyset$. Also we are given that $X_1 \cap A, X_2 \cap D  \not = \emptyset$.
Since $a_1\d f_1\c f_k\d c_1\d a_1$ is a hole (where $a_1\in A\cap X_1$ and $c_1 \in C\cap X_2$)
it follows that $k$ is even.
Since $f_1$ is minor, $X_1 \cap A$ is complete to $X_1 \cap B$, and so $A,B$ are not subsets of $X_1$;
and similarly $C\cap X_2$ is complete to $D \cap X_2$ and therefore
$C,D$ are not subsets of $X_2$. So all the eight sets $A\cap X_1, A \setminus X_1$ etc.
are nonempty. Choose a square $a_1\d b_1\d b_2\d a_2\d a_1$ such that $f_1$ is adjacent to $a_1$ and
not to $a_2$;
and choose an antisquare $c_1\d d_1\d d_2\d c_2\d c_1$ such that $f_k$ is adjacent to $d_1$ and not to $d_2$.
It follows that $f_1$ is nonadjacent to $b_2$, since $X_1 \cap A$ is complete to $X_1 \cap B$, and
$f_k$ is not adjacent to $c_1$ since $X_2 \cap C$ is complete to $X_2 \cap D$. But then
$a_1\d f_1\c f_k \d d_1 \d b_2\d d_2 \d c_1 \d a_1$ is an odd hole, a contradiction. This proves (3).
\\
\\
(4) {\em Not both  $X_1 \subseteq A \cup C$ and $X_2 \subseteq B\cup D$ .}
\\
\\
For assume both these hold. Then again, the only edges between $V(K)$ and $F$ are between $X_1$ and $f_1$
and between $X_2$ and $f_k$. By (1) and (2), again all four of the sets $A \cap X_1, B \cap X_2,
C \cap X_1, D \cap X_2$ are nonempty. There are two cases, depending on the parity of $k$. First assume $k$
is odd. Then $A \cap X_1$ is anticomplete to $B \cap X_2$ (for if $ab$ were an edge there, then
$a\d f_1\d \cdots\d f_k\d b\d a$ would be an odd hole), and so $A \setminus X_1, B\setminus X_2$ are nonempty;
and similarly $C \cap X_1$ is anticomplete to $D \cap X_2$, and therefore
$C \setminus X_1, D\setminus X_2$ are nonempty. Choose
a square $a_1\d b_1\d b_2\d a_2\d a_1$ such that $f_1$ is adjacent to $a_1$ and not to $a_2$,
and choose an antisquare $c_1\d d_1\d d_2\d c_2\d c_1$ such that $f_k$ is adjacent to $d_1$ and not to $d_2$.
Since $a_1\d f_1 \c f_k\d d_1\d b_2\d a_2\d a_1$ is not an odd hole it follows that $b_2 \in X_2$,
and therefore $a_1 \notin X_1$ (since  $A \cap X_1$ is anticomplete to $B \cap X_2$), and
$c_2 \notin X_1$ similarly. But then
the three paths $a_2\d b_2,c_2\d d_1,a_1\d f_1 \c f_k$ form a long prism, contrary to $G\in \mathcal{F}_5$.
Now assume $k$ is even. Then $A \cap X_1$ is anticomplete to $B\setminus X_2$ (for if $a \in A \cap X_1$
is adjacent to $b \in B\setminus X_2$ then $a\d f_1\d \cdots\d f_k\d d\d b\d a$ is an odd hole, where
$d \in X_2 \cap D$) . Similarly $A \setminus X_1$ is anticomplete to $B \cap X_2$,
$C \cap X_1$ is anticomplete to $D \setminus X_2$, and $C \setminus X_1$ is anticomplete to $D \cap X_2$.
Choose $a \in A \cap X_1$ and a neighbour $b$ of $a$ in $B$; then $b \in X_2$. Similarly choose
$c \in C \cap X_1$ and $d \in D \cap X_2$, adjacent. Then the three paths $a\d b,c\d d,f_1\d \cdots\d f_k$
form a prism, and so $k = 2$. If $f_1$ is $C$-complete then since $C \cap X_1 = C$
is anticomplete to $D \setminus X_2$, it follows that $f_2$ is $D$-complete; and then we can add $f_1$
to $A$ and $f_2$ to $B$, contrary to the maximality of the cube. So $C \not \subseteq X_1$. Choose
an antisquare  $c_1\d d_1\d d_2\d c_2\d c_1$ such that $f_1$ is adjacent to $c_1$ and not to $c_2$. It follows
that $f_2$ is adjacent to $d_2$ and not to $d_1$. If $f_1$ is $A$-complete, then as before $f_2$ is
$B$-complete, and we can add $f_1$ to $C$ and $f_2$ to $D$ (because $f_1\d d_1\d f_2\d c_2\d f_1$ is a new
antisquare), a contradiction. So $f_1$ has a nonneighbour in $A$, and we can
choose  a square $a_1\d b_1\d b_2\d a_2\d a_1$ such that $f_1$ is adjacent to $a_1$ and not to $a_2$.
It follows that $f_2$ is adjacent to $b_1$ and not to $b_2$. But then $a_1\d f_1\d f_2\d d_2\d b_2\d d_1\d c_2\d a_1$
is an odd hole, a contradiction. This proves (4).

\bigskip

From (1)-(4), the first assertion of the theorem follows. Now let us prove the second assertion.
We may assume $X$ meets both $A \cup C$ and $B \cup D$, and so from what we just proved,
either $X \subseteq C\cup D$ or
$X \subseteq A\cup B$. Suppose first that $X \subseteq C\cup D$. If possible, choose
$c \in C \cap X$ and $d \in D \cap X$,
nonadjacent, and choose a path $P$ joining them with interior in $F$. Let $a_1\d b_1 \d b_2 \d a_1 \d a_1$
be a square; then the three paths $a_1\d b_1,a_2 \d b_2, c\d P\d d$ form a long prism, a contradiction.
So there are no such $c,d$, and the theorem holds.

Now assume that $X \subseteq A\cup B$. Assume $X \cap A$ is not complete to $X \cap B$, and
choose a path $a\d f_1\d \cdots\d f_k\d b$, where $a \in A,b\in B$ are nonadjacent and
$f_1,\ldots,f_k \in F$, with $k$ minimum. Since $f_1$ is minor, its neighbours in $A$ are complete
to its neighbours in $B$, and so $k \ge 2$.
Let $A'$ be the set of all vertices $a \in A$ such that
$a$ is adjacent to $f_1$ and there is a nonneighbour $b$ of $a$ in $B$ adjacent to $f_k$. By assumption
$A' \not = \emptyset$. Define $B'$ similarly in $B$. If $A' = A$ and $B' = B$, then $f_1$ is $A$-complete,
and so there are no edges between $\{f_1,\ldots,f_{k-1}\}$ and $B$, from the minimality of $k$; and similarly
$f_k$ is $B$-complete and there are no edges between $\{f_2,\ldots,f_k\}$ and $A$. Choose a square
$a_1\d b_1\d b_2\d a_2\d a_1$; then $a_1\d b_1,a_2\d b_2, f_1 \c f_k$ form a prism, so $k = 2$, and we can add
$f_1$ to $C$ and $f_2$ to $D$, contrary to the maximality  of the cube. So we may assume that
$A' \not = A$. Choose a square $a_1\d b_1\d b_2\d a_2\d a_1$ so that $a_1 \in A'$ and $a_2 \notin A'$.
Choose $c \in C$ and $d \in D$, adjacent. Choose $b \in B'$ nonadjacent to $a_1$ (this exists from the
definition of $A'$). From the minimality of $k$, $a_1\d f_1\d \cdots\d f_k\d b$ is a path. From the hole
$a_1\d f_1\d \cdots\d f_k\d b\d d\d c\d a_1$ we deduce that $k$ is even. Since $b$ is not adjacent to $a_1$, $b$ is
different from $b_1$.
Suppose that $f_k$ is adjacent to $b_2$. Then the set of attachments of $\{f_1,\ldots,f_k\}$ with respect
to the prism formed by $a_1\d b_1,a_2\d b_2,c\d d$ is not local, and yet it has no attachment in $c\d d$, so
by \ref{legaljump}, both $a_2$ and $b_1$ are attachments. Since $a_2,b_1$ are
nonadjacent, it follows from
the minimality of $k$ and \ref{prismjumps} that $a_2$ is adjacent to $f_1$ and $b_1$ to $f_k$, contradicting
that $a_2 \notin A'$.

So $f_k$ is not adjacent to $b_2$. Then $b$ is different from $b_2$.
Since $c$ has no neighbour in the connected set $F' = \{f_1,\ldots,f_k,b\}$, and the set of attachments
of $F'$ is not local with respect to the prism formed by $a_1\d b_1,a_2\d b_2,c\d d$, it follows from
\ref{legaljump} that $F'$ has an attachment in $a_2\d b_2$. If $a_2$ is not an attachment
then $b_2$ is,
and from the minimality of $k$ it follows that $b$ is the unique neighbour of $b_2$ in $F'$; but then
$a_1\d f_1\d \cdots\d f_k\d b$, $a_2\d b_2$,$c\d d$ form a long prism, a contradiction. So $a_2$ is an attachment
of $F'$. Since $a_2\d a_1\d f_1\c f_k\d b\d a_2$ is not an odd hole, $a_2$ has a neighbour in
$\{f_1\l f_k\}$.
If $b_1$ also has a neighbour in $\{f_1\l f_k\}$, then (since $a_2,b_1$ are nonadjacent) from the
minimality of $k$ and \ref{prismjumps} it follows that $a_2$ is adjacent to $f_1$ and $b_1$ to $f_k$,
and hence $a_2 \in A'$, a contradiction. So $b_1$ has no neighbour in $\{f_1\l f_k\}$.
Since $a_1\d f_1\d \cdots\d f_k\d b\d b_1\d a_1$ is not an odd hole it follows that $b_1$ is not adjacent to $b$,
and therefore has no neighbours in $F'$.
Let $P$ be the path between $a_2$ and $b$ with interior in $F'$. From \ref{legaljump},
$a_1$ has a neighbour in $P\setminus a_2$.  But the
only neighbour of $a_1$ in $F'$ is $f_1$, so $f_1$ is in $P\setminus a_2$, and therefore $f_1$ is adjacent
to $a_2$, and there are no other edges between $a_2$ and $F'$. Since $a_2 \notin A'$ it follows that
$a_2$ is adjacent to $b$. But then the set of neighbours of $b$ in the prism formed by
$a_1\d b_1, a_2 \d b_2, c\d d$ is not local, and yet none are in the path $a_1\d b_1$, contrary to
\ref{legaljump}.  This proves \ref{cubeFnbrs}.\bbox

The main result of this section is \ref{summary}.6, which we restate, the following:

\begin{thm}\label{cube}
Let $G\in \mathcal{F}_5$. If $G$ contains a double diamond as an induced subgraph, then either one of
$G,\overline{G}$ admits a 2-join, or $G$ admits a balanced skew partition. In particular, every
recalcitrant graph belongs to $\mathcal{F}_6$.
\end{thm}
\Proof We may assume that $G,\overline{G}$ do not admit 2-joins, and $G$ does not admit a balanced
skew partition. Suppose for a contradiction that $G$ contains a double diamond; then it contains
a cube, and so there is a maximal cube $(A,B,C,D)$ in $G$, forming $K$. Let $F$ be the set of all
minor vertices in $V(G) \setminus V(K)$, and $Y$ the set of all major ones.
\\
\\
(1) {\em Every anticomponent $Y_1$ of $Y$ is complete to one of $A \cup B, C\cup D, A \cup D, B \cup C$,
and every edge from $A\cup D$ to $B \cup C$ has a $Y_1$-complete end.}
\\
\\
This is immediate from \ref{cubeFnbrs} by taking complements.
\\
\\
(2) {\em There is no anticomponent of $Y$ that is complete to $ A \cup D$ or $B \cup C$.}
\\
\\
For suppose such a component exists, say $Y_1$. From the symmetry we may assume it is complete to $ A \cup D$.
Define $L$ to be the union of $C$ and all components of $F$ with an attachment in $C$, and $M$ to
be the union of $B$ and all other components of $F$; and define $X$ to be the set of all
$Y_1$-complete vertices of $G$ not in $L\cup M$. So all major vertices belong to $Y_1 \cup X$, and
the four sets $L,M,X\cup A \cup D,Y_1$ are nonempty and partition $V(G)$; and since $Y_1$ is complete
to $X\cup A \cup D$, and there are no edges between $L,M$ by \ref{cubeFnbrs}, it follows that
$(L \cup M, X\cup A \cup D \cup Y_1)$ is a skew partition of $G$. By \ref{geteven} it is not loose.
We claim it is balanced. For by \ref{balancev}, $(L,D)$ is balanced, since any vertex in $B$ is $D$-complete
and $L$-anticomplete. Let $u,v \in L$ be adjacent, and suppose they are joined by an odd antipath
$Q_1$ with interior in $Y_1$. If they both have nonneighbours in $D$, then since $D$ is anticonnected
they are also joined by an antipath $Q_2$ with interior in $D$, which is also odd since its union with
$Q_1$ is an antihole, contradicting that $(L,D)$ is balanced. So we may assume that $u$ is $D$-complete.
Hence $u \notin C$, and so $u$ belongs to some component $F_1$ of $F$ with an attachment in $C$. Since
$u$ is minor, all its neighbours in $C$ are adjacent to all its neighbours in $D$, and hence it has
no neighbours in $C$; so $v \in F_1$. Since $F_1$ has an attachment in $C$ and in $D$ (because
$u$ has neighbours in $D$) it follows that $F$ has no attachments in $A$, and so $u,v$ have no neighbours
in $A$. But then $a\d u\d Q_1\d v\d a$ is an odd antihole (where $a \in A$), a contradiction.
Next suppose there exist nonadjacent $u,v \in Y_1$, joined by an odd path $P$ with interior in
$L$. By what we just proved about odd antipaths, it follows that $P$ has length $\ge 5$.
Now $A \cup D$ is anticonnected, and there is no $A \cup D$-complete vertex in $L$, since every vertex
in $L$ is minor or belongs to $C$. Hence the ends of $P$ are $A\cup D$-complete and its internal vertices
are not. But this contradicts \ref{RRR}. By \ref{onepair}, $G$ admits a balanced skew partition, a
contradiction. This proves (2).
\\
\\
(3) {\em There is no component of $F$ such that its set of attachments in $K$ is a subset of one
of $A \cup C, B \cup D$.}
\\
\\
This follows from (2) by taking complements.
\\
\\
(4) {\em There do not exist both a component $F_1$ of $F$ with set of attachments contained in $A \cup B$
and an anticomponent $Y_1$ of $Y$ complete to $A \cup B$.}
\\
\\
For assume that such $F_1,Y_1$ exist. Define $M =  C \cup D \cup (F \setminus F_1)$, and $X$ to
be the set of all $Y_1$-complete vertices in $V(G) \setminus (M \cup F_1)$. So  $A \cup B \subseteq X$,
and the four sets $F_1, M, Y_1,X$ are all nonempty and form a partition of $V(G)$. Since
$Y_1$ is complete to $X$ and there are no edges between $F_1$ and $M$, it follows that
$(F_1 \cup M, Y_1 \cup X)$ is a skew partition of $G$.
Suppose that $u,v \in F_1$ are adjacent and are joined by an odd antipath $Q$ with interior in $Y_1$.
Since $a\d u\d Q\d v\d a$ is not an odd antihole (where $a \in A$), it follows that one of $u,v$ has a neighbour
in $A$, say $u$. Since by (3) there is an attachment of $F_1$ in $B$, and no vertex in $B$ is $A$-complete,
it follows from \ref{cubeFnbrs} that not every vertex in $A$ is an attachment of $F_1$, and so
$u$ is not $A$-complete.  Choose a square $a_1\d b_1\d b_2\d a_2\d a_1$ such that
$u$ is adjacent to $a_1$ and not to $a_2$. Since $a_2\d u\d Q\d v\d a_2$ is not
an antihole, and $u$ is not adjacent to $a_2$, it follows that $va_2$ is an edge.
Since $u$ is minor, it is not adjacent to $b_2$, by\ref{cubevnbrs}; and
since $b_2\d u\d Q\d v\d b_2$ is not an antihole, $v$ is
adjacent to $b_2$. Similarly $b_1$ is adjacent to $u$ and not $v$. But then $G|\{a_1,a_2,b_1,b_2,u,v,c,d\}$
(where $c \in C$ and $d \in D$ are adjacent) is $L(K_{3,3}\setminus e)$, a contradiction. So
there is no such edge $uv$. By taking complements it follows that there do not exist nonadjacent
vertices in $Y_1$ joined by an odd path with interior in $F_1$. Hence by \ref{onepair}, $G$ admits
a balanced skew partition, a contradiction.  This proves (4).
\\
\\
(5) {\em There do not exist both a component $F_1$ of $F$ with set of attachments contained in $C \cup D$
and an anticomponent $Y_1$ of $Y$ complete to $C \cup D$.}
\\
\\
For assume such $F_1,Y_1$ exist. Define $M =  A \cup B \cup (F \setminus F_1)$, and $X$ to
be the set of all $Y_1$-complete vertices in $V(G) \setminus (M \cup F_1)$. So  $C \cup D \subseteq X$,
and the four sets $F_1, M, Y_1,X$ are all nonempty and form a partition of $V(G)$. Since
$Y_1$ is complete to $X$ and there are no edges between $F_1$ and $M$, it follows that
$(F_1 \cup M, Y_1 \cup X)$ is a skew partition of $G$. Suppose that $u,v \in F_1$ are adjacent
and joined by an odd antipath $Q$ with interior in $Y_1$. Choose $c \in C$ and $d \in D$, nonadjacent.
Since $c\d u\d Q\d v\d c$ is not an odd antihole, $c$ is adjacent to one of $u,v$, and similarly so is $d$.
So $u,v$ are both attachments of $F_1$, contrary to \ref{cubeFnbrs}. Hence there are no such $u,v$. It follows
by taking complements that there are no two nonadjacent $u,v \in Y_1$ joined by an odd path with interior
in $F_1$; and so by \ref{onepair}, $G$ admits a balanced skew partition, a contradiction. This proves (5).

\bigskip

Now if $Y = \emptyset$, then by (3) it follows that $G$ admits a 2-join, a contradiction. So $Y$ is nonempty,
and by taking complements,  $F$ is nonempty.
By (4), passing to the complement if necessary, we may assume that there is no anticomponent
of $Y$ that is complete to $A \cup B$. Hence $Y$ is complete to $C \cup D$, by (1) and (2). Since $Y$ is nonempty, it
follows from (5) that there is no component $F_1$ of $F$ with set of attachments contained in $C \cup D$;
so by (3), all attachments of $F$ belong to $A \cup B$. Choose an anticomponent
$Y_1$ of $Y$. By (3) and \ref{cubeFnbrs}, $Y_1$ is not $A$-complete or $B$-complete.
Let $X$ be the set of $Y_1$-complete vertices in $A \cup B \cup C \cup D$. Let $L$ be the union of
$A \setminus X$ and all components of $F$ that have an attachment in $A \setminus X$; and let $M$ be the
union of $B \setminus X$ and all other components of $F$. By (1) there are no edges between
$A \setminus X$ and $B \setminus X$; and therefore by \ref{cubeFnbrs}, no component of $F$ has
attachments in both $A \setminus X$ and $B \setminus X$. Hence
there is no edge between $L$ and $M$. Since $L,M,X\cup (Y \setminus Y_1), Y_1$ is a partition
of $V(G)$, and $Y_1$ is complete to $X\cup (Y \setminus Y_1)$, it follows that
$(L \cup M, X \cup (Y \setminus Y_1) \cup Y_1)$ is a skew partition of $G$. No vertex of
$D$ has a neighbour in $L$, and so it is loose, contrary to \ref{geteven}.
Hence there is no such graph $G$. This proves \ref{cube}.

\section{Consequences}

Disposal of the long prism and double diamond has a number of consequences that we develop in
this section. First, we
can now prove a form of Chv\'{a}tal's skew partition conjecture \cite{starcut}, that no minimum
imperfect graph $G$ admits a skew partition, because it is a consequence of \ref{cube} and the following.

\begin{thm}\label{oddskew}
Let $G \in \mathcal{F}_6$. If $G$ admits a skew partition, then it admits a balanced skew partition.
\end{thm}
\Proof Suppose $G$ admits a skew partition. Neither of $G,\overline{G}$ contains as
an induced subgraph either a long prism, or a double diamond, or $L(K_{3,3}\setminus e)$, so the
result follows from \ref{findprism}.  This proves \ref{oddskew}. \bbox

Consequently we have the following:
\begin{thm}\label{bipattach}
Let $G \in \mathcal{F}_6$, and assume that $G$ admits no balanced skew partition.
Let $X,Y \subseteq V(G)$ be nonempty, disjoint, and complete to each other.
\begin{itemize}
\item If $X \cup Y = V(G)$, then either $G$ is complete, or $\overline{G}$ has exactly two components,
both with $\le 2$ vertices (and hence $|V(G)| \le 4$).
\item If $X \cup Y \not = V(G)$, then $V(G) \setminus (X \cup Y)$ is connected, and if in addition
$|X| > 1$, then every vertex in $X$ has a neighbour in $V(G) \setminus (X \cup Y)$.
\end{itemize}
\end{thm}
\Proof By \ref{oddskew}, $G$ admits no skew partition.
Assume first that $X \cup Y = V(G)$. Then $\overline{G}$ is not connected; let its components be
$B_1\l B_k$ say. We may assume that $G$ is not complete, and therefore we may assume that some
$B_i$, say $B_1$, has cardinality $>1$. Choose $x,y \in B_1$, adjacent. Then
$(V(G)\setminus\{x,y\},\{x,y\})$ is not a skew partition, and so
$G \setminus\{x,y\}$ is connected. Hence $k = 2$ and $B_1 = \{x,y\}$. Similarly $B_2$ has cardinality
$\le 2$, and so $|V(G) \le 4$ and the theorem holds.
Now assume that $G \setminus (X \cup Y)$ is nonnull.
Suppose that $V(G) \setminus (X \cup Y)$  is not connected; then $(V(G) \setminus (X \cup Y), X \cup Y)$
is a skew partition, a contradiction. So $V(G) \setminus (X \cup Y)$ is connected.
Now suppose some $x \in X$ has no neighbour in $V(G) \setminus (X \cup Y)$. Hence
$V(G) \setminus ((X \setminus \{x\}) \cup Y)$ is not connected, and since $G$ admits no skew partition it
follows that $X = \{x\}$. This proves \ref{bipattach}.\bbox

Here is another consequence:

\begin{thm}\label{nobanister}
Let $G  \in \mathcal{F}_6$. Let $C$ be a cycle in $G$ of length $\ge 6$, with vertices $p_1 \l p_n$ in order,
and let $1 < h < i$ and $i+1 < j < n$. Let $C$ be induced except possibly for an edge $p_hp_j$.
Let $Y \subseteq V(G) \setminus V(C)$ be anticonnected, such that the only $Y$-complete vertices in
$C$ are $p_n,p_1,p_i,p_{i+1}$. Suppose there is a path
$F$ of $G \setminus Y$ from $p_h$ to $p_j$ (possibly of length 1), such that there are no edges between
its interior and $V(C) \setminus \{p_h,p_j\}$. Then some vertex of $F$ is $Y$-complete.
\end{thm}
\Proof Assume no vertex of $F$ is $Y$-complete.
Since the hole \[p_1 \c p_h \d F \d p_j \c p_n \d p_1\] is even, and the path
$p_1 \c p_h \c p_i$ is even (by \ref{greentouch}), it follows that the path
\[p_i \d p_{i-1} \c p_h \d F \d p_j \c p_n\] is odd, and therefore has length 3 by
\ref{RRR}. So $F$ has length 1, and $i = h+1$ and $n = j+1$. Similarly $h = 2$ and $j = i+2$, and
so $n = 6$.  Then $p_2,p_5$ are adjacent, so there is an antipath $Q$
joining them with interior in $Y$. But then in $\overline{G}$, the three paths
$p_1\d p_4$,$p_5\d p_2$, $p_3\d Q \d p_6$ form a long prism, a contradiction.
This proves \ref{nobanister}. \bbox

There is a variant of \ref{evenantipath2}, the following.

\begin{thm}\label{evenantipath3}
Let $G  \in \mathcal{F}_6$, and let $p_1 \c p_m$ be a path in $G$.
Let $2 \le s\le m-2$, and let $p_s\d q_1 \c q_n\d p_{s+1}$ be an antipath,
where $n \ge 2$. Assume that $p_1,p_m$ are both adjacent to all of
$q_1,\ldots,q_n$. Then $n$ is even and $m = 4$.
\end{thm}
\Proof If $n$ is even then $p_s\d q_1\c q_n\d p_{s+1}$ is an odd antipath, and
$p_1,p_m$ are complete to its interior; and hence $p_1,p_m$ are both adjacent to one of
$p_s,p_{s+1}$. So $s = 2$ and $m = s+2$, and therefore $m = 4$. Now assume $n$ is odd; then
$p_s\d q_1 \c q_n\d p_{s+1}$ is an even antipath of length $\ge 4$, contrary to \ref{doubleRR2}
applied in $\overline{G}$ to this antipath and the sets $\{p_1\l p_{s-1}\}$,$\{p_{s+2}\l p_n\}$.
This proves  \ref{evenantipath3}.\bbox

There is a strengthening of \ref{evengap}:
\begin{thm}\label{evengap2}
Let $G  \in \mathcal{F}_6$, let $C$ be a hole in $G$, and let
$Q\subseteq V(G) \setminus V(C)$ be anticonnected. Let $P$ be a path in $C$ of length $>1$
so that its ends are $X$-complete and its internal vertices are not. Then $P$ has even length.
\end{thm}
\Proof  The claim is trivial if $C$ has length 4, so we assume it has length $\ge 6$.
Let the vertices of $C$ be $p_1,\ldots,p_n$ in order, and let $P$ be $p_1 \d \cdots \d p_k$
say, where $3 \le k < n$. Assume $k$ is even. Then by \ref{RRR} applied to $P$ we deduce that
$P$ has length 3, so $k= 4$. By \ref{greentouch} every $X$-complete vertex is adjacent to one of
$p_2,p_3$, so there are none in the interior of the odd path $p_4 \d p_5 \d \cdots \d p_n \d p_1$.
By \ref{RRR} this path also has length 3, so $n = 6$. Let $Q$ be the shortest antipath with
interior in $X$, joining either $p_2,p_3$ or $p_5,p_6$. From the symmetry we may assume its vertices
are $p_2\d q_1 \c q_m \d p_3$ say. Then $Q$ is odd since it can be completed to an antihole via
$p_3 \d p_1 \d p_4 \d p_2$; and since $p_5 \d p_2 \d Q \d p_3 \d p_5$ is therefore not an antihole,
it follows that $p_5$ (and similarly $p_6$) has a nonneighbour in the interior of $Q$. From the
choice of $Q$ it follows that $p_5,p_6$ both have exactly one nonneighbour in the interior of $Q$;
one is nonadjacent to $q_1$ and the other to $q_m$. Suppose that $m > 2$. If $p_5$ is nonadjacent
to $q_1$ then the three antipaths $q_1\c q_m, p_5\d p_3, p_2\d p_6$ for a long prism in $\overline{G}$,
contrary to $G  \in \mathcal{F}_6$; while if $p_5$ is nonadjacent to $q_m$ then
$q_1\c q_m, p_6\d p_3, p_2\d p_5$ form a long prism, again a contradiction. So $m= 2$. But then
$G|\{p_1\c p_6,q_1,q_2\}$ is $L(K_{3,3}\setminus e)$ if $p_5$ is nonadjacent to $q_1$, and
a double diamond if $p_5$ is nonadjacent to $q_2$, again contrary to $G  \in \mathcal{F}_6$. This
proves \ref{evengap2}. \bbox

There is also a strengthening of \ref{hole&antipath}; we no longer need the vertex $z$.

\begin{thm}\label{hole&antipath2}
Let $G  \in \mathcal{F}_6$, let $C$ be a hole in $G$ of length $\ge 6$, with vertices $p_1,\ldots,p_m$ in
order, and let $Q$ be an antipath with vertices $p_1,q_1,\ldots,q_n,p_2$, with length $\ge 4$
and even.
There is at most one vertex in $\{p_3,\ldots,p_m\}$ complete to either $\{q_1,\ldots,q_{n-1}\}$
or $\{q_2,\ldots,q_n\}$, and any such vertex is one of $p_3,p_m$.
\end{thm}
\Proof Suppose first that one of $q_1,\ldots,q_n$ belongs to the hole. Since it is adjacent to at least
one of $p_1,p_2$ (since $Q$ is an antipath), we may assume that it is $p_m$; and since it is nonadjacent
to $p_2$, it follows that $p_m = q_n$. So $p_3 \not = q_1$ (since $q_1$ is adjacent to $q_n$), and therefore
no more of $q_1,\ldots,q_n$ belong to $C$. Suppose that there exists $i$ with $3 \le i < m$
such that $p_i$ is complete
to either  $\{q_1,\ldots,q_{n-1}\}$ or $\{q_2,\ldots,q_n\}$. If $i < m-1$ then $p_i$ is not adjacent
to $p_m = q_n$, so $p_i$ is complete to $\{q_1,\ldots,q_{n-1}\}$; but then $p_i \d p_1 \d q_1 \d \cdots \d
q_n \d p_i$ is an odd antihole. So $i = m-1$. By \ref{evengap2} applied to the path
$p_{m-1} \d p_m \d p_1 \d p_2$ it follows that $p_{m-1}$ is not complete to $\{q_1,\ldots,q_{n-1}\}$, and
therefore it is complete to $\{q_2,\ldots,q_n\}$ and nonadjacent to $q_1$. But then
$p_2 \d p_{m-1} \d q_1 \d \cdots \d q_n \d p_2$ is an odd antihole, a contradiction. So there is
no such $i$, and therefore the theorem holds in this case.

So we may assume that none of $q_1,\ldots,q_n$ belong to $C$.
Let $X = \{q_1,\ldots,q_n\}$, and let $Y_1,Y_2$ be
the sets of vertices in $\{p_3,\ldots,p_m\}$ complete to $X\setminus q_n$, $X\setminus q_1$
respectively.
\\
\\
(1) $Y_1 \subseteq Y_2 \cup \{p_m\}$, {\em and} $Y_2 \subseteq Y_1 \cup \{p_3\}$.
\\
\\
This is proved as in the proof of \ref{hole&antipath}.
\\
\\
(2) {\em If $Y_1 \not \subseteq \{p_m\}$ then $p_3 \in Y_1 \cap Y_2$, and if
$Y_2 \not \subseteq \{p_3\}$ then $p_m \in Y_1 \cap Y_2$.}
\\
\\
For assume $Y_1 \not \subseteq \{p_m\}$, and choose $i$ with $3 \le i \le m-1$ minimum so
that $p_i \in Y_1$. By (1), $p_i \in Y_2$, so we may assume $i>3$, for otherwise the claim
holds. By \ref{evengap2} applied to the anticonnected set $X \setminus q_n$,
$i$ is even.  The path $p_1\d \cdots\d p_i$ is odd, and between $X\setminus q_1$-complete
vertices, so by \ref{evengap2} it contains another in its interior, say $p_h$. From the
minimality of $i$, $p_h \notin Y_1$, so by (1) $ h = 3$, and \ref{evengap2} applied to the
path $p_3 \c p_i$ implies that $i = 4$. Choose $j$ with $4 \le j \le m$ maximum so that
$p_j \in Y_2$. By (1), $p_j$ is $X$-complete. By \ref{evenantipath3} applied to
$p_j \c p_m \d p_1 \c p_4$  we deduce that $j \le 5$, and so $j \not = m$.
By \ref{evengap2} applied to the path $p_j \c p_m \d p_1$ and anticonnected set $X \setminus q_1$,
it follows that $j$ is odd, and so $j = 5$. From \ref{evengap2} applied to the path
 $p_5 \c p_m \d p_1 \d p_2$ and anticonnected set $X \setminus q_n$, we deduce that there exists
$k$ with $6 \le k \le m$ such that $p_k \in Y_1$. Since it is not in $Y_2$, it follows from
(1) that $k = m$, and so $p_m \in Y_1 \setminus Y_2$. But then $p_3 \d q_1 \c q_n \d p_m \d p_3$
is an odd antihole, a contradiction. This proves (2).

\bigskip

Now not both $p_3,p_m$ are in $Y_1 \cap Y_2$, for otherwise $Q$ could be completed to an
odd antihole via $p_2\d p_m\d p_3\d p_1$. Hence we may assume $p_3  \notin Y_1 \cap Y_2$, and
so from (2), $Y_1 \subseteq \{p_m\}$. By (1), $Y_2 \subseteq \{p_3\} \cup Y_1$, and so
$Y_1 \cup Y_2 \subseteq \{p_3,p_m\}$. We may therefore assume that $Y_1 \cup Y_2 = \{p_3,p_m\}$,
for otherwise the theorem holds. In particular, $p_3 \in Y_2$. If also $p_m \in Y_2$, then
$p_3\d p_4\d \cdots\d p_m$ is an odd path between $X\setminus q_1$-complete vertices, and none of
its internal vertices are $X\setminus q_1$-complete, contrary to \ref{evengap2}. So
$p_m \notin Y_2$, and so $p_m \in Y_1$; but then $p_3\d q_1\d q_2\d \cdots\d q_n\d p_m\d p_3$ is an
odd antihole, a contradiction. This proves \ref{hole&antipath2}. \bbox

This implies a strengthening of \ref{bighole&antihole}:
\begin{thm}\label{hole&antihole}
Let $G  \in \mathcal{F}_6$.  Let $C$ be a hole of length $>4$ and $D$ an antihole of length $>4$.
Then $|V(C) \cap V(D)| \le 2$.
\end{thm}
\Proof Assume that $|V(C) \cap V(D)| \ge 3$; then by taking complements if necessary, we may assume that
there are three vertices in $V(C) \cap V(D)$ such that exactly one pair of them is adjacent. Hence we can number
the vertices of $C$ as $p_1,\ldots,p_m$ in order, and the vertices of $D$ as
$p_1,q_1,\ldots,q_n,p_2,p_k$ for some $k$ with $4 \le k \le m-1$. (Possibly the hole and antihole also
share some fourth vertex.) Hence the antipath $p_1\d q_1\c q_n\d p_2$ has length $\ge 4$ and even. The
vertex $p_k$ is complete to $\{q_1,\ldots,q_n\}$, and different from $p_3,p_m$, contrary to
\ref{hole&antipath2}. This proves \ref{hole&antihole}.\bbox

\section{Odd wheels}

A {\em wheel} in a graph $G$ is a pair $(C,Y)$, satisfying:
\begin{itemize}
\item $C$ is a hole of length $\ge 6$
\item $Y$ is a non-null anticonnected set disjoint from $C$
\item there are two disjoint edges of $C$, both $Y$-complete.
\end{itemize}
We need to study how the remainder of a nonconforming graph
can attach onto a wheel. Conforti and Cornu\'{e}jols also made such a study - see \cite{C&C}
and \cite{CCVZ} for several theorems related to the results of this section.
We call $C$ the {\em rim} and $Y$ the {\em hub} of the wheel.
A maximal path in a path or hole $H$
whose vertices are all $Y$-complete is called a {\em segment} or {\em $Y$-segment} of $H$.
A wheel $(C,Y)$ is {\em odd} if some segment has odd length.
In this section we prove that there are no odd wheels in a recalcitrant graph.
(G\'{e}rard Cornu\'{e}jols informs us that he and his co-workers proved the same result, independently,
but, like us, assuming the truth of \ref{longprism} - see \cite{CCVZ}.)

Let us say that distinct vertices $u,v$ of the rim of a wheel $(C,Y)$ have the {\em same wheel-parity}
if there is a path of the rim joining them containing an even number of $Y$-complete edges
(and hence by \ref{evengap}, so does the second path, if $u,v$ are nonadjacent); and
{\em opposite wheel-parity} otherwise.

\begin{thm}\label{wheelvnbrs}
Let $G  \in \mathcal{F}_6$, and let $(C,Y)$ be a wheel in $G$. Let $v \in V(G) \setminus (V(C)\cup Y)$, such
that $v$ is not $Y$-complete.  Suppose that there exist neighbours of $v$ in $C$ with
opposite wheel-parity. Then in every path of $C$ between them there is a $Y \cup \{v\}$-complete edge.
Moreover, either:
\begin{itemize}
\item $v$ has only two neighbours in $C$, and they are adjacent and both $Y$-complete, or
\item there is a $3$-vertex path $p_1\d p_2 \d p_3$ in $C$, so that $p_1,p_2,p_3$ are all
$Y\cup \{v\}$-complete, and every other neighbour of $v$ in $C$ has the same wheel-parity as $p_1$, or
\item $(C, Y \cup \{v\})$ is a wheel.
\end{itemize}
\end{thm}
\Proof
\\
\\
(1) {\em Let $P$ be a path in $C$ of length $\ge 1$, such that its ends are adjacent to $v$ and
have opposite wheel-parity. Then either some internal vertex of $P$ is a neighbour of $v$, or
$P$ has length 1.}
\\
\\
Let $C$ have vertices $p_1,\ldots,p_n$ in order, and let $P$ be the path $p_1\c p_j$ say, where $j < n$.
We assume no internal vertex of $P$ is a neighbour of $v$, and that $j \ge 3$. From the hole
$v \d p_1 \c p_j \d v$ it follows that $j$ is odd.
Since $p_1,p_j$ have opposite wheel-parity with respect to $(C,Y)$, there are an odd number
of $Y$-complete edges in $P$.  Choose $Y' \subseteq Y$
minimal such that $Y'$ is anticonnected and there are an odd number of $Y'$-complete edges in $P$.
From \ref{evengap} applied to the hole $v \d p_1 \c p_j \d v$, it contains
just one $Y'$-complete edge and only two $Y'$-complete vertices. Hence there exists $i$ with $1 \le i < j$
so that $p_i,p_{i+1}$ are the only $Y'$-complete vertices in $P$. Since $j$ is odd, it follows
that exactly one of $i-1, j-i$ is even; so (by replacing $P$ by its reverse if necessary) we may
assume that $i$ is odd. So $p_j$ is different from $p_{i+1}$, and hence $p_j$ is not $Y'$-complete.
There are two disjoint $Y'$-complete edges in $C$, so one of them does not use $p_i$; and therefore
it does not use $p_1$ either (for $p_1$ is not $Y'$-complete unless $i = 1$). Hence both its ends are
in $\{p_{j+1}\l p_n\}$. Consequently $n \ge j+2$, and since $n$ is even and $j$ is odd it follows that
$n \ge j+3$. Therefore there is a $Y'$-complete vertex in $\{p_{j+2} \l p_{n-1}\}$.

Suppose that $v$ has a neighbour in $\{p_{j+2}\l p_{n-1}\}$. Then there is a path $Q$ from $v$ to a
$Y'$-complete vertex $u$ say, with $V(Q) \subseteq \{v, p_{j+2}\l p_{n-1}\}$, such that no internal vertex
of $Q$ is $Y'$-complete. The path $p_i\c p_1 \d v \d Q \d u$ has both ends $Y'$-complete, and no internal
vertex $Y'$-complete, and the $Y'$-complete vertex $p_{i+1}$ has no neighbour in its interior; so this
path is even, that is, $Q$ is odd. Hence the path $p_{i+1}\c p_j \d v \d Q \d u$ is odd, and so by
\ref{RRR} has length 3; and hence $j = i+2$ and $Q$ has length 1. Also, every $Y'$-green vertex is adjacent
to one of $p_j,v$, by \ref{greentouch}; and so $p_i$ is adjacent to $v$, and so $i = 1, j = 3$; and
$v$ is adjacent to every $Y'$-complete vertex in $C$ except $p_2$ and possibly $p_4$ (for no others are
adjacent to $p_3$). In particular, there are two nonadjacent $Y' \cup \{v\}$-complete vertices in $C$,
and so by \ref{evengap} there are an even number of $Y' \cup \{v\}$-complete edges in $C$. But
all $Y'$-complete edges of $C$ are $Y' \cup \{v\}$-complete except $p_1p_2$ and possibly
$p_4p_5$; and since there are also an even number of $Y'$-complete edges in $C$, it follows that
$p_4, p_5$ are $Y'$-complete, and $v$ is adjacent to $p_5$ and not to $p_4$.
But then the vertices $v,p_1,p_2,p_3,p_4,p_5$ violate
\ref{nobanister}.

This proves that $v$ has no neighbour in  $\{p_{j+2}\l p_{n-1}\}$. Choose $k$ with $j \le k \le n$ minimum
such that $p_k$ is $Y'$-complete. Since there is a $Y'$-complete vertex in $\{p_{j+2} \l p_{n-1}\}$,
it follows that $k < n$. From \ref{evengap} it follows that the path $p_{i+1}\c p_k$ is even, and
so $k$ is even. Suppose that $v$ is not adjacent to $p_{j+1}$. Since $v \d p_j \c p_n \d v$ is not
an odd hole, it follows that $v$ is not adjacent to $p_n$, so $p_1,p_j$ are its only neighbours
in $C$. But $p_i \c p_1 \d v \d p_j \c p_k$ is odd, and therefore has length 3 by \ref{RRR}; and by
\ref{greentouch}, every $Y'$-complete vertex in $C$ is adjacent to $v$ except possibly $p_{j-1},p_{j+1}$,
a contradiction since there is a $Y'$-complete vertex in $\{p_{j+2}\l p_{n-1}\}$. So $v$ is
adjacent to $p_{j+1}$. Since $v \d p_{j+1} \c p_n \d p_1\d v$ is not an odd hole, it follows that
$v$ is also adjacent to $p_n$, so it has exactly four neighbours
in $C$. Choose $m$ with $k \le m \le n$ maximum so that $p_m$ is $Y'$-complete.
It follows that $m \ge j+2$. If $m = n$ then $k$ has no neighbours in the interior of the odd path
$p_{i+1} \c p_j \d v \d p_n$, and the ends of this path are $Y'$-complete and its internal vertices are
not, contrary to \ref{greentouch}. So $m < n$.
Then \ref{evengap} applied to the path $p_m \c p_n \d p_1 \c p_i$ implies
that $m$ is odd, and therefore $m > k$. Suppose that $m > k+1$. Then $p_m \c p_n \d v \d p_{j+1} \c p_k$
is an odd path, and $p_{i+1}$ has no neighbour in its interior, contrary to \ref{greentouch}. So
$m = k+1$, and there is symmetry between the paths $p_1\c p_j$ and $p_{j+1} \c p_n$. Both these
paths have length $\ge 2$; suppose they both have length 2. Then $n = 6$, and the only $Y' \cup \{v\}$-complete
vertices in $C$ are $p_1,p_4$, contrary to \ref{evengap2}. So one of the paths has length $>2$, and from
the symmetry we may assume that $j \ge 4$. Hence the hole $H = v \d p_1 \c p_j \d v$ has length $\ge 6$,
and the only $Y'$-complete vertices in it are $p_i,p_{i+1}$. By \ref{RRC}, $Y'$ contains a hat or a leap.
But $p_{k+1}$ has no neighbour in this hole, so the pair $(V(H),Y')$ is balanced by \ref{balancev}, and
hence there is no leap. So there is a hat; that is, there exists $y \in Y'$ with no neighbours in
$H$ except $p_i,p_{i+1}$. From the minimality of $Y'$ it follows that $Y' = \{y\}$. But then
$G |(V(C) \cup \{v,y\})$ is the line graph of a bipartite subdivision of $K_4$, a contradiction.
This proves (1).

\bigskip

From (1) the first assertion of the theorem follows. Now we prove the second assertion.
Suppose that $v$ has at least four neighbours in $C$, two with the same wheel-parity, and
two others with the opposite wheel-parity. Then there are two disjoint paths as in (1), and therefore from
(1) there are two disjoint $Y\cup \{v\}$-complete edges in $C$, and so $(C, Y\cup \{v\})$ is a wheel
and the theorem holds. So we may assume that $C$ has vertices $p_1,\ldots,p_n$ in order, and
$v$ is adjacent to $p_1$, and $v$ has no other neighbour in $C$ with the same wheel-parity as $p_1$.
Since $v$ has at least one other neighbour, we may assume it has a neighbour in
$V(C) \setminus \{p_1,p_n\}$.
Choose $i >1$ minimum so that $v$ is adjacent to $p_i$; then $i < n$, so by (1), $ i = 2$. So $p_2$ is
$Y \cup \{v\}$-complete. If $v$ has a third neighbour in $C$ then similarly $p_n$ is $Y \cup \{v\}$-complete
and the theorem holds; and if not then again the theorem holds.  This proves \ref{wheelvnbrs}.\bbox

\begin{thm}\label{wheelFnbrs}
Let $G  \in \mathcal{F}_6$, and let $(C,Y)$ be a wheel in $G$. Let $F \subseteq V(G) \setminus (V(C)\cup Y)$
be connected, such that no vertex in $F$ is $Y$-complete. Let $X \subseteq V(C)$
be the set of attachments of $F$ in $C$. Suppose that there exist vertices in $X$ with
opposite wheel-parity, and there are two vertices in $X$ that are nonadjacent. Then either:
\begin{itemize}
\item there is a vertex $v \in F$ so that $(C, Y \cup \{v\})$ is a wheel, or
\item there is a vertex $v \in F$ with at least four neighbours in $C$, and a $3$-vertex path
$p_1\d p_2 \d p_3$ in $C$, so that $p_1,p_2,p_3$ are all
$Y\cup \{v\}$-complete, and every other neighbour of $v$ in $C$ has the same wheel-parity as $p_1$,
or
\item there is a $3$-vertex path $p_1\d p_2 \d p_3$ in $C$, all $Y$-complete, and a path
$p_1 \d f_1\c f_k \d p_3$  with interior in $F$, such that there no edges between $\{f_1\l f_k\}$
and $\{p_4\l p_n \}$.
\end{itemize}
\end{thm}
\Proof We may assume that $F$ is minimal. If $|F| = 1$ then the result follows from \ref{wheelvnbrs},
so we assume $|F| \ge 2$.
\\
\\
(1) {\em If there do not exist nonadjacent vertices in $X$ with different wheel-parity, then the theorem
holds.}
\\
\\
For there exist vertices in $X$ with different wheel-parity, which are therefore adjacent; say
$p_1,p_2$, where $C$ has vertices $p_1\l p_n$ in order. So $p_1,p_2$ are both $Y$-complete, since they
have different wheel-parity. There is a third attachment of $F$, since there are
two that are nonadjacent, say $p_i$ where $3 \le i \le n$. Since $p_1,p_2$ have different wheel-parity,
we may assume that $p_2,p_i$ have different wheel-parity; and therefore $p_2,p_i$ are adjacent,
that is, $i = 3$, and $p_3$ is $Y$-complete.
Suppose $F$ has a fourth attachment $p_j$ say, where $4 \le j \le n$. From the
symmetry we may assume $j \not = n$; and so $p_j$ is nonadjacent to both $p_1,p_2$, and one of
these has different wheel-parity from $p_j$, a contradiction. So $p_1,p_2,p_3$ are the only
attachments of $F$, and then the theorem holds. This proves (1).

\bigskip
From (1) we may assume there are nonadjacent vertices in $X$ with opposite wheel-parity, say $x_1,x_2$,
and therefore $F$ is the interior of a path between $x_1,x_2$, from the minimality of $F$. Let
$C$ have vertices $p_1\l p_n$ in order; then we may assume that there exists $m$ with
$3 \le m \le n-1$ such that $p_1,p_m$ have opposite wheel-parity,
and there is a path $p_1\d f_1\c f_k \d p_m$ where $F = \{f_1\l f_k\}$.
Let $X_1$ be the set of attachments in $C$ of $F \setminus f_k$, and $X_2$ the set of attachments
of $F \setminus f_1$. From the minimality of $F$, for $i = 1,2$ either all members of $X_i$ have the
same wheel-parity, or there are at most two members of $X_i$, adjacent if there are two. Since $k \ge 2$
it follows that $X_1 \cup X_2 = X$.
\\
\\
(2) {\em $X_1$ and $X_k$ do not both have members of opposite wheel-parity.}
\\
\\
For suppose they do; then $X_1,X_2$ both consist of exactly two adjacent vertices of opposite wheel-parity,
say $X_1 = \{p_1,p_2\}$ and $X_2 = \{p_{m'},p_{m'+1}\}$. So $p_1,p_2,p_{m'},p_{m'+1}$ are all $Y$-complete, and
all distinct since two of them are nonadjacent and of
opposite wheel-parity. So the only edges between $F$ and $\{p_1,p_2\}$ are incident with $f_1$, and
similarly for $f_k$. But then $G$ contains a long prism since $n \ge 6$, a contradiction. This proves
(2).
\\
\\
(3) {\em If $X_1$ has members of opposite wheel-parity then the theorem holds.}
\\
\\
For assume $X_1$ has members of opposite wheel-parity.  Then we may assume its only members are
$p_1,p_2$, and they are both
$Y$-complete. From (1) we may assume that all members of $X_2$ have the same wheel-parity
as $p_2$. In particular, $p_1$ has no neighbour in $F \setminus f_1$.
So the only edges between $F$ and $C$ are $f_1p_1$, edges incident with $p_2$, and edges incident
with $f_k$.  Suppose that $p_2$ also has no neighbour in $F\setminus f_1$, and therefore
$p_2$ is adjacent to $f_1$. If $f_k$ has
a unique neighbour $x$ in $C$, then $x$ can be linked onto the triangle $\{p_1,p_2,f_1\}$; if
$f_k$ has two nonadjacent neighbours in $C$ then $f_k$ can be linked onto the same triangle; and
if it has exactly two neighbours and they are adjacent, then $G$ contains a long prism, in
each case a contradiction. So $p_2$ has a neighbour in $F\setminus f_1$. Let $R_1$ be the path
$p_1 \d f_1 \c f_k$, and let $R_2$ be the path from $p_2$ to $f_k$ with interior in $F\setminus f_1$.
Then $p_1$ has no neighbours in $R_2 \setminus p_2$.
Let $Q_1$ be the path from $f_k$ to $p_n$ with interior in $C \setminus p_1$. Now
$p_1 \d R_1 \d f_k \d Q_1 \d p_n \d p_1$ is a hole, so $R_1$ and $Q_1$ have lengths of opposite parity; and
since this hole contains an odd number of $Y$-complete edges (since all neighbours of $f_k$ have wheel-parity
opposite from that of $p_1$) it follows that it contains exactly one such edge and only two
$Y$-complete vertices. Since $p_1$ is $Y$-complete, the other is therefore $p_n$.
The path $p_2 \d R_2 \d f_k \d Q_1 \d p_n$ is between $Y$-complete vertices, and no internal
vertex is $Y$-complete, and the $Y$-complete vertex $p_1$ has no neighbour in its interior; so it
is even by \ref{greentouch}, that is, $R_1,R_2$ have opposite parity.  Now
there is a $Y$-complete vertex in $\{p_4 \l p_{n-1}\}$; for there are two disjoint $Y$-complete
edges in $C$, and an even number of $Y$-complete edges in $C$. Let $p_s$ be such a vertex, where
$4 \le s \le n-1$. We claim that
$f_k$ has a neighbour in $\{p_4 \l p_{n-1}\}$. For if not, then since $X \not = \{p_n,p_1,p_2\}$
(because there are nonadjacent vertices in $X$ of opposite wheel-parity), it follows that $f_k$ is adjacent
to $p_3$. Since $p_s $ is not in $Q_1$, it follows that $p_3$ is not in $Q_1$, and so $f_k$ has
another neighbour, which must be $p_n$; but then $f_k \d p_3 \d p_4 \c p_n \d f_k$ is an odd hole.
So $f_k$ has a neighbour in $\{p_4 \l p_{n-1}\}$; and therefore there is a path $Q_2$ from
$f_k$ to some $x$, such that $x$ is the unique $Y$-complete vertex in $Q_2$, and $V(Q_2 \setminus f_k)
\subseteq  \{p_4 \l p_{n-1}\}$. Now the path $p_2 \d R_2 \d f_k \d Q_2$ has both ends $Y$-complete,
and no internal vertex $Y$-complete, and the $Y$-complete vertex $p_1$ has no neighbour in its
interior, so it is even by \ref{greentouch}. Therefore the path $p_1 \d R_1 \d f_k \d Q_2$ is odd,
since $R_1,R_2$ have opposite parity; and again its ends are $Y$-complete and its internal
vertices are not. So it has length 3, by \ref{RRR}, and so $k = 2$; and every $Y$-complete
vertex is adjacent to one of $f_1,f_2$. Consequently there is no $Y$-complete vertex in $C$
different from $p_1$ with the same wheel-parity as $p_1$, a contradiction. This proves (3).

\bigskip

From (3) we may assume that all members of $X_1$ have the same wheel-parity, and
all members of $X_2$ have the opposite wheel-parity. It follows that $X_1 \cap X_2 = \emptyset$, and
so there are no edges between the interior of $F$ and $C$. So $X_1$ is the set of neighbours of $f_1$
in $C$, and $X_2$ is the set of neighbours of $f_k$ in $C$.
\\
\\
(4) {\em At least one of $f_1,f_k$ has only one neighbour in $C$.}
\\
\\
For assume they both have at least two. Then there are disjoint paths $Q,R$ of $C$, and
both containing neighbours of both $f_1,f_k$. Choose $Q,R$ minimal, and let $Q$ have ends
$q_1,q_2$; then from the minimality of $Q$, $q_1$ is the unique neighbour of one of $f_1,f_k$ in $Q$,
and $q_2$ is the unique neighbour of the other. Let $f_1q_1$ and $f_kq_2$ be edges say. Similarly
let $R$ have ends $r_1,r_2$, where $f_1r_1,f_kr_2$ are edges. Since $q_1,q_2$ have opposite
wheel-parity, it follows that there are an odd number of $Y$-complete edges in the hole
$f_1\c f_k \d q_2 \d Q \d q_1 \d f_1$; so by \ref{evengap} there is exactly one, and just two
$Y$-complete vertices. If there are no edges between $Q$ and $R$ this contradicts \ref{nobanister}. Since
$Q,R$ are disjoint subpaths of $C$, all the edges between them join their ends; so we may assume that
$q_1$ is adjacent to one of $r_1,r_2$.  From the hole $f_1 \c f_r \d q_2 \d Q \d q_1 \d f_1$ it follows
that $Q$ has parity $k-1$, and similarly so does $R$. Suppose first that $q_1$ is adjacent to $r_1$. Since
$q_1 \d Q \d q_2 \d f_r \d r_2 \d R \d r_1 \d q_1$ is not an odd hole, it follows that $q_2$ is
adjacent to $r_2$, and hence $G$ contains a long prism, since $C$ has length $\ge 6$, a contradiction.
So $q_1$ is adjacent to $r_2$. Since $q_1$ is a neighbour of $f_1$ and $r_2$ of $f_k$, it follows
that $q_1,r_2$ have opposite wheel-parity, and since they are adjacent, they are both $Y$-complete.
Let $q'$ be the neighbour of $q_1$ in $Q$, let $Q' = Q \setminus q_1$, let $r'$ be the neighbour
of $r_2$ in $R$, and let $R' = R \setminus r_2$. Since
in the hole $f_1\c f_k \d q_2 \d Q \d q_1 \d f_1$ there are only two $Y$-complete vertices and they are
adjacent, it follows that the second is $q'$, and similarly $r'$ is $Y$-complete. If $q_2$ is adjacent
to $r_1$ then not both $q_2,r_1$ are $Y$-complete since $C$ has length $\ge 6$; and so
there are exactly three $Y$-complete edges in $C$, contrary to \ref{evengap}. It follows that
$q_2$ is not adjacent to $r_1$. From the hole $q_1 \d Q \d q_2 \d f_r \d r_2 \d q_1$ it follows that
$Q$ has odd length, and therefore so does $R$ and $k$ is even. But then the path
$q' \d Q'\d q_2 \d f_r \c f_1 \d r_1 \d R' \d r'$ has odd length, its ends are $Y$-complete and its
internal vertices are not, and so by \ref{RRR} it has length 3; that is, $Q$,$R$ have length 1 and $k = 2$.
Hence the path $r_1 \d f_1 \d f_2 \d q_2$ is odd, its ends are $Y$-complete, and its internal vertices
are not, so every $Y$-complete vertex is adjacent to one of $f_1,f_2$. Let $ab,a'b'$ be two
$Y$-complete edges of $C$, disjoint and so that there are no edges from $\{a,b\}$ to $\{a',b'\}$.
Then each of $a,b,a',b'$ is adjacent
to one of $f_1,f_2$, and since all neighbours of $f_1$ in $C$ have opposite wheel-parity from all
neighbours of $f_2$ in $C$, we may assume that $a,a'$ are adjacent to $f_1$ and $b,b'$ to $f_2$.
But this contradicts \ref{nobanister}. This proves (4).

\bigskip

From (4) we may assume that $X_1$ has only one member, say $p_1$. Choose $i,j$ with
$2 \le i , j \le n$, such that $p_i,p_j$ are adjacent to $f_k$, with $i$ minimum and $j$ maximum.
From the hole  $p_1\d f_1 \c f_k \d p_i \d p_{i-1} \c p_1$ ($= H_1$ say) we deduce that $i,k$ have the
same parity, and from the hole $p_1\d f_1 \c f_k \d p_i \d p_{i+1} \c p_n \d p_1$ (= $H_2$ say) that
$j,k$ have the same parity. (So either $p_i = p_j$ or $p_i,p_j$ are nonadjacent.)
Since $p_1,p_i$ have different wheel-parity, and so do $p_1,p_j$,
there is an odd number of $Y$-complete edges in each of $H_1,H_2$; and therefore there is exactly
one $Y$-complete edge and exactly two $Y$-complete vertices in each of the holes, by
\ref{evengap}. Suppose that $i = j$. Then there are only two $Y$-complete edges in $C$, and therefore
they are disjoint, and $p_1,p_i$ are not $Y$-complete (since $H_1,H_2$ both have only two $Y$-complete
vertices), contrary to \ref{nobanister}. So $j > i$, and hence $j \ge i+2$.
If $p_1$ is not $Y$-complete, then the $Y$-complete edge in $H_1$ is disjoint from
the path $p_1\d f_1 \c f_k$, and so is the one in $H_2$; but this contradicts \ref{nobanister}.
So $p_1$ is $Y$-complete. Since $H_1$ contains only two $Y$-complete vertices and they are adjacent,
the other is $p_2$, and similarly $p_n$ is $Y$-complete.
\\
\\
(5) {\em $f_k$ has no neighbour in $\{p_3\l p_{j-2}\}$.}
\\
\\
For assume it does. We claim there is also a $Y$-complete vertex in this set; for otherwise the only
$Y$-complete vertices in $C$ are $p_n,p_1,p_2$ and possibly $p_{j-1}$, which is impossible since there
are two disjoint $Y$-complete edges and an even number of $Y$-complete edges in $C$. Hence there is a
path $P$ say from $f_k$ to some $x$ so that $x$ is the unique $Y$-complete vertex in $P$ and
$V(P\setminus f_k) \subseteq \{p_3\l p_{j-2}\}$. The path
$p_n \d p_{n-1}\c p_j \d f_k \d P \d x$ is even, since
its ends are $Y$-complete, no internal vertex is $Y$-complete, and the $Y$-complete vertex $p_1$
has no neighbour in its interior. The path $p_1 \d f_1 \c f_k \d P \d x$ is therefore odd
(since $k,j$ have opposite parity), and also its ends are $Y$-complete and no internal vertex is
$Y$-complete; so it has length 3 by \ref{RRR}, and hence $k = 2$ and every $Y$-complete vertex is adjacent
to one of $f_1,f_2$, by \ref{greentouch}. So there is no $Y$-complete vertex in $C\setminus p_1$ with
the same wheel-parity as $p_1$, a contradiction. This proves (5).

\bigskip

Since $f_k$ is adjacent to $p_i$, and $i < j$ and $j-i$  is even, it follows from (5) that $i = 2$,
and similarly $f_k$ has no neighbours in $\{p_{i+2}\l p_{n-1}\}$ and $j = n$. So $f_k$ has no
neighbours in $\{p_3\l p_{j-2}\} \cup \{p_{i+2}\l p_{n-1}\} = \{p_3 \l p_{n-1}\}$, and therefore
$p_2,p_n$ are its only neighbours, contradicting that there are nonadjacent vertices in $X$ of
opposite wheel-parity. This proves \ref{wheelFnbrs}.\bbox

The main result of this section is \ref{summary}.7, which we restate, the following.

\begin{thm}\label{oddwheel}
Let $G  \in \mathcal{F}_6$. If there is an odd wheel in $G$ then $G$ admits a balanced skew partition.
In particular, every recalcitrant graph belongs to $\mathcal{F}_7$.
\end{thm}
\Proof
Suppose $(C,Y)$ is an odd wheel with $Y$ maximal, and subject to that, such that
the number of $Y$-complete edges in $C$ is minimum. (We refer to these conditions as the ``optimality'' of
$(C,Y)$.)
\\
\\
(1) {\em There is no vertex $v \in V(G) \setminus (V(C) \cup Y)$ such that $v$ is not $Y$-complete
and has nonadjacent neighbours in $C$ of opposite wheel-parity.}
\\
\\
Suppose there is such a vertex $v$.  Suppose first that there is an odd $Y \cup \{v\}$-segment in $C$.
From the maximality of $Y$, $(C,Y \cup \{v\})$ is therefore not a wheel, and so there is a unique
$Y \cup \{v\}$-complete edge in $C$. By \ref{RRC}, either $v$ has only two neighbours in $C$, or some
vertex of $Y$ has only three, in either case a contradiction. So there is no odd
 $Y \cup \{v\}$-segment in $C$.
Define a ``line'' to be a maximal subpath of $C$ with no internal
vertex adjacent to $v$. It follows that every edge of $C$ is in a unique line.
Let $C$ have vertices $p_1\l p_n$ in order, and let $S$ be an odd $Y$-segment.

Since there are no odd $Y \cup \{v\}$-segments, it follows that an even
number of edges of $S$ are $Y \cup \{v\}$-complete. Hence an odd number are not, and therefore there is a
line $L$ containing an odd number of edges of $S$ that are not $Y \cup \{v\}$-complete. In particular
it contains at least one edge that is $Y$-complete and not $Y \cup \{v\}$-complete, so $L$ has
length $>1$. Let the ends of $L$ be $p,q$. By \ref{wheelvnbrs}, $p$ and $q$ have the same wheel-parity
with respect to $(C,Y)$, and so $L$ contains an odd number of edges of some other $Y$-segment $S' \not = S$.
In particular, there are two disjoint $Y$-complete edges in the hole $v \d p \d L \d q \d v$ ( = $H$ say);
so $H$ has length $\ge 6$ (because $v$ is not $Y$-complete) and so $(H,Y)$ is a wheel. Moreover
it is an odd wheel, for it contains an odd number of edges of $S$, and those edges form either one
or two $Y$-segments in $H$, and one of these segments is odd. Since there is a $Y \cup \{v\}$-complete
edge in $C$ (by \ref{wheelvnbrs}, since $v$ has neighbours in $C$ of opposite wheel-parity)
which therefore does not belong to $L$, this contradicts the optimality of $(C,Y)$. This
proves (1).

\bigskip

Since $(C,Y)$ is an odd wheel, $C$ has at least two segments, and therefore there are vertices
$u,v$ in $C$ with different wheel-parity and neither of them $Y$-complete. Let $X$ be the set
of all $Y$-complete vertices in $V(G)$. Then $|X| >1$, since there at least four in $C$; so by
\ref{bipattach}, we may assume that
$V(G) \setminus (X \cup Y)$ is nonempty and connected ( = $Z$ say), and every
vertex in $X$ has a neighbour in it, for otherwise $G$ admits a balanced skew partition and the theorem
holds. In particular $u,v \in Z$, so there is a minimal
connected subset $F$ of $Z$ such that there are two vertices of $C \setminus X$ (say $p,q$) of opposite
wheel-parity, both with neighbours in $F$. Since $p,q$ have opposite wheel-parity and are not
$Y$-complete, they are not adjacent. From the minimality of $F$, $F$ is a path, and no vertex
of $F$ is in $C$. By \ref{wheelFnbrs} and (1),
there is a 3-vertex path $p_1\d p_2 \d p_3$ in $C$, all $Y$-complete, and a path
$p_1 \d f_1\c f_k \d p_3$  with interior in $F$, such that there no edges between $\{f_1\l f_k\}$
and $\{p_4\l p_n \}$.  But then $C \setminus p_2$ can be completed
to a hole $C'$ say, via $p_1 \d f_1\c f_k \d p_3$i; and $C'$ has length $\ge 6$.
For every odd segment $S$ of $(C,Y)$, either
it contained both or neither of the edges $p_1p_2,p_2p_3$; and so in either case an odd number
of edges of $S$ belong to $C'$. Since $(C,Y)$ has an odd segment and there are an even number
of $Y$-complete edges in $C$, it has at least two odd segments. It follows that there are
two disjoint $Y$-complete edges in $C'$, and so $(C',Y)$ is a wheel. Since an odd number
of edges of the odd segment $S$ belong to $C'$, it follows that  $(C',Y)$ is an odd wheel, contrary
to the optimality of $(C,Y)$. This proves \ref{oddwheel}. \bbox

\section{Another extension of the Roussel-Rubio lemma}

Let $\{a_1,a_2,a_3\}$ be a triangle in $G$. A {\em reflection} of this triangle is another triangle
$\{b_1,b_2,b_3\}$ of $G$, disjoint from the first, so that the only edges between it and the first triangle
are $a_1b_1,a_2b_2,a_3b_3$. Hence these six vertices induce a prism. A subset $F$ of $V(G)$ is
said to {\em catch} the triangle $\{a_1,a_2,a_3\}$ if it is connected and disjoint from that triangle and
$a_1,a_2,a_3$ all have neighbours in $F$.
We begin with the following extremely useful little fact.

\begin{thm}\label{tricatch}
Let $A$ be a triangle in a graph $G \in \mathcal{F}_7$, and let $F \subseteq V(G)\setminus A $
catch $A$. Then either $F$ contains a reflection of $A$, or some vertex of $F$ has
$\ge 2$ neighbours in $A$.
\end{thm}
\Proof Suppose not, and choose $F$ minimal such that it catches $A$. Let
$A = \{a_1,a_2,a_3\}$ say, and for $i = 1,2,3$, let $B_i$ be
the set of neighbours of $a_i$ in $F$. Thus the three sets $B_1,B_2,B_3$ are pairwise disjoint
and nonempty.
\\
\\
(1) {\em There is no path in $F$ meeting all of $B_1,B_2,B_3$.}
\\
\\
For assume there is, and choose it minimal. So then we may assume there is a path $P$ from $b_1 \in B_1$
to $b_2 \in B_2$, such that some vertex of $P$ is in $B_3$, and for $i = 1,2$, $b_i$ is the only
vertex of $P$ in $B_i$. Since $B_3$ is disjoint from $B_1 \cup B_2$, every vertex of $B_3$ in $P$
is an internal vertex of $P$; and so $P$ has length $\ge 2$. But then $(C,\{b_3\})$ is an odd wheel,
where $C$ is the hole $a_1 \d b_1 \d P \d b_2 \d a_2 \d a_1$, contrary to $G \in \mathcal{F}_7$. This proves (1).

\bigskip

Choose $b_1\in F$ so that $F \setminus b_1$ is connected; then from the minimality of $F$,
$F\setminus x$ does not catch $A$, and so we may assume that $B_1 = \{b_1\}$. Since $F$ is connected and
$|F| \ge 2$, there is a second vertex $b_2 \not = b_1$ in $F$ so that $F \setminus b_2$ is connected,
and so similarly we may assume $B_2 = \{b_2\}$. Let $P$ be a path in $F$ between $b_1,b_2$.
By (1) no vertex of $P$ is in $B_3$, so $F$ contains a connected subset $F'$ including $V(P)$ which contains
exactly one vertex of $B_3$. From the minimality of $F$, $|B_3| = 1$; let $B_3 = \{b_3\}$ say.
Let $Q$ be a minimal path in $F$ such that $b_3 \in V(Q)$ and some vertex of $P$ has a neighbour
in $Q$. From the minimality of $Q$ it follows that $Q$ is vertex-disjoint from $P$, and $Q$ has
ends $b_3,x$ say, where $x$ is the unique vertex of $Q$ with neighbours in $P$. From the minimality
of $F$, $x$ either has one neighbour in $P$, or just two neighbours and they are adjacent; for if it
has two nonadjacent neighbours, any vertex of $P$ between them could be deleted from $F$,
contrary to the minimality of $F$. If $x$ has just one neighbours $y$ say in $P$, then $y$ can be
linked onto the triangle $A$, contrary to \ref{trianglev}; so it has two adjacent. Since $G$
does not contain a long prism it follows that $Q$ has length 0 and $P$ has length 1, and so $F$ contains
a reflection of $A$, as required. This proves \ref{tricatch}.\bbox

We did not use the full strength of $G \in \mathcal{F}_7$ in proving \ref{tricatch}; we just used that
there were no odd wheels with hubs of cardinality 1. This suggest that there should be some generalization
of \ref{tricatch} whose proof does use the full strength of the hypothesis that there are no odd wheels,
and that is true, but not easy - it will be a consequence of the main result of this section.

Before we start on that, let us give a strengthening of \ref{RRC} for graphs in $\mathcal{F}_7$.

\begin{thm}\label{RRstrip}
Let $G \in \mathcal{F}_7$, and let $F,Y \subseteq V(G)$ be disjoint, such that $F$ is connected and
$Y$ is anticonnected. Let $a_0,b_0 \in V(G) \setminus (F \cup Y)$ and $a,b \in F$ such that
$a \d a_0 \d b_0 \d b$ is a $3$-edge path in $G$. Suppose that:
\begin{itemize}
\item $a_0,b_0$ are both $Y$-complete, and $a,b$ are not $Y$-complete,
\item the only neighbours of $a_0,b_0$ in $F$ are $a$ and $b$ respectively,
\item $F \setminus a$ and $F \setminus b$ are both connected.
\end{itemize}
Then either:
\begin{enumerate}
\item there is a vertex in $Y$ with no neighbour in $F$, or
\item there are two nonadjacent vertices $y_1,y_2 \in Y$, such that $a$ is the only neighbour
of $y_1$ in $F$, and $b$ is the only neighbour of $y_2$ in $F$.
\end{enumerate}
\end{thm}
\Proof We may assume that every vertex in $Y$ has a neighbour in $F$, for otherwise statement 1 of the
theorem holds.
\\
\\
(1) {\em There exist nonadjacent $y_1,y_2$ in $Y$, such that $y_1$ is adjacent to $a$ and not $b$, and
$y_2$ is adjacent to $b$ and not $a$.}
\\
\\
For choose a path $P$ between $a$ and $b$. Then the hole $a_0 \d a \d P_0 \d b \d b_0 \d a_0$
($= C$, say) has length $\ge 6$. If there are any $Y$-complete vertices in $P$, then they belong to
the interior of $P$ since $a,b$ are not $Y$-complete, and there is an odd number of $Y$-complete
edges in $P$, by \ref{evengap}; but then $(C,Y)$ is an odd wheel (the path $a_0 \d b_0$ is an odd
$Y$-segment), a contradiction. So there are no $Y$-complete vertices in $P$. By \ref{RRC} applied to $C$,
$Y$ contains either
a hat or a leap. Suppose first it contains a hat, that is, there is a vertex $y \in Y$ with no
neighbour in $P$. By the assumption above, $y$ has a neighbour in $F$. Consequently $F$ catches the triangle
$\{a_0,b_0,y\}$. But $y$ is not adjacent to $a$ or $b$ since it has no neighbour in $P$, and $a$
is the only vertex in $F$ adjacent to $a_0$, and the same for $b,b_0$; and $a,b$ are
nonadjacent, so $F$ contains no reflection of the triangle.  This contradicts \ref{tricatch}.
Hence there is no such $y$, and so there is a leap. This proves (1).
\\
\\
(2) {\em There is no path in $F$ between $a$ and $b$ such that $y_1$ or $y_2$ has a neighbour in
its interior.}
\\
\\
For suppose there is such a path, $P'$ say. Then the set $\{y_1,y_2\}$ contains neither a leap not a hat
for the hole $a_0 \d a \d P' \d b \d b_0 \d a_0$ ( = $C$ say), and so by \ref{RRC} there is a vertex
in $P$ adjacent to both
$y_1,y_2$. By \ref{evengap} there is an even number of $\{y_1,y_2\}$-complete edges in this
hole, and since $a,b$ are not $\{y_1,y_2\}$-complete, $(C,\{y_1,y_2\})$ is an odd wheel,
a contradiction. This proves (2).

\bigskip

Now if neither of $y_1,y_2$ has any more neighbours in $F$ then statement 2 of the theorem holds;
so we assume at least one of them has another neighbour in $F$. Since $F \setminus a$, $F \setminus b$
are both connected, there is a connected subset $F'$ of $F \setminus \{a,b\}$, so that both $a$ and $b$
have neighbours in $F'$, and at least one of $y_1,y_2$ has a neighbour in $F'$. Choose $F'$ minimal
with these properties. At least one of $y_1,y_2$ has a neighbour (say $x$) in $F'$.
We claim that $F' \setminus x$ is connected. For if not, let
$L$ be a component of it, and $M$ the union of the other components. From the minimality of $F$, not
both $a,b$ have neighbours in $L \cup \{x\}$, and not both have neighbours in $M \cup \{x\}$; so we may
assume all neighbours of $a$ in $F'$ are in $L$, and all neighbours of $b$ are in $M$. But then there
is a path from $a$ to $b$ with interior in $F$ and with $x$ in its interior, contrary to (2). This
proves that $F' \setminus x$ is connected. There is a path from $a$ to $b$ with interior in $F'$, and
$x$ is not in it, by (2), and it has length $>1$ since $a,b$ are nonadjacent. So $a$,$b$ both have
neighbours in $F' \setminus x$. From the minimality of $F'$, $y_1$ and $y_2$ both have no neighbours
in $F' \setminus x$. We claim that $x$ is adjacent to both $y_1$ and $y_2$. For it is adjacent to
at least one, say $y_1$; let $Q$ be a path from $x$ to $b$ with interior in $F'$. Then $y_1 \d x \d Q \d b$
is a path, since $y_1$ has no more neighbours in $F'$. Since $b_0 \d y_1 \d x \d Q \d b \d b_0$
is a hole it follows that $Q$ is odd. Therefore $a_0 \d y_1 \d x \d Q \d b \d y_2 \d a_0$ is not a hole,
and so $y_2$ has neighbours in $Q$. Since it has no neighbours in $F' \setminus x$, this proves our claim
that $x$ is adjacent to both $y_1,y_2$.

With $Q$ as before, and therefore odd, it follows that
$y_2 \d x \d Q \d b \d y_2$ is not a hole, and therefore $Q$ has length 1, that is, $x$ is adjacent
to $b$. Similarly $x$ is adjacent to $a$; but then $x \d a \d a_0 \d b_0 \d b \d x$ is an odd hole,
a contradiction. This proves \ref{RRstrip}.\bbox

The following is a variant of \ref{RRstrip}, not so symmetrical, but more useful.
\begin{thm}\label{RRstrip2}
Let $G \in \mathcal{F}_7$, and let $F,Y \subseteq V(G)$ be disjoint, such that $F$ is connected and
$Y$ is anticonnected. Let $a_0,b_0 \in V(G) \setminus (F \cup Y)$ and $a,b \in F$ such that
$a \d a_0 \d b_0 \d b$ is a $3$-edge path in $G$. Suppose that:
\begin{itemize}
\item $a_0,b_0$ are both $Y$-complete, and $a,b$ are not $Y$-complete,
\item the only neighbours of $a_0,b_0$ in $F$ are $a$ and $b$ respectively,
\item $F \setminus a$ is connected.
\end{itemize}
Then  there is a vertex $y \in Y$ with no neighbour in $F \setminus a$.
\end{thm}
\Proof
If $F \setminus b$ is connected, the result follows from \ref{RRstrip}. So assume it is not, and
let $F_1'$ be the component of $F \setminus b$ that contains $a$, and $F_2'$ the union of all the other
components. For $i = 1,2$ let $F_i = F_i' \cup \{b\}$. Then $F_1 \setminus a$, $F_1 \setminus b$ are
both connected, so by \ref{RRstrip} either there exists $y \in Y$ with no neighbour in $F_1$, or there
exist nonadjacent $y_1,y_2 \in Y$ with no neighbours in $F_1$ except $a,b$ respectively.
Suppose the first. If $y$ has a neighbour in $F_2$ then $b$ can be linked onto the triangle
$\{y,a_0,b_0\}$, a contradiction to \ref{trianglev}; and if not then $y$ satsifies the theorem.
Now suppose the second. If $y_1$ has neighbours in $F_2$ then $(F\setminus\{a\}) \cup \{y_2\}$
catches the triangle
$\{a,a_0, y_1\}$; the only neighbours of $a,a_0,y_1$ belong to the disjoint sets $F_1',\{y_2\},F2'$;
and there is no reflection since there are no edges between $y_1$ and $F_1'$, contrary to \ref{tricatch}.
So $y_1$ has no neighbours in $F_2$.
This proves \ref{RRstrip2}.\bbox

The next result is just a technical lemma for use in proving the main result of this section, which is
\ref{tripleRR}.
\begin{thm}\label{tripleRR1}
Let $G \in \mathcal{F}_7$ and let $P$ be a path in $G$ with length $>1$, with vertices
$p_1,\ldots,p_n$ in order. Let $X,Y \subseteq V(G) \setminus V(P)$ be anticonnected sets, so that
$X \cup Y$ is anticonnected, $p_1$ is $X$-complete, and $p_n$ is the unique $Y$-complete vertex in $P$.
(Note that $X,Y$ need not be disjoint.) Let $z \in V(G) \setminus (X \cup Y \cup V(P))$, complete
to $X \cup Y$ and with no neighbours in $P$.  Assume that $p_n$ is not $X$-complete. Let
$p_n \d x_1 \c x_k \d y$ be an antipath with interior in $X$ from $p_n$ to some $y \in Y$. Then
$p_{n-1}$ is nonadjacent to $x_1$.
\end{thm}
\Proof Let $F = \{p_{n-1},x_1\l x_k\} \cup Y$. Since $p_{n-1}$ is not $Y$-complete it follows that $F$ is
anticonnected, and both $F \setminus p_{n-1}$, $F \setminus x_1$ are anticonnected.
The only nonneighbour of $z$ in $F$ is $p_{n-1}$, and the only nonneighbour of $p_n$ in $F$ is $x_1$;
and we may assume that $p_{n-1}$ is adjacent to $x_1$. Now $p_{n-1} \d z \d p_n \d x_1$
is a path in $\overline{G}$, and $F$ is connected in  $\overline{G}$, and $\{p_1 \l p_{n-2}\}$ is anticonnected
in $\overline{G}$. Also, $z$ and $p_n$ are $\{p_1 \l p_{n-2}\}$-complete in  $\overline{G}$, and
$p_{n-1},x_1$ are not. We may therefore apply \ref{RRstrip} in $\overline{G}$, and deduce that there is
a vertex in $\{p_1 \l p_{n-2}\}$ which is complete (in $G$) to $F \setminus p_{n-1}$. But then this vertex is
$Y$-complete, a contradiction. This proves \ref{tripleRR1}. \bbox

We gave in \ref{doubleRR} an extension of the Roussel-Rubio lemma to two anticonnected sets
instead of one (we haven't had much use of that theorem yet, but its time is coming.) In that extension the
two sets had to be complete to each other. Now we prove a similar result where
the two sets are not complete to each other. Incidentally, unlike \ref{doubleRR},
what we are going to prove here is
not true for general Berge graphs --- we need the hypothesis that $G \in \mathcal{F}_7$.

\begin{thm}\label{tripleRR}
Let $G \in \mathcal{F}_7$ and let $P$ be an odd path in $G$ with length $>1$, with vertices
$p_1,\ldots,p_n$ in order. Let $X,Y \subseteq V(G) \setminus V(P)$ be anticonnected sets, so that
$X \cup Y$ is anticonnected, $p_1$ is $X$-complete, and $p_n$ is the unique $Y$-complete vertex in $P$.
(Note that $X,Y$ need not be disjoint.) Let $z \in V(G) \setminus (X \cup Y \cup V(P))$, complete to
$X \cup Y$ and with no neighbours in $P$.  Then an odd number of edges of $P$ are $X$-complete.
\end{thm}
\Proof If possible choose a counterexample $P,X,Y$ such that
\begin{enumerate}
\item $P$ is minimal
\item subject to condition 1, $X \cup Y$ is minimal, and
\item subject to conditions 1 and 2, $|X| + |Y|$ is minimum.
\end{enumerate}
We refer to this property as the ``optimality'' of $P,X,Y$.
\\
\\
(1) {\em No vertex of $P\setminus p_1$ is $X$-complete.}
\\
\\
If $p_n$ is $X$-complete, then since $P$ has odd length $>1$, and the $X$-complete vertex $z$
has no neighbour in $P$, it follows from \ref{greentouch} and \ref{evengap} that there are an odd
number of $X$-complete edges in $P$, and the theorem holds, a contradiction. So $p_n$ is not
$X$-complete. By \ref{tripleRR1}, $p_{n-1}$ is not $X$-complete. Since $p_1$ is $X$-complete, we
can choose $i$ with $1 \le i \le n $ maximum such that $p_i$ is $X$-complete. So $i \le n-2$.
Since $z$ has no neighbour in the path $p_1 \c p_i$, if $i$ is even then there is an odd number
of $X$-complete edges in this path and hence in $P$, by \ref{greentouch} and \ref{evengap}.
So we may assume that $i$ is odd. Hence the theorem is also false for $X,Y$ and the
path $p_i \l p_n$. From the optimality of $P,X,Y$ it follows that $i = 1$. This proves (1).
\\
\\
(2) {\em Suppose that $x_1,x_2 \in X$ are distinct and such that $X \setminus x_i$ is anticonnected
for $i = 1,2$. Then $X \cap Y = \emptyset$, and one of $x_1,x_2$ is the unique vertex of
$X$ with a nonneighbour in $Y$.}
\\
\\
For if $(X \setminus x_i) \cup Y$ is not anticonnected for some $i$, then
$Y$ is disjoint from $X \setminus x_i$ (since both these sets are anticonnected), and $Y$ is complete
to $X \setminus x_i$; and therefore $x_i \notin Y$ (since $x_i$ has a nonneighbour in $X \setminus x_i$),
so $X \cap Y = \emptyset$. But then the statement of (2) holds. So we may assume that
$(X\setminus x_i) \cup Y$ is anticonnected for $i = 1,2$. From the optimality of $P,X,Y$ it follows
that the theorem holds for $X\setminus x_i,Y,P$; and so, since $p_n$ is the unique $Y$-complete vertex
in $P$, it follows that there are an odd number of $X\setminus x_i$-complete edges in $P$, for $i = 1,2$.
For $i = 1,2$ let $W_i$ be the set of $X \setminus x_i$-complete vertices in $P$. So
$W_1 \cap W_2 = \{p_1\}$.  Let $Q$
be an antipath in $X$ between $x_1$ and $x_2$. We claim that $Q$ is odd. For since
$W_1 \cap W_2 = \{p_1\}$, there are nonadjacent vertices $p_i,p_j$ of $P$, such that $p_i \in W_1
\setminus W_2$ and $p_j \in W_2 \setminus W_1$; and since $p_i \d x_1 \d Q \d x_2 \d p_j \d p_i$ is an
antihole it follows that $Q$ is odd. Let us say a {\em line} is a minimal subpath of $P\setminus p_1$ meeting
both $W_1$ and $W_2$. So every line has length $\ge 1$, and has one end in $W_1$ and the other in $W_2$,
and has no more vertices in either $W_1$ or $W_2$. If some line $L$ has odd length $>1$, then
$(L,X \setminus x_1, X \setminus x_2)$ is another counterexample to the theorem, contrary to the
optimality of $P,X,Y$; and if some line has length 1, say $p_i\d p_{i+1}$ where $p_i \in W_1$, then
$z \d p_i \d x_1 \d Q \d x_2 \d p_{i+1} \d z$ is an odd antihole, a contradiction. Hence every line
is even.  Choose $i$ minimum with $2 \le i \le n$ such that $\{p_2\l p_i\}$
includes a line. (This is possible since both $W_1,W_2$ meet $P\setminus p_1$.) Since all lines have
length $\ge 2$ it follows that $i \ge 4$. From the minimality of $i$,
$\{p_2\l p_{i-1}\}$ does not include a line, and so for some $k \in \{1,2\}$,
the path $p_1\c p_i$ has both ends $Y \setminus x_k$-complete and no internal vertex
$X \setminus x_k$-complete. But this path has length $\ge 2$, and $z$ has no neighbour in it, so by
\ref{greentouch}
it is even, that is, $i$ is odd. Choose $j$ with $j \ge 2$ maximum so that $\{p_j\l p_n\}$ includes a line.
Since every line has length $\ge 2$ it follows that $2 \le j \le n-2$. From the maximality of
$j$ it follows that for some $k \in \{1,2\}$, $W_k \cap \{p_j\l p_n\} = \{p_j\}$.
If the path $p_j \l p_n$ has odd length, then $p_j \l p_n, X \setminus x_k,Y$ is a counterexample to
the theorem, contrary to the optimality of $P,X,Y$. So $n-j$ is even, and hence $j$ is even.
Now $i$ is odd, so if $i \ge j$ then $p_j\c p_i$ is an odd line, a contradiction. Hence $i < j$,
and $j-i$ is odd.  Now the edges
$p_{i-1}p_i, p_jp_{j+1}$ are in lines. Consequently we may choose $r,s$ with $i \le r < s \le j$
such that $p_r,p_s \in W_1 \cup W_2$, and
the edges $p_{r-1}p_r,p_sp_{s+1}$ are in lines, and $s-r$ is odd; and therefore we may
choose such $r,s$ with $s-r$ minimum. If there is a line contained in the path $p_r \c p_s$,
say $p_h\c p_k$, then since $k-h$ is even, one of the paths $p_r \c p_h$ and $p_k \c p_s$
is odd, contrary to the minimality of $s-r$. So we may assume that $p_r\l p_s$ all belong to $W_1$.
Since $p_{r-1}p_r,p_sp_{s+1}$ are in lines, and $p_r,p_s \in W_1$, there exist $q,t$
with $2 \le q < r < s < t \le n$ such that $p_q\c p_r$ and $p_s \c p_t$ are lines.
All lines are even, so $r-q$ and $t-s$ are even, and therefore $t-q$ is odd.
Moreover $p_q,p_t \in W_2$, and none of $p_{q+1} \l p_{t-1}$ belongs to $W_2$, so
and the path $p_q\c p_t$ is odd, and $z$ has no neighbour in it, contrary to
\ref{greentouch}. This proves (2).
\\
\\
(3) {\em There is an antipath $x_1\c x_s \d y_1 \c y_t$ such that $s, t>1$ and $X = \{x_1\l x_s\}$,
and $Y = \{y_1\l y_t\}$.}
\\
\\
For if $|X| = 1$, $X = \{x\}$ say, then $z\d x\d p_1\c p_n$ is an odd path of length $\ge 5$
between $Y$-complete vertices, and none of its internal vertices are $Y$-complete, contrary to
\ref{RRR}. So $|X| \ge 2$, and similarly $|Y| \ge 2$. Hence there are at least two vertices $x \in X$
so that $X \setminus x$ is anticonnected, and from (2) , $X \cap Y = \emptyset$, and there is
a unique vertex $x\in X$ with nonneighbours in $Y$. By (2), there do not exist two vertices
$x' \in X\setminus x$ so that $X\setminus x'$ is anticonnected; and therefore $X$ is an antipath with one
end $x'$. The same applies for $Y$, and this proves (3).

\bigskip

Choose $t'$ with $1 \le t' \le t$, minimum so that $p_1$ is nonadjacent to $y_{t'}$. (This is possible
since $p_1$ is not $Y$-complete.) So $x_1\c x_s \d y_1 \c y_{t'} \d p_1$ is an antipath. Define
$W = (X \setminus x_1) \cup \{y_1\l y_{t'-1}\}$.
\\
\\
(4) {\em For every subpath $P'$ of $P$, if the ends of $P'$ are adjacent to $x_1$, then there are
an even number of $W$-complete edges in $P'$.}
\\
\\
For suppose not; then we may choose $P'$ so that no internal vertex of $P'$ is adjacent to $x_1$.
Let $P'$ be $p_h\c p_k$ say, where $1 \le h < k \le n$. Choose $i,j$ with $h \le i \le j \le k$ such
that $p_i,p_j$ are $W$-complete, with $i$ minimum and $j$ maximum. Since $p_k$ is not $X$-complete
it follows that $p_k$ is not $W$-complete (because it is adjacent to $x_1$), and so $j < k$.
Since there are an odd number of $W$-complete edges in $p_h\c p_k$, it follows that $k \ge h+2$,
and $x_1 \d p_h\c p_k \d x_1$ is a hole (so $k-h$ is even), containing an odd number of $W$-complete
edges. By \ref{evengap} it contains exactly one, and only two $W$-complete vertices; so $j = i+1$.
The path $z\d x_1 \d p_h \c p_i$ has both ends $W$-complete, and no internal vertex $W$-complete,
and the $W$-complete vertex $p_{j}$ has no neighbour in its interior (since $j < k$); so it is even,
by \ref{greentouch}, and hence $i-h$ is even. Since $k-h$ is even, it follows that
$p_j\c p_k \d x_1 \d z$ is an odd path; and again its ends are $W$-complete and its internal vertices
are not. By \ref{RRR} it has length 3, so $k = j+1$; and by \ref{greentouch}, every $W$-complete
vertex is adjacent to one of $p_k,x_1$. But no $W$-complete vertex in $P$ is adjacent to $x_1$ except
$p_1$, since no other vertex of $P$ is $X$-complete. So every $W$-complete vertex in $P\setminus p_1$
is adjacent to $p_k$, and so must be one of $p_{k-1},p_{k+1}$.
In particular, since $i < k-1$ it follows that $i = 1$, and so $j = 2, k = 3$, and the $W$-complete
vertices in $P$ are $p_1,p_2$ and possibly $p_4$.

By \ref{tripleRR1} (with $X$ and $Y$ exchanged), $p_2$ is nonadjacent
to $y_{t'}$. Choose $d$ with $1 \le d \le n$ minimum so that $y_{t'}$ is adjacent to $p_d$; then
$d \ge 3$. Then the path $p_1 \c p_d \d y_{t'} \d z$ has length $\ge 4$, and its ends are
$W\cup\{x_1\}$-complete,
and its internal vertices are not, so it is even by \ref{RRR}. Hence $d$ is odd, and the path
$p_1 \c p_d \d y_{t'}$ is odd. None of its internal vertices are $X$-complete, and the $X$-complete
vertex $z$ has no neighbour in its interior, and one end $p_1$ is $X$-complete, so the other end $y_{t'}$
is not; and hence $t' = 1$, since all other vertices of $Y$ are $X$-complete. So $W = X \setminus x_1$.
Let $V = X \setminus x_s$.  Now the path $p_1\c p_d \d y_1$ is between $V$-complete vertices, and is odd
and has length $>1$, and the $V$-complete vertex $z$ has no neighbour in its interior; so by
\ref{greentouch}, there is a $V$-complete edge in its interior.  Choose $c$ with $2 \le c \le d$
minimum such that $p_c$ is $V$-complete. Since $p_2$ is nonadjacent to $x_1$ it follows that $c \ge 3$.
Since $p_1\c p_c$ is between $V$-complete vertices and its internal vertices are not $V$-complete
and $z$ has no neighbour in it, it is even by \ref{greentouch}, and so $c$ is odd. We already saw
that $p_1,p_2$ and possibly $p_4$ are $W$-complete, and $c\ge 3$, so we may choose $b$ with
$2 \le b \le c$ maximum such that $p_b$ is $W$-complete. Hence $b = 2$ or 4.
The path $p_b \c p_c$ is odd, and $p_b$ is $W$-complete, and $p_c$ is
$V$-complete, and no other vertices of the path are either $W$- or $V$-complete. If
$c - b > 1$ then $p_b \c p_c, W,V$ is a counterexample to the theorem, contradicting the optimality
of $X,Y,P$ . So $c = b+1$. Then $z \d p_b \d x_1 \c x_s \d p_c \d z$ is an
antihole , so $s$ is odd. But then $p_2 \d x_1 \c x_s \d y_1 \d p_2$ is an odd antihole, a contradiction.
This proves (4).

\bigskip

Choose $h$ with $1 \le h \le n$ maximum such that $x_1$ is adjacent to $p_h$. Since
$x_1\d p_h \c p_n$ is between $Y$-complete vertices (since $s \ge 2$) and none of its internal
vertices are $Y$-complete, and the $Y$-complete vertex $z$ has no neighbour in its interior, this
path either has length 1 or even length by \ref{greentouch}. So either $h = n$ or $h$ is odd.
From the optimality of $P,X,Y$, it follows that $P, W,Y$ is not a counterexample to the theorem, and
so there are an odd number of $W$-complete edges in $P$.  Since $x_1$ is adjacent to $p_1$, from (4)
there are an even number of $W$-complete edges between $p_1$ and $p_h$, so there are an odd number in
the path $p_h \c p_n$, and in particular $h < n$, so $h$ is odd. Choose $i,j$ with
$h \le i \le j \le n$ such that $p_i,p_j$ are $W$-complete, with $i$ minimum and $j$ maximum.
Hence $j > i$.
Since $z \d x_1 \d p_h \c p_i$ is a path of length $\ge 2$ between $W$-complete vertices, and its
internal vertices are not $W$-complete, and the $W$-complete vertex $p_j$ has no neighbour in its interior,
it follows from \ref{greentouch} that $i-h$ is even.
\\
\\
(5) $h > 1$.
\\
\\
For assume $h = 1$; so $p_1$ is the only neighbour of $x_1$ in $P$. Let $S$ be the antipath
\[x_1 \c x_s \d y_1 \c y_t' \d p_1.\] Now $x_1 \d S \d p_1 \d z$ is an antipath, of length $\ge 4$;
all its internal vertices have neighbours in $P \setminus p_1$, and its ends do not. By
\ref{RRR} applied in $\overline{G}$, it follows that this antipath has even length and so $S$ has odd
length. Its ends have no neighbours in $P \setminus\{p_1,p_2\}$, and $z$ is complete to its interior
and also has no neighbours in $P \setminus\{p_1,p_2\}$; so by \ref{greentouch} applied in
$\overline{G}$, some internal vertex of $S$ has no neighbour in $P \setminus\{p_1,p_2\}$. But they
are all adjacent to $p_j$ or to $p_n$, so $j = 2$. By \ref{tripleRR1}, $p_2$ is nonadjacent to $y_{t'}$,
and also to $x_1$ since it is not $X$-complete. Therefore $p_2 \d x_1 \c x_s \d y_1 \c y_{t'} \d p_2$
is an antihole $D$ say. Choose $d$ with $1 \le d \le n$ minimum so that $y_{t'}$ is adjacent to
$p_d$; then $d \ge 3$, and so $x_1 \d p_1 \c p_d \d y_{t'} \d x_1$ is a hole of length $\ge 6$, with
three vertices in common with $D$, namely $p_2, x_1, y_{t'}$. From \ref{hole&antihole}, $D$ has
length 4, and so $t' = 1$ and $s = 2$. Since $W = \{x_2\}$ and $j = 2$, it follows that the only edges
between $x_1,x_2$ and $P$ are $x_1p_1,x_2p_1,x_2p_2$. But then the three paths $p_1\d x_1, x_2 \d z,
p_2 \c p_d \d y_1$ form a long prism, a contradiction. This proves (5).

\bigskip

From (5), since $p_h$ is adjacent to $x_1$, it follows that $p_h$ is not complete to $X \setminus x_1$,
and therefore $h < i < j$. Choose $s'$ with $1 \le s' \le s$ minimum such that $p_h$ is nonadjacent
to $x_{s'}$. So $p_j \d x_1 \c x_{s'} \d p_h \d p_j$ is an antihole, and so $s'$ is even.
Hence $x_1 \c x_{s'} \d p_h \d z$ is an odd antipath; all its internal vertices have neighbours
in $\{p_{h+1} \l p_n\}$, and its ends do not, so by \ref{RRR} it has length 3, that is,
$s' = 2$. The set $F = \{x_2, p_h \l p_n\}$ is connected; the only neighbour of $x_1$ in $F$ is
$p_h$; the only neighbour of $z$ in $F$ is $x_2$. Since $x_1,z$ are
$(X \setminus \{x_1,x_2\}) \cup Y$-complete, and $p_h,x_2$ are not (for $p_h$ is not $Y$-complete), it
follows from \ref{RRstrip} that there is a vertex in $(X \setminus \{x_1,x_2\}) \cup Y$ with no
neighbour in $F$ except possibly $x_2$. But every vertex in $(X \setminus \{x_1,x_2\}) \cup Y$ is adjacent
to either $p_j$ or to $p_n$, a contradiction. This proves \ref{tripleRR}. \bbox

\section{Pseudowheels}

Our next goal is to prove a version of \ref{oddwheel} for ``pseudowheels'', in which one of the
``vertices'' is really an anticonnected set. Fortunately we don't need to generalize \ref{oddwheel}
completely; it is enough to generalize the case when there is a segment of the wheel of length 1,
and one of its vertices has blown up to become the anticonnected set. (We did try to generalize
\ref{oddwheel} completely, but were unable to do it and it gave us a lot of trouble; so eventually
we found a way to make do with this special case.)

We begin with an even more special case, a form of \ref{nobanister} when one vertex is replaced by an
anticonnected set.
\begin{thm}\label{pseudohex}
Let $G \in \mathcal{F}_7$, and let  $X, Y$ be disjoint nonempty anticonnected subsets of $V(G)$,
complete to each other.
Let $p_1\d p_2 \d p_3 \d p_4 \d p_5$ be a track in $G\setminus (X \cup Y)$, induced except possibly
for the edge $p_2p_5$. Let $X$ be complete to $p_1,p_5$ and not to $p_2,p_3,p_4$. If $p_1,p_3,p_4$
are $Y$-complete then so is one of $p_2,p_5$.
\end{thm}
\Proof Assume not. Then in $\overline{G}$, $\{p_1,p_3,p_5\}$ is a triangle, and
the connected set $F = X \cup Y \cup \{p_2,p_4\}$ catches it. In $\overline{G}$,
the  only neighbours of $p_5$ in $F$ are in $Y \cup \{p_2\}$, the only neighbours of $p_3$ in $F$
are in $X$, and the only neighbour of $p_1$ in $F$ is $p_4$. Hence no vertex of $F$ has two neighbours
in the triangle, so by \ref{tricatch}, $F$ contains a reflection of the triangle. So (back in $G$)
there are vertices $b_1 \in X$ and $b_2 \in Y \cup \{p_2\}$ such that $b_1,b_2,p_4$ are pairwise
nonadjacent, and $b_1$ is adjacent to $p_1,p_5$ and not $p_3$, and $b_2$ is adjacent to $p_1,p_3$
and not $p_5$. Since $p_4$ is $Y$-complete and $b_2,p_4$ are nonadjacent it follows that
$b_2 \notin Y$, and so $b_2 = p_2$.  If $p_2,p_5$ are adjacent then $Y$ and the six vertices
$p_1\l p_5,b_1$ contradict \ref{nobanister}, and otherwise they form an odd wheel, a contradiction.
This proves \ref{pseudohex}. \bbox

There is a reformulation of \ref{doubleRR2} that we sometimes need:
\begin{thm}\label{doubleRRodd}
Let $G \in \mathcal{F}_7$, and let $X, Y$ be disjoint nonempty anticonnected subsets of $V(G)$,
complete to each other.  Let $P$ be a path in $G$ with even length $>0$, with vertices
$p_1,\ldots,p_n$ in order, so that $p_1$ is
$X$-complete, $p_n$ is not $X$-complete and $p_n$ is the unique $Y$-complete vertex of $P$.
Suppose that there is a $Y$-complete vertex in $G$ nonadjacent
to both $p_{n-1},p_{n-2}$.  Then either:
\begin{itemize}
\item  there is an odd number of $X$-complete edges in $P$, or
\item $n = 3$ and there is an odd antipath joining $p_{n-1}$ and $p_n$ with interior in $X$.
\end{itemize}
\end{thm}
\Proof
Choose an $X$-complete vertex $p_i$ in $P$ with $i$ maximum.
Suppose first that $i$ is even. Then the path $p_1 \c p_i$ is odd, and we may
assume that an even number of its edges are $X$-complete. So it has length $>1$; by \ref{evengap}, none
of its internal vertices are $X$-complete; and by \ref{RRR} it has length 3, and there is an odd
antipath $Q$ joining $p_2,p_3$ with interior in $X$. But there is another, $R$ say, with interior in $Y$,
which also
must be odd. Hence $R$ cannot be completed to an antihole via $p_3 \d p_n \d p_2$,
and so $n \le 4$. Since $n$ is odd it follows that $n = 3$; but then by hypothesis there is an
$Y$-complete vertex $v$ nonadjacent to $p_1,p_2$, and then $v \d p_1 \d R \d p_2 \d v$ is an odd antihole,
a contradiction.

Now we assume that $i$ is odd.
Hence the path $p_i \c p_n$ is even, and by \ref{doubleRR2} it has length 2.
Let $Q$ be the antipath between $p_{n-2},p_{n-1}$ with interior in $Y$, and $R$ the antipath between
$p_{n-1},p_n$ with interior in $X$.  By hypothesis there is a $Y$-complete vertex nonadjacent to
$p_{n-1},p_{n-2}$, and therefore $Q$ is even, so $R$ is odd by \ref{doubleRR2}. Hence $R$ cannot be
completed to an antihole via $p_n \d p_1 \d p_{n-1}$; and so $n = 3$ and the theorem holds.
This proves \ref{doubleRRodd}.\bbox

We need the following extension of \ref{evengap}.

\begin{thm}\label{pseudogap}
Let $G \in \mathcal{F}_7$, and let $X, Y$ be disjoint nonempty anticonnected subsets of $V(G)$,
complete to each other. Let $P$ be a path $p_1\c p_n$ of $G \setminus (X \cup Y)$, where $n \ge 5$,
such that $p_1,p_n$ are the only $X$-complete vertices of $P$. Then $P$ has even length. Assume that
at least two vertices of $P$ are $Y$-complete,
and let $P'$ be a maximal subpath of $P$ such that none of its internal vertices are $Y$-complete. Then
the length of $P'$ has the same parity as the number of ends of $P'$ that belong to $\{p_1,p_n\}$ and
are not $Y$-complete. Moreover, the number of $Y$-complete edges
of $P$ has the same parity as the number of elements of $\{p_1,p_n\}$ that are $Y$-complete.
\end{thm}
\Proof
Since $P$ is a path of
length $\ge 4$, and its ends are $X$-complete and its internal vertices are not, it follows that
$P$ has even length.  Let $P'$ be a maximal subpath of $P$ of length $\ge 2$ in which no internal
vertex is $Y$-complete, and assume first that both ends of $P'$ are $Y$-complete.
Suppose $P'$ has odd length, and let its ends be $p_i,p_j$ where $i < j$.  Then \ref{RRR} implies
that $j-i$ = 3, and
there is an odd antipath $Q$ joining $p_{i+1},p_{i+2}$ with interior in $Y$. Since $n\ge 5$, either
$n >j$ or $1<i$, and from the symmetry
between $p_1$ and $p_n$ we may assume the latter. Since $p_{i+1},p_{i+2}$
are not $X$-complete, they are joined by an antipath $Q'$ with interior in $X$. Since
$Q \cup Q'$ is an antihole it follows that $Q'$ is odd.  But then $p_n\d p_{i+1}\d Q' \d p_{i+2} \d p_n$
is an odd antihole, a contradiction. So in this case $P'$ has even length.
We may therefore assume that an end of $P'$ is not $Y$-complete, and from the maximality of $P'$,
any such end is either $p_1$ or $p_n$, and we may assume it is $p_n$ from the symmetry. The other
end of $P'$ is therefore not $p_1$ since at least two vertices of $P$ are $Y$-complete, and so
it is $p_i$, where $i$ is maximum with $2 \le i \le n$ such that $p_i$ is $Y$-complete.
Since $i > 1$, no vertex of $P'$
is $X$-complete except $p_n$. Suppose that $P'$ is even; then we may apply \ref{doubleRR2}.
We deduce that $P'$ has length 2,
and so $i = n-2$. Now the antipath joining $p_{n-2},p_{n-1}$ with interior in $X$ is even since it
can be completed to an antihole via $p_{n-1} \d p_1 \d p_{n-2}$; and the antipath joining
$p_{n-1},p_n$ with interior in $Y$ is even since it can be completed to an antihole via
$p_n \d p_h \d p_{n-1}$, where $p_h$ is some $Y$-complete vertex with $1 \le h <i$.
But this contradicts \ref{doubleRR2}. Consequently $P'$ is odd, as required.

By counting we see that that the number of $Y$-complete edges in $P$ has the same parity as
the number of ends of $P$ that are $Y$-complete. This proves \ref{pseudogap}.\bbox

Let us say a {\em pseudowheel} in a graph $G$ is a triple $(X,Y,P)$, satisfying:
\begin{itemize}
\item $X, Y$ are disjoint nonempty anticonnected subsets of $V(G)$, complete to each other
\item $P$ is a path $p_1\c p_n$ of $G \setminus (X \cup Y)$, where $n \ge 5$
\item $p_1,p_n$ are the only $X$-complete vertices of $P$
\item $p_1$ is $Y$-complete, and so is at least one other vertex of $P$; and $p_2,p_n$ are not $Y$-complete.
\end{itemize}

\begin{thm}\label{pseudogap2}
Let $(X,Y, P)$ be a pseudowheel in a graph $G \in \mathcal{F}_7$, where $P$ is $p_1\c p_n$.
Then $P$ contains an odd number, at least $3$, of $Y$-complete edges, and $P$ has length $\ge 6$.
\end{thm}
\Proof By \ref{pseudogap}, $P$ contains an odd number of $Y$-complete edges, since an odd number
of ends of $P$ are $Y$-complete. Suppose it only contains one, say $p_ip_{i+1}$. Since $p_2,p_n$ are
not $Y$-complete it follows that $3 \le i \le n-2$. So there is an antipath
joining  $p_i,p_{i+1}$ with interior in $X$, and by \ref{evenantipath3} applied to
the path $P$, this antipath has length 2, that is, there exists $x \in X$  nonadjacent to both
$p_i,p_{i+1}$. Let $C$ be a hole containing
$x, p_i,p_{i+1}$ and with $C \setminus x \subseteq P$. Then $(C,Y)$ is an odd wheel, since $C$ contains
the $Y$-complete vertices $x, p_i,p_{i+1}$ and it also contains $p_{i-1},p_{i+2}$ which are not
$Y$-complete, contrary to $G \in \mathcal{F}_7$. So at least three edges of $P$ are $Y$-complete, and therefore
$P$ has length $\ge 6$.  This proves \ref{pseudogap2}. \bbox

A pseudowheel $(X,Y,P)$ in $G$ is {\em optimal} if
\begin{itemize}
\item there is no pseudowheel $(X',Y',P')$ in $G$ such that the number of $Y'$-complete vertices in $P'$
is less than the number of $Y$-complete vertices in $P$, and
\item there is no pseudowheel $(X,Y',P)$ in $G$ such that $Y \subset Y'$.
\end{itemize}

\begin{thm}\label{pseudovnbrs}
Let $G \in \mathcal{F}_7$, and let $(X,Y, P)$ be an optimal pseudowheel in $G$, where $P$ is $p_1\c p_n$.
Let $v \in V(G) \setminus (X \cup Y \cup V(P))$, not $Y$-complete.  Then there is a subpath
$P'$ of $P$ such that
\begin{itemize}
\item $V(P')$ contains all the neighbours of $v$ in $P$,
\item there is no $Y$-complete vertex in the interior of $P'$, and
\item if $v$ is $X$-complete, then either $V(P') = \{p_1\}$, or $p_n \in V(P')$.
\end{itemize}
\end{thm}
\Proof Choose $h,k$ with $1 \le h \le k \le n$ such that $v$ is adjacent to $p_h,p_k$, with $h$ minimum
and $k$ maximum. (If this is impossible then the theorem holds.) Choose $i,j$
with $2 \le i \le j \le n$ such that $p_i,p_j$ are $Y$-complete, with $i$ minimum and $j$ maximum.
By \ref{pseudogap} it follows that $i$ is odd and $j$ is even, and $j-i \ge 3$ by \ref{pseudogap2},
since all $Y$-complete edges in $P$ lie in the path $p_i \c p_j$.
\\
\\
(1) {\em If $v$ is both adjacent to $p_1$ and $X$-complete then the theorem holds.}
\\
\\
For from the optimality of
$(X,Y,P)$ it follows that $(X, Y \cup \{v\},P)$ is not a pseudowheel, and so $p_1$ is the only
$Y \cup \{v\}$-complete vertex in $P$. By \ref{doubleRRC} (with $X,Y$ replaced by $Y \cup \{v\}, X$)
we deduce that either there exists $y \in Y \cup \{v\}$ nonadjacent to all $p_2,\ldots,p_n$, or
there exist nonadjacent $y_1,y_2 \in Y \cup \{v\}$ so that
$y_1\d p_2\d \cdots\d p_n\d y_2$ is a path. But $p_i$ is $Y$-complete and
$3 \le i \le n-1$, so the second statement does not hold; and the first holds only of $y = v$. This
proves (1).
\\
\\
(2) {\em We may assume that there is a path $Q$ from $v$ to some vertex $q$, such that
$q$ is the only $Y$-complete vertex in $Q$, and $V(Q \setminus v) \subseteq \{p_{i+1}\l p_{j-1}\}$.}
\\
\\
For by \ref{pseudogap} and the fact that there is a $Y$-complete edge in $P$, it follows that there is a
$Y$-complete vertex in $\{p_{i+1}\l p_{j-1}\}$. If $v$ has a neighbour in this set then the claim
holds, so suppose it does not.  We may assume $v$ has a neighbour in $\{p_1\l p_i\}$, for otherwise
the theorem holds. Suppose it also has a neighbour
in $\{p_j\l p_n\}$. Then there is a hole $C$ containing $v$, with $C \setminus v \subseteq P$,
such that $p_i \c p_j$ is a path of $C$. Since all $Y$-complete edges in $P$ belong to this path,
and there are an odd number of them, it follows that there is an odd number ($\ge 3$) of $Y$-complete
edges in $C$, contrary to \ref{evengap}. So $v$ has no neighbours in  $\{p_j\l p_n\}$, and hence
$k \le i$. We may therefore assume that $v$ is $X$-complete, so $k>1$ by (1).
The path $v\d p_k \c p_n$ has length $\ge 4$, and its ends are $X$-complete, and its internal
vertices are not, so by \ref{RRR} it has even length, and therefore the path $v \d p_k \c p_i$ is even.
But $v$ is the only $X$-complete vertex in $v \d p_k \c p_i$, and $p_i$ is its only $Y$-complete
vertex (since $k>1$), so by \ref{doubleRR2}, this path has length 2, and so $k = i-1$. There is no
odd antipath joining $v,p_k$ with interior in $Y$, since the $Y$-complete vertex $p_j$ is nonadjacent
to $v,p_k$; and there is no odd antipath joining $p_k,p_i$ with interior in $X$, since the $X$-complete
vertex $p_n$ is nonadjacent to  $p_k,p_i$, contrary to \ref{doubleRR2}.  This proves (2).
\\
\\
(3) {\em If $v$ is $X$-complete then the theorem holds.}
\\
\\
For then we may assume that $v$ is nonadjacent to $p_1$ by (1).
If $h$ is odd then $p_1\c p_h \d v$ is an odd path with ends
$X$-complete and its internal vertices not, so it has length 3 by \ref{RRR}; but the
$X$-complete vertex $p_n$ has no neighbour in its interior (since
$n \ge 5$), contrary to \ref{greentouch}. So $h$ is even.
Suppose that one of $p_2\l p_h$ is $Y$-complete. Then $h \not = 2$ since $p_2$ is not $Y$-complete, so
$h \ge 4$, and $h <j $ by (2). Hence $(X,Y,p_1\c p_h\d v)$ is a pseudowheel, not containing $p_j$,
contrary to the optimality of $(P,X,Y)$.  So there are no $Y$-complete vertices in  $\{p_2\l p_h\}$, and so
$i > h$.  Let $Q,q$ be as in (2).  Since the $X\cup Y$-complete vertex $p_1$ has no
neighbours in $Q$, the pairs $(V(Q),X),(V(Q),Y)$ are balanced by \ref{balancev}; so by \ref{doubleRR},
$Q$ has odd length. Hence the path $p_1 \c p_h \d v \d Q\d q$ has odd length, and its ends are
$Y$-complete, and its internal vertices are not. By \ref{RRR} it has length 3; so $h = 2$ and
$v$ is adjacent to $q$. Also every $Y$-complete vertex is adjacent to one of $v,p_2$, by \ref{greentouch},
so they are all adjacent to $v$ except $p_1$ and possibly $p_3$.
Suppose $p_i$ is adjacent to $v$, and is
therefore $Y \cup \{v\}$-complete. Then the path $p_1 \c p_i$ has even length; the only
$X$-complete vertex in it is $p_1$; and the only $Y \cup \{v\}$-complete vertex in it is $p_i$.
By \ref{doubleRR2} it has length 2. But there is an $X$-complete vertex nonadjacent to both
$p_2,p_3$, namely $p_n$, and there is a $Y\cup \{v\}$-complete vertex nonadjacent to both
$p_1,p_2$, since there is a $Y$-complete vertex in $\{p_4 \l p_n\}$ by \ref{pseudogap}, and it
is necessarily adjacent to $v$ as we already saw. Hence both pairs $(\{p_1,p_2\},Y \cup \{v\})$
and $(\{p_2,p_3\},X)$ are balanced by \ref{balancev}, contrary to \ref{doubleRR2}. This proves
that $p_i$ is not adjacent to $v$, and therefore $i = 3$. Choose $h' >i$ minimum so that $v$ is
adjacent to $p_{h'}$. From the hole $v\d p_2 \c p_{h'} \d v$ it follows that $h'$ is even.
From \ref{doubleRRodd} applied to the even path $p_3 \c p_{h'} \d v$, and using the fact that the
$X \cup Y$-complete vertex $p_1$ has no neighbour in this path, we deduce that there is a
$Y$-complete edge in $p_3 \c p_{h'} \d v$. Since $v$ is adjacent to every $Y$-complete vertex in $P$
except $p_1,p_3$, it follows that the only such edge is $p_3p_4$, and therefore $h' = 4$. But then
the track $p_1\c p_4 \d v$ violates \ref{pseudohex}. This proves (3).

\bigskip

Henceforth we may therefore assume that $v$ is not $X$-complete. If $k \le h+1$ then the theorem holds,
so we assume $k \ge h+2$.
\\
\\
(4) {\em If $v$ is not adjacent to $p_1$ then the theorem holds.}
\\
\\
For let $P'$ be the path
$p_1\c p_h\d v \d p_k \c p_n$. Suppose that any of $p_2\l p_h,p_k\l p_n$ is $Y$-complete. Then $P'$ has
length $\ge 4$, since $h > 1$ and $p_2,p_n$ are not $Y$-complete, and so $(X,Y,P')$ is a pseudowheel.
By the optimality of $(X,Y,P)$ it follows that there are no $Y$-complete vertices among
$\{p_{h+1}\l p_{k-1}\}$; but then the claim holds. So we may assume that none of
$p_2\l p_h,p_k\l p_n$ is $Y$-complete, and therefore $h < i \le j < k$, and since $j-i \ge 3$
it follows that $k -h \ge 5$.  Let $Q,q$ be as in (2).  Then
$q\d Q \d v \d p_k \c p_n$ is a path, $R$ say; the only $Y$-complete vertex in $R$ is $q$; the only
$X$-complete vertex in $R$ is $p_n$; and the $X\cup Y$-complete vertex $p_1$ has no neighbour in its
interior. By \ref{doubleRR}, $R$ is odd.
Therefore the paths $p_1\c p_h \d v \d Q \d q$ and $p_1 \c p_h \d v \d p_k \c p_n$ have
lengths of opposite parity. For the first path, its ends are $Y$-complete and its
internal vertices are not. For the second, its ends are $X$-complete and its internal vertices are not.
So one of them has length 3, and so $h = 2$, and there is an odd antipath joining $v,p_2$ with
interior in one of $X,Y$. Since $v,p_2$ are joined by an antipath with interior in $X$ and by
another with interior in $Y$, and all such pairs of antipaths have the same parity (since their
union is an antihole), it follows that $v,p_2$ are joined by an odd antipath with interior in each
of $X$,$Y$. Hence every $X$-complete vertex is adjacent to one of $v,p_2$, and so is every $Y$-complete
vertex. In particular $k = n$, and $v$ is adjacent to every $Y$-complete vertex in $P$ except
$p_1$ and possibly $p_3$. But then $R$ has length 2, contradicting
that it has odd length. This proves (4).

\bigskip

Henceforth then we assume that $v$ is adjacent to $p_1$ and not $X$-complete.
\\
\\
(5) {\em $p_{n-1}$ is not $Y \cup \{v\}$-complete.}
\\
\\
For suppose it is. Since
$n \ge 7$, it follows from \ref{RRR} applied to $P \setminus p_n$ and $Y \cup \{v\}$ that there is
a $Y \cup \{v\}$-complete vertex in $\{p_2\l p_{n-2}\}$; choose such a vertex, $p_{j'}$ say, with $j'$
maximum. Now $j = n-1$. If $j' < j-1$ then $j-j'$ is even from \ref{greentouch} applied to $p_j' \c p_j$; but then
the odd path $p_{j'} \c p_n$ contains no $Y \cup \{v\}$-complete edges, and $p_1$ is $X$-complete,
$Y \cup \{v\}$-complete and has no neighbours in the path $p_{j'} \c p_n$, contrary to \ref{tripleRR}.
So $j' = j-1$. Let $F = X \cup Y \cup \{v,p_{n-1}\}$. Then $F$ is anticonnected, and
each of $p_1,p_{n-2},p_n$ has a nonneighbour in $F$; the only nonneighbour of $p_1$ in $F$ is $p_{n-1}$;
all nonneighbours of $p_{n-2}$ in $F$ belong to $X$; and all nonneighbours of $p_n$ in $F$ belong
to $Y \cup \{v\}$. So in $\overline{G}$, the connected set $F$ catches the triangle
$\{p_1,p_{n-2},p_n\}$, and by \ref{tricatch} it contains a reflection of the triangle,
which is impossible since $p_{n-1}$ is complete (in $G$) to $Y \cup \{v\}$. This proves (5).
\\
\\
(6) {\em $v$ is not adjacent to $p_n$.}
\\
\\
For suppose it is. By \ref{pseudogap2} there are at least three $Y$-complete edges in $P$, and so
there is a $Y$-complete vertex $p_a$ in $P$ with $a\ge 3$, even and
different from $p_{n-1}$.
Thus $j-a$ is even, and so by \ref{evengap} there is an even number of $Y$-complete edges in the even path
$p_a\c p_j$, and hence in the odd path $p_a \c p_n$. But $p_a$ is $Y$-complete, and $p_n$ is the
unique $X \cup \{v\}$-complete vertex in this path, contrary to \ref{tripleRR}. This proves (6).
\\
\\
(7) {\em There is no neighbour $p_m$ of $v$ in $P$ with $1 \le m \le n$ so that $v, p_m$ are joined by
an odd antipath with interior in $Y$.}
\\
\\
For suppose such a neighbour exists. So $1 < m < n$ by (6), and there is an antipath
joining $v,p_m$ with interior in $X$, which therefore is also odd, since its union with
the antipath through $Y$ is an antihole. Since it cannot be completed to an odd antihole via
$p_m\d p_n\d v$, it follows that $m = n-1$, and in particular $m$ is even.
Since $j$ is even, either $p_j = p_m$ or $p_j$ is nonadjacent to $p_m$; and in
either case it follows that $p_j$ is adjacent to $v$, since every $Y$-complete vertex is adjacent
to one of $v, p_m$.  By (5) $n-j \ge 3$ and odd,  and the path $p_j\c p_n$
(with anticonnected sets $X$ and $Y \cup \{v\}$) violates \ref{tripleRR}.  This proves (7).

\bigskip

Suppose that $j \ge k$, and let $P'$ be the path $p_1\d v \d p_k \c p_n$. Then $P'$ has length $\ge 4$,
since $p_{n-1}$ is not $Y \cup \{v\}$-complete, and so $(X,Y,P')$ is a pseudowheel; and by the
optimality of $(X,Y,P)$ it follows that there are
no $Y$-complete vertices in $p_2 \c p_{k-1}$, contrary to (2). So $j < k$.
Let $Q,q$ be as in (2), and assume first that $Q$ is even. Then the path $p_1\d v \d Q \d q$ has odd
length; its ends are $Y$-complete, and its internal vertices are not, so by \ref{RRR} it has length 3,
and its internal vertices are joined by an odd antipath with interior in $Y$, contrary to (7).
So $Q$ is odd.

\bigskip

Next assume that $k$ is even. Then the path $p_1\d v \d p_k \c p_n$ is odd, and its ends
are $X$-complete, and its internal vertices are not, so by \ref{RRR} it has length 3, and $k = n-1$, and
its internal vertices $v, p_{n-1}$ are joined by an odd antipath with interior in $X$. Since $p_{n-1}$
is not $Y \cup \{v\}$-complete, they are also joined by an odd antipath with interior in $Y$,
contrary to (7).
This proves that $k$ is odd. Hence the path $q \d Q\d v \d p_k \c p_n$ is even, and by (6) it
has length $>2$ contrary to \ref{doubleRR2}. This proves \ref{pseudovnbrs}.\bbox

\begin{thm}\label{pseudoFnbrs}
Let $G \in \mathcal{F}_7$, and let $(X,Y, P)$ be an optimal pseudowheel in $G$, where $P$ is $p_1\c p_n$.
Let $F \subseteq V(G) \setminus (X \cup Y \cup V(P))$ be connected, such that no vertex in $F$
is $Y$-complete. Then there is a subpath
$P'$ of $P$ such that
\begin{itemize}
\item $V(P')$ contains all the attachments of $F$ in $P$,
\item there is no $Y$-complete vertex in the interior of $P'$, and
\item if some vertex of $F$ is $X$-complete then either $V(P') = \{p_1\}$ or $p_n \in V(P')$.
\end{itemize}
\end{thm}
\Proof
Suppose the theorem is false, and choose a minimal counterexample $F$. From \ref{pseudovnbrs}
$|F| \ge 2$.
\\
\\
(1) {\em Some vertex in $F$ is $X$-complete.}
\\
\\
For suppose not. Since $F$ is a counterexample, it has attachments $p_a,p_c$ such that there is
a $Y$-complete vertex $p_b$ with $a < b < c$. From the minimality of $F$, $F$ is the interior
of a path $p_a \d f_1 \c f_k \d p_c$. Let
$W_1$ be the set of attachments of $F \setminus f_k$ in $P$, and $W_2$ the set of attachments of
$F \setminus f_1$ in $P$. From the minimality of $F$, for $i = 1,2$ there is a subpath
$p_{a_i}\c p_{b_i}$ of $P$ ( $ = P_i$ say), so that no internal vertex of $P_i$ is $Y$-complete,
and $W_i \subseteq V(P_i)$. Choose $P_1,P_2$ minimal; then $p_{a_1}$ is a neighbour of some member of
$F \setminus f_k$, and therefore of $f_1$ from the minimality of $F$, and similarly
$p_{b_2}$ is a neighbour of $f_k$, and $p_1 \c p_{a_1}\d f_1 \c f_k \d p_{b_2} \c p_n$ is a path
$P'$ say. Suppose that there is a $Y$-complete vertex in $P'$ different from $p_1$. Then $P'$
has length $\ge 4$, and $(X,Y,P')$ is a pseudowheel, contrary to the optimality of $(X,Y,P)$.
So there are no $Y$-complete vertices in $P'$. But also there are none in $\{p_{a_1+1}\l p_{b_1 -1}\}$
and none in $\{p_{a_2+1}\l p_{b_2 -1}\}$, so all the $Y$-complete vertices of $P$ belong to
$\{p_{b_1}\l p_{a_2}\}$, except for $p_1$. By \ref{pseudogap2} there are an odd number, at least 3, of
$Y$-complete edges in this path. From the minimality of $F$,
$f_1\c f_k \d p_{a_2} \d p_{a_2-1} \c p_{b_1} \d f_1$ is a hole, which therefore also contains
an odd number $\ge 3$ of $Y$-complete edges. But this
contradicts \ref{evengap}. This proves (1).

\bigskip

From (1), since $F$ is a counterexample, there is a $Y$-complete vertex in $F$ and there exists
$a$ with $1 < a<b$ such that $p_a$ is an attachment of $F$ and $p_b$ is $Y$-complete. From the minimality
of $F$, there is a path $p_a\d f_1\c f_k$ such that $F = \{f_1 \l f_k\}$ and $f_k$ is the unique
$X$-complete vertex in $F$.  Let
$W_1$ be the set of attachments of $F \setminus f_k$ in $P$, and $W_2$ the set of attachments of
$F \setminus f_1$ in $P$. From the minimality of $F$, for $i = 1,2$ there is a subpath
$p_{a_i}\c p_{b_i}$ of $P$ ( $ = P_i$ say), so that no internal vertex of $P_i$ is $Y$-complete,
and $W_i \subseteq V(P_i)$, and either $b_2 = n$ or $a_2 = b_2 = 1$.

First assume that $b_2 = n$.
Choose $P_1,P_2$ minimal; then $p_{a_1}$ is a neighbour
of $f_1$, and $p_1 \c p_{a_1}\d f_1 \c f_k$ is a path
$P'$ say. Suppose that there is a $Y$-complete vertex in $P'$ different from $p_1$. Then $P'$
has length $\ge 4$, and $(X,Y,P')$ is a pseudowheel, contrary to the optimality of $(X,Y,P)$.
So there are no $Y$-complete vertices in $P'$. But also there are none in $\{p_{a_1+1}\l p_{b_1 -1}\}$
and none in $\{p_{a_2+1}\l p_{b_2 -1}\}$, so all the $Y$-complete vertices of $P$ belong to
$\{p_{b_1}\l p_{a_2}\}$, except for $p_1$. By \ref{pseudogap2} there are an odd number, at least 3, of
$Y$-complete edges in this path. From the minimality of $F$,
$f_1\c f_k \d p_{a_2} \d p_{a_2-1} \c p_{b_1} \d f_1$ is a hole, which therefore also contains
an odd number $\ge 3$ of $Y$-complete edges. But this contradicts \ref{evengap}.

So we may assume that $a_2 = b_2 = 1$, and that $p_1 \in W_2$, and therefore $b_1 > 1$.
From the minimality of $F$ there are
no edges between $F \setminus f_1$ and $V(P) \setminus p_1$.  Choose $P_1$ minimal. So
$p_{b_1}$ is adjacent to $f_1$, and either $a_1 = 1$ or $p_{a_1}$ is adjacent to $f_1$. Suppose
first that an odd number of edges of the path $p_1\c p_{a_1}$ are $Y$-complete. Hence $p_1$
has no neighbours in $F \setminus f_k$, and so $f_1\c f_k \d p_1 \c p_{a_1} \d f_1$ is a hole.
It contains an odd number of $Y$-complete edges, and at least three $Y$-complete vertices, because
$p_1$ is $Y$-complete and $p_2$ is not, a contradiction to \ref{evengap}. So there are an even
number of $Y$-complete edges in the path $p_1\c p_{a_1}$, and therefore an odd number in
$p_{b_1}\c p_n$, since there are an odd number in $P$, and none in $P_1$.  Therefore there are an
odd number in the path $f_k \c f_1 \d p_{b_1} \c p_n$ ($ = R$ say). But an edge of $p_{b_1}\c p_n$
is $Y$-complete and $p_n$ is not, so $b_2 \le n-2$; and since $k \ge 2$, it follows that $R$
has length $\ge 4$. Also, at least two vertices of $R$ are $Y$-complete, and its ends are not
$Y$-complete, and its ends are its only $X$-complete vertices. This contradicts \ref{pseudogap},
So there is no such $F$. This proves \ref{pseudoFnbrs}.\bbox

Now we come to the main result of this section, \ref{summary}.8, which we restate, the following.

\begin{thm}\label{pseudowheel}
Let $G \in \mathcal{F}_7$. If it contains a pseudowheel then it admits a balanced skew partition. In
particular, every recalcitrant graph belongs to $\mathcal{F}_8$.
\end{thm}
\Proof
Suppose $G$ contains a pseudowheel; then it contains an optimal pseudowheel, say $(X,Y, P)$,
where $P$ is $p_1\c p_n$. Let $Z$ be the set of all $Y$-complete vertices in $G$. So
$Y,Z$ are disjoint, nonempty, and complete to each other, and $|Z| \ge 2$. Let $F_0 = V(G) \setminus (Y \cup Z)$.
 By \ref{bipattach},  we may assume that $F_0$ is connected and every vertex in $Z$ has a neighbour in $F_0$,
for otherwise the theorem holds. Choose
$i>1$ so that $p_ip_{i+1}$ is $Y$-complete, and let $A$,$B$ be the two components of $V(P\setminus p_i)$.
Since $p_1,p_{i+1}$ both have neighbours in $F_0$, it follows that $F_0$ contains a minimal
connected set so that there are vertices in $A$ and in $B$ with neighbours in $F$. From the minimality
of $F$ it is disjoint from $V(P)$; and disjoint from $X \cup Y$ since $X \subseteq Z$, contrary to
\ref{pseudoFnbrs}. This proves \ref{pseudowheel}.\bbox

\section{Wheel systems - the base case}

Our next goal is to show that if there is a wheel in a member
of $\mathcal{F}_8$ then the graph admits a balanced skew partition, and in particular that
there is no wheel in a recalcitrant graph. Assuming there is no balanced skew partition, the strategy is
to show that
there is no anticonnected set which is maximal such that there is a wheel of which it is a hub. In other words,
we want to show that given any wheel, there is a second wheel whose hub is a proper superset of the hub
of the first wheel. To prove this, we need to generalize wheels into what we call ``wheel systems'', and prove
a long and complicated theorem about wheel systems, by induction on some parameter.  The base case of this
inductive proof is the theorem of this section. (Incidentally, we will not need the
hypothesis that there is no pseudowheel in $G$ for several more sections, so what we are proving here
is true also for graphs in $\mathcal{F}_7$, and we formulate it that way, although we only need
it for graphs in $\mathcal{F}_8$.)

\begin{thm}\label{wheelbase}
Let $G \in \mathcal{F}_7$, and let $Y,A \subseteq V(G)$, so that $Y$ is anticonnected and
$A$ is connected. Let $x_0,x_1,v,z \in V(G) \setminus (Y \cup A)$ be distinct, satisfying:
\begin{itemize}
\item $x_0,x_1$ are nonadjacent and are both $Y$-complete,
\item $z$ is adjacent to $x_0,x_1,v$ and has no neighbour in $A$,
\item no vertex in $A$ is adjacent to both $x_0,x_1$,
\item $v$ is nonadjacent to one of $x_0,x_1$, and is not $Y$-complete,
\item every vertex in $Y$ that is nonadjacent to $v$ has a neighbour in $A$ and is adjacent to $z$, and
\item $x_0,x_1,v$ all have neighbours in $A$.
\end{itemize}
Then $z$ is $Y$-complete and there is a wheel $(C,Y)$ in $G$ with
$x_0,x_1,z \in V(C) \subseteq \{x_0,x_1,z\} \cup A$.
\end{thm}
\Proof We assume by induction that the result holds for all smaller sets $Y$; and for fixed $Y$, we
assume that $|A|$ is minimum satisfying the hypotheses.
\\
\\
(1) {\em There exists $y \in Y$ adjacent to $z$ and with a neighbour in $A$,
so that $Y \setminus y$ is empty or anticonnected.}
\\
\\
For if $|Y| = 1$, let $Y = \{y\}$; then since $v$ is not $Y$-complete it follows that $y$ is nonadjacent
to $v$, and therefore is adjacent to $z$ and has a neighbour in $A$
and the claim holds. So assume $|Y|>1$, and choose distinct
$y_1,y_2 \in Y$ so that $Y \setminus y_i$ is anticonnected (i = 1,2). Not both $y_1,y_2$ is the unique
nonneighbour of $v$ in $Y$; so we may assume that $v$ is not $Y\setminus y_2$-complete.
By the minimality of $|A|+|Y|$, $z$ is $Y\setminus y_2$-complete and there is a
$Y\setminus y_2$-complete vertex in $A$; and in particular, $y_1$ is adjacent to $z$ and has a
neighbour in $A$, so we may set $y = y_1$. This proves (1).

\bigskip

Let $y$ be as in (1), and let $Y' = Y \setminus y$.
\\
\\
(2) {\em Either $v$ is $Y'$-complete and nonadjacent to $y$, or $z$ is $Y$-complete and there is
a path $x_0 \d p_1 \c p_n \d x_1$ from $x_0$ to $x_1$ with interior in $A$, containing at least two
$Y'$-complete edges.}
\\
\\
For if $v$ is $Y'$-complete the first assertion holds, so we assume not; and in particular
$Y'$ is nonempty. By induction $z$ is $Y'$-complete, and therefore $Y$-complete, and there is a path
as in the claim. This proves (2).
\\
\\
(3) {\em There is no connected $F \subseteq A$ containing neighbours of all of $x_0,x_1,v,y$ except $A$ itself.}
\\
\\
For suppose there is. From the minimality of $A$, some member of $Y$ has no neighbour in $F$ and is
nonadjacent to $v$. In particular, $v$ is not $Y'$-complete, so $Y'$ is nonempty and by (2), at
least two vertices of $A$ are $Y'$-complete. Since $F \not = A$, there exists $f \in A \setminus F$ so that
$A \setminus f$ is connected. But every vertex in $Y \cup \{x_0,x_1,v\}$ has a neighbour in $A \setminus f$;
for all members of $Y'$ have at least two neighbours in $A$ (since $A$ contains two $Y'$-complete vertices),
and $x_0,x_1,v,y$ have neighbours in $F$. This contradicts the minimality of $A$, and therefore proves (3).

\bigskip

Let $x_0 \d p_1 \c p_n \d x_1$ be a path from $x_0$ to $x_1$ with interior in $A$, and let $C$ be the hole
$z\d x_0\d p_1 \c p_n \d x_1 \d z$.
\\
\\
(4) {\em If any vertex of $p_1\l p_n$ is $Y\cup \{v\}$-complete then $z$ is $Y$-complete; and if $z$
is $Y$-complete then we may assume that no edge of $x_0 \d p_1 \c p_n \d x_1$ is $Y$-complete.
In particular neither of $p_1,p_n$ is $Y\cup\{v\}$-complete.}
\\
\\
For let $p_i$ be $Y\cup\{v\}$-complete, say, and suppose $z$ is not $Y$-complete. By (2), $v$ is $Y'$-complete
and nonadjacent to $y$. Let $Q$ be an antipath between $z,y$ with interior in $Y'$, and let $R$ be an
antipath between $v,p_i$ with interior in $\{x_0,x_1\}$. Then $z\d Q \d y \d v \d R \d p_i \d z$ is an
antihole, meeting the hole $C$ in at least three vertices, contrary to \ref{hole&antihole}. This proves
the first assertion. The second is immediate, for otherwise $(C,Y)$ satisfies the theorem. For the third,
note that if say $p_n$ is $Y\cup\{v\}$-complete, then $p_nx_1$ is a $Y$-complete edge, a contradiction.
This proves (4).
\\
\\
(5) {\em If $x_0$ is adjacent to $v$, then $v$ is nonadjacent to all of $p_2\l p_n$.}
\\
\\
For suppose $v$ is adjacent to one of
$p_2\l p_n$, and choose $i$ with $2 \le i \le n$ maximum such that $v$ is adjacent to $p_i$.
Suppose first that $i = n$. Since $x_0,x_1,p_n$ belong to $C$, there is no
antihole of length $\ge 5$ containing them by \ref{hole&antihole}. By (4), $p_n$ is not
$Y$-complete, and hence there is an antipath between $p_n,v$ with
interior in this set, and it can be completed via $v\d x_1 \d x_0 \d p_n$ to an antihole of length $\ge 5$
containing $x_0,x_1,p_n$, a contradiction.  So $i < n$.

Since the hole $C$ is even, it follows that $n$ is odd.
From the hole $z\d v \d p_i \c p_n \d x_1 \d z$ it follows that $i$ is odd. Since $i>1$,
$x_0 \d v \d p_i \c p_n \d x_1$ is an odd path of length $\ge 5$. Its ends are
$Y \cup \{z\}$-complete, and its internal vertices are not, so by \ref{RRR}, $Y \cup \{z\}$
is not anticonnected. Hence $z$ is $Y$-complete.  The ends of the same path are both
$Y$-complete, so by \ref{RRR}, some edge
of the path is $Y$-complete.  Since $v$ is not $Y$-complete, this edge belongs to $C$, contrary to (4).
This proves (5).

\bigskip

Let us choose $C$ so that either $v$ is $Y'$-complete or $(C,Y')$ is a wheel (this is possible by (2)).
\\
\\
(6) {\em If $x_0$ is adjacent to $v$, then not both $v,y$ have neighbours in $\{p_1\l p_n\}$.}
\\
\\
For if they do, then by (5) $p_1$ is the only neighbour of $v$ in $\{p_1\l p_n\}$.
Suppose first that $v$ is adjacent to $y$. By (2), $z$ is $Y$-complete, and $(C,Y')$
is a wheel, and so every vertex in
$Y'$ has a neighbour in $\{p_2\l p_n\}$. By (4) $p_1$ is not $Y$-complete.
Therefore $z,x_0$ are the only $Y\cup \{v\}$-complete vertices in $C$,
and by \ref{RRC} there is a hat or a leap. Since all vertices in $Y'$ have a neighbour
in $\{p_2\l p_n\}$, and $y$ is adjacent to $x_1$, it follows that there is no hat, and so
$y,v$ form a leap, a contradiction since they are adjacent. So $v$ is nonadjacent to $y$.
Choose $j$ with $1 \le j \le n$ minimum so that $y$ is adjacent to $p_j$. From the hole
$z\d v \d p_1 \c p_j \d y \d z$ we deduce that $j$ is odd, and
therefore $x_0\d p_1 \c p_j\d y \d x_0$ is not a hole, that is, $j = 1$, and hence $p_1$ is adjacent to $y$.
By (4) $p_1$ is not $Y'$-complete. If $v$ is $Y'$-complete,
then an antipath between $p_1$ and $y$ with interior in $Y'$ can be extended to an antihole via
$y\d v \d x_1 \d p_1$, and this antihole shares the vertices $p_1,x_1,v$ with the hole
$z\d v \d p_1 \c p_n \d x_1 \d z$, contrary to \ref{hole&antihole}.
So $v$ is not $Y'$-complete. By (2), $z$ is $Y$-complete, and $(C,Y')$ is a wheel. By \ref{wheelvnbrs}
applied to the wheel $(C,Y')$ and vertex $v$, it follows that $p_1$ is $Y'$-complete and therefore
$Y\cup \{v\}$-complete, contrary to (4). This proves (6).
\\
\\
(7) {\em Not both $v,y$ have neighbours in $\{p_1\l p_n\}$.}
\\
\\
For by (6) we may assume that $v$ is nonadjacent to $x_0$, and similarly nonadjacent to $x_1$.
Choose $i$ with $1 \le i \le n$ maximum so that $v$ is adjacent to $p_i$.  From the hole
$z\d v \d p_i \c p_n \d x_1 \d z$ it follows that $i$ is odd.  Suppose first that $v$ is not
$Y'$-complete. By (2), $z$ is $Y$-complete and $(C,Y')$ is a wheel.
By \ref{wheelvnbrs}, $p_i,z$ have the same wheel-parity, and so there are an odd number of $Y'$-complete
edges in $p_i\c p_n\d x_1$. By (4) no edge of the path $x_0\d p_1 \c p_n \d x_1$ is $Y$-complete.
Consequently $zx_1$ is the unique $Y$-complete edge of the hole
$z\d v\d p_i \c p_n \d x_1 \d z$ ($ = C_1$ say). Suppose that $y$ is
nonadjacent to all $v,p_i\l p_n$.  Now $y$ has a neighbour in
$\{p_1\l p_n\}$ by hypothesis, so $\{p_1\l p_n,v\}$ ($= F$ say)
catches the triangle $\{z,x_1,y\}$. The only neighbour of $z$ in $F$ is $v$; the only neighbour of
$x_1$ in $F$ is $p_n$; and $y$ is nonadjacent to both $v,p_n$ by assumption. By \ref{tricatch},
$F$ includes a reflection of the triangle; but then $i= n$ and there is an antihole of length 6 using
$z,x_1,p_n$, contrary to \ref{hole&antihole}. This proves that $y$ is adjacent to one of $v,p_i \l p_n$.
Since there is an odd number of $Y'$-complete edges in the path $p_i\c p_n\d x_1$, it follows
that every member of $Y$ is adjacent to one of $v,p_i \l p_n$.  Consequently
$Y$ contains no hat for $C_1$.  Assume that $C_1$ has length $\ge 6$. By \ref{RRC},
$Y$ contains a leap, so there are nonadjacent $y_1,y_2 \in Y$ so that
$y_1\d v \d p_i \c p_n \d y_2$ is a path, of odd length $\ge 5$. But the ends of this path are
$\{x_0,x_1\}$-complete and its internal vertices are not, contrary to \ref{RRR}. So $C_1$ has length 4,
that is, $i = n$, and $p_n$ is $Y'$-complete.
By (4) it follows that $p_n$ is nonadjacent to $y$, and therefore
$y$ is adjacent to $v$ (since $Y$ contains no hat for $C_1$).
Let $Q$ be an antipath between $v,y$ with interior in $Y'$; then $x_1 \d v \d Q \d y \d p_n$
is an antipath, and $z$ is complete to its interior, and $x_0$ is complete to all its interior
except $x_1$, contrary to \ref{hole&antipath2} applied to this antipath and the hole $C$.
This proves (7) assuming that $v$ is not $Y'$-complete.

We therefore assume that $v$ is $Y'$-complete, and consequently nonadjacent to $y$.
Now $\{v,p_1\l p_n\}$ is connected and catches the triangle $\{z,x_1,y\}$. By
\ref{hole&antihole}, it contains no reflection of the triangle, since as before that would give an antihole
of length 6 with 3 vertices in $C$.  So by \ref{tricatch}, there is a vertex in $\{v,p_1\l p_n\}$
with two neighbours in the triangle. The only neighbour of $z$ in it is $v$, which is nonadjacent
to both $x_1,y$. The only neighbour of $x_1$ in it is $p_n$, and therefore $y$ is adjacent
to $p_n$.  We recall that $i$ is maximum such that $v$ is adjacent to $p_i$. Since $y$ is
adjacent to $p_n$, we may choose $j$ with $i \le j \le n$ minimum such that $y$ is adjacent to $p_j$.
From the hole $z \d v \d p_i \c p_j \d y \d z$ we see that $j$ is odd.
Suppose $j \not = i$. Then the path $v \d p_i \c p_j \d y$ is even and has length $\ge 4$.
By \ref{doubleRR2} with anticonnected sets $\{x_0,x_1\}, Y' \cup \{z\}$ we deduce that $Y' \cup \{z\}$ is not
anticonnected, and hence $z$ is $Y$-complete. By \ref{doubleRRodd} with sets $\{x_0,x_1\},Y',$ since the
$\{x_0,x_1\}\cup Y'$-complete vertex $z$ has no neighbours in the interior of $P$, it follows that
there are an odd number of $Y'$-complete edges in the path $v \d p_i \c p_j \d y$. Since $y$ is not
$Y'$-complete, they all belong to the path $v \d p_i \c p_j$. Let $C_1$ be the hole
$z\d v\d p_i \c p_n \d x_1 \d z$; then $(C_1,Y')$ is a wheel. Now $p_n$ is not $Y$-complete by
(4), and therefore not $Y'$-complete. Since there is no $Y$-complete edge in the odd path
$p_j \c p_n$, and the $Y$-complete vertex $z$ has no neighbour in its interior, it follows from
\ref{greentouch} that $p_j$ is not $Y$-complete and hence not $Y'$-complete. Consequently $p_j,p_n$
have opposite wheel-parity (with respect to the wheel $(C_1,Y')$) and yet are both not $Y'$-complete,
and so this wheel is an odd wheel, contrary to $ G \in \mathcal{F}_7$.
This proves that $j = i$,  that is, $y$ is adjacent to $p_i$.

Now suppose $i < n$. If $p_i$ is not $Y$-complete then an antipath between $p_i$ and $y$ with
interior in $Y$ can be extended via $y\d v \d x_1 \d p_i$ to an antihole sharing the vertices
 $p_i,x_1,v$ with the hole $z\d v\d p_i\d p_n \d x_1 \d z$ ($= C_1$ say), contrary to \ref{hole&antihole}.
So $p_i$ is $Y$-complete, and therefore so is $z$, by (4). But then $(C_1,Y)$ is an odd wheel,
since $z,x_1,p_i$ are $Y$-complete and $v,p_n$ are not (by (4)), contrary to $ G \in \mathcal{F}_7$.
So $i = n$, and hence $p_n$ is adjacent to both $v,y$. From the symmetry
between $x_0,x_1$ it follows that $p_1$ is adjacent to both $v,y$. By (4), $p_1,p_n$
are not $Y$-complete.  So in $\overline{G}$, the connected set $Y \cup \{p_1,p_n\}$
catches the triangle $\{x_0,x_1,v\}$; $x_0,x_1,v$ all have unique neighbours in it, namely
$p_n,p_1,y$ respectively; and these three vertices do not form a triangle since $yp_1$ is
not an edge (of $\overline{G}$), contrary to \ref{tricatch}. This proves (7).
\\
\\
(8) {\em If $v$ is nonadjacent to $y$ and adjacent to one of $x_0,x_1$ then the theorem holds.}
\\
\\
For assume $v$ is nonadjacent to $y$ and adjacent to $x_0$ say.
Now $A \cup \{x_1\}$ catches the triangle $\{z,x_0,v\}$; it contains no reflection of this triangle,
since $x_0,x_1$ have no common neighbour in $A$; and the unique neighbour of $z$ in this set
is nonadjacent to both $x_0,v$. So by \ref{tricatch} it follows that there is a vertex
in $A$ adjacent to both $x_0,v$.
Also, $A\cup {v}$ catches the triangle $\{z,x_1,y\}$. Suppose that $A\cup \{v\}$ contains a
reflection of this triangle; then there exist $f\in A$ adjacent to $x_1,v$ and not to $y$.
Since $f \in A$ it follows that $f$ is nonadjacent to $x_0$; but then
$f \d v \d x_0 \d y \d x_1 \d f$
is an odd hole, a contradiction. Hence  by \ref{tricatch} there is a vertex
in $A$ adjacent to both $x_1,y$.
Consequently from (3), $A$ is the vertex set of a path $f_1\c f_k$, where
$f_1$ is adjacent to $x_0,v$, and $f_k$ to $x_1,y$. Since $f_1 \in A$ it follows that
$f_1$ is not adjacent to $x_1$.

Now assume that $f_1$ is not the unique neighbour of $v$ in $A$. From (3),
$f_1$ is the unique neighbour of $x_0$ in $A$. By (7), $f_k$ is not the unique
neighbour of $x_1$ in $A$, and so from (3) it is the unique neighbour of $y$ in $A$.
In particular $y$ is not adjacent to $f_1$.
Both $x_0,z$ have unique neighbours in $A \cup \{x_1\} = F$ say, namely
$f_1,x_1$ respectively. Now $x_0,z$ are both $\{v,y\}$-complete, and $f_1,x_1$ are not.
Since $F \setminus x_1$ is connected,
this contradicts \ref{RRstrip2}. So  $f_1$ is the unique neighbour of $v$ in $A$.
Suppose that $f_k$ is the unique neighbour
of $y$ in $A$. Then both $z,y$ have unique neighbours in $A \cup \{v\}$, namely $v,f_k$ respectively;
and $z,y$ are $\{x_0,x_1\}$-complete, and $v,f_k$ are not. Once again this contradicts \ref{RRstrip2}.
So $f_k$ is not the unique neighbour of $y$ in $A$, and therefore it is the unique neighbour of $x_1$
in $F$.

Suppose that $f_k$ is $Y$-complete. Since $f_k = p_n$, it follows from (4) that $z$ is not $Y$-complete;
and so $v$ is $Y'$-complete by (2), and an antipath between
$z,y$ with interior in $Y'$ can be extended to an antihole via $y \d v \d f_k \d z$, which shares the
vertices $z,v,f_k$ with the hole $z\d v \d f_1 \c f_k \d x_1 \d z$ ($=C_1$ say), contrary to \ref{hole&antihole}.
So $f_k$ is not $Y$-complete and therefore not $Y'$-complete (and in particular, $Y'$ is nonempty).

Suppose that $z$ is not $Y$-complete; and therefore $Y' \cup \{z\}$ is anticonnected, and $v$ is
$Y'$-complete.
Choose $h$ with $1 \le h < k$ minimum so that $f_h$ is adjacent to $y$ (this exists since $f_k$
is not the unique neighbour of $y$ in $A$). The path $v \d f_1 \c f_h \d y$ is even, since it
can be completed to a hole via $y \d z \d v$, and therefore the path $v \d f_1 \c f_h \d y \d x_1$ is odd
(this is a path since $f_k$ is the unique neighbour of $x_1$ in $A$); and the ends of this path are
$Y' \cup \{z\}$-complete, and its internal vertices are not. By \ref{RRR} it has length 3. So $f_1$
is adjacent to $y$ and $v$. If
$f_1$ is not $Y'$-complete, then an antipath between $f_1,y$ with
interior in $Y'$ can be completed to an antihole via $y \d v\d x_1 \d f_1$, which shares the vertices
$v,x_1,f_1$ with the hole $C_1$, contrary to \ref{hole&antihole}; while if $f_1$ is $Y$-complete, then
an antipath between $z,y$ with interior in $Y'$ can be completed to an antihole via $y \d v \d x_1 \d
f_1 \d z$, again contrary to \ref{hole&antihole}. This proves that $z$ is $Y$-complete.

In the hole $C_1$, $z,x_1$ are $Y$-complete and $v,f_k$ are not; so since $ G \in \mathcal{F}_7$, no other vertex
of $C_1$ is $Y$-complete. By \ref{RRC}, $Y$ contains a leap or hat for $C_1$. From a hypothesis of
the theorem, every vertex in $Y$
has a neighbour in $A \cup\{v\}$, so there is no hat, and hence there exist nonadjacent $y_1,y_2$
in $Y$ so that $y_1\d v\d f_1 \c f_k \d y_2$ is a path. Since both ends of this path are
$\{x_0,x_1\}$-complete, and no internal vertex is $\{x_0,x_1\}$-complete, this contradicts \ref{RRR}.
This proves (8).
\\
\\
(9) {\em There is no connected $F \subseteq A$ containing neighbours of $x_0,x_1,v$ except $A$ itself.}
\\
\\
For suppose that such a set $F$ exists with $F \not = A$, and choose $f \in A\setminus F$ so that
$A \setminus f$ is connected. From the minimality of $A$, there exists $y' \in Y$ nonadjacent to $v$ with
no neighbour in $A \setminus f$, and therefore $f$ is the unique neighbour of $y'$ in $A$. If $y' \in Y'$,
then $v$ is not $Y'$-complete, and therefore by (2) there are two $Y'$-complete vertices in $A$, a contradiction.
So $y' = y$, and therefore $y$ is not adjacent to $v$.  Suppose that $v$ is not
adjacent to $f$. Then both $z,y$ have unique neighbours in $A \cup \{v\}$, namely $v,f$; $z,y$
are $\{x_0,x_1\}$-complete, and $v,f$ are not; $f\d y\d z \d v$ is a path; and $x_0,x_1$ both have
neighbours in $A$, contrary to \ref{RRstrip2}. So $v$ is adjacent to $f$. By (8) we may assume that
$v$ is nonadjacent to both $x_0,x_1$. Since $f$ is not $\{x_0,x_1\}$-complete, we may assume from
the symmetry that $f$ is nonadjacent to $x_1$. Now $A \cup \{v\}$ catches the triangle
$\{z,y,x_1\}$; the only neighbour of $z$ in $A \cup \{v\}$ is $v$; the only neighbour of $y$ in
$A \cup \{v\}$ is $f$; $v,f$ are both nonadjacent to $x_1$; and so by \ref{tricatch}, $A \cup \{v\}$
contains a reflection of the triangle. Hence there exists $f_1 \in A\setminus f$, adjacent to $v,f,x_1$ and
not to $y$ (and therefore not to $x_0$).
Since every path between $x_0,x_1$ with interior in $A$ has length $\ge 4$ it follows that
$x_0$ is nonadjacent to $f,f_1$, and this restores the symmetry between $x_0,x_1$; and consequently
by the same argument there exists $f_0 \in A \setminus f$ adjacent to $v,f,x_0$ and not to $y,x_1$.
Since $z\d x_0 \d f_0 \d f_1 \d x_1 \d z$ is not an odd hole, $f_0$ is nonadjacent to $f_1$; but then
$x_0\d f_0\d f \d f_1 \d x_1$ violates (7).  This proves (9).

\bigskip

From (7) and (9), it follows that there exists $f \in A$ such that $A \setminus f$ is
connected, $f$ does not belong to $C$, and $f$ is the unique neighbour of $v$ in $A$.
\\
\\
(10) {\em If $v$ is adjacent to one of $x_0,x_1$ then the theorem holds.}
\\
\\
For let $v$ be adjacent to $x_0$ say. Suppose first that $x_0$ is not adjacent to $f$.
Then $A \cup \{x_1\}$ catches the triangle $\{z,v,x_0\}$; the only neighbour
of $z$ in $A \cup \{x_1\}$ is $x_1$; the only neighbour of $v$ in $A \cup \{x_1\}$ is $f$;
$x_1,f$ are both nonadjacent to $x_0$; and $A \cup \{x_1\}$ contains no reflection of the triangle since
that would give a 6-antihole with
3 vertices in common with $C$, contradicting \ref{tricatch}. So $x_0$ is adjacent to $f$, and
therefore $x_1$ is nonadjacent to both $v,f$.
By (8) we may assume that $v$ is adjacent to $y$, and therefore not $Y'$-complete. By (2)
$z$ is $Y'$-complete and $(C,Y')$ is a wheel. Let $v\d q_1\c q_k\d x_1$
be a path between $v,x_1$ with interior in $A$ (so $f = q_1$) and let $C_1$ be the hole
$z\d v\d q_1\c q_k\d x_1 \d z$. From (9), $A = \{q_1\l q_k\}$. Since $q_k = p_n$ and
$z$ is $Y$-complete, it follows from (4) that $q_k$ is not $Y$-complete.  Since $(C_1,Y)$ is not
an odd wheel, it follows that $(C_1,Y)$ is not a wheel, and so $z,x_1$ are the only $Y$-complete
vertices in $C_1$, by \ref{evengap}. By \ref{RRC}, $Y$ contains a leap or hat for $C_1$. But $y$ is
adjacent to $v$, and all other vertices of $Y$ have at least two neighbours in $\{p_1\l p_n\}$, which
is a subset of $\{q_1\l q_k\}$, a contradiction. This proves (10).

\bigskip

We may therefore assume that $v$ is nonadjacent to both $x_0,x_1$. From
(9) one of $x_0,x_1$ has a unique neighbour in $A$, and from the
symmetry we may assume it is $x_1$. Let its neighbour be $f_1$. By (7) and (9), $v$ has
no neighbour in $\{p_1\l p_n\}$, and in particular $f \not = f_1$.
Let $Q$ be a path in $A$ between $f,f_1$, say $f = q_1\c q_k = f_1$, so
$z\d v \d q_1 \c q_k \d x_1 \d z$ is a hole ($C_1$ say).
\\
\\
(11) {\em $z$ is not $Y'$-complete, and $v$ is $Y'$-complete and nonadjacent to $y$.}
\\
\\
For assume $z$ is $Y'$-complete.  So $z,x_1$ both have unique neighbours in $A \cup \{v\}$, namely $v,f_1$.
By (4), $f_1$ is not $Y$-complete. So $z,x_1$ are $Y$-complete,
and $v,f_1$ are not. By \ref{RRstrip2}, it follows that some vertex in $Y$
has no neighbour in $A$. But $y$ has a neighbour in $A$ by (1), and so some vertex in $Y'$ has no neighbour
in $A$. In particular, there is no $Y'$-complete vertex in $A$, and so by (2), $v$ is $Y'$-complete
and nonadjacent to $y$. From \ref{RRstrip2} applied to the path $v\d z \d x_1 \d f$ and the anticonnected set
$\{y\}$, it follows that $y$ is adjacent to $f_1$.  Since $(C_1,Y)$ is not an odd wheel, it follows from \ref{RRC}
that $Y$ contains a leap or a hat for $C_1$. Since all members
of $Y'$ are adjacent to $v$ and $y$ is adjacent to $f_1$, there is no hat, and the leap must use
$y$; so we may assume $y,y'$ is a leap for some $y' \in Y'$. Hence $y\d f_1 \d Q \d f \d v \d y'$ is a path.
Since this path has odd length $\ge 5$, and its ends are $\{x_0,x_1\}$-complete and its internal
vertices are not, this contradicts \ref{RRR}. So $z$ is not $Y'$-complete. The claim follows from
(2). This proves (11).
\\
\\
(12) {\em $y$ is nonadjacent to all of $q_1\l q_{k-1}$.}
\\
\\
For suppose not, and choose $i$ with $1 \le i \le k$ minimum so that $y$ is adjacent to $q_i$.
From the hole $z\d v \d q_1 \c q_i \d y \d z$ it follows that $i$ is odd. So by (10) we may assume that
$v \d q_1 \c q_i \d y \d x_1$ is an odd path. Its ends are $Y' \cup \{z\}$-complete,
its internal vertices are not, and $Y' \cup \{z\}$ is anticonnected by (11); so it has
length 3 by \ref{RRR}, that is, $i = k$ and $y$ is adjacent to $f$.
If $f$ is not $Y'$-complete, an antipath between $f,y$ with interior in $Y'$
can be completed to an antihole via $y\d v\d x_1 \d f$, sharing the vertices $v,x_1,f$
with $C_1$, contrary to \ref{hole&antihole}. So $f$ is $Y'$-complete. Since $z$ is not, an
antipath between $z,y$ with interior in $Y'$ can be completed to an antihole via $y\d v\d x_1 \d f \d z$,
again contrary to \ref{hole&antihole}. This proves (12).

\bigskip

To conclude, $A \cup \{v\}$ catches $\{y,z,x_1\}$, and so by \ref{tricatch}, $y$ is adjacent to $f_1 = q_k$.
Suppose that $x_0$ is adjacent to one of $q_1\l q_k$. Then  $\{p_1\l p_n\} \subseteq \{q_1\l q_k\}$
from the minimality of $A$, and so
the neighbours of $y$ in $C$ are precisely $x_0,z,x_1,q_k=p_n$, contrary to \ref{evengap} applied
to $C$ and $y$. So  $x_0$ is nonadjacent to all of $q_1\l q_k$; but then
$v\d q_1\c q_k \d y \d x_0$ is an odd path of length $\ge 5$, its ends are $Y'\cup\{z\}$-complete,
and its internal vertices are not, contrary to \ref{RRR}.
This proves \ref{wheelbase}.\bbox

\section {Wheel systems}

We continue with the proof that there is no wheel in a recalcitrant graph. The remainder of the argument
breaks into two parts. First, we invent an object called
a ``$Y$-diamond wheel system'', where $Y$ is an anticonnected subset of $V(G)$; and in this section we
prove that if $G$ contains a $Y$-diamond, then it also contains a wheel $(C,Y)$ (for the same set $Y$).
And we prove the same thing for another special kind of wheel system, a ``$Y$-square''.
The second half of the
argument, in the next few sections, is roughly as follows. Suppose there is a wheel, and choose
one, say $(C,Y')$, with $Y'$ maximal. One can argue that, since $G$ has no skew partition, there is
a vertex $y$ that could be added to $Y'$, forming a new anticonnected set $Y$ say, so that there are
precisely two $Y$-complete edges in $C$ and they are consecutive; and moreover there is a ``tail'',
a path from $y$
to a vertex of $C$ that is not $Y'$-complete, with certain special properties. We grow $C$ into
what we call a ``wheel system'' relative to $Y$, and make it maximal, and analyze how the remainder
of $G$ attaches to this wheel system. Since there is no skew partition in $G$, we can show that
$G$ contains a $Y$-diamond or $Y$-square. But therefore $G$ also contains a wheel with hub $Y$, contrary
to the maximality of $Y_0$.

Let $G$ be a graph. A {\em frame} in $G$ is a pair $(z,A_0)$, where $z \in V(G)$, and $A_0$ is
a nonnull connected subset of $V(G)\setminus z$, containing no neighbours of $z$.
With respect to a given frame $(z,A_0)$,  a {\em wheel system} in $G$ of {\em height $t\ge 1$}
is a sequence $x_0\l x_t$ of distinct vertices of $G \setminus (A_0 \cup \{z\})$,
satisfying the following conditions:
\begin{enumerate}
\item $A_0$ contains neighbours of $x_0$ and of $x_1$, and no vertex in $A_0$ is $\{x_0,x_1\}$-complete.
\item For $2 \le i \le t$, there is a connected subset of $V(G)$ including $A_0$, containing a
neighbour of $x_i$, containing no neighbour of $z$, and containing no $\{x_0\l x_{i-1}\}$-complete vertex.
\item For $1 \le i \le t$, $x_i$ is not $\{x_0\l  x_{i-1}\}$-complete.
\item $z$ is adjacent to all of $x_0\l x_t$.
\end{enumerate}

Note that this definition is symmetric between $x_0,x_1$, so $x_1,x_0,x_2\l x_t$ is another wheel system.
Let $x_0\l x_t$ be a wheel system of height $t$. For $1 \le i \le t$ we define
$X_i = \{x_0\l x_i\}$, and we define $A_i$ to be the maximal connected subset of $V(G)$
that includes $A_0$, contains no neighbour of $z$, and contains no $X_i$-complete vertex.
So for each $i$, $A_{i-1} \subseteq A_i$. Note that condition 2 above just says
that $x_i$ has a neighbour in $A_{i-1}$.

We need three special kinds of wheel systems. Let  $x_0\l x_t$ be a wheel system,
and define $X_i,A_i$ as above. Let $Y \subseteq V(G)$ be nonempty and anticonnected, such that
$Y$ is disjoint from $\{z,x_0\l x_t\}$, and $x_0\l x_{t-1}$ are all $Y$-complete and $x_t$ is not.
We say $x_0\l x_t$ is a
\begin{itemize}
\item  {\em $Y$-diamond} if
$t \ge 3$,
$x_t$ is $X_{t-2}$-complete, and
$x_t$ has a neighbour in $A_{t-2}$
\item {\em $Y$-square} if
$t \ge 3$,
$x_t$ is adjacent to $x_{t-1}$,
$x_t$ has no neighbour in $A_{t-2}$, and
there is a vertex in $A_{t-1}$ adjacent to $x_t$ with a neighbour in $A_{t-2}$
\item {\em $Y$-square-diamond} if
$t \ge 4$,
$x_t$ is $X_{t-2}$-complete, $x_{t-1}$ is not $X_{t-3}$-complete,
$x_t$ has no neighbour in $A_{t-3}$, $x_{t-1}$ has a neighbour in $A_{t-3}$, and
there is a vertex in $A_{t-2}$ adjacent to both $x_t,x_{t-1}$ with a neighbour in $A_{t-3}$.
\end{itemize}
(Note that $Y$-square-diamonds are $Y$-diamonds.)

The main result of this section is the following.
\begin{thm}\label{betterwheel}
Let $ G \in \mathcal{F}_7$, let  $(z,A_0)$ be a frame, and let $Y \subseteq V(G)\setminus (A_0\cup \{z\})$
be nonempty and
anticonnected. Suppose that there is either a $Y$-diamond, a $Y$-square, or a $Y$-square-diamond in $G$.
Then $z$ is $Y$-complete and $G$ contains a wheel $(C,Y)$.
\end{thm}
\noindent{\bf Proof of \ref{betterwheel}, assuming \ref{t=3}, \ref{betterwheel1}, \ref{betterwheel2}, and
\ref{betterwheel3}.}
\\
\\
(1) {\em If $Y_1 \subseteq V(G)$ is anticonnected with $Y \subseteq Y_1$, and there is a $Y_1$-square-diamond
in $G$ of height $t+1$, we may assume that there is an anticonnected $Y_2$ with
$Y_1 \subseteq Y_2 \subseteq V(G)$
such that there is a $Y_2$-square or $Y_2$-diamond in $G$ of height $\le t$.}
\\
\\
We proceed by induction on $t$.
By \ref{t=3} it follows that $t \ge 4$, and from \ref{betterwheel3}, we may assume that
there is an anticonnected superset
$Y_2$ of $Y_1$ such that $G$ contains a $Y_2$-square-diamond of height $t$. From the inductive
hypothesis, we may assume there is an anticonnected superset $Y_3$ of $Y_2$ such that $G$ contains a
$Y_3$-diamond or $Y_3$-square of height $t-1$; but then the claim holds. This proves (1).

\bigskip

From (1) we may assume that there is an anticonnected superset $Y_1$ of $Y$ such that for some
$t$, $G$ contains a $Y_1$-square or $Y_1$-diamond, of height $t$ say. Choose $t$ minimum.
If $t=3$ the theorem follows from \ref{t=3}, so we may assume $t \ge 4$. From (1) there is no
anticonnected superset $Y_2$ of $Y_1$ such that $G$ contains a $Y_2$-square-diamond of height $\le t$;
and the result follows from \ref{betterwheel1} and \ref{betterwheel2}. This proves
\ref{betterwheel}.\bbox

In \ref{betterwheel} we do not require that $z$ is $Y$-complete. This is included in the statement
of \ref{betterwheel}
merely because in the
inductive proof of \ref{betterwheel}, it is helpful to be proving the stronger statement. When we apply
the result later, we know that $z$ is $Y$-complete anyway.

First we show:

\begin{thm}\label{t=3}
Let $ G \in \mathcal{F}_7$, let  $(z,A_0)$ be a frame, and let $Y \subseteq V(G)\setminus (A_0\cup \{z\})$
be nonempty and anticonnected. There is no $Y$-square of height $3$ or $Y$-square-diamond of height $4$ in $G$;
and if $x_0\l x_3$ is a $Y$-diamond of height $3$, then $z$ is $Y$-complete and $G$ contains a wheel
$(C,Y\cup \{x_3\})$.
\end{thm}
\Proof
Let $x_0\l x_t$ be a wheel system in $G$, and let $X_i,A_i$ be defined as before.
Suppose first that $x_0\l x_t$ is a $Y$-square of height 3. So $t = 3$,
$x_3$ is adjacent to $x_2$, $x_3$ has no neighbour in $A_1$, and
there is a vertex $q$ in $A_2$ adjacent to $x_3$ with a neighbour in $A_1$.
From the maximality of $A_1$ it follows that $q$ is $X_1$-complete, and therefore nonadjacent to $x_2$
(since it belongs to $A_2$ and so is not $X_2$-complete). Let
$Q$ be a path from $q$ to $x_2$ with interior in $A_1$; so $Q$ has length $\ge 2$.
But $Q$ is even since it can be completed to a hole via $x_2\d x_3 \d q$, and so $q \d Q \d x_2 \d z$
is an odd path; its ends are $X_1$-complete, and its internal vertices are not.
By \ref{RRR} it has length 3, and there is an antipath with interior in $X_1$, joining its
middle vertices ($x_2$ and $r$ say). This antipath can be completed via $r\d z \d q \d x_2$ to an antihole
of length $\ge 6$, containing $x_0,x_1$ and $z$. But let $P$ be a path from $x_0$ to $x_1$ with
interior in $A_0$; then it has length $\ge 3$ since $A_0$ contains no vertex adjacent to both
$x_0,x_1$, and hence $z \d x_0 \d P \d x_1 \d z$ is a hole of length $\ge 6$ containing $x_0,x_1$ and $z$.
But this contradicts \ref{hole&antihole}, as required.

Now suppose $x_0\l x_t$ is a $Y$-square-diamond of height 4. So $t = 4$, $x_4$ is $X_2$-complete,
$x_3$ is not $X_1$-complete, $x_4$ has no neighbour in $A_1$, $x_3$ has a neighbour in $A_1$, and
there is a vertex $q$ in $A_2$ adjacent to both $x_4,x_3$ with a neighbour in $A_1$.
As before $q$ is $X_1$-complete, and therefore not adjacent to $x_2$; let
$Q$ be a path from $q$ to $x_2$ with interior in $A_1$. The proof is completed exactly as in the
previous paragraph.

So now we may assume that $x_0\l x_t$ is a $Y$-diamond of height 3. So $t = 3$,  $x_3$ is
$X_1$-complete (and therefore nonadjacent to $x_2$), and $x_3$ has a neighbour in $A_1$.
But then from \ref{wheelbase} with $A = A_1$,$v = x_1$ and anticonnected set $Y \cup \{x_3\}$, the
result follows. This proves \ref{t=3}. \bbox

We remark that the pieces of this jigsaw do not seem to fit well together. There is some annoying wastage
in \ref{t=3}; we produce a wheel with hub $Y\cup\{x_3\}$, and all we use in proving \ref{betterwheel}
is that there is a wheel with hub $Y$. Perhaps there is a better way to organize it, but so far it eludes us.

\begin{thm}\label{betterwheel1}
Let $ G \in \mathcal{F}_7$, let $(z,A_0)$ be a frame, and let $Y \subseteq V(G)\setminus (A_0\cup \{z\})$
be nonempty and anticonnected. Let $x_0\l x_t$ be a $Y$-diamond in $G$ of height $t \ge 4$. Suppose that there
is no anticonnected set $Y'$ with $Y \subseteq Y' \subseteq V(G)$ such that either:
\begin{itemize}
\item there is a $Y'$-diamond in $G$ of height $t-1$, or
\item there is a $Y'$-square in $G$ of height $t-1$, or
\item there is a $Y'$-square-diamond in $G$ of height $t$.
\end{itemize}
Then $z$ is $Y$-complete and $G$ contains a wheel $(C,Y)$.
\end{thm}
\Proof
Assume that either $z$ is not $Y$-complete or $G$ contains no wheel $(C,Y)$. Define $X_i,A_i$ as usual. So
$x_t$ is $X_{t-2}$-complete, and $x_t$ has a neighbour in $A_{t-2}$, and $Y$ is complete to $X_{t-1}$
and not to $x_t$.
\\
\\
(1) {\em Not both $x_t$ and $x_{t-1}$ have neighbours in $A_{t-3}$.}
\\
\\
For suppose they do. If $x_{t-1}$ is $X_{t-3}$-complete, then
\[x_0\l x_{t-1}\]
is a $Y \cup \{x_t\}$-diamond of height $t-1$, while if $x_{t-1}$ is not $X_{t-3}$-complete, then
\[x_0\l x_{t-3},x_{t-1},x_t\]
is a $Y$-diamond of height $t-1$, in both cases a contradiction. This proves (1).
\\
\\
(2) {\em There is a vertex $q$ in $A_{t-2}$ adjacent to both $x_t$ and $x_{t-1}$,
and a path $R$ in $A_{t-2}$ from $q$ to $A_{t-3}$ so that not both $x_t$ and $x_{t-1}$
have neighbours in $A_{t-3} \cup V(R\setminus q)$.}
\\
\\
For let $F$ be a minimal connected subgraph of $A_{t-2}$ including $A_{t-3}$ and containing neighbours
of both $x_t$ and $x_{t-1}$. If $x_t,x_{t-1}$ have a common neighbour in $F$, then the claim is
satisfied (from the minimality of $F$), so we assume not.
Let $P$ be a path between $x_t$ and $x_{t-1}$ with interior in $F$, say
$x_t \d p_1 \c p_n \d x_{t-1}$. Hence $P$ has length $>2$, and
the hole $z \d x_1 \d P \d x_2 \d z$ (= $C$ say) it follows that $P$ is even.
The only $X_{t-2}$-complete vertices in $C$ are $z$ and $x_t$, so by \ref{RRC}, $X_{t-2}$
contains a leap or a hat for $C$. Suppose it contains a leap; then there are nonadjacent
$x_i,x_j \in X_{t-2}$ so that $x_i \d p_1 \c p_n \d x_{t-1} \d x_j$ is an odd path. Since
$x_i,x_j$ are $Y \cup \{x_t\}$-complete, it follows from \ref{RRR} that this path contains another
$Y \cup \{x_t\}$-complete vertex, which must be $p_1$ since no others are adjacent to $x_t$.
Its ends are also $Y \cup \{x_t,z\}$-complete, and no internal vertex is $Y \cup \{x_t,z\}$-complete, so
by \ref{RRR}, $Y \cup \{x_t,z\}$ is not anticonnected, that is, $z$ is $Y$-complete. But then
let $C_1$ be the hole $z\d x_i \d p_1 \c p_n \d x_{t-1} \d z$; then $(C_1,Y)$ is a wheel,
a contradiction.

So $X_{t-2}$ contains a hat for $C$; that is, there exists $x_i \in X_{t-2}$
with no neighbours in $C$ except $x_t,z$. Hence the path $x_i \d x_t \d p_1 \c p_n \d x_{t-1}$ is odd
and has length $\ge 5$, and its ends are $Y \cup \{z\}$-complete, and no internal vertex is
$Y \cup \{z\}$-complete, so by \ref{RRR}, $z$ is $Y$-complete.  Let $S$ be a path between
$x_i$ and $x_{t-1}$ with interior in $F$. Then $V(S \cup P) \setminus \{x_i,x_t\}$ ( = $F'$ say)
is connected and catches the triangle $\{z,x_i,x_t\}$. The only neighbour of $z$ in $F'$ is
$x_{t-1}$, which is nonadjacent to both $x_i,x_t$. If $F'$ contains a reflection of the triangle,
there is an antihole of length 6 containing $z,x_{t-1},x_t$, which is impossible by \ref{hole&antihole}
since these three vertices belong to $C$. So by \ref{tricatch}, there is a vertex in $F'$ adjacent to
both $x_i,x_t$. Since $x_i$ has no neighbour in $P \setminus x_t$, it follows that
both $x_t,x_{t-1}$ have neighbours in the interior of $S$, and so there is a path $P'$ between
$x_t,x_{t-1}$ with $P' \setminus x_t$ a subpath of $S \setminus x_i$. As before $P'$ has length $\ge 4$,
and so $S$ has length $\ge 4$, and $P',S$ both have even length since they can be completed to holes
through $z$.  Since the $X_{t-2}$-complete vertex $z$ has no neighbours in the interior of $P'$, from
\ref{doubleRR2} (applied to $P'$ with anticonnected sets $Y$ and $X_{t-2}$) it follows that
there is a $Y$-complete edge in $P'$, and since $x_t$ is not $Y$-complete, there is
therefore one in $S$.
But since the edges $zx_{t-1},zx_i$ are also $Y$-complete, we deduce that there are at least three
$Y$-complete edges in the hole $z \d x_i \d S \d x_{t-1} \d z$, and so that hole is the rim
of a wheel with hub $Y$, a contradiction. This proves (2).

\bigskip

Choose $q, R$ as in (2) with $R$ minimal, and let $R$ be $r_1\c r_n$, where $r_1 = q$ and $r_n$ is the only
vertex of $R$ in $A_{t-3}$ .
\\
\\
(3) {\em $x_{t-1}$ has neighbours in $A_{t-3}$.}
\\
\\
For assume not.
Since $x_{t-1}$ has no neighbours in $A_{t-3}$ it follows that $q \notin A_{t-3}$, and so $R$ has
length $>0$.
Suppose first that every antipath between $x_{t-1}$ and $q$ with interior in $X_{t-2}$ is odd, and
let $Q$ be such an antipath. Since all internal vertices of $Q$ have neighbours in $A_{t-3}$, and
$z$ is complete to its interior and anticomplete to $A_{t-3}$, it follows from \ref{greentouch}
applied in $\overline{G}$ that one end of $Q$ has a neighbour in  $A_{t-3}$. By hypothesis,
$x_{t-1}$ does not, so $q$ does. From the maximality of $A_{t-3}$ it follows that $q$ is $X_{t-3}$-complete;
and since $q \in A_{t-2}$ and is therefore not $X_{t-2}$-complete, $q$ is nonadjacent to $x_{t-2}$.
Now by assumption, every every antipath between $x_{t-1}$ and $q$ with interior in $X_{t-2}$ is odd,
and so $x_{t-2}$ is adjacent to $x_{t-1}$. But then
\[x_0\l x_{t-1}\]
is a $Y \cup \{x_t\}$-square of height $t-1$, a contradiction.
So we may assume some antipath $Q$ between between $x_{t-1}$ and $q$ with interior in $X_{t-2}$ is
even.

From (2), not both $x_t,x_{t-1}$ have neighbours in $A_{t-3} \cup V(R \setminus q)$.
Suppose that $x_{t-1}$ has such a neighbour, and so $x_t$ does not. Since by assumption
$x_{t-1}$ has no neighbours in $A_{t-3}$, it follows that all neighbours of $x_{t-1}$ in
$A_{t-3} \cup V(R \setminus q)$ lie in the interior of $R$, and in particular $R$ has length $\ge 2$.
The antipath $x_t \d x_{t-1} \d Q \d q$ is odd, and its ends have no neighbours in the connected set
$A_{t-3} \cup \{r_3\l r_n\}$. Since $z$ is complete to its interior and anticomplete to
$A_{t-3} \cup \{r_3\l r_n\}$, it follows from \ref{greentouch} applied in $\overline{G}$
that some internal vertex of this antipath has no neighbours in $A_{t-3} \cup \{r_3\l r_n\}$.
But all internal vertices of $Q$ lie in $X_{t-2}$ and therefore have neighbours in $A_{t-3}$;
so $x_{t-1}$ has no neighbour in $A_{t-3} \cup \{r_3\l r_n\}$. Hence $r_2$ is its only neighbour in
$A_{t-3} \cup V(R \setminus q)$. Suppose that every antipath between $x_{t-1}$ and $r_2$ with interior
in $X_{t-2}$ is odd, and let $Q'$ be such an antipath. All internal vertices of $Q'$ have neighbours
in the connected set $A_{t-3}$, and $z$ is complete to the interior of $Q'$ and anticomplete to
 $A_{t-3}$; so by \ref{greentouch} applied in $\overline{G}$, it follows that $r_2$ has neighbours
in $A_{t-3}$. From the maximality of $A_{t-3}$, $r$ is $X_{t-3}$-complete, and therefore not adjacent
to $x_2$. Since by assumption every antipath between $x_{t-1}$ and $r_2$ with interior
in $X_{t-2}$ is odd, it follows that $x_{t-1}$ is adjacent to $x_{t-2}$. But then
\[x_0\l x_{t-1}\]
is a $Y \cup \{x_t\}$-square of height $t-1$, a contradiction.
So some antipath $Q'$ between $x_{t-1}$ and $r_2$ with interior
in $X_{t-2}$ is even. Hence the antipath $x_{t-1} \d Q' \d r_2 \d z$ is odd. All its internal vertices
have neighbours in the connected set $A_{t-3}\cup \{r_3\l r_n\}$ and its ends do not, so by
\ref{RRR} this antipath has length 3, that is, $Q'$ has length 2. Let $x_i$ be its middle vertex.
Then the connected set $A_{t-3} \cup V(R\setminus\{r_1,r_2\}) \cup \{x_i,x_t,z\}$  ( $ = F$ say)
catches the triangle $\{r_1,r_2,x_{t-1}\}$; the only neighbours of $r_1$ in $F$ are $x_t$ and
possibly $x_i$; the neighbours of $r_2$ in $F$ lie in $A_{t-3} \cup \{r_3\}$; and the
only neighbour of $x_{t-1}$ in $F$ is $z$. This contradicts \ref{tricatch}, since $z$ has no
neighbour in $A_{t-3} \cup \{r_3\}$.

So $x_{t-1}$ has no neighbours in  $A_{t-3} \cup V(R \setminus q)$.  Now the antipath
$z \d q \d Q \d x_{t-1}$ is odd, and all its internal vertices have neighbours in
$A_{t-3} \cup V(R \setminus q)$, and its ends do not, so by \ref{RRR} it has length 3, that is,
$Q$ has length 2 (let its middle vertex be $x_i$); and there is an odd path $P$ between $q,x_i$
with interior in $A_{t-3} \cup V(R \setminus q)$. Let $C$ be the hole $z \d x_{t-1} \d q \d P \d x_i \d z$;
then $C$ has length $\ge 6$. By \ref{hole&antihole} there is no antihole of length $\ge 6$ containing
$q,x_i,x_{t-1}$. If $q$ is not $Y$-complete then an antipath between $q,x_t$ with interior in $Y$
can be completed to such
an antihole via $x_t\d x_{t-1} \d x_i \d q$, so $q$ is $Y$-complete; and if $z$ is not $Y$-complete,
an antipath between $z$ and $x_t$ with interior in $Y$ can be extended to such an antihole, via
$x_t \d x_{t-1} \d x_i \d q \d z$. So $z$ is also $Y$-complete.  Now both $x_i,q$ are
$Y \cup \{x_t\}$-complete, and therefore by \ref{greentouch} applied to $P$, it follows
that $P$ contains a $Y \cup \{x_t\}$-complete edge. Hence the hole $C$ contains at least three
$Y$-complete edges, a contradiction. This proves (3).

\bigskip

From (3) and the choice of $R$ it follows that $x_t$ has no neighbours in  $A_{t-3} \cup V(R\setminus q)$.
Let $Q$ be an antipath between $q$ and $x_{t-1}$ with interior in $X_{t-2}$. Then
$z \d q \d Q \d x_{t-1} \d x_t$ is an antipath of length $\ge 4$, and its ends have no neighbours in
the connected set $A_{t-3} \cup V(R \setminus q)$, and its internal vertices do, so by \ref{RRR}
it has even length, that is, $Q$ is even. The antipath $x_t \d x_{t-1} \d Q \d q$ is therefore odd,
and its internal vertices have neighbours in $A_{t-3}$, and $z$ is complete to its interior and anticomplete
to $A_{t-3}$, so by \ref{greentouch} applied in $\overline {G}$, it follows that one of its ends, and hence
$q$ has a neighbour in $A_{t-3}$. From the maximality of $A_{t-3}$ it follows that $q$ is
$X_{t-3}$-complete and therefore nonadjacent to $x_{t-2}$. If $x_{t-1}$ is not $X_{t-3}$-complete, then
\[x_0\l x_t\]
is a $Y$-square-diamond of height $t$; while if $x_{t-1}$ is $X_{t-3}$-complete, then
\[x_0\l x_{t-1}\]
is a $Y \cup \{x_t\}$-diamond of height $t-1$, in both cases a contradiction.
This proves \ref{betterwheel1}.\bbox

\begin{thm}\label{betterwheel2}
Let $ G \in \mathcal{F}_7$, let $(z,A_0)$ be a frame, and let $Y \subseteq V(G)\setminus (A_0\cup \{z\})$
be nonempty and
anticonnected. Let $x_0\l x_t$ be a $Y$-square in $G$ of height $t \ge 4$. Then there is an anticonnected
set $Y'$ with  $Y \subseteq Y' \subseteq V(G)\setminus (A_0 \cup \{z\})$ such that either:
\begin{itemize}
\item there is a $Y'$-diamond in $G$ of height $t-1$, or
\item there is a $Y'$-square in $G$ of height $t-1$, or
\item there is a $Y'$-square-diamond in $G$ of height $t$.
\end{itemize}
\end{thm}
\Proof
Assume that no such $Y'$ exists.  Define $X_i,A_i$ as usual. So $x_t$ is adjacent to $x_{t-1}$,
$x_t$ has no neighbour in $A_{t-2}$, there is a vertex $q$ in $A_{t-1}$ adjacent to $x_t$ with a
neighbour in $A_{t-2}$, and $Y$ is complete to $X_{t-1}$ and not to $x_t$.
From the maximality of $A_{t-2}$ it follows that $q$ is $X_{t-2}$-complete. Since $q \in A_{t-1}$, it is not
$X_{t-1}$-complete, and so $q$ is nonadjacent to $x_{t-1}$.
\\
\\
(1) {\em $x_{t-1}$ has neighbours in $A_{t-3}$.}
\\
\\
For suppose not.  Let $R$ be a path between $q$ and $x_{t-1}$
with interior in $A_{t-2}$. Then $R$ has length $\ge 2$, and from the hole $q\d R \d x_{t-1} \d x_t \d q$
it follows that $R$ has even length. So the path $q\d R \d x_{t-1} \d z$ is odd, and its ends are
$X_{t-2}$-complete, and its interior vertices are not, so by \ref{RRR} it has length 3, that is,
$R$ has length 2. Let its middle vertex be $r$. Since $x_{t-1}$ has no neighbour in $A_{t-3}$, it
follows that $r \in A_{t-2} \setminus A_{t-3}$.  Let $Q$ be an antipath
between $r$ and $x_{t-1}$ with interior in $X_{t-2}$. Since $r \d Q \d x_{t-1} \d q \d z \d r$
is an antihole, it follows that $Q$ is odd. All its internal vertices have neighbours in
$A_{t-3}$, and one end $x_{t-1}$ does not, and $z$ is complete to its interior and
anticomplete to $A_{t-3}$. By \ref{greentouch} applied in $\overline{G}$, it follows that $r$
has neighbours in $A_{t-3}$.  Hence $r$ is $X_{t-3}$-complete, and nonadjacent to $x_{t-2}$.
Since $z \d x_{t-1} \d r \d q \d x_{t-2} \d z$ is not an odd hole it follows that
$x_{t-2}$ is adjacent to $x_{t-1}$. But then
\[x_0\l x_{t-1}\]
is a $\{q\}$-square of height $t-1$, and yet $z$ is not $\{q\}$-complete, a contradiction. This proves (1).
\\
\\
(2) {\em $q$ has neighbours in $A_{t-3}$.}
\\
\\
For suppose not. Let $S$ be an antipath between $x_t$ and $x_{t-1}$ with $V(S) \subseteq X_t$, that is,
with interior in $X_{t-2}$. Then $x_t\d S \d x_{t-1} \d q$ is an antipath with length $\ge 3$; by (1),
all its internal vertices have neighbours in $A_{t-3}$, and its ends do not, and $z$ is complete
to its interior and anticomplete to $A_{t-3}$; so by \ref{greentouch} applied in $\overline{G}$ it
follows that $S$ has odd length. But then $x_t\d S \d x_{t-1} \d q \d z$ has odd length $\ge 5$, and its
internal vertices have neighbours in $A_{t-2}$ and its ends do not, contrary to \ref{RRR} applied
in $\overline{G}$. This proves (2).

\bigskip

If $x_{t-1}$ is $X_{t-3}$-complete, then
\[x_0\l x_{t-1}\]
is a $\{q\}$-diamond of height $t-1$, and yet $z$ is not $\{q\}$-complete, a contradiction.
So $x_{t-1}$ is not $X_{t-3}$-complete. It follows from (2) that if $x_t$ is $X_{t-3}$-complete then
\[x_0\l x_{t-3},x_{t-1},x_{t-2},x_t\]
is a $Y$-square-diamond of height $t$, while if $x_t$ is not $X_{t-3}$-complete then
\[x_0\l x_{t-3},x_{t-1},x_t\] is a $Y$-square of height $t-1$, in either case a contradiction. This
proves \ref{betterwheel2}.\bbox

\begin{thm}\label{betterwheel3}
Let $ G \in \mathcal{F}_7$, let $(z,A_0)$ be a frame, and let $Y \subseteq V(G)\setminus (A_0\cup \{z\})$
be nonempty and anticonnected. Let $x_0\l x_{t+1}$ be a $Y$-square-diamond in $G$ of height $t+1 \ge 5$.
Then there is a nonempty
anticonnected set $Y'$ with $Y' \subseteq V(G)\setminus (A_0 \cup \{z\})$ such that either
$Y \subseteq Y'$ or $z$ is not $Y'$-complete, and such that either:
\begin{itemize}
\item there is a $Y'$-diamond in $G$ of height $t-1$, or
\item there is a $Y'$-square in $G$ of height $t-1$, or
\item there is a $Y'$-square-diamond in $G$ of height $t$.
\end{itemize}
\end{thm}
\Proof
Suppose that no such $Y'$ exists.
Let $x_0\l x_{t+1}$ be a $Y$-square-diamond in $G$, and define $X_i,A_i$ as usual. So
$x_{t+1}$ is $X_{t-1}$-complete, $x_t$ is not $X_{t-2}$-complete,
$x_{t+1}$ has no neighbour in $A_{t-2}$, $x_t$ has a neighbour in $A_{t-2}$,
there is a vertex $q$ in $A_{t-1}$ adjacent to both $x_{t+1},x_t$ with a neighbour in $A_{t-2}$,
and $Y$ is complete to $X_t$ and not to $x_{t+1}$. From the maximality of $A_{t-2}$ it follows
that $q$ is $X_{t-2}$-complete, and therefore nonadjacent to $x_{t-1}$.

Choose a path $v_1\c v_s$ with $s$ minimum such that $v_1\l v_s \in A_{t-2}$,
and $v_1$ is adjacent to $q$, and $v_s \in A_{t-3}$. (If $q$ has a neighbour in $A_{t-3}$ then $s = 1$.)
Let $R$ be a path between $q$ and $x_{t-1}$ with interior in $A_{t-2}$, and if possible
with interior in $A_{t-3}\cup \{v_1\l v_s\}$.
Then $R$ has length $\ge 2$, and from the hole $q\d R \d x_{t-1} \d x_{t+1} \d q$
it follows that $R$ has even length. So the path $q\d R \d x_{t-1} \d z$ is odd, and its ends are
$X_{t-2}$-complete, and its internal vertices are not, so by \ref{RRR} it has length 3, that is,
$R$ has length 2. Let its middle vertex be $r$.
\\
\\
(1) {\em $x_{t-1}$ has neighbours in $A_{t-3}$.}
\\
\\
For suppose not. It follows that $r\in A_{t-2} \setminus A_{t-3}$.  Let $Q$ be an antipath
between $r$ and $x_{t-1}$ with interior in $X_{t-2}$. Since $r \d Q \d x_{t-1} \d q \d z \d r$
is an antihole, it follows that $Q$ is odd. All its internal vertices have neighbours in
$A_{t-3}$, and one end $x_{t-1}$ does not, and $z$ is complete to its interior and
anticomplete to $A_{t-3}$. By \ref{greentouch} applied in $\overline{G}$, it follows that $r$
has neighbours in $A_{t-3}$.  Hence $r$ is $X_{t-3}$-complete, and nonadjacent to $x_{t-2}$.
Since $z \d x_{t-1} \d r \d q \d x_{t-2} \d z$ is not an odd hole it follows that
$x_{t-2}$ is adjacent to $x_{t-1}$. But then
\[x_0\l x_{t-1}\]
is a $\{q\}$-square of height $t-1$, and yet $z$ is not $\{q\}$-complete, a contradiction. This proves (1).

\bigskip

From (1) it follows that it is possible to choose $R$ with interior in $A_{t-3}\cup \{v_1\l v_s\}$,
and therefore we have done so.
\\
\\
(2) {\em $q$ has neighbours in $A_{t-3}$, and therefore $r \in A_{t-3}$.}
\\
\\
For suppose it does not.  Let $Q$ be an antipath between $x_{t-1}$ and $r$ with interior in $X_{t-2}$.
From the antihole $x_{t-1}\d Q \d r \d z \d q \d x_{t-1}$ it follows that $Q$ is odd. Hence
the antipath $q\d x_{t-1} \d Q \d r \d x_{t+1}$ is odd with length $\ge 5$; and its internal vertices
have neighbours in $A_{t-3} \cup \{v_2\l v_s\}$, and its ends do not, contrary to \ref{RRR}
applied in $\overline{G}$. This proves (2).
\\
\\
(3) {\em $x_{t-1}$ is not $X_{t-3}$-complete.}
\\
\\
For if it is, then
\[x_0\l x_{t-1}\]
is a $\{q\}$-diamond of height $t-1$, and yet $z$ is not $\{q\}$-complete, a contradiction.
This proves (3).
\\
\\
(4) {\em $x_t$ has no neighbour in $A_{t-3}$.}
\\
\\
For suppose $x_t$ has a neighbour in $A_{t-3}$.  If $x_t$ is $X_{t-3}$-complete then since it is
not $X_{t-2}$-complete, it is nonadjacent to $x_{t-2}$, and therefore
\[x_0\l x_{t-2},x_t\]
is a $Y \cup \{x_{t+1}\}$-diamond of height $t-1$; while if $x_t$ is not $X_{t-3}$-complete then
\[x_0\l x_{t-3},x_{t-1},x_t,x_{t+1}\]
is a $Y$-square-diamond of height $t$, in either case a contradiction. This proves (4).

\bigskip

In particular, $x_t$ is not adjacent to $r$. Since $z \d x_t \d q \d r \d x_{t-1} \d z$ is not an
odd hole it follows that $x_t$ is adjacent to $x_{t-1}$. If $x_t$ is $X_{t-3}$-complete, then
\[x_0\l x_{t-3},x_{t-1},x_{t-2},x_t\]
is a $Y \cup \{x_{t+1}\}$-square-diamond of height $t$; while if $x_t$ is not $X_{t-3}$-complete, then
\[x_0\l x_{t-3},x_{t-1},x_t\]
is a $Y \cup \{x_{t+1}\}$-square of height $t-1$, in either case a contradiction. This proves
\ref{betterwheel3}.\bbox

\section{Finding a wheel system}

In this section we apply the results of the previous two sections to prove a powerful result about wheel
systems that will be the engine behind almost all the remainder of the paper. First we need a lemma.

\begin{thm}\label{pseudohole}
Let $ G \in \mathcal{F}_7$, and let $X,Y$ be disjoint nonempty anticonnected subsets of $V(G)$, complete
to each other. Let $p_1\l p_n$ be a path in $G\setminus (X \cup Y)$ of length $\ge 4$, such that
$p_1,p_n$ are $X$-complete and $p_2\l p_{n-1}$ are not. Suppose that either:
\begin{enumerate}
\item $p_1,p_2,p_3$ are $Y$-complete, or
\item there exists $i$ with $1 \le i \le n-3$ such that $p_i,p_{i+1},p_{i+2},p_{i+3}$ are all $Y$-complete, or
\item there exists  $i$ with $1 \le i \le n-3$ such that $p_{i+1},p_{i+2}$ are $Y$-complete and
$p_i, p_{i+3}$ are not.
\end{enumerate}
Then there is a wheel in $G$ with hub $Y$.
\end{thm}
\Proof In the second and third case let $i$ be as given, and in the first case let $i = 1$.
Let $Q$ be an antipath joining $p_{i+1},p_{i+2}$ with interior in $X$.  Since $1 <i+1, i+2 < n$, and $n\ge 5$,
and $p_1,p_n$ are both complete
to the interior of $Q$, it follows from \ref{evenantipath3} that $Q$ has length 2, that is, there
exists $x \in X$ nonadjacent to both $p_{i+1},p_{i+2}$. Choose $h$ with $1 \le h \le i$ maximum so that
$x$ is adjacent to $p_h$, and choose $j$ with $i+3 \le j \le n$ minimum so that $x$ is adjacent to
$p_j$. Then $x\d p_h \c p_j \d x$ is a hole of length $\ge 6$, say $C$, and $x,p_i,p_{i+1},p_{i+2},p_{i+3}$
are all vertices of it, and $x, p_{i+1},p_{i+2}$ are $Y$-complete. In the first case $xp_1,p_1p_2,p_2p_3$
are all $Y$-complete edges of the hole, so $(C,Y)$ is a wheel. In the second case, the three edges of
$p_i\d p_{i+1} \d p_{i+2}\d p_{i+3}$ are all $Y$-complete edges of $C$, so again $(C,Y)$ is a wheel.
In the third case, \ref{evengap} implies that $(C,Y)$ is a wheel (and in this case it is in fact an odd wheel,
a contradiction).  This proves \ref{pseudohole}.\bbox

Let $Y$ be a nonempty anticonnected subset of $V(G)$, let $(z,A_0)$ be a frame with $A_0 \cup \{z\}$
disjoint from $Y$, and let $x_0\l x_{t+1}$ be a wheel system with respect to this frame.
We say $Y$ is a {\em hub}
for the wheel system if $t \ge 1$, $z,x_0\l x_t$ are all $Y$-complete and $x_{t+1}$ is not.

The main result of this section is the following. (Now we need to use that there are no psuedowheels,
so we are back in $\mathcal{F}_8$).

\begin{thm}\label{drivingwheel}
Let $G \in \mathcal{F}_8$, let $Y \subseteq V(G)$ be nonempty and anticonnected, let $(z,A_0)$ be a frame
with $Y \cap (A_0 \cup \{z\}) = \emptyset$, and let $x_0\l x_{t+1}$ be a wheel system with hub $Y$,
and with $t \ge 2$.  Define $X_i, A_i$ as usual.  Then either
\begin{itemize}
\item $x_{t+1}$ has a neighbour in $A_{t-1}$, or
\item some member of $Y$ is nonadjacent to $x_{t+1}$ and has no neighbour in $A_t$, or
\item there are $\ge 2$ members of $Y$ that are nonadjacent to $x_{t+1}$ and have no neighbour
in $A_{t-1}$, or
\item there is a wheel in $G$ with hub $Y$.
\end{itemize}
\end{thm}
\Proof We may assume that $x_{t+1}$ has no neighbour in $A_{t-1}$.
\\
\\
(1) {\em There do not exist $x_i,x_j \in X_t$ joined by an odd path $x_i\d x_{t+1}\d P \d x_j$
such that $x_i,x_j \in X_t$ and $P$ has interior in $A_t$.}
\\
\\
For assume such a path exists, and let $P$ have vertices $x_{t+1}\d p_1 \c p_n \d x_j$. There is an
even path $S$ between $x_i$ and $x_j$ with interior in $A_{t-1}$. Since
$x_i\d x_{t+1}\d P\d x_j\d S\d x_i$ is not an odd hole, and $x_{t+1}$ has no neighbours in $A_{t-1}$,
it follows that $\{p_1\l p_n\} \cup A_{t-1}$ is connected.
Since $p_1 \notin A_{t-1}$, there exists $k$ such that $p_k \notin A_{t-1}$ and $p_k$ has a neighbour
in $A_{t-1}$; and since $p_k$ is not adjacent to $z$, it follows from the maximality of $A_{t-1}$ that
$p_k$ is $X_{t-1}$-complete. Since at least one of $x_i,x_j$ is in $X_{t-1}$, it follows that
$k = n$ and $i = t$. But  $\{p_1\l p_n, x_j\} \cup A_{t-1}$ ($=F$ say) catches the triangle $\{z,x_{t+1},x_t\}$;
the only neighbour of $z$ in $F$ is $x_j$; the only neighbour of $x_{t+1}$ in $F$ is $p_1$; and $x_j,p_1$
are nonadjacent, and are both nonadjacent to $x_t$, contrary to \ref{tricatch}. This proves (1).

\bigskip

Since $x_{t+1}$ has a neighbour in $A_t$ and none in $A_{t-1}$, there is a path from $x_t$ to $A_{t-1}$
with interior in $A_t \setminus A_{t-1}$. Hence there is a path
$x_{t+1}\d p_1 \c p_m$ such that $p_1\l p_m \in A_t\setminus A_{t-1}$ and $p_m$ has a neighbour in $A_{t-1}$,
and $p_1\l p_{m-1}$ have no neighbours in $A_{t-1}$.  (Hence $m \ge 1$, and $p_m$ is $X_{t-1}$-complete.)
Choose such a path so that if possible, every member of $Y$ has a neighbour
in $A_{t-1} \cup \{x_{t+1},p_1\l p_m\}$.
\\
\\
(2) {\em We may assume that one of $x_0\l x_t$ is nonadjacent to both $x_{t+1},p_1$.}
\\
\\
For certainly there is an antipath $Q$ joining $x_{t+1},p_1$ with interior in $X_t$, since $x_{t+1},p_1$
are not $X_t$-complete. Suppose that $Q$ is odd. Every vertex of the
interior of $Q$ has neighbours in the connected set $A_{t-1}$, and $x_{t+1}$ does not, and $z$ is complete
to the interior of $Q$ and anticomplete to $A_{t-1}$; so by \ref{greentouch} applied in $\overline{G}$
it follows that $p_1$ has a neighbour in $A_{t-1}$. Hence $m = 1$, and $p_1$ is $X_{t-1}$-complete, and
therefore not adjacent to $x_t$. If $x_t$ is also nonadjacent to $x_{t+1}$ then the claim holds, and if
$x_t$ is adjacent to $x_{t+1}$, then
\[x_0\l x_{t+1}\]
is a $Y$-square, and the theorem holds by \ref{betterwheel}.

Now assume that $Q$ is even.
The antipath $z\d p_1 \d Q \d x_{t+1}$ is therefore odd and has length $\ge 3$; all its internal vertices
have neighbours in the connected set $A_{t-1} \cup \{p_2\l p_m\}$, and its ends do not. So it has length 3,
by \ref{RRR} applied in $\overline{G}$, and hence $Q$ has length 2. This proves (2).
\\
\\
(3) {\em We may assume that every vertex in $Y$ has a neighbour in  $A_{t-1} \cup \{x_{t+1},p_1\l p_m\}$.}
\\
\\
For suppose some $y \in Y$ has no such neighbour. If $y$ has no neighbour in $A_t$ then
the theorem holds, so we assume it does. Consequently
there is a connected subset $F$ of $A_t$ including $A_{t-1} \cup \{p_1\l p_m\}$ which contains a
neighbour of $y$, and we may choose $F$ minimal with this property. Since $y$ has no neighbour
in  $A_{t-1} \cup \{p_1\l p_m\}$, it follows from the minimality of $F$ that $y$ has a
unique neighbour in $F$, say $f$, and therefore $f \in A_t \setminus A_{t-1}$.
There is a path $R$ between $y$ and $x_{t+1}$ with interior in $F$,
and therefore $z\d x_{t+1} \d R \d y \d z$ is a hole ($C$ say), and so
$R$ has even length. Suppose it has length $\ge 4$. The only $X_t$-complete vertices in $C$ are $z,y$,
so by \ref{RRC}, $X_t$ contains a hat or leap. By (1) there is no leap, so there exists
$x \in X_t$ with no neighbours in $C$ except $y, z$. But $F\cup \{x_{t+1}\}$ catches the triangle
$\{x,y,z\}$; the only neighbour of $z$ in $F\cup \{x_{t+1}\}$ is $x_{t+1}$; the only neighbour of $y$
in $F\cup \{x_{t+1}\}$ is $f$; and $x_{t+1},f$ are nonadjacent, and both nonadjacent to $x$, contrary to
\ref{tricatch}. So $R$ has length 2, and therefore $x_{t+1}$ is adjacent to $f$.

Since $y$ has no neighbour in $\{x_{t+1}\} \cup A_{t-1}$, we may assume that all other members of $Y$
have neighbours in $\{x_{t+1}\}\cup A_{t-1}$, for otherwise the theorem holds.  We recall that
initially we chose the path $x_{t+1} \d p_1 \c p_m$ so that $p_m$ has a neighbour in $A_{t-1}$,
and $p_1\l p_{m-1}$ do not, and if possible every member of $Y$ has a neighbour in
$A_{t-1}\cup\{x_{t+1},p_1\l p_m\}$.  Since $f$ is adjacent to $y,x_{t+1}$, it follows
that $f$ has no neighbours in $A_{t-1}$, and $f$ is nonadjacent to
$p_2\l p_m$, since otherwise there would be a better choice of path using $f$.
Let $Q$ be an antipath between $f,x_{t+1}$ with interior in $X_t$.
Every internal vertex of $Q$ has a neighbour in $A_{t-1}$, and its ends do not, and $z$
is complete to the interior of $Q$ and anticomplete to $A_{t-1}$; so by \ref{greentouch} applied in
$\overline{G}$, it follows that $Q$ is even.  So the
antipath $y\d x_{t+1} \d Q \d f$ is odd, and all its internal vertices have neighbours in
$A_{t-1}\cup \{p_1\l p_m\}$, and $y$ does not; and $z$ is complete to the interior of the antipath and
anticomplete to $\{p_1\l p_i\}$. By \ref{greentouch} applied in $\overline{G}$, it follows that
$f$ has a neighbour in $A_{t-1}\cup \{p_1\l p_m\}$, and therefore $f$ is adjacent to $p_1$.
By (2) there exists $x \in X_t$ nonadjacent to $x_{t+1},p_1$. Consequently,
$\{z,y,x,p_2\l p_m\} \cup A_{t-1}$
(= $F'$ say) catches the triangle $\{x_{t+1},f,p_1\}$. The only neighbour of
$x_{t+1}$ in $F'$ is $z$; the only neighbour of $f$ in $F'$ is $y$; and $z,y$ are both
nonadjacent to $p_1$. By \ref{tricatch}, $F'$ contains a reflection of the triangle, and hence there is a
vertex in $F'$ adjacent to $y,z,p_1$. But the only neighbours of $z$ in $F'$ are $x,y$, and $x$
is nonadjacent to $p_1$, a contradiction. This proves (3).

\bigskip
Since $p_m$ is $X_{t-1}$-complete it follows that $x_0\l x_{t-1}$ all have neighbours in $p_1\l p_m$.
Since $x_t,p_m$ have neighbours in $A_{t-1}$ and none of $x_{t+1},p_1\l p_{m-1}$ have neighbours in $A_{t-1}$,
we can extend the path $x_{t+1}\d p_1 \c p_m$ to a path $x_{t+1}\d p_1 \c p_m\d p_{m+1} \c p_{n}$ containing
neighbours of all members of $X_t$.
By (2), we can choose $i$ with $2 \le i \le n$ maximum so that some vertex of $X_t$ is nonadjacent
to all of $x_{t+1},p_1\l p_{i-1}$; and choose $s$ with $0 \le s \le t$ such that $x_s$ is nonadjacent
to all of $x_{t+1},p_1\l p_{i-1}$. Since every vertex in $X_t$ has a neighbour in $\{x_{t+1},p_1\l
p_{n}\}$, it follows from the maximality of $i$ that every vertex in
$X_t$ is adjacent to one of $x_{t+1},p_1\l p_i$, and particular, $x_s$ is adjacent to $p_i$. Note that
if $i > m$ then $s = t$, since $p_n$ is $X_{t-1}$-complete.
\\
\\
(4) {\em $i$ is odd, and $p_i$ is $Y$-complete.}
\\
\\
For $z\d x_{t+1} \d p_1 \c p_i \d x_s \d z$ is a hole $C$ say, and so $i$ is odd. Suppose $p_i$ is not
$Y$-complete. Now $C$ has length $\ge 6$, and $z,x_s$ are
$Y$-complete, and $x_{t+1},p_i$ are not. Since $(C,Y)$ is not an odd wheel, \ref{RRC} implies that
$Y$ contains a leap or hat
for $C$. Suppose it contains a leap; then there are nonadjacent
$y_1,y_2 \in Y$ so that $y_1 \d x_{t+1} \d p_1 \c p_i \d y_2$ is a path. Since this path is odd and has
length $\ge 5$, and its ends are $X_t$-complete and its internal vertices are not, this contradicts
\ref{RRR}. So $Y$ contains a hat, that is, there exists $y \in Y$ nonadjacent to $x_{t+1},p_1\l p_i$.
By (3), $y$ has a neighbour in $A_{t-1}\cup \{p_j: i+1 \le j \le m\}$.

Suppose first that $i \le m$, and let $p_i\d r_1 \c r_k\d y$
be a path from $p_i$ to $y$ with interior in $A_{t-1} \cup \{p_{i+1}\l p_m\}$.
Then $z\d x_{t+1} \d p_1 \c p_i \d r_1 \c r_k \d y \d z$ is a hole of length $\ge 6$, and the
only $X_t$-complete vertices in this hole are $z,y$. Since this hole is not the rim of an odd wheel,
\ref{RRC} implies that $X$ contains a hat or leap, and so
some $x \in X_t$ has no neighbour in $\{x_{t+1},p_1\l p_i,r_1\l r_{k-1}\}$, contrary to the choice of $x_s$.

Now suppose that $i > m$, and so $s = t$. Let $p_m\d r_1 \c r_k\d y$ be a path from $p_m$ to $y$ with interior
in $A_{t-1}$. Again, $z\d x_{t+1} \d p_1 \c p_m \d r_1 \c r_k \d y \d z$ is a hole of length $\ge 6$,
and its only $X_t$-complete vertices are $z,y$. By \ref{RRC} $X_t$ contains a hat or leap. By (1)
it contains no leap, so there exists
$x \in X_t$ nonadjacent to all $x_{t+1},p_1\l p_m, r_1\l r_k$. Since $p_m$ is $X_{t-1}$-complete,
it follows that $x = x_t$. Now $\{x_{t+1},p_1\l p_i,r_1\l r_k\}$ ($=F$ say) is connected, and catches
the triangle $\{y,z,x_t\}$; the only neighbour of $z$ in $F$ is $x_{t+1}$; the only neighbour of $y$ in $F$
is $r_k$ (because $y$ is nonadjacent to $x_{t+1},p_1\l p_i$); and the only neighbour of $x_t$ in $F$ is
$p_i$ (because $x_t$ is a hat). Since $x_{t+1}$ is not adjacent to $p_i$, this contradicts \ref{tricatch}.
This proves (4).
\\
\\
(5) {\em Let $R$ be a path from $x_t$ to some vertex $r$, such that $r$ is the unique $X_{t-1}$-complete
vertex in $R$, and $V(R\setminus x_t) \subseteq A_{t-1} \cup \{p_1\l p_m\}$. Then $R$ is odd, and has
length $\ge 3$. In particular, $x_t$ is nonadjacent to $p_m,p_{m-1}$.}
\\
\\
For assume that $R$ is even. Then the path
$z \d x_t \d R \d r$ is odd, and its ends are $X_{t-1}$-complete, and its internal vertices are not,
so by \ref{RRR}, it has length 3, that is, $R$ has length 2. Let $q$ be the middle vertex of $R$.
By \ref{RRR} there is an odd antipath $Q$ joining $q,x_t$ with interior in $X_{t-1}$. Now $p_m$ is
$X_{t-1}$-complete and nonadjacent to $x_t$, and since $Q$ cannot be completed to an antihole
via $x_t \d p_m \d q$, it follows that $p_m$ is adjacent to $q$. Suppose first that
$q \in \{p_1\l p_m\}$; then it follows
that $q = p_{m-1}$. Hence $q \d Q \d x_t \d p_m$ is an even antipath of length $\ge 4$; $q$ is
its only vertex that is anticomplete to $A_{t-1}$, and $p_m$ is its only vertex that is anticomplete
to $\{z,x_{t+1},p_1\l p_{m-2}\}$. Since the sets $A_{t-1}$, $\{z,x_{t+1},p_1\l p_{m-2}\}$ are each connected
and anticomplete to each other, this contradicts \ref{doubleRR2} applied in $\overline{G}$. So
$q \in A_{t-1}$, and in particular $x_t$ is nonadjacent to $p_m, p_{m-1}$.
Let $R'$ be a path between $x_t,p_m$ with interior in $\{z,x_{t+1},p_1\l p_m\}$;
then $x_t \d R \d p_m \d R' \d x_t$ is a hole of length $\ge 6$ sharing the vertices $x_t,q,p_m$
with the antihole $q \d Q \d x_t \d p_m \d z \d q$, contrary to \ref{hole&antihole}. So $R$ is odd. Since
$r$ is not $X_t$-complete, it follows that $R$ has length $\ge 3$. The last assertion of
the claim is immediate. This proves (5).
\\
\\
(6) {\em We may assume that none of $x_{t+1},p_1\l p_{i-1}$ is $X_{t-1}$-complete, and in particular
$i \le m$.}
\\
\\
For suppose first that one of $p_1\l p_{i-1}$ is $X_{t-1}$-complete, and choose $h$ with
$1 \le h < i$ maximum so that $p_h$ is $X_{t-1}$-complete. Since $p_h$ is not adjacent to $p_s$
it follows that $s = t$, and therefore $p_i$ is not $X_{t-1}$-complete.
By (5), $i-h$ is even, and so the
path $p_h \c p_i \d x_t \d z$ is even and has length $\ge 4$. Since its only $X_{t-1}$-complete
vertices are its ends, and since $z,x_t,p_i$ are $Y$-complete by (4), it follows from \ref{pseudohole}
that there is a wheel with hub $Y$, and the theorem holds. So we may assume that none of
$p_1\l p_{i-1}$ is $X_{t-1}$-complete, and in particular $i \le m$, since $p_m$ is $X_{t-1}$-complete.
Now assume that $x_{t+1}$ is $X_{t-1}$-complete. Since $x_{t+1}$ is nonadjacent to $x_s$ it follows
that $s = t$. Let $R$ be a path between $x_t,p_m$ with interior
in $A_{t-1}$. By (5), $R$ is odd, and so the path $x_{t+1}\d p_1 \c p_i \d x_t \d R \d p_m$ is odd, of
length $\ge 5$, its ends are $X_{t-1}$-complete, and its internal vertices are not, contrary to
\ref{RRR}. This proves (6).

\bigskip

Choose $k$ with $ i \le k \le m$ minimum such that $p_k$ is $X_{t-1}$-complete.
\\
\\
(7) {\em None of $x_{t+1},p_1\l p_{k-1}$ is $X_{t-1}$-complete, and $k$ is odd.}
\\
\\
The first assertion follows from (6). Hence the path $z\d x_{t+1} \d p_1 \c p_k$ has length $\ge 4$,
and its ends are $X_{t-1}$-complete, and its internal vertices are not; so by \ref{RRR}, it has
even length. This proves (7).
\\
\\
(8) {\em $x_t$ is adjacent to one of $p_1\l p_k$.}
\\
\\
For suppose $x_t$ is nonadjacent to all of $p_1\l p_k$. From the definition of $i$ it follows that
$x_t$ is adjacent to $x_{t+1}$. Let $S$ be a path between $x_t,p_k$ with interior in
$A_{t-1} \cup \{p_{k+1}\l p_m\}$, and let $C$ be the hole $x_t \d x_{t+1} \d p_1 \c p_k \d S \d x_t$.
Since $C$ is even and $k$ is odd, it follows that $S$ is even, and so by (5), some internal vertex of $S$ is
$X_{t-1}$-complete. The path $z \d x_t \d S \d p_k$ is odd, and its
ends are $X_{t-1}$-complete, so by \ref{evengap} it contains an odd number of $X_{t-1}$-complete
edges. Since $x_t$ is not $X_{t-1}$-complete, all these $X_{t-1}$-complete edges belong to $S$ and
hence to $C$, and there are no further $X_{t-1}$-complete edges in $C$. Thus an odd number of edges
of $C$ are $X_{t-1}$-complete, and so by \ref{evengap} there is exactly one, and exactly two
$X_{t-1}$-complete vertices. Since $p_k$ is $X_{t-1}$-complete, the second such vertex is the neighbour
of $p_k$ in $S$. This therefore does not belong to $A_{t-1}$, and so $k<m$, and $p_{k+1}$ is the second
$X_{t-1}$-complete vertex of $C$. By \ref{RRC} applied to $C$, $X_{t-1}$ contains a leap or hat, and in
either case some $x \in X_{t-1}$ is nonadjacent to all of $x_t,x_{t+1},p_1$, and adjacent to
$p_k$. Hence $(V(C) \setminus \{x_t,x_{t+1}\})\cup \{x\}$ ($=F$ say) catches the triangle
$\{z,x_t,x_{t+1}\}$; the only neighbour
of $z$ in $F$ is $x$; the only neighbour of $x_{t+1}$ in $F$ is $p_1$; and $x,p_1$ are nonadjacent,
and are both nonadjacent to $x_t$, contrary to \ref{tricatch}. This proves (8).
\\
\\
(9) {\em $p_k$ is $Y$-complete.}
\\
\\
For suppose not. Then $i < k$, by (4). But then $z,X_{t-1}$ are $Y$-complete and $x_{t+1},p_k$ are not, and
some vertex of the path $x_{t+1} \d p_1 \c p_k$ is $Y$-complete (namely $p_i$); and so
$(X_{t-1},Y,z\d x_{t+1} \d p_1 \c p_k)$ is a pseudowheel, contrary to $G \in \mathcal{F}_8$. This proves (9).

\bigskip

By (8), we may choose $j$ with $1 \le j \le k$ maximum such that $x_t$ is adjacent to $p_j$.
By (5), $k-j$ is even and $\ge 2$. Suppose that $p_j$ is $Y$-complete. The path $z\d x_t \d p_j \c p_k$ has
even length $\ge 4$, and its only $X_{t-1}$-complete vertices are its ends, and $z,x_t,p_j,X_{t-1}$ are all
$Y$-complete, so by \ref{pseudohole}, there is a wheel with hub $Y$ and the theorem holds. So we may
assume that $p_j$ is not $Y$-complete. Now the path $x_t \d p_j \c p_k$ has odd length $\ge 3$, and
both its ends are $Y$-complete, and the $Y$-complete vertex $z$ has no neighbour in its interior, so by
\ref{RR} and \ref{evengap}, an odd number of its edges are $Y$-complete. Since $p_j$ is not $Y$-complete,
an odd number of edges of $p_j \c p_k$ are $Y$-complete. The path $z\d x_{t+1} \d p_1 \c p_k$ ($= P$ say) is even,
by (7), and since its ends are $Y$-complete, it follows that an even number of its edges are $Y$-complete,
by \ref{evengap}. We deduce that an odd number of edges of $z\d x_{t+1}\d p_1\c p_j$ are $Y$-complete.
There is therefore a $Y$-segment $P'$ of this path that has odd length. Since $p_j$ is not $Y$-complete,
it follows that $P'$ is also a $Y$-segment of $P$.
If $P'$ has length $>1$ then \ref{pseudohole} applied to $P$ implies
that there is a wheel with hub $Y$, and the theorem holds. So we may assume that $P'$ has length 1. But
both vertices of $P'$ are internal vertices of $P$, since $x_{t+1},p_j$ are
not $Y$-complete, and again \ref{pseudohole} applied to $P$
implies there is a wheel with hub $Y$. This proves \ref{drivingwheel}.\bbox

Before we apply this result, it is helpful to transform it a little.

\begin{thm}\label{drivingwheel2}
Let $G \in \mathcal{F}_8$, let $Y \subseteq V(G)$ be nonempty and anticonnected, let $(z,A_0)$ be a frame
with $Y \cap (A_0 \cup \{z\}) = \emptyset$, and let $x_0\l x_{t+1}$ be a wheel system with hub $Y$,
where $t \ge 1$.  Define $A_i, X_i$ as usual, and assume that at most $y \in Y$
has no neighbour in $A_1$.  Suppose there is no wheel with hub $Y$. Then there exists $r$
with $1 \le r \le t$, and a member $y \in Y$, with the following properties:
\begin{itemize}
\item $y$ is nonadjacent to $x_{t+1}$ and has no neighbour in $A_r$
\item $x_{t+1}$ has a neighbour in $A_r$, and a non-neighbour in $X_r$.
\end{itemize}
\end{thm}
\Proof We proceed by induction on $t$. If $t = 1$ then \ref{wheelbase} implies that
there exists $y \in Y$ nonadjacent to $x_{t+1}$ and with no neighbour in $A_t$,
and the theorem holds. So we may assume $t \ge 2$.
If $x_{t+1}$ has no neighbour in $A_{t-1}$, then the result follows from \ref{drivingwheel}, since
at most one member of $Y$ has no neighbour in $A_{t-1}$. So we assume that  $x_{t+1}$ has a
neighbour in $A_{t-1}$.  If $x_{t+1}$ is $X_{t-1}$-complete then
\[x_0\l x_{t+1}\]
is a $Y$-diamond, and the theorem holds by \ref{betterwheel}.
If $x_{t+1}$ is not $X_{t-1}$-complete, then
\[x_0\l x_{t-1},x_{t+1}\]
is a wheel system with hub $Y$, and the result follows from the inductive hypothesis. This proves
\ref{drivingwheel2}.\bbox

We also need a further transformation, which uses the following lemma.

\begin{thm}\label{gethop}
Let $G \in \mathcal{F}_8$, not admitting a balanced skew partition, let $Y \subseteq V(G)$ be nonempty
and anticonnected, let $(z,A_0)$ be a frame
with $Y \cap (A_0 \cup \{z\}) = \emptyset$, and let $x_0\l x_s$ be a wheel system with $s \ge 1$, such
that $z,x_0\l x_s$ are $Y$-complete. Then there is a sequence $x_{s+1}\l x_{t+1}$ with $t \ge s$
such that $x_0\l x_{t+1}$ is a wheel system with respect to the frame $(z,A_0)$, with hub $Y$.
\end{thm}
\Proof
For choose a sequence $x_{s+1}\l x_t$, all $Y$-complete and such that $x_0\l x_t$ is a wheel system
with respect to $(z,A_0)$, with $t$ maximum. So $t \ge 1$. Define $X_i$ and $A_i$ as usual.
From \ref{bipattach}, there is a path $P$ from $z$ to
$A_t$, disjoint from $X_t$ and containing no $X_t$-complete vertex
except $z$. Let $v$ be the neighbour of $z$ in this path.
From the maximality of $A_t$, it follows that $P$ has length 2. So $v$ has a neighbour in $A_t$, and
therefore $x_0\l x_t,v$ is a wheel system. From the maximality of $s$ it follows
that $v$ is not $Y$-complete, and therefore $Y$ is a hub for this wheel system. This proves \ref{gethop}.
\bbox

\begin{thm}\label{drivingwheel3}
Let $G \in \mathcal{F}_8$, not admitting a balanced skew partition, let $Y \subseteq V(G)$ be nonempty and
anticonnected, let $(z,A_0)$ be a frame
with $Y \cap (A_0 \cup \{z\}) = \emptyset$, and let $x_0\l x_s$ be a wheel system with $s \ge 1$, such
that $z,x_0\l x_s$ are $Y$-complete. Define $A_i, X_i$ as usual, and assume that every member of $Y$
has a neighbour in $A_s$, and at most one member
of $Y$ has no neighbour in $A_1$. Suppose there is no wheel with hub $Y$. Then there exists $r$
with $1 \le r < s$, and a member $y \in Y$, and a vertex $v$ with the following properties:
\begin{itemize}
\item $y$ is nonadjacent to $v$ and has no neighbour in $A_r$
\item $v$ is adjacent to $z$, and has a neighbour in $A_r$, and a non-neighbour in $X_r$.
\end{itemize}
\end{thm}
\Proof  By \ref{gethop},
there is a sequence $x_{s+1}\l x_{t+1}$ with $t \ge s$
such that $x_0\l x_{t+1}$ is a wheel system with respect to the frame $(z,A_0)$, with hub $Y$.
By \ref{drivingwheel2}, there exists $r$
with $1 \le r \le t$, and a member $y \in Y$, such that
$y$ is nonadjacent to $x_{t+1}$ and has no neighbour in $A_r$, and
$x_{t+1}$ has a neighbour in $A_r$, and a non-neighbour in $X_r$. Since every member of $Y$ has a neighbour
in $A_s$, it follows that $r < s$, and the result follows. This proves \ref{drivingwheel3}. \bbox

\section{Wheels with tails}

We continue with the proof that recalcitrant graphs do not contain wheels. The idea is, suppose $G$
contains a wheel, and look at one, $(C,Y)$ say, with the largest hub. Then we find another vertex to
add to the hub, which almost makes a larger hub, but not quite, because where there should be edges joining
it to $C$ there is instead a path with certain properties (a ``tail''). Then we use the machinery of the
previous section to show that there is a wheel with a larger hub anyway. So now we need some lemmas
about tails.

If $(C,Y)$ is a wheel in $G$, and there is no wheel $(C', Y')$ with $Y \subset Y'$, we say $(C,Y)$
is an {\em optimal} wheel.

Let $(C,Y)$ be a wheel in $G$. A {\em kite} for  $(C,Y)$ is a vertex $y \in V(G)\setminus (Y \cup V(C))$,
not $Y$-complete,
that has at least four neighbours in $C$, three of which are consecutive
and $Y$-complete.

\begin{thm}\label{nokite}
Let $G \in \mathcal{F}_8$, not admitting a balanced skew partition, and let $(C,Y)$ be an optimal
wheel in $G$. Then there is no kite.
\end{thm}
\Proof
Assume $y$ is a kite. Let $x_0\d z \d x_1$ be a subpath of $C$, all $Y$-complete and adjacent to $y$.
Let $A_0 = V(C) \setminus \{z,x_0,x_1\}$, so $x_0,x_1$ is a wheel system with respect to $(z,A_0)$, and
$x_0,x_1$ are $Y \cup \{y\}$-complete. By \ref{gethop},
there is a sequence $x_2\l x_{t+1}$ with $t \ge 1$
such that $x_0\l x_{t+1}$ is a wheel system with respect to the frame $(z,A_0)$, with hub $Y\cup \{y\}$.
Since every vertex of $Y\cup \{y\}$ has a neighbour in $A_0$, and there is no wheel with hub
$Y\cup \{y\}$ from the optimality of $(C,Y)$, this contradicts \ref{drivingwheel2}.
This proves \ref{nokite}. \bbox

Let $(C,Y)$ be a wheel in $G$,
let $z\in V(C)$, and let $x_0,x_1$ be the neighbours of $z$ in $C$. A path $T$ of
$G \setminus \{x_0,x_1\}$ from $z$ to $V(C)\setminus \{z,x_0,x_1\}$ is called a {\em tail} for $z$ (with
respect to the wheel $(C,Y)$) if
\begin{itemize}
\item $x_0,z,x_1$ are all $Y$-complete,
\item there is a $Y$-complete edge in $C\setminus \{x_0,z,x_1\}$
\item the neighbour of $z$ in $T$ is adjacent to $x_0,x_1$,
\item no vertex of $T$ is in $Y$,
\item no internal vertex of $T$ is $Y$-complete, and
\item no vertex of $G$ is a kite for $(C,Y)$.
\end{itemize}

\begin{thm}\label{docktail}
Let $G \in \mathcal{F}_8$, let $(C,Y)$ be an optimal wheel, let $z\in V(C)$, and let $x_0,x_1$ be the neighbours of
$z$ in $C$. Let $T$ be a tail for $z$, and let $y$ be the neighbour of $z$ in $T$. Let
$A_0 = V(C) \setminus\{z,x_0,x_1\}$, and let $x_0\l x_{t+1}$ be a wheel system with respect to the
frame $(z,A_0)$, with hub $Y\cup \{y\}$. Define $A_1\l A_{t+1}$ as usual.
Then either $y$ is adjacent to $x_{t+1}$, or $y$ has a neighbour in $A_t$.
\end{thm}
\Proof We assume for a contradiction that $y$ has no neighbour in $A_t$, and $y$ is not adjacent
to $x_{t+1}$.  Let $y \d u_1 \c u_n$ be a minimal
subpath of $T\setminus z$ so that $u_n$ has a neighbour in $A_t$; so  $n>0$. From
the maximality of $A_t$ it follows that $u_n$ is $X_t$-complete and therefore $X_1$-complete since $t \ge 1$;
and since $T$ is a tail it follows that none of $u_1\l u_n$ are $Y$-complete. Let $P$ be a path
from $u_n$ to some $Y$-complete vertex $p$ say, so that no vertex of $P \setminus p$ is $Y$-complete.
\\
\\
(1) {\em $P$ is odd.}
\\
\\
For $P$ has length $\ge 1$ since no vertex of $T\setminus z$ is $Y$-complete; and
the only $X_t$-complete vertex of $P$ is $u_n$, and the only $Y$-complete
vertex of $P$ is $p$. Since $z$ is complete to $X_t$ and to $Y$, and anticomplete to $V(P)$, it follows
from \ref{doubleRR} that $P$ has odd length. This proves (1).

\bigskip
Since $y,u_1\l u_{n-1}$ have no neighbours in $A_t$ it follows that $z\d y\d u_1 \c u_n \d P \d p$ is a path,
$Q$ say.
\\
\\
(2) {\em We may assume that $Q$ has even length $\ge 4$, and so $n$ is even.}
\\
\\
For the ends of $Q$ are $Y$-complete, and since
none of $y,u_1\l u_n$ are $Y$-complete, it follows that no internal vertex of $Q$ is $Y$-complete.
Suppose that $Q$ has length 3. So $n = 1$, and there is an odd antipath joining
$y,u_1$ with interior in $Y$. Hence every $Y$-complete vertex in $G$ is adjacent to one of $y,u_1$.
In particular, since $y$ has no neighbour in $A_t$, it follows that $u_1$ is adjacent to all the $Y$-complete
vertices in $C$ except $z$ (for we already showed that it is $X_t$-complete and therefore adjacent to
$x_0,x_1$). But then $u_1$ is a kite for $(C,Y)$, a contradiction.
So we may assume that $Q$ does not have length 3. Hence by \ref{RRR},
$Q$ has even length.  From (1), it follows that $n$ is even. This proves (2).
\\
\\
(3) {\em $x_{t+1}$ is adjacent to one of $u_1\l u_{n-1}$.}
\\
\\
For suppose not. Choose a path $N$ from $x_{t+1}$ to $u_n$ with interior in $A_t$ (possibly of length 1).
Then $z\d y \d u_1 \c u_n \d N \d x_{t+1} \d z$ is a hole, and since $n$ is even
it follows that $N$ is even. Hence $z\d x_{t+1} \d N \d u_n$ is an odd path; its ends are $X_t$-complete,
its internal vertices are not, and the $X_t$-complete vertex $y$ has no neighbour in its interior,
contrary to \ref{greentouch}. This proves (3).
\\
\\
(4) {\em $x_{t+1}$ is not $Y$-complete.}
\\
\\
For suppose it is. Since $G \in \mathcal{F}_8$, the triple $(Y,X_{t+1},Q)$ is not a pseudowheel. Since
$y,p$ are not $X_{t+1}$-complete, it follows that no internal vertex of $S$ is
$X_{t+1}$-complete.  By \ref{doubleRRC} applied to $Q$,
$Y$ and $X_{t+1}$, it follows that there exists $x \in X_{t+1}$ with no neighbour
in $Q\setminus z$ except possibly $p$. But $x_{t+1}$ is adjacent to one of $u_1\l u_{n-1}$ by (3), and
all other members of $X_{t+1}$ are adjacent to $y$, a contradiction. This proves (4).

\bigskip

Since $x_{t+1}$ has a neighbour in $A_t$, there is a path $R$ from $x_{t+1}$ to some $Y$-complete vertex $r$
in $A_t$ with $V(R \setminus x_{t+1})\subseteq A_t$
so that no vertex of $R \setminus r$ is $Y$-complete.
\\
\\
(5) {\em $R$ has odd length.}
\\
\\
For certainly $R$ has length $\ge 1$; suppose it
has length 2, and let its middle vertex be $a$ say. There is an antipath joining $x_{t+1},a$ with interior
in $Y$, and it is odd since it can be complete to an antihole via $a\d z \d r \d x_{t+1}$. Now
$x_{t+1},a$ are not $X_t$-complete (since $a \in A_t$) and so there is an antipath joining $x_{t+1},a$
with interior in $X_t$, which is therefore also odd, since its union with the antipath with interior
in $Y$ is an antihole. But $y$ is $X_t$-complete and nonadjacent to both $x_{t+1}$ and $a$ (since it
has no neighbour in $A_t$), and so this antipath can be completed to an odd antihole via $a\d y \d x_{t+1}$,
a contradiction. This proves that $R$ does not have length 2. Hence the path $z\d x_{t+1} \d R \d r$ does
not have length 3; its ends are $Y$-complete and its internal vertices are not, and it has
length $>1$, so by \ref{RRR} it has even length, that is, $R$ has odd length. This proves (5).
\\
\\
(6) {\em If $x_{t+1}$ is adjacent to $u_1$ then $u_1$ is $X_t$-complete.}
\\
\\
For suppose not; then there is an antipath $L$ say joining $x_{t+1},u_1$ with interior
in $X_t$. So $z\d u_1 \d L \d x_{t+1} \d y$ is an antipath of length $\ge 4$; all its internal
vertices have neighbours in $A_t \cup \{u_2\l u_n\}$, and its ends do not. By \ref{RRR} applied in
$\overline{G}$, it has even length, and so $u_1 \d L \d x_{t+1} \d y$ is an odd antipath. Since
$u_1,y$ are not $Y$-complete, they are joined by an antipath with interior in $Y$, which is therefore
also odd; and yet it can be completed to an odd antihole via $y\d p \d u_1$ (we recall that $p \in A_t$ is
$Y$-complete). This proves (6).
\\
\\
(7) {\em None of $u_1\l u_{n-1}$ is $X_t$-complete.}
\\
\\
For suppose that one of $u_1\l u_{n-1}$ is $X_t$-complete, and let $S$ be a path from $x_{t+1}$ to some
$X_t$-complete vertex $s$ say, with $V(S \setminus x_{t+1}) \subseteq \{u_1\l u_{n-1}\}$, such that
$s$ is the only $X_t$-complete vertex in $S$. Certainly $S$ has length $\ge 1$. Suppose it has
even length. Then the path $z\d x_{t+1} \d S \d s$ is odd, and its ends are $X_t$-complete, and its
internal vertices are not; so by \ref{greentouch}, the $X_t$-complete vertex $y$ has a neighbour in
its interior, contrary to (6).  So $S$ has odd length. The path
$s\d S \d x_{t+1} \d R \d r$ therefore has even length; its
only $X_t$-complete vertex is $s$, and its only $Y$-complete vertex is $r$, so by \ref{doubleRR2}, the
path has length 2, that is, both $R,S$ have length 1. Moreover, either $x_{t+1},r$ are joined by an
odd antipath with interior in $X_t$, or $x_{t+1},s$ are joined by an odd antipath with interior in $Y$.
The first is impossible since the antipath could be completed to an odd antihole via $r \d y \d x_{t+1}$,
so the second
holds. In particular, every $Y$-complete vertex is adjacent to one of $x_{t+1},s$, and therefore all
such vertices in $A_t$ are adjacent to $x_{t+1}$. In particular, $x_{t+1}$ is adjacent to all the
$Y$-complete vertices in $C$ except possible $x_0,x_1$. Since there is a $Y$-complete edge in
$C \setminus\{x_0,z,x_1\}$ from the definition of a tail, it follows that
$x_{t+1}$ has two adjacent neighbours in $C$ of opposite wheel-parity, and at least one
other neighbour in $C$; but it is not a kite, and the wheel is optimal,
contrary to \ref{wheelvnbrs}. This proves (7).

\bigskip

By (3) we may choose $i$ with
$1 \le i \le n-1$ minimum such that  $x_{t+1}$ is adjacent to $u_i$. By (7), the only
$X_t$-complete vertices in the hole $z \d y\d u_1 \c u_i \d x_{t+1} \d z$ are $z$,$y$, and therefore
by (6) this hole has length $\ge 6$.  By \ref{RRC} $X_t$ contains a leap or a hat. If it contains a leap,
there are nonadjacent vertices in $X_t$, joined by an odd path of length $\ge 5$ with interior in
$\{u_1\l u_i,x_{t+1}\}$, and consequently with no internal vertex $Y$-complete. Since both its ends
are $Y$-complete, this contradicts \ref{RRR}. So there is a hat, that is, there exists $x \in X_t$
with no neighbours in $\{u_1\l u_i,x_{t+1}\}$. Then
$A_t \cup \{u_1\l u_n,x_{t+1}\}$ ($= F$ say) catches the triangle $\{z,y,x\}$; the only neighbour of $z$
in $F$ is $x_{t+1}$; the only neighbour of $y$ in $F$ is $u_1$; and both $x_{t+1},u_1$ are nonadjacent to $x$.
Moreover $x_{t+1}$ is nonadjacent to $u_1$, and so $F$ contains no reflection of the triangle. This
contradicts \ref{tricatch}, and therefore proves \ref{docktail}.\bbox

We combine the previous result with \ref{drivingwheel2} to prove the following.

\begin{thm}\label{notail}
Let $G \in \mathcal{F}_8$, not admitting a balanced skew partition, and let $(C,Y)$ be an optimal wheel
in $G$. Then no vertex of $C$ has a tail.
\end{thm}
\Proof
Suppose $z \in V(C)$ has a tail $T$; let $y$ be the neighbour of $z$ in $T$, and let $x_0,x_1$ be the
neighbours of $z$ in $C$. Let $A_0 = V(C) \setminus \{z,x_0,x_1\}$, so $x_0,x_1$ is a wheel system
with respect to $(z,A_0)$, and $x_0,x_1$ are $Y\cup \{y\}$-complete. By \ref{gethop} there exist
$x_2\l x_{t+1}$ with $t \ge 1$ such that $x_0 \l x_{t+1}$ is a wheel system with respect to $(z,A_0)$,
with hub $Y\cup \{y\}$.  Define $A_i, X_i$ as usual. From the construction, all members of $Y$ have
a neighbour in $A_0$. By \ref{drivingwheel2}, there exists $r$ with $1 \le r \le t$, such that
$y$ is nonadjacent to $x_{t+1}$ and has no neighbour in $A_r$, and
$x_{t+1}$ has a neighbour in $A_r$, and a non-neighbour in $X_r$. Now $x_0\l x_r,x_{t+1}$ is a wheel system
with hub $Y\cup \{y\}$, and $T$ is a tail for $z$, contrary to \ref{docktail}. This proves \ref{notail}.
\bbox

\section{The end of a wheel}

In this section we complete the proof that there is no wheel in a recalcitrant graph.

\begin{thm}\label{wheeltohalf}
Let $G \in \mathcal{F}_8$, not admitting a balanced skew partition, and let $(C,Y)$ be an optimal wheel in
$G$. Then there is a subpath $c_1\d c_2 \d c_3$ of $C$ such that $c_1,c_2,c_3$ are
all $Y$-complete, and a path $c_1\d p_1 \c p_k \d c_3$ such that none of $p_1\l p_k$ are in $V(C) \cup Y$,
none of them is $Y$-complete, and none of them has a neighbour in $V(C)\setminus \{c_1,c_2,c_3\}$.
\end{thm}
\Proof
There are
two nonadjacent $Y$-complete vertices in $C$ with opposite wheel-parity, say $a,b$, and by \ref{bipattach},
there is a path $P$ in $G$ joining them so that none of its interior vertices is in $Y$ or is $Y$-complete.
There may be internal vertices of $P$ that belong to $C$, but we may choose a subpath $P'$ of $P$, with ends
$a',b'$ say, so that $a',b'\in V(C)$ have opposite wheel-parity and $P'$ has minimum length. It follows
that no vertex of the interior of $P'$ is in $C$. Suppose $a',b'$ are adjacent; then since they are in $C$
and have opposite wheel-parity, they are both $Y$-complete, and therefore neither is in the interior of
$P$, and so
$a,b$ are adjacent, a contradiction. So $a',b'$ are nonadjacent. Let $F$ be the interior of $P'$;
then no vertex of $F$ is in $Y \cup V(C)$, no vertex of $F$ is $Y$-complete, and there are attachments
of $F$ in $C$ which are nonadjacent and have opposite wheel-parity. The result follows from \ref{nokite} and
\ref{wheelFnbrs} applied to $F$. This proves \ref{wheeltohalf}.

Now we can prove \ref{summary}.9, which we restate.

\begin{thm}\label{nowheel}
Let $G \in \mathcal{F}_8$, not admitting a balanced skew partition; then there is no wheel in $G$.
In particular, every recalcitrant graph belongs to $\mathcal{F}_9$.
\end{thm}
\Proof
Suppose there is a wheel in $G$, and choose an optimal wheel $(C,Y)$ so that $C$ contains as few
$Y$-complete edges as possible.
\\
\\
(1) {\em Exactly $4$ edges of $C$ are $Y$-complete.}
\\
\\
For by \ref{wheeltohalf} there is a subpath $c_1\d c_2 \d c_3$ of $C$ such that $c_1,c_2,c_3$ are
all $Y$-complete, and a path $c_1\d p_1 \c p_k \d c_3$ such that none of $p_1\l p_k$ are in $V(C) \cup Y$,
none of them is $Y$-complete, and none of them has a neighbour in $V(C)\setminus \{c_1,c_2,c_3\}$.
Let $C'$ be the hole formed by the union of the paths $C \setminus c_2$, $c_1\d p_1 \c p_k \d c_3$.
Then it has length $\ge 6$, and it contains fewer $Y$-complete edges than $C$. From the choice
of $(C,Y)$ it follows that $(C',Y)$ is not a wheel, and since $C$ has at least 4 $Y$-complete edges,
and $C'$ has only two fewer, it follows that exactly 4 edges of $C$ are $Y$-complete. This
proves (1).

\bigskip

Since $(C,Y)$ is not an odd wheel, there are vertices  $x_0,z,x_1,c_1,c_2,c_3$ of $C$, in order, and
all distinct except possibly $x_1 = c_1$ or $c_3 = x_0$, so that
the $Y$-complete edges in $C$ are $x_0z,zx_1, c_1c_2, c_2c_3$.
\\
\\
(2) {\em There is no tail for $z$ with respect to $(C,Y)$.}
\\
\\
For suppose there is a tail; then by \ref{drivingwheel3} there is a wheel with hub a
proper superset of $Y$, a contradiction. This proves (2).

\bigskip

Let $A_0 = V(C) \setminus \{z,x_0,x_1\}$.  Since $G$
does not admit a skew partition, there is a path $T$ of $G \setminus \{x_0,x_1\}$ from $z$ to $A_0$,
such that no vertex in its interior is in $Y$ or $Y$-complete.  Let $y$ be the neighbour of $z$ in $T$.
\\
\\
(3) {\em $y$ is not adjacent to both $x_0,x_1$.}
\\
\\
For assume it is. By \ref{nokite} there is no kite for $(C,Y)$, so with respect to the wheel $(C,Y)$,
$T$ is a tail for $z$ (because at least one of the $Y$-complete edges $c_1c_2,c_2c_3$
belongs to $C \setminus \{x_0,z,x_1\}$). This contradicts (2), and therefore proves (3).
\\
\\
(4) {\em $y$ has no neighbour in $A_0$.}
\\
\\
For suppose first that it has a neighbour in $A_0 \setminus c_2$, say $c$.
Then $c,z$ are nonadjacent and have opposite wheel-parity in the wheel $(C,Y)$, and not both neighbours
of $c$ in $C$ are $Y$-complete, by (1), and not both neighbours of $z$ in $C$ are adjacent to $y$, by (2); so
\ref{wheelvnbrs} implies that $(C,Y\cup\{y\})$ is a wheel, a contradiction. So $y$ has no neighbour
in $A_0 \setminus c_2$. Next suppose that $y$ is adjacent to $c_2$.  From the symmetry we may
assume that $x_0 \not = c_3$. Let $Q$ be the path of $C \setminus z$
between $x_0,c_3$; so $Q$ has length $>0$, and even length by \ref{evengap}.
Since $x_0 \d Q \d c_3 \d c_2 \d y \d x_0$ is not an odd hole, it follows that $y$ is not adjacent to $x_0$.
But then the hole $x_0 \d Q \d c_3 \d c_2 \d y \d z \d x_0$ is the rim of an odd wheel with hub $Y$,
contrary to $G \in \mathcal{F}_8$. So $y$ is not adjacent to $c_2$. This proves (4).

\bigskip

Let $T$ have vertices $z\d y\d v_1\c v_{n+1}$, where $v_{n+1} \in A_0$.  From (4), $n \ge 1$.
By choosing $T$ of minimum length we may assume that
none of $y,v_1\l v_{n-1}$ have neighbours in $A_0$.
\\
\\
(5) {\em If $n = 1$ then no neighbour of $v_1$ in $A_0$ is $Y$-complete.}
\\
\\
For otherwise we may assume $v_2$ is $Y$-complete. From the symmetry we may assume that $x_0 \not = c_3$.
Let $Q$ be the path of $C \setminus z$ between $x_0,c_3$; so $Q$ has length $>0$, and even length
by \ref{evengap}.
Since $y,v_1$ are not $Y$-complete, there is an antipath
joining them with interior in $Y$, and it is odd since it can be completed to an antihole via
$v_1\d z\d v_2 \d y$. Hence every $Y$-complete vertex is adjacent to one of $y,v_1$, and since
$c_2,c_3$ are $Y$-complete and not adjacent to $y$ by (4), it follows that $v_1$ is adjacent to
$c_2,c_3$. By (3), $v_1$ is adjacent to one of $x_0,x_1$, and so it has two nonadjacent neighbours in
$C$ of opposite wheel-parity. By \ref{wheelvnbrs}, there are 3 consecutive vertices in $C$, all $Y$-complete
and adjacent to $v_1$. By (2), $v_1$ has no other neighbours in $C$.
Hence $x_1 = c_1$ and the neighbours
of $v_1$ in $C$ are $c_1, c_2, c_3$. Consequently $x_0$ is adjacent to $y$; but then
$x_0 \d Q \d c_3 \d v_1 \d y \d x_0$ is an odd hole, a contradiction. This proves (5).
\\
\\
(6) {\em One of $x_0,x_1$ has no neighbours in $\{y,v_1\l v_n\}$.}
\\
\\
For let $P$ be a path $y \d p_1 \c p_k$ from $y$ to some $Y$-complete vertex $p_k \in A_0$, with interior in
$A_0 \cup \{v_1\l v_n\}$, such that $p_k$ is the only $Y$-complete vertex in $P$. Since none of
$y,v_1\l v_{n-1}$ have neighbours in $A_0$ it follows that
$\{y,v_1\l v_n\} \subseteq \{y,p_1\l p_{k-1}\}$. From (5), $k \ge 3$.
Since $G \in \mathcal{F}_8$, $(Y,\{x_0,x_1\}, z\d y \d p_1 \c p_k)$ is not a pseudowheel. But the ends of
the path $z\d y \d p_1 \c p_k$ are $Y$-complete and its internal vertices are not; the path has length
$\ge 4$; $Y, z$ are $\{x_0,x_1\}$-complete, and $y,p_k$ are not. So no other vertices of the path are
$\{x_0,x_1\}$-complete. By \ref{doubleRRC}, applied to the same path and same anticonnected sets,
it follows that one of $x_0,x_1$ is nonadjacent to all of $y,p_1\l p_{k-1}$. Since
$\{y,v_1\l v_n\} \subseteq \{y,p_1\l p_{k-1}\}$, this proves (6).

\bigskip

Let $F = \{y,v_1\l v_n\}$.
From the symmetry we may assume that $x_0$ has no neighbours in $F$. Let $S$ be a
path from $y$ to $x_0$ with interior in $F \cup A_0$. It follows that $S$ has length
$\ge 3$. Let $C'$ be the hole $z\d y \d S \d x_0 \d z$; so $C'$ has length $\ge 6$. Suppose that $x_0$
is different from $c_3$.  Since $(C',Y)$ is not an odd wheel, it follows that $(C',Y)$ is not a wheel, and so
$x_0,z$ are the only $Y$-complete vertex in $C'$. By \ref{RRC}, $Y$ contains a leap or a hat.
A leap would imply there are two vertices in $Y$, joined by an odd path of length $\ge 5$ with interior
in $F \cup A_0$. Hence its ends are $\{x_0,x_1\}$-complete, and its internal vertices
are not, contrary to \ref{RRR}. So $Y$ contains a hat, that is, there exists $y' \in Y$ with no neighbour
in $C'$ except $z,x_0$. But $F \cup A_0$ catches the triangle $\{x_0,y',z\}$;
the only neighbour of $x_0$ in $F\cup A_0$ is its neighbour in $S$, say $s$; the only neighbour
of $z$ in $F\cup A_0$
is $y$; and $s,y$ are nonadjacent, and both nonadjacent to $y'$, contrary to \ref{tricatch}. This proves that
$x_0 = c_3$, and therefore $x_1 \not = c_1$. By exchanging $x_0,x_1$, we deduce that $x_1$ has
a neighbour in $F$. There are therefore two attachments of $F$ in $C$ with
opposite wheel-parity, and two that are nonadjacent. By (1), \ref{wheelFnbrs}, \ref{nokite} and the
optimality of the wheel, and since $x_0 = c_3$ has no neighbour in $F$, it follows that there is a path $R$
between $z,c_2$ with interior in $F$, and no vertex of $C$ has neighbours in the interior of $R$
except $z,c_2$. But then the hole formed by the union of $R$ and the path $C \setminus x_0$ is the rim
of an odd wheel with hub $Y$, a contradiction. This proves \ref{nowheel}.\bbox

\begin{thm}\label{notailsystem}
Let $G \in \mathcal{F}_9$, admitting no balanced skew partition, let  $(z,A_0)$ be a frame and
$x_0\l x_s$ a wheel system with respect to it,
and define $X_i, A_i$ as usual. Then there is no vertex $y \in V(G) \setminus \{z,x_0\l x_s\}$
that is $\{z,x_0\l x_s\}$-complete and has a neighbour in $A_s$.
\end{thm}
\Proof For suppose there is such a frame, wheel system, and $y$, and choose them with $s$ minimum
(it is important here that we minimize over all choices of the frame, not just of the wheel system);
say $(z,A_0)$, $x_0\l x_s$ and $y$ respectively. By \ref{drivingwheel3}, there exists $r$
with $1 \le r < s$, and a vertex $v$ such that $y$ is nonadjacent to $v$ and has no neighbour in $A_r$, and
$v$ is adjacent to $z$, and has a neighbour in $A_r$, and a non-neighbour in $X_r$.
Then $(y,A_0)$ is a frame, and $x_0\l x_r$ is a wheel system with respect to it, and $z$ is
$\{y,x_0\l x_r\}$-complete, and has a neighbour in $A_r'$ (namely $v$), where $A_r'$ is the maximal
connected subset of $V(G)$ including $A_0$ and containing no neighbour of $y$ and no $X_r$-complete
vertex. But this contradicts the minimality of $s$. This proves \ref{notailsystem}.\bbox

Now we can prove \ref{summary}.10, the following.

\begin{thm}\label{notriad}
Let $G \in \mathcal{F}_9$, admitting no balanced skew partition, and let $C$ be a hole in $G$ of
length $\ge 6$. Then there is no vertex
of $G \setminus V(C)$ with three consecutive neighbours in $C$. In particular, every recalcitrant graph
belongs to $\mathcal{F}_{10}$.
\end{thm}
\Proof
Suppose that there is such a vertex, say $y$, and let it be adjacent to $x_0,z,x_1 \in V(C)$, where
$x_0\d z \d x_1$ is a path.
Let $A_0 = V(C) \setminus \{z,x_0,x_1\}$. By \ref{notailsystem} applied to $(z,A_0)$ and $x_0,x_1$,
it follows that $y$ has no other neighbour in $C$. Choose $t$ maximum so that
there is a sequence $x_2\l x_t$ with the following properties:
\begin{itemize}
\item for $2 \le i \le t$, there is a connected subset $A_{i-1}$ of $V(G)$ including $A_{i-2}$, containing a
neighbour of $x_i$, containing no neighbour of $z$ or $y$, and containing no
$\{x_0\l x_{i-1}\}$-complete vertex,
\item for $1 \le i \le t$, $x_i$ is not $\{x_0\l  x_{i-1}\}$-complete, and
\item $x_0\l x_t$ are $\{y,z\}$-complete.
\end{itemize}
Since $G$ admits no skew partition by \ref{oddskew}, there is a path $P$ from $\{z,y\}$ to $A_0$, disjoint
from $\{x_0\l x_t\}$
and containing no $\{x_0\l x_t\}$-complete vertex in its interior. Choose such a path of minimum
length. From the symmetry between
$z,y$ we may assume its first vertex is $y$; say the path is $y\d p_1 \c p_{k+1}$, where
$p_{k+1} \in A_0$.  From the minimality of the length of $P$ it follows that $z$ is not adjacent to
any of $p_2\l p_k$. If $z$ is adjacent to $p_1$ then we may set $x_{t+1} = p_1$, contrary to
the maximality of $t$. So $p_1\l p_{k+1}$ are all nonadjacent to $z$. Hence
$(z,A_0)$ is a frame, and $x_0\l x_t$ is a wheel system with respect to it, and $y$ is adjacent to
all of $z,x_0\l x_t$, and there is a connected subset of $V(G)$ including $A_0$, containing a
neighbour of $y$, containing no neighbour of $z$, and containing no
$\{x_0\l x_t\}$-complete vertex. But this contradicts \ref{notailsystem}. This proves \ref{notriad}.\bbox

This has the following useful corollary, which is \ref{summary}.11.

\begin{thm}\label{noantihole}
Let $G \in \mathcal{F}_{10}$; then $G$ does not contain both a hole of length $\ge 6$ and an antihole
 of length $\ge 6$. In particular, for every recalcitrant graph $G$, one of $G, \overline{G}$ belongs
to $\mathcal{F}_{11}$.
\end{thm}
\Proof
Let $C$ be a hole and $D$ an antihole, both of length $\ge 6$. Let $W = V(C) \cap V(D)$,
$A = V(C) \setminus W$, and $B = V(D) \setminus W$. Let $W,A,B$ have cardinality $w,a,b$ respectively.
Let there be $p$ edges between $A$ and $W$, $q$ edges between $B$ and $W$, $r$ edges between
$A$ and $B$, and $s$ edges with both ends in $W$.
Let there be $p'$ nonedges between $A$ and $W$, $q'$ nonedges between $B$ and $W$,
$r'$ nonedges between $A$ and $B$, and $s'$ nonedges with both ends in $W$.
By \ref{evengap}, and since $G \in \mathcal{F}_{10}$, every vertex in $B$ has at most $\frac{1}{2}(a+w)$
neighbours in $C$, so $q+r \le \frac{1}{2}(a+w)b$. Also, every vertex in $W$ has at most two
neighbours in $A\cup W$, so $p+2s \le 2w$. Summing, we obtain
\[ p + q + r + 2s \le \frac{1}{2}ab + \frac{1}{2}bw + 2w.\]
By the same argument in the complement we
deduce that
\[ p' + q' + r' + 2s' \le \frac{1}{2}ab + \frac{1}{2}aw + 2w.\]
But
\[ p+p' + q + q' + r + r' +  2s + 2s' = ab + aw + bw + w(w-1),\]
so
\[ 4w \ge \frac{1}{2}aw + \frac{1}{2}bw + w(w-1),\]
that is,
\[ w(a + b + 2w - 10) \le 0.\]
Since $a+w,b+w \ge 6$, it follows that $w = 0$, and so $C,D$ are disjoint.
Moreover, equality holds throughout this calculation, so every vertex in $D$ is adjacent to exactly half
the vertices of $C$ and vice versa. By \ref{evengap}, and since $G \in \mathcal{F}_{10}$, it follows that
for each $v \in D$, its neighbours in $C$ are pairwise nonadjacent.
Let $C$ have vertices $c_1\l c_m$ in order, and let $D$ have vertices
$d_1\l d_n$. So for every vertex of $D$, its set of neighbours in $V(C)$ is either the set of all $c_i$
with $i$ even, or the set with $i$ odd, and the same with $C,D$ exchanged. We may assume that
$c_1$ is adjacent to $d_1$. Hence the edges between $\{c_1,c_2,c_4,c_5\}$ and $\{d_1,d_2,d_4,d_5\}$ are
$c_1d_1,c_1d_5, c_2d_2, c_2d_4, c_4d_2, c_4d_4,c_5d_1,c_5d_5$; and so the subgraph induced on these
eight vertices is a double diamond, contrary to $G \in \mathcal{F}_{10}$. This proves \ref{noantihole}.

Let us mention a theorem of \cite{CCZ}, which could be applied at this stage as an alternative to
the next section, the following (and see also \cite{C&C} for some related material):

\begin{thm}\label{CCZthm}
Let $G\in \mathcal{F}_5$. Suppose that for every hole $C$ in $G$ of length $\ge 6$, and every vertex
$V \in V(G) \setminus V(C)$, either:
\begin{itemize}
\item $v$ has $\le 3$ neighbours in $C$, or
\item $v$ has exactly $4$ neighbours in $C$, say $a,b,c,d$, where $ab$ and $cd$ are edges, or
\item $v$ is $V(C)$-complete, or
\item no two neighbours of $v$ in $C$ are adjacent.
\end{itemize}
Suppose also that the same holds in $\overline{G}$. Then either one of $G, \overline{G}$ is
bipartite or a line graph of a bipartite graph, or $G$ admits a loose skew partition.
\end{thm}

The method we give below is somewhat shorter, however.

\section{The end}

We recall that we are trying to prove \ref{recalcitrant}. In view of \ref{noantihole}, it suffices
to show the following, which is \ref{summary}.12, and the objective of the remainder of the paper:

\begin{thm}\label{bipartite}
Let $G \in  \mathcal{F}_{11}$; then either $G$ is complete, or $G$ is bipartite, or $G$ admits a balanced skew partition.
\end{thm}

We begin with a further strengthening of \ref{RRR}, as follows.

\begin{thm}\label{RRRR}
Let $G \in  \mathcal{F}_{11}$, and let $P$ be a path in $G$ with odd length.
Let $X \subseteq V(G)$ be anticonnected, so that both ends of $P$ are $X$-complete. Then
some edge of $P$ is $X$-complete.
\end{thm}
\Proof
Suppose not; then from \ref{RRR}, $P$ has length 3 (let its vertices be $p_1,p_2,p_3,p_4$ in order)
 and $p_2,p_3$ are joined by an antipath $Q$ with interior in $X$. But then
$p_2\d Q \d p_3 \d p_1 \d p_4 \d p_2$ is an antihole of length $>4$, a contradiction. This proves
\ref{RRRR}.\bbox

\begin{thm}\label{pseudohole2}
Let $G \in  \mathcal{F}_{11}$. Let $X \subseteq V(G)$ be nonempty
and anticonnected, and let $p_1\c p_n$ be a path of $G\setminus X$ with $n \ge 4$, such that
$p_1,p_n$ are $X$-complete and $p_2\l p_{n-1}$ are not. There is no vertex
$y \in V(G)\setminus (X \cup \{p_1\l p_n\})$
such that $y$ is $X$-complete and adjacent to $p_1,p_2$.
\end{thm}
\Proof Suppose such a vertex $y$ exists.
By \ref{RRRR}, $n$ is odd, and therefore $n \ge 5$. Let $Q$ be an antipath joining $p_2,p_3$
with interior in $X$. Since $Q$ can be completed to an antihole via $p_3\d p_n \d p_2$, it follows
that $Q$ has length 2, and so there exists $x \in X$ nonadjacent to $p_2,p_3$. Since $x$ is adjacent
to $p_n$, we may choose $i$ with $2 \le i \le n$ minimum so that $x$ is adjacent to $p_i$. Hence
$x\d p_1 \c p_i \d x$ is a hole of length $\ge 6$, and $y$ has three consecutive neighbours in it, contrary
to $G \in  \mathcal{F}_{11}$. This proves \ref{pseudohole2}.\bbox

The next is a strengthening of \ref{tricatch}.

\begin{thm}\label{triplecatch}
Let $G \in  \mathcal{F}_{11}$. Let $X_1,X_2,X_3$ be disjoint nonempty
anticonnected sets, complete to each other. Let $F \subseteq V(G) \setminus (X_1 \cup X_2 \cup X_3)$
be connected, such that for $i = 1,2,3$ there is an $X_i$-complete vertex in $F$. Then there is a vertex
in $F$ complete to two of $X_1,X_2,X_3$.
\end{thm}
\Proof Suppose not; then we may assume $F$ is minimal with this property.
\\
\\
(1) {\em If $p_1\l p_n$ is a path in $F$, and $p_1$ is its unique $X_1$-complete vertex and
$p_n$ is its unique $X_2$-complete vertex then $n$ is even.}
\\
\\
For $n>1$, since no vertex is both $X_1$-complete and $X_2$-complete. Assume $n$ is odd; then
by \ref{doubleRR2}, $n = 3$. But there is an antipath $Q_1$ between $p_2,p_3$ with interior in
$X_1$, and an antipath $Q_2$ between $p_1,p_2$ with interior in $X_2$; and then
$p_2\d Q_1 \d p_3 \d p_1 \d Q_2 \d p_2$ is an antihole of length $>4$, a contradiction. This proves (1).

\bigskip

From the minimality of $F$, there are  (up to symmetry) three cases:
\begin{enumerate}
\item for $i = 1,2,3$ there is a unique $X_i$-complete vertex $v_i \in F$; there is a vertex $u \in F$
different from $v_1,v_2,v_3$, and three paths $P_1,P_2,P_3$ in $F$, all of length $\ge 1$, such that
each $P_i$ is from $v_i$ to $u$, and for $1 \le i < j \le 3$, $V(P_i\setminus u)$ is disjoint from
$V(P_j \setminus u)$ and there is no edge between them, or
\item for $i = 1,2,3$ there is a unique $X_i$-complete vertex $v_i \in F$; there are three paths
$P_1,P_2,P_3$ in $F$, where each $P_i$ is from $v_i$ to some $u_i$ say, possibly of length 0;
and for $1 \le i < j \le 3$, $V(P_i)$ is disjoint from $V(P_j)$ and the only edge between $V(P_i),V(P_j)$
is $u_iu_j$, or
\item for $i = 1,2$ there is a unique $X_i$-complete vertex $v_i \in F$, and there is a path $P$ in
$F$ between $v_1,v_2$ containing at least one $X_3$-complete vertex.
\end{enumerate}

Suppose that the first holds, and let $P_1,P_2,P_3$ be as in the first case. Then some two of $P_1,P_2,P_3$
have lengths of the same parity, and their union violates (1).

Now suppose the second holds, and for $i = 1,2,3$ let $u_i,v_i,P_i$ be as in the second case.
Let $Q_1$ be an antipath joining $u_2,u_3$ with interior in $X_1$,
and define $Q_2,Q_3$ similarly.  If $P_1,P_2,P_3$ all have length 0,
then the union of $Q_1,Q_2,Q_3$ is an antihole of length $>4$, a contradiction. So
we may assume that $P_1$ has length $>0$, and hence $u_1 \not = v_1$. Since $v_1\d u_2 \d Q_1 \d u_3\d v_1$
is an antihole, $Q_1$ has length 1. Since $u_1,u_3$ are not
$X_1$-complete, they are joined by an antipath with interior in $X_1$, and its union with $Q_2$ is an
antihole; so $Q_2$ has length 2, and similarly so does $Q_3$. For $i = 1,2,3$ let $x_i$ be the
middle vertex of $Q_i$. So
\[V(P_1\setminus u_1) \cup V(P_2\setminus u_2) \cup V(P_3\setminus u_3) \cup \{x_1,x_2,x_3\}\]
is connected, and catches the triangle $\{u_1,u_2,u_3\}$; and none of its vertices have two neighbours
in the triangle, and it contains no reflection of the triangle since there is no antihole of length 6.
This is contrary to \ref{tricatch}.

Now suppose the third holds, and let $v_1,v_2,P$ be as in the third case. Let $P$ have
vertices $p_1\l p_n$ where $v_1 = p_1$ and $v_2 = p_n$. Since one of its vertices is $X_3$-complete
and $p_1,p_n$ are not, it follows that $n \ge 3$; and by (1), $n$ is odd, so $n \ge 4$. Choose $i$
minimum and $j$ maximum with $1 \le i,j \le n$ such that $p_i,p_j$ are $X_3$-complete. So $i>1$,
and $i$ is even by (1), and similarly $j<n$ and $j$ is odd. So the path $p_i \c p_j$ has odd length,
and so by \ref{RRRR} one of its edges is $X_3$-complete, say $p_kp_{k+1}$ where $2 \le k \le n-2$.
Now $p_k,p_{k+1}$ are joined by an antipath with interior in $X_1$, and by another with interior
in $X_2$, and the union of these is an antihole; so they both have length 2. Hence for $i = 1,2$ there
exist $x_i \in X_i$ nonadjacent to both $p_k,p_{k+1}$. Let $R$ be a path between $p_{k+2},p_{k-1}$
with interior in $(V(P)\setminus\{p_k,p_{k+1}\}) \cup\{x_1,x_2\}$. Then $R$ can be completed to
a hole $C$ via $p_{k-1}\d  p_k \d p_{k+1} \d p_{k+2}$, and $C$ has length $\ge 6$, and at least one edge of
$C$ is $X_3$-complete, namely $p_kp_{k+1}$, and at least one more vertex of it is $X_3$-complete, since
$R$ uses at least one of $x_1,x_2$. But this contradicts \ref{evengap}, and the hypothesis that
$G \in  \mathcal{F}_{11}$.

This proves \ref{triplecatch}.\bbox

\begin{thm}\label{pseudoholewithhat}
Let $G \in  \mathcal{F}_{11}$, admitting no balanced skew partition. Let $X,Y$ be disjoint anticonnected subsets
of $V(G)$, complete to each other, and let $p_1\c p_n$ be a path of $G \setminus (X \cup Y)$, with
$n \ge 2$, such that $p_1$ is the unique $X$-complete vertex in the path, and $p_n$ is the unique
$Y$-complete vertex. Then there is no $z \in V(G) \setminus (X \cup Y \cup \{p_1\l p_n\})$, complete
to $X \cup Y$ and nonadjacent to $p_1,p_n$.
\end{thm}
\Proof Suppose that $z$ exists, and choose $X$ maximal.
By \ref{bipattach}, there is a path $Q$ in $G$ from $z$ to $p_1$, so that none of its internal
vertices is in $X$ or is $X$-complete. Since no vertex of $\{p_2\l p_n\}$ is $X$-complete, we may
choose $Q$ so that if $z$ has a neighbour in $\{p_2\l p_n\}$ then $V(Q) \subseteq \{z,p_1\l p_n\}$.
The connected subset $V(Q\setminus z)\cup \{p_1\l p_n\}$ ($=F$ say)
contains an $X$-complete vertex, a $Y$-complete vertex, and a $\{z\}$-complete vertex.
The only $X$-complete vertex in $F$ is $p_1$, and that is not $Y$-complete or $\{z\}$-complete; so
by \ref{triplecatch} some vertex in $F$ is $Y$-complete and
adjacent to $z$. If $z$ has a neighbour in $\{p_1\l p_n\}$, then $V(Q) \subseteq \{z,p_1\l p_n\}$,
and so $p_n$ is the only vertex of $F$ that
is $Y$-complete; and it is not adjacent to $z$, a contradiction. So $z$ has no neighbour in $\{p_1\l p_n\}$,
and therefore
only one vertex in $F$ is adjacent to $z$, the neighbour of $z$ in $Q$, say $q$.
Hence $q$ is nonadjacent to $p_1$, for otherwise we could add $q$ to $X$, contrary to the maximality
of $X$. Consequently $Q$ has length $>2$. This contradicts \ref{pseudohole2} applied to $Q$,$X$ and
any vertex $y \in Y$. This proves \ref{pseudoholewithhat}.\bbox

We deduce

\begin{thm}\label{holewithhat}
Let $G \in  \mathcal{F}_{11}$, admitting no balanced skew partition, and let $C$ be a hole.
If $z \in V(G) \setminus V(C)$
has two neighbours in $C$ that are adjacent, then $C$ has length $4$ and $z$ has a third neighbour in $C$.
In particular, $G$ has no antipath of length $4$.
\end{thm}
\Proof Let $C$ be the hole with vertices $p_1\l p_{n+2}$ in order, and assume some $z \in V(G) \setminus V(C)$
is adjacent to $p_{n+1},p_{n+2}$. By \ref{pseudoholewithhat}, taking $X = \{p_{n+1}\}$ and
$Y = \{p_{n+2}\}$ we deduce that $z$ is adjacent to at least one of $p_1, p_n$. Since $G \in  \mathcal{F}_{11}$
it follows that $C$ has length $4$. This proves \ref{holewithhat}.\bbox

\begin{thm}\label{doublecatch}
Let $G \in  \mathcal{F}_{11}$, admitting no balanced skew partition. Let $X_1,X_2,X_3$ be pairwise disjoint,
nonempty, anticonnected subsets of $V(G)$, complete to each other. Let
$F \subseteq V(G) \setminus (X_1 \cup X_2 \cup X_3)$ be connected, so that for at least two
values of $i \in \{1,2,3\}$, every member of $X_i$
has a neighbour in $F$. Then some vertex of $F$ is complete to two of $X_1,X_2,X_3$.
\end{thm}
\Proof Assume not, and choose a counterexample with $X_1\cup X_2 \cup X_3\cup F$ minimal.
Suppose $F$ contains an $X_i$-complete vertex for two values of $i \in \{1,2,3\}$, say $i = 1,2$; and choose
a path $p_1\c p_n$ of $F$ such that $p_1$ is $X_1$-complete and $p_n$ is $X_2$-complete, with $n$
minimum. So $n \ge 2$. From the minimality of $F$, $F = V(P)$, and there is a vertex $x_1\in X_1$ such that
$p_1$ is its only neighbour in $F$, and there exists $x_2 \in X_2$ so that $p_n$ is its only neighbour
in $F$. By \ref{holewithhat} applied to the hole $x_2\d x_1 \d p_1 \c p_n \d x_2$ and any $x_3\in X_3$,
it follows that $n = 2$. Let $Q$ be an antipath between $p_1,p_2$
with interior in $X_3$; since $p_1$ has a nonneighbour $x_2\in X_2$, and $p_2$ has a nonneighbour
$x_1 \in X_1$, it follows that $x_1\d p_2 \d Q \d p_1 \d x_2$ is an antipath of length $\ge 5$,
contrary to \ref{pseudoholewithhat}.

So there is at most one $i$ so that $F$ contains $X_i$-complete vertices, and
from the symmetry we may assume that $F$ contains no  $X_1$- or $X_2$-complete vertices.
We may also assume that all members of $X_1$ have neighbours in $F$, and therefore
$|X_1| \ge 2$; choose distinct $x_1,x_1' \in X_1$ so that
$X_1 \setminus x_1,X_1 \setminus x_1'$ are both anticonnected. From the minimality of $X_1$, there is a vertex
$f$ of $F$ complete to two of $X_1 \setminus x_1,X_2,X_3$, and therefore complete to $X_1 \setminus x_1$
and $X_3$, and similarly a vertex $f'$ of $F$ complete to $X_1 \setminus x_1'$ and $X_3$. Let $P$
be a path in $F$ between $f,f'$. Since all vertices of $X_1\cup X_3$ have neighbours in $V(P)$, the
minimality of $F$ implies that $F = V(P)$; and moreover, since all vertices of $(X_1\setminus x_1)\cup X_3$
are adjacent to $f$, the minimality of $F$ implies that $f'$ is the unique neighbour of $x_1$ in $F$.
Similarly $f$ is the unique neighbour of $x_1'$ in $F$.  Let $Q$ be an
antipath in $X_1$ joining $x_1,x_1'$. Since $f$ has a nonneighbour $x \in X_2$,
$x\d f \d x_1 \d Q \d x_1'$ is an antipath, and so $Q$ has length 1, and hence $x_1,x_1'$ are nonadjacent.
From the minimality of $F$, there exists $x_2 \in X_2$ with no neighbour in $F \setminus f$. If $x_2$
is also nonadjacent to $f$, then $x_2\d x_1 \d f' \d P \d f \d x_1' \d x_2$ is a hole of length $\ge 6$,
and any member of $X_3$ has three consecutive neighbours on it, contrary to $G \in  \mathcal{F}_{11}$.
But then $x_1$ has two consecutive neighbours on the hole $x_1' \d f \d P \d f' \d x_1 \d x_2 \d x_1'$,
and this hole has length $>4$, contrary to \ref{holewithhat}.
This proves \ref{doublecatch}.\bbox

Now we can complete the proof of \ref{bipartite}, and hence of \ref{recalcitrant} and therefore
of \ref{decomp} and \ref{SPGC}, as follows.

\noindent{\bf Proof of \ref{bipartite}.}
Let $G \in  \mathcal{F}_{11}$, admitting no balanced skew partition.
Suppose that $G$ is not bipartite; hence it has a triangle, and so we may choose disjoint nonempty
anticonnected sets $X_1\l X_k$, complete to each other, with $k \ge 3$. Choose these with maximal union.
Let $N$ be the set of all $X_k$-complete vertices in $G$. Suppose first that $X_k \cup N = V(G)$.
By \ref{bipattach}, either $G$ is complete (and the theorem holds), or $\overline{G}$ has exactly
two components, both with $\le 2$ vertices (and therefore $G$ is bipartite and again the theorem holds).

So we may assume that $X_k \cup N \not = V(G)$.  By
\ref{bipattach}, the set $V(G) \setminus (X_k\cup N)$ ($=F$ say) is connected and every vertex of $N$ has
a neighbour in it, and in particular, all vertices of $X_1\cup X_2$ have a neighbour in it. By
\ref{doublecatch}, some vertex $v \in F$ is complete to two of $X_1,X_2,X_3$. We may assume that
$v$ is complete to $X_1\l X_i$ say where $2 \le i \le k$, and not complete to $X_{i+1} \l X_k$. Define
\[X_{i+1}' = X_{i+1} \cup \c \cup X_k \cup \{v\};\]
then the sets $X_1\l X_i, X_{i+1}'$ violate the optimality of the choice of $X_1\l X_k$. This proves
\ref{bipartite}.\bbox

\section{Acknowledgements}

We worked on pieces of this with several other
people, and we would particularly like to thank Jim Geelen, Bruce Reed, Chunwei Song and
Carsten Thomassen for their help.

We would also like to acknowledge our debt to Michele Conforti, G\'{e}rard Cornu\'{e}jols
and Kristina Vu\v{s}kovi\'{c}; their pioneering work was very important to us, and
several of their ideas were seminal for this paper. In particular, they made the basic
conjecture \ref{decomp}, and suggested that Berge graphs containing wheels should have skew partitions
induced by the wheels, and Vu\v{s}kovi\'{c}
convinced us that prisms were interesting, and a theorem like \ref{longprism} should exist.
In addition, they proved a sequence of steadily improving theorems, saying that any
counterexample to the strong perfect graph conjecture must contain subgraphs of certain kinds
(wheels and parachutes), and those gave us a great incentive to work from the other end, proving that
minimum counterexamples cannot contain subgraphs of certain kinds.

Finally, we would like to acknowledge the American Institute of Mathematics, who
supported two of us full-time for six months to work on this project.


\begin{thebibliography}{99}
\bibitem{Berge} C. Berge, ``F\"{a}rbung von Graphen, deren s\"{a}mtliche bzw. deren ungerade Kreise
starr sind'', {\em Wiss. Z. Martin-Luther-Univ. Halle-Wittenberg Math.-Natur. Reihe}
10 (1961), 114.
\bibitem{starcut} V. Chv\'{a}tal, ``Star-cutsets and perfect graphs'', {\em J. Combinatorial Theory,
Ser. B}, 39 (1985), 189-199.
\bibitem{M-join} V. Chv\'{a}tal and N. Sbihi, ``Bull-free Berge graphs are perfect'', {\em
Graphs and Combinatorics}, 3 (1987), 127-139.
\bibitem{C&C} M. Conforti and G. Cornu\'{e}jols, ``Graphs without odd holes, parachutes or proper
wheels: a generalization of Meyniel graphs and of line graphs of bipartite graphs'',
manuscript, November 1999.
\bibitem{CCV} M. Conforti, G. Cornu\'{e}jols and K. Vu\v{s}kovi\'{c}, ``Square-free perfect graphs'',
manuscript, Feb. 2001.
\bibitem{CCVZ} M. Conforti, G. Cornu\'{e}jols, K. Vu\v{s}kovi\'{c} and G. Zambelli,
``Decomposing Berge graphs containing proper wheels'', manuscript, March 2002.
\bibitem{CCZ}  M. Conforti, G. Cornu\'{e}jols and G. Zambelli, ``Decomposing Berge graphs containing
no proper wheels, stretchers or their complements'', manuscript, May 2002.
\bibitem{2-join} G. Cornu\'{e}jols and W.H. Cunningham, ``Compositions for perfect graphs'',
{\em Discrete Math. } 55 (1985), 245-254.
\bibitem{Lovasz} L. Lov\'{a}sz, ``A characterization of perfect graphs'', {\em J. Combinatorial Theory,
Ser. B}, 13 (1972), 95-98.
\bibitem{Reed} J. Ramirez-Alfonsin and B. Reed (eds.), {\em Perfect Graphs}, Wiley, Chichester, 2001.
\bibitem{R&R} F. Roussel and P. Rubio, ``About skew partitions in minimal imperfect graphs'',
{\em J. Combinatorial Theory, Ser. B}, 83 (2001), 171-190.
\bibitem{Tutte} N. Robertson, P. Seymour and R. Thomas, ``Tutte's edge-colouring conjecture'',
        {\em J. Combinatorial Theory, Ser. B}, 70 (1997), 166 - 183.
\bibitem{tu} P. Seymour, ``Decomposition of regular matroids'',
                   {\em J. Combinatorial Theory, Ser. B,} 28 (1980), 305-359.
\bibitem{Pfaffian} N. Robertson, P. Seymour and R. Thomas, ``Permanents, Pfaffian
        orientations, and even directed circuits'', {\em Annals of Math.}, 150 (1999), 929-975.
\bibitem{Wagner} K. Wagner, ``\"{U}ber eine Eigenschaft der ebenen Komplexe'', {\em Math. Ann.}
114 (1937), 570-590.
\end{thebibliography}
\end{document}